\documentclass[11pt]{amsart}
\usepackage{amsmath,amsfonts,amssymb,amscd,amsthm,amsbsy}
\numberwithin{equation}{section}
\usepackage{setspace}

\usepackage{setspace}
\usepackage{amssymb,amsfonts,amsthm,amscd,stmaryrd,dsfont,esint,upgreek}
\reversemarginpar
\def\nl{\newline}
\newcommand{\beq}{\begin{equation}}
\newcommand{\Eq}{\end{equation}}
\newtheorem{thm}{Theorem}[section]
\newtheorem{lem}{Lemma}[section]

\newtheorem{cor}{Corollary}[section]
\newtheorem{prop}{Proposition}[section]
\newtheorem{defn}{Definition}[section]

\def\Binom#1#2{\begin{matrix} #1\\ #2\end{matrix}}
\def\un#1{\underline{#1}}
\def\sbullet{\ {\scriptstyle\bullet}\ }
\def\bcite{[~\quad~]}
\def\ds{\displaystyle}
\def\bone{\text{\bf 1}}
\def\ep{\varepsilon}
\def\C{{\mathbb C}}
\def\F{{\mathbb F}}
\def\N{{\mathbb N}}
\def\nat{{\mathbb N}}
\def\R{{\mathbb R}}
\def\1{{\mathbb 1}}
\def\tee{{\mathbb T}}
\newcommand{\T}{{\mathbb T}}
\def\Z{{\mathbb Z}}
\def\A{{\mathcal A}}
\def\FF{{\mathcal F}}
\def\M{{\mathcal M}}
\def\P{{\mathbb P}}
\def\E{{\mathbb E}}
\def\loc{\operatorname{loc}}
\def\tr{\operatorname{tr}}
\def\esssup{\mathop{\rm ess\, sup}\limits}
\def\circW{{\buildrel \bullet\over W}}
\def\circX{{\buildrel \bullet\over X}}
	\def\circW{{\dot W}}
	\def\circX{{\dot X}}
\def\Circ#1{{\buildrel \circ\over #1}}
\def\CCirc#1{{\buildrel {\circ\circ}\over #1}}
\def\CCCirc#1{{\buildrel {\circ\circ\circ}\over #1}}
\def\Inf{\mathop{\rm inf\vphantom{p}}\limits}
\def\Partial#1#2{\frac{\partial #1}{\partial #2}}
\def\circfrown#1{\mathop{#1}\limits^{\dotfrown}}
\def\dotfrown{{\buildrel \circ\over{\frown}}}
\def\Circfrown#1{\mathop{#1}\limits^{\Dotfrown}}
\def\Dotfrown{{\dot{\frown}}}
\def\Bigfrown#1{\mathop{#1}\limits^{\dot{\text{\LARGE$\frown$}}}}
\def\dothat#1{{\buildrel \circ\over{\widehat{#1}}}}
\def\Frown#1{\mathop{#1}\limits^{\frown}}

\newcommand{\on}{\operatorname}

\newcommand{\I}{\on{I}}
\newcommand{\II}{\on{II}}

\newcommand{\abs}[1]{ \left| #1 \right|}
\newcommand{\dist}{\on{dist}}
\newcommand{\eps}{\varepsilon}
\newcommand{\ind}{\mathbf{1}}
\newcommand{\ip}[1]{\langle #1 \rangle}
\newcommand{\Lip}{\on{Lip}}
\newcommand{\norm}[1]{ \left\| #1 \right\| }
\newcommand{\nor}[2]{\left\|#1\right\|_{#2}}
\newcommand{\oline}[1]{\overline{#1}}
\newcommand{\oo}{\infty}
\newcommand{\pars}[1]{\left(#1\right)}
\newcommand{\uline}[1]{\underline{#1}}

\newcommand{\mfk}{\mathfrak}
\newcommand{\mcl}{\mathcal}
\newcommand{\mbb}{\mathbb}
\newcommand{\mbf}{\mathbf}
\title[Pathwise solutions nonlinear equations rough time dependence]{Pathwise solutions for fully nonlinear first- and second-order partial differential
equations with multiplicative rough time dependence}
\author[Panagiotis E. Souganidis]
{Panagiotis E. Souganidis}
\address{Department of Mathematics, University of Chicago, Chicago, Illinois 60637, USA}
\email{souganidis@math.uchicago.edu}
\thanks{  Partially supported by the National Science
Foundation Grants DMS-1266383 and DMS-1600129, the Office for Naval Research grant N000141712095 and the Air Force Office for Scientific Research grant FA9550-18-1-0494.}
\dedicatory{Version: \today}

\begin{document}

%\author{Panagiotis E. Souganidis$^{(1,2)}$}
%\thanks{$^{(1)}$ Department of Mathematics, The University of Chicago,
%Chicago, IL 60637, USA,\newline e-mail: souganidis@math.uchicago.edu}
%\thanks{$^{(2)}$ Partially supported by the NSF grants DMS-1266383 and DMS-160012
%\dedicatory}%{Version: \today}

%%%%\rangle p%%%%%%%%%%%
\begin{abstract}
The notes are an overview of the theory of pathwise  weak solutions of   two classes of scalar fully nonlinear first- and second-order degenerate  parabolic partial differential  equations with multiplicative rough time dependence, a special case being Brownian. These  are Hamilton-Jacobi,   Hamilton-Jacobi-Isaacs-Bellman  and quasilinear divergence form equations including  multidimensional scalar conservation laws. If the time dependence is ``regular'', the weak solutions are respectively the viscosity and  entropy/kinetic solutions. The main results are the well-posedness and qualitative properties of the solutions. 
Some concrete applications are also discussed. 
%The results presented here are part of the ongoing development of the theory in collaboration with  P.-L. Lions.
%The development of the theory for the first class of equations  is based on joint work with P.-L. Lions. 
%The results about  quasilinear divergence form equations  are based on joint work with  P.-L. Lions, B. Perthame and B. Gess.
\end{abstract}

\maketitle      %% SHORT TITLE FOR RUNNING HEAD
%\markboth{Fully nonlinear first- and 
%second-order stochastic pde}

\baselineskip=14pt      %% DRAFT double-space mode

\setcounter{section}{-1}
\section{Introduction}

I present  an overview of the theory of pathwise  weak solutions of  two classes of scalar fully nonlinear first- and second-order degenerate  parabolic (stochastic) partial differential  equations (spde for short) with multiplicative rough  time dependence, a special case being Brownian. These  are Hamilton-Jacobi,  Hamilton-Jacobi-Isaacs-Bellman   and quasilinear divergence form partial differential equations (pde for short) including  multidimensional scalar conservation laws. If the time dependence is ``regular'', the weak solutions are respectively the viscosity and  entropy/kinetic solutions. The main results are the well-posedness and qualitative properties of the solutions. Some concrete applications are also discussed both to motivate as well as to show the scope of the theory. 
Most of the results presented here are part of the ongoing development of the theory in collaboration with  P.-L. Lions \cite{lionssouganidis1, lionssouganidis2, lionssouganidis3, lionssouganidis4, lionssouganidis5, lionssouganidis6, lionssouganidisbook, lionssouganidis7}. 
%The development of the theory for the first class of equations  is based on joint work with P.-L. Lions. 
The results about  quasilinear divergence form equations  are based on joint work with  P.-L. Lions, B. Perthame and B. Gess \cite{lionsperthamesouganidis1, lionsperthamesouganidis2, lionsperthamesouganidis3, gesssouganidis1, gesssouganidis2, gesssouganidis3, gessperthamesouganidis}. 

\smallskip

Problems of the type discussed here arise in several applied contexts and models  for a wide variety of phenomena and applications
including mean field games, turbulence, phase transitions and front propagation with 
random velocity, nucleations in physics, macroscopic limits of particle systems, pathwise stochastic control theory,
stochastic optimization with partial observations, stochastic selection, etc.. 
\smallskip

The general classes of evolution equations considered in these notes are 
\begin{equation}\label{eq0.1}
du = 
F (D^2 u, D u, u,x,t) dt + \sum_{i=1}^m H^i(Du,u, x,t) \cdot dB_i \ 
\text{in}  \ Q_T:=\R^d\times (0,T], % \R^d \times (0,T]\ ,
\end{equation}
and 
\begin{equation}\label{quasilinear}
du + \sum_{i=1}^{d}\partial_{x_{i}}(A^{i}(u,x,t))\cdot dB_{i}-\text{div}(A(u,x,t)Du)dt=0 \ \text{ in }   \ Q_T,%\R^{N}\times(0,T],
\end{equation}
with initial condition 
\begin{equation}\label{eq0.2}
u(\cdot,0) = u_0  \  \text{on} \ \R^d.
\end{equation}

Here $F=F(X,p,u,x,t), H^1=H^1(p,u,x,t), \ldots,  H^m=H^m(p,u,x,t), A^1=A^1(u,x,t),\dots, A^d=A^d(u,x,t)$ and  $A=A(u,x,t)$ are (at least) continuous functions of their arguments (exact assumptions will be shown later),   $F$ and $A$ are respectively degenerate elliptic in $X$ and monotone in $u$, $B:=(B_1,\ldots, B_m)$ and $B=(B_1,\ldots, B_d)$ are, for example,    continuous geometric rough  in time  and ``$\cdot$'' simply denotes the way $B$ acts on the $H^i$ and $A^i$.  When $B$  is a Brownian path,  ``$\cdot$''  becomes the usual  Stratonovich differential ``$\circ$'', something justified by the fact  that  the pathwise solutions  may be obtained  as the limit of solutions of equations with smooth signals. %$B$ replaced by smooth approximations.
% and  $Q_T:=\R^d \times (0,T]$ and $Q_T:=\R^d \times (0,\infty)$. 
The $B_i$'s  can be taken to be  approximations of ``colored white noise.''  For simplicity, below we assume that any spatial dependence on the signal $B_i$ is part of 
%that is
%\begin{equation*}
%B(x,t) = \sum_{i=1}^m c_i (x)B_i,
%\end{equation*}
%where the $c_i$'s are smooth in $x$. Given the notation used here, however, the $c_i$'s are 
$H^i$ and the  $A^i$. Finally, $Q_\oo:=\R^d \times (0,\oo).$

\smallskip

When $B$ is either smooth or has bounded variation, then ``$d$'' is the regular time derivative and \eqref{eq0.1} and \eqref{quasilinear} 
are ``regular''  equations,  which have been  studied using respectively the viscosity  and entropy/kinetic theories. When the driving signals are regular (``non rough''),  I refer to the equations as   ``deterministic'' or ``non-rough''. If the signals are ``rough'', the equations will be called ``rough'' or ``stochastic'' when the path is Brownian.
% In the rest of the paper, I will be using the terms deterministic and/or non-rough when the driving signals are smooth. 
\smallskip 

The theory presented in these notes is a pathwise one and  simply treats $B$ as the time derivative of a continuous function.
When the $H^i$'s and $A^i$'s  are respectively  independent of $(u,x)$ and $x$, the general qualitative theory does not need any other assumption but continuity. When there is spatial dependence,  then it is necessary to argue differently. %  and the theory of rough paths comes in handy. % as it will become clear later.

\smallskip

There is a vast literature for linear and quasilinear versions of \eqref{eq0.1} as well as work for some versions of \eqref{quasilinear}. Listing all the references is not possible in this introduction.  Some  connections are made in main the body of the notes.

%\smallskip

\subsection*{Organization of the notes} Concrete examples where \eqref{eq0.1} and \eqref{quasilinear} arise are presented in Section~1. 
Section~2 discusses the main difficulties  and explains why the Stratonovich formulation is more appropriate. Sections 3 to Section 13 are devoted to the pathwise solutions of Hamilton-Jacobi and Hamilton-Jacobi-Isaacs equations. In Section 3, I present new results about nonlinear equations with linear rough path dependence,  I introduce the system of characteristics, and I discuss a short time classical result about stochastic Hamilton-Jacobi equations in the smooth regime. Section 4 is about fully nonlinear equations with semilinear rough path dependence. Section 5 is about formulae or the lack thereof for Hamilton-Jacobi equations with time dependence. Section 6  discusses the simplest possible nonlinear pde with rough time signals as the limit of regular approximations. Section 7 is about pathwise solutions of  nonlinear first-order pde with nonsmooth Hamiltonians and rough signals. In Section 8, I present new results about the qualitative properties of the pathwise solutions.  Section 9 is devoted to the well-posedness theory of the pathwise solutions with  spatially depended  $H^i$'s. Section 10 is about  Perron's method, while Section 11 discusses the convergence of approximation schemes with error estimates. In Section 12 I present new results about the homogenization of pathwise solutions. Section 13 is about  the asymptotics of stochastically perturbed reaction-diffusion equations.  The results  about quasilinear divergence form equations including multi-dimensional stochastic conservation laws are presented in Section 14. Finally, the Appendix  summarizes few  basic things from the classical theory of viscosity solutions that are used in the notes.
%and present the proofs of some results used earlier.  

%In Section 2 we discuss the main difficulties associated 
%with equations like \eqref{eq0.1} with $\zeta$ as in \eqref{eq0.4}. 
%Here we only remark that the assumptions on $\zeta$ and $Z$ 
%put \eqref{eq0.1} outside the reach of the viscosity theory so far, 
%which up to now  considered equations with 
%at most $L^1$-dependence in time, i.e., $Z$ of bounded variation. 
%This work provides  an extension of 
%the theory of viscosity solutions to more general time dependence and 
%suggests a possible strategy to consider even more general equations.

%%%%%%%%%%%%%%%%%%%%%%%%
\section{Motivation and some examples}\label{sec:motivations}
\smallskip

%I  present here examples that in different contexts lead to 
%equations like \eqref{eq0.1} and \eqref{quasilinear}.
%\smallskip 
A discussion follows about a number of results that have been or may be solved using the theory presented in here. %It should be noted, that there are many open questions related to the problems discussed here. 
In several places, to keep  the discussion simple,  the presentation  is informal. 
%n order to keep the presentation simple.
%the discussion a bit informal.
% the discussion in this section what follows to keep the presentation simple in many places the discussion is somehow informal.
%\smallskip

\subsection*{Motion of interfaces} An important  question in pde and geometry as well as applications like phase transitions is the understanding of the long 
time behavior of solutions of reaction-diffusion equations and the properties of the developing interfaces, which separate the regions where the solutions approach the different equilibria of the equation. 
\smallskip

A classical  and well studied problem in this context is the asymptotic behavior of the solution $u^\ep$ to 
the so called  Allen-Cahn equation 
\begin{equation*}%\label{AC}
u^\ep_t-\Delta u^\ep + \frac{1}{\ep^2}W'(u^\ep)=0 \  \text{in}  \ Q_T,
\end{equation*} 
where $W:\R\to \R$ is a double-well potential with wells of equal depth located at, for example, at $\pm 1$. It is well known that as, $\ep \to 0$, $u^\ep \to \pm1$  inside and outside an interface moving with normal velocity 
$V=-\kappa$, where $\kappa$ is the mean curvature. The interface is the zero-level set of the solution of the level-set pde
\begin{equation}\label{ls}
v_t =\left(I-\frac{Dv}{|Dv|}\otimes \frac{Dv}{|Dv|}\right):D^2v   \  \text{in} \ Q_T,
%v_t - \text{tr}\left[(I-\frac{Dv}{|Dv|}\otimes \frac{Dv}{|Dv|})D^2v\right]  =0 \ \text{in} \ Q_T,
\end{equation}
where for $A,B \in {\mathcal S}^d$, the space of symmetric $d\times d$ matrices,  $A:B:= \text{tr} (AB)$ and $I$ is the identity matrix in $\R^d$.
\smallskip

For the applications, however, it is  interesting to consider potentials with wells at locations which change with the scale $\ep$ and to identify the exact scaling at which something nontrivial comes up. An example of such a problem is 
\begin{equation*}%\label{AC}
u^\ep_t-\Delta u^\ep+\frac{1}{\ep^2}(W'(u^\ep) + \ep c(t))=0  \  \text{in} \  Q_T,
\end{equation*} 
for some smooth function $c=c(t)$, which leads, as $\ep \to 0$, to an interface moving with normal velocity  $V=-\kappa + \alpha c(t) $, where   $\alpha \in \R$  is a ``universal'' constant which is independent of $c$. 
\smallskip

A natural question is what happens if $c$ is irregular and, in particular, if $c=dB$, where $B$ is a Brownian path. Note that such perturbations often appear  in the hydrodynamic limit of interacting particle systems.
It turns out that in this case the oscillations of the wells due to $dB$ are too strong for the system to stabilize.  However, as it was it was shown by Lions and Souganidis \cite{lsallencahn},
%without any regularity restrictions by Lions and Souganidis \cite{lsallencahn},  
if $B$ is replaced by a ``mild'' approximation $B^\ep$, then the asymptotic interface moves with normal velocity
\[V=-\kappa + \alpha dB,\]
and is characterized as a level set of the solution of the ``stochastic'' level-set pde
\begin{equation}\label{stochastic level set}
dv = \left[\left(I-\frac{Dv}{|Dv|}\otimes \frac{Dv}{|Dv|}\right):D^2v\right]dt + \alpha |Dv|\cdot dB \ \text{in} \ Q_T.
\end{equation}

More details including references as well as a sketch of the proof of the result in \cite{lsallencahn} are presented in Section 13. 

%The asymptotic behavior of the solutions to the perturbed Allen-Cahn equation with $c=dB$ and $B$ a space time Brownian motion was conjectured by Otha, Jasnow and Kawasaki \cite{othajasnowkawasaki} , while the unperturbed equation was proposed by Allen and Cahn \cite{allencahn}  as a model to study phase transitions. The rigorous justification of the conjectured behavior for the latter as well as its perturbed version with $c$ smooth was obtained by Evans, Soner and Souganidis \cite{evanssonersouganidis} and Barles and Souganidis \cite{barlessouganidis};  see Souganidis \cite{bardicime, montrealfronts} for a  comprehensive overview of the theory. For the former problem with $c=dB^\ep$ and $B^\ep$ a mild approximation of  a time Brownian motion, Yip \cite{yip} and Funaki \cite{funaki} obtained results for short time and convex initial interfaces. Lions and Souganidis \cite{lionssouganidisbook} proved the full global in time result in this case and they also showed that the behavior conjectured in \cite{othajasnowkawasaki} 
%cannot be correct if $c = dB$ and $B$ a Brownian motion in time--(see Section The explanation for the latter is that the oscillations due to the presence of $dB$ are so strong that they interfere with the stability properties of the equilibria of the potential. 

\subsection*{ A stochastic selection principle}  A classical question in the theory of level set  interfacial motions  is whether there is ``fattening'', that is, if there are configurations (initial data) such that the zero level set of the solution $v$ to \eqref{ls} develops interior. For the motion by mean curvature, it is known that, if the initial configuration  is two touching balls, then, for positive times, the evolving front is  a ``surface'' that looks like the boundary of either two separated shrinking balls or some connected open set which moves in time, and there are well defined minimal and maximal moving boundaries. %The reason behind such behavior is that the curvature is not well defined at the touching point and the initial surface can be thought as the boundary of either two balls which are infinitesimally separated or an open set which is ``almost'' pinched. 
%\smallskip 
\vskip.125in

As it is often the case the introduction of stochasticity resolves this ambiguity and provides a definitive selection principle. Indeed, it was proved by Souganidis and Yip \cite{souganidisyip} without any regularity restrictions on the evolving set (see also  Dirr,  Luckhaus and Novaga \cite{DLN} for a short time result), that the zero level sets of the  solutions $v^{\pm \ep}$ of the stochastically perturbed level set pde 
%\begin{equation}\label{st}
\[dv^{\pm \ep} =\left[ \left(I-\frac{Dv^{\pm \ep}}{|Dv^{\pm \ep}|}\otimes \frac{Dv^{\pm \ep}}{|Dv^{\pm \ep}|}\right):D^2v^{\pm \ep}\right] dt  \pm\ep |Dv^{\pm \ep}|\circ dB  \  \text{in}  \  Q_T,\]
with initial data two touching balls, never develop interior and, as $\ep \to 0$, converge in the Hausdorff distance to the maximal interface of the unperturbed problem. % described above.
\smallskip

\subsection*{Pathwise stochastic control theory}
To keep the notation simple I assume here that $d=1$. A typical stochastic control  problem with finite horizon $T>0$ consists of 

(i)~a controlled stochastic 
differential equation (sde for short)
\begin{equation*}
%\begin{cases}
dX_s =b(X_s,\alpha_s)ds  + \sqrt 2 \sigma_1(X_s,\alpha_s)dB_{1,s} +  \sqrt{2} \sigma_2(X_s)\circ dB_{2,s}  \quad (  0\leq t\leq s\leq  T) \quad X_t=x, %\\[1.2mm]
%\end{cases}
\end{equation*}
where $(B_{1,t})_{t\geq 0}$ and $(B_{2,t})_{t\geq 0}$ are two independent Brownian motions with respective filtrations $ ({\mathcal F}^{B_1}_t)_{t\geq 0}$ and   $({\mathcal F}^{B_2}_t)_{t\geq 0}$, $(\alpha_t)_{t \geq 0} \in \mathcal A$, 
%\begin{equation*}
%dX_t = b(X_t,\alpha_t)dt  +  \sqrt{2} \sigma (X_t) \circ dB \ \quad \ \  X_0=x,
%\end{equation*}
the set of admissible ${\mathcal F}^{B_1}_t$-progressively measurable controls with values in $A$ a subset of 
some $\R^k$, 
%satisfying appropriate measurability properties and 
and 
\smallskip

(ii)~a pay-off functional, which, to simplify the presentation, here is taken to be
\begin{equation*}
J(x,t; \alpha)=E_{x,t} [g(X_T)  | {\mathcal F}^{B_2}_T],
\end{equation*}
the goal being  to minimize the pay-off over $\mathcal A.$
\smallskip

The associated value function, which is defined by 
\[u(x,t) = \text{essinf}_{\alpha \in \mathcal A} J(x,t;\alpha),  
\]
has been shown  in Lions and Souganidis \cite{lionssouganidis2, lionssouganidisbook} (see also Buckdahn and Ma \cite{buckdahnma3} for a special  case) to be the pathwise solution of the stochastic associated Bellman equation 
\[%\begin{cases}  
du + \underset{\alpha \in A}\inf
% _{\pi \in A}
 \left[ \sigma_1^2 (x,\alpha) u_{xx} + b(x,\alpha) u_x  \right] \ dt + \sqrt 2 \sigma_2(x) u_x \circ dB_2=0 \ \text{in} \ Q_T \quad  u(\cdot, T)= g, 
\]
%
%\begin{equation*}
%du +  \sup_{\alpha\in A} (-b(x,\alpha) \cdot Du)dt  - \sqrt{2}  \sigma (x)Du\circ dB=0 \  \text{in}   \  Q_T  \quad  u(\cdot,0)=g,
%\end{equation*}
which is a special case of \eqref{eq0.1} with  $F$ nonlinear and $H$  linear; notice that to be consistent with control theoretic formulation of the problem the equation is written backwards in time.  
 \smallskip

The aim of the  classical stochastic control theory with the stochastic dynamics above, is to minimize over $\mathcal A$ the ``averaged'' payoff 
\[\overline J(x,t; \alpha)=E_{x,t} [g(X_T)].
\]
It is a classical fact that the value function 
\[\overline u(x,t) =\text{essinf}_{\alpha \in \mathcal A} \overline J(x,t;\alpha )\]
is the unique viscosity solution of the deterministic Bellman terminal valued problem
\[  \overline u_t + \inf_{\alpha \in A}\left[  \left(\sigma_1^2(x,\alpha) + \sigma_2^2(x) \right) {\overline u}_{xx} +  \left (b(x,\alpha) + \sigma_{2,x} \sigma_2(x) \right ) {\overline u}_x\right]=0 \ \text{in} \ Q_T \quad 
 u(\cdot, T)= g.
 \]

\subsection*{Mean field games}  A typical example of the Lasry-Lions mean field theory \cite{lasrylions1, lasrylions2, lasrylions3}  is the study of the asymptotic  behavior, as $L\to \infty$, of the law $\mathcal{L}(X^1_t, \dots, X^L_t)$ of the solution of the sde
\[dX^i = \sigma \left(X^i, \frac{1}{L-1} \sum_{j\not=i} \delta_{X^j}\right) \circ dB \qquad \qquad (i=1,...,L).\]
Here  $ \delta_y$ is the  Dirac mass at $y$ and  $\sigma \in C^{0,1}(\R^d \times \mathcal{P}(\R^d); \mathcal{S}^d)$, %$%\mathcal{S}^N$ and 
 $\mathcal{P}(X)$ being the set   of 
% symmetric $N \times N$ matrices and  
 probability  measures on $X$. %\ \ $\sum^* \text{adjoint of $\sum$}$}%$}
\smallskip

The result (see Lions \cite{lionscollege}) is that,  as $L\to\infty$, in the sense of measures and for all $t >0$, 
\[\mathcal{L}(X^1_t, \dots, X^L_t) \to 
  \pi_t \in \mathcal{P}(\mathcal{P}(\R^d)),\]
where the  density $(m_t)_{t\geq 0}$ of the evolution in time of $(\pi_t)_{t\geq 0}$, which is  defined, for all  $U \in C(\mathcal{P}(\R^d))$, by 
\[\int U(m)d\pi_t(m) = {E}[U(m_t)],\]
solves the stochastic conservation law 
\[dm + \text{div}_x (\sigma^T(m,x) \circ dB)=0 \ \text{in}  \ Q_T,\]
which is a special case of \eqref{quasilinear}. Here $\sigma^T$ is the transpose of the matrix $\sigma$. 

%\end{document}
\section{The Main Difficulties and the Choice of Stochastic Calculus.}

\subsection*{Difficulties} Given that, in general and without rough signals,  \eqref{eq0.1} and \eqref{quasilinear} do not have global smooth solutions,  it is natural to expect that this is the case in the presence of  rough  time dependence.  
\smallskip
%\vskip.075in 

It is also not possible to use directly  the standard viscosity and entropy solutions  of  the ``deterministic'' theory, since they depend on inequalities satisfied either at some special points or after integration. Consider, for example, \eqref{quasilinear} with  $d=1$ and $A\equiv 0$. An entropy solution  must satisfy, in the sense of distributions,  the  weak entropy inequality 
$dS(u) + Q(u)_x\cdot dB\leq 0$   for all pairs $(S,Q)$ of convex entropy $S$ and entropy flux $Q$.  The inequality does not  make sense if $B$ is a rough path.  A similar difficulty arises when dealing with viscosity inequalities.
\smallskip

Moreover, the lack of regularity does not allow to express the solutions in any form involving time integration as is the case for sde, that is to say, for example, that $u$ solves $du=H(Du)\cdot dB \ \text{in} \ Q_T$ if,  for all $x\in \R^d$ and $s,t \in [0,T]$ with $s>0$, 
\[u(x,t)=u(x,s) +\int_s^tH(Du(x,\tau))\cdot dB(\tau).\]

\smallskip

Another possibility, at least when $m=1$, is to take advantage of the multiplicative noise to change time and obtain an equation without rough parts. For example,  formally, if $du+H(Du)\cdot dB=0$, the change of time   $u(x,t)=U(x,B(t))$  yields that  $U$ must be a global smooth solution to the forward-backward time homogeneous Hamilton-Jacobi equation $U_t + H(DU)=0$ in $\R^d\times (-\infty, \infty)$. It is, of course,  well known that such solutions do not exist in general.  Behind this difficulty is the basic fact that the nonlinear problems develop shocks which are not reversible, while  the changing sign of the rough signals, in some sense, forces the solutions to move forward and backward in time.  Note that  the time change works in intervals where $dB$ does not change sign. More details about this are given later in the notes.%the latter problem cannot have a smooth solution.
\smallskip

A natural question is whether it is possible to solve the equations  in law. 
Recall that solving the sde $dX = \sqrt2 \sigma (X_t) dB$
in law is equivalent to understanding,
for all smooth $\phi$ and $T >0$, the solutions $u$ 
of the initial  value problem 
$$
u_t = \ds  \sigma \sigma^T:D^2u		
\  \text{in} \  Q_T \ \quad   
u(\cdot, 0)=\phi. $$

For the equations here the state variable must belong to  a suitable function space and the corresponding 
spde is set in infinite dimensions.  For example, the infinite dimensional pde describing the law of 
%Indeed let $U\in C(\overline Q_T$ and write $U(f,t)$. 
%infinite dimensional equation for the law would bcorresponding to 
$ du = \sqrt{2}\, H(Du)\circ dB$
is, formally,  
$$U_t  = D^2 U (H(Df),H(Df)).$$
The problem is that the Hessian $D^2U$ is an 
unbounded operator independently of the choice of the base space. 
Such pdes are far away from the theory of viscosity 
solutions in infinite dimensions  developed by Crandall and Lions \cite{crandalllionsinfinite1, crandalllionsinfinite2}.
%\cite{CL2,CL3}, etc..
\smallskip

Solving linear stochastic pde in law is related to the martingale approach which has been used 
successfully in linear and some quasilinear settings. A partial list of references is Chueshov and Vuillermot \cite{chueshovvuillermot1, chueshovvuillermot2}, Da Prato, Ianelli and Tubaro \cite{dapratotubaro}, Gerencs\'{e}r, Gy\"{o}ngy and Krylov \cite{GGK},  Huang and Kushner \cite{huangkushner}, Krylov \cite{krylov, krylov2}, Krylov and R{\"o}ckner \cite{krylovrockner},  Rozovski{\u\i} \cite{rozovski1, rozovski2}, Pardoux \cite{pardoux1, pardoux2, pardoux3}, Watanabe \cite{watanabe}. 
The methodology requires some tightness (compactness) which typically follows from estimates on the derivatives 
of the solutions.  In general, the latter are not available for  nonlinear problems. 
%The latter is a problem in infinite dimensions, since it requires estimates 
%on  derivatives of the solutions, which are not, in general, available.  
%This difficulty is related to the nonlinear character of \eqref{eq0.1}.
%This program has been carried out successfully for linear and quasilinear 
%uniformly elliptic pde. 
%When dealing with degenerate fully nonlinear spde which do not have, 
%in general, smooth solutions, there are various serious difficulties. 
\smallskip

\subsection*{The choice of stochastic calculus; Stratonovich vs It\^o} When studying sdes,  it is important to decide if they are written in Stratonovich or It\^o form, each of which having advantages and disadvantages; for example, more regularity and chain rule for the former and less regularity but no chain rule for the latter.
\smallskip

At first glance, the choice of calculus does not seem to be relevant for the nonlinear problems discussed here due to the lack of regularity.   This is, however, not the case. The actual formulation plays an important role in the interpretation, well-posedness, stability  and construction of the solutions, which, typically,  are obtained as limits of solutions with regular time dependence. The discussion below touches upon some of these issues.

%are almost never written explicitly. 
%the choice of calculus   This is closely 
%related to the way solutions are constructed,  the nature of the nonlinearity, and well-posedness of the equations.
\smallskip

The advantage of the Stratonovich formulation can be seen in the following rather simple example. Consider, for $\lambda\geqq 0$, the  It\^o-form spde 
\begin{equation*}%\label{takis1}
du = 
\lambda u_{xx} dt + \sqrt{2}u_x dB
\ \text{in} \  Q_T.
\end{equation*}
The change of variables 
%\begin{equation*}
$u(x,t) = v (x+\sqrt 2 B(t),t)$
%\end{equation*}
yields  that $v$ satisfies the (deterministic) pde
\begin{equation*}
v_t = (\lambda -1) v_{xx} \  \text{in} \  Q_T,
\end{equation*}
which is well-posed if and only if $\lambda \geq 1$.
\smallskip

 Of course this is not an issue if the spde was in Stratonovich form to begin with. In that case the change of variables yields the equation 
 $$v_t = \lambda v_{xx} \  \text{in} \  Q_T,$$
which is well posed if and only if $\lambda\geq 0$, as is this case when  $B$ is a smooth path.
\smallskip

%Staying with the same example, it is worth remarking that $u$ is  given in terms of $v$ if $B$ is smooth. 
Consider, for example, a family $(B^\ep)_{\ep >0}$ of smooth approximations of  the Brownian motion $B$ and the  solution $u^\ep$ of the equation
$$u^\ep_t = 
 u^\ep_{xx} + u^\ep_x \dot B^\ep.$$ 
It is immediate  that $u^\ep(x,t) = v (x+B^\ep(t),t)$ with $v$ solving $v_t =  v_{xx}$. Letting $\ep \to 0$ then  yields that $u^\ep \to u$, which solves 
% corresponding solutions, it is immediate from the transformed equation that the solutions will converge, for $\lambda> 0$, to a solution of the stochastic pde 
\begin{equation*} 
du = 
 u_{xx} dt + u_x \circ dB.
\end{equation*}
Another example, where the use of Stratonovich appears to be necessary, is the application  to front propagation via the level set pde. One of the important 
elements of the theory is that the moving interfaces depend only on the initial one and not the particular choice of the initial datum of the pde. This is equivalent to the requirement that   the equations are invariant under increasing changes of the unknown.
\smallskip

%difficulty is encountered, concerns the level
%set pde in front propagation. 
%One of the very essential elements of the level set theory is that the front   
%depends only on the (zero) level sets and not any other properties 
%of the function used to define it initially.
%Analytically this  follows from the fact that the level pde is geometric, 
%i.e., it satisfies a structural
%property implying that, if $u$ is a solution and $\phi :\R\to\R$ 
%is nondecreasing, then $\phi (u)$ is also a solution. 
Consider, for example, the pde 
\begin{equation*}
u_t + 
%\beta |Du| \text{ div}(\frac{Du}{|Du|}) +
 |Du|=0  \ .
\end{equation*}
Arguing as if the solution $u$ were smooth (the argument can be made 
rigorous using viscosity solutions), it is straightforward to check that,  
for nondecreasing $\phi$, $\phi (u)$ is also a solution; 
note that the monotonicity of $\phi$ is important when 
dealing  with viscosity solutions.
\smallskip

The next example shows that the It\^o formulation is the wrong one.  Assume that  level set pde of the interfacial motion  $V=dB$ with $B$ a Brownian motion is  %Tcorresponding level set stochastic pde in It\^o's form is 
\begin{equation*}
du =|Du| dB.
\end{equation*}
If $u$ is a smooth solution and $\phi :\R\to \R$ is smooth and nondecreasing,   
It\^o's formula yields that 
\begin{equation*}
d\phi (u) = |D\phi (u)| dB + \frac12 \phi'' (u) |Du|^2,
\end{equation*}
which  is not the same equation as the one satisfied by $u$. This is of course not the case if the level set pde was written in the Stratonovich form, which, however, requires a priori additional regularity 
which is not available here. Indeed, if $du=H(Du)\circ dB$, then, in It\^o's form
$$du= H(Du)dB+ \frac 12 \left \langle D^2u DH(Du), DH(Du) \right \rangle dt.  $$
where, for $x,y \in \R^d$, $\langle x,y \rangle$ is the usual inner product.  To make, however,  sense of this last equation, it is necessary to have information about $D^2u$ which, in general, is not available.
 
\smallskip

In the context of second- and first-order (deterministic) pde the 
difficulties due to the lack of regularity are  overcome using viscosity solutions. Their definition is based on inequalities which, as mentioned earlier, cannot be expected to make sense in the presence of rough signals. 
\smallskip

There is, however, a reformulation of the definition for viscosity solutions, which, at first glance, appears to be more conducive to stochastic calculus. %However,  its Statonovich formulation fails  since it requires regularity which is not available, while, given in It\^o's form, it  leads to the wrong conclusions. 

% which are formulated below in a  way that is very convenient to explain the issue of the need of more regularity that arises  when using  Stratonovich calculus.
%The usual definition of the latter,
%however, cannot be used here, since, with It\^o's formulation,  
%the inequalities go the wrong way.
\smallskip 

Indeed, for $B$ smooth,  consider again the equation
$u_t = H(Du,x)\dot B.$
The definition of viscosity subsolutions 
is equivalent to the requirement that,
for any smooth $\phi : \R^d \to \R$, the map $t \to 
\max (u-\phi)$ satisfies, in the viscosity sense,  the differential inequality
$$\frac{d}{dt} \sup (u(\cdot,t)-\phi) 
\leqq \sup_{\bar x(t) \in \text{argmax} (u(\cdot,t)-\phi)} (H(D\phi (\bar x(t)), \bar x(t))\dot B).$$
%where $\bar x(t)$ denotes a point where $\max (u(\cdot,t)-\phi)$ 
%is achieved --- there may, of course, exist several such points, etc..
\smallskip
%
%For general $x$-dependent Hamiltonian $H$ the requirement is that 
%$$\frac{d}{dt} \max (u (\cdot,t)-\phi) \leqq \sup_{\bar x} \left(H(D\phi (\bar x(t)),\bar x(t)) \dot B\right)\ .$$
%\smallskip
%

If $B$ is a Brownian motion, then,  
assuming that there exists a unique maximum 
point $\bar x(t)$ of $u(\cdot,t)-\phi$,
the Stratonovich formulation should be 
$$\frac{d}{dt}\max (u(\cdot,t)-\phi)
\leqq H(D\phi(\bar x(t)),\bar x(t) ) \circ dB,$$
a fact which, however, breaks down due to the lack of regularity 
in $t$ of the map $t\mapsto \bar x(t)$. 
\smallskip

If $\dot B\in L^1((0,T))$, then the above inequality is meaningful and   
%makes sense, and % it has been used 
has been used by Lions and Perthame \cite{lionsperthame} and Ishii \cite{ishiitime} 
to study viscosity solutions of Hamilton-Jacobi equations 
with $L^1$-time dependence. % and by \cite{N} for second-order equations.
%\smallskip
\vskip.05in

The regularity concerns can, of course, be relaxed, if the inequality above is required to hold in It\^o's sense.  This, however,  leads to a contradiction to  the classical fact that 
 the maximum of two subsolutions is a subsolution. 
%Another way to see, again formally, that It\^o's formulation is not
%appropriate is to compute the It\^o's  differential of the max of two 
%functions and to compare it to the individual differentials.
\smallskip

Recall that, if $u$ and $v$ are actually differentiable with 
respect to $t$, then
$$\frac{d}{d t} (\max (u,v)) 
= \mathds{1}_{\{u(\cdot,t) > v(\cdot,t)\}} 
u_t + \mathds{1}_{\{u(\cdot,t) \leqq v(\cdot,t)\}}  v_t,$$
where $\mathds{1}_A$ denotes the characteristic function of the set $A$.
\smallskip

If 
$$u_t = H(Du)\ ,\quad v_t = H(Dv)\quad\text{and}\quad  H(0)=0\ ,$$
it follows that 
$$ \frac{d}{d t} \max (u,v) \leqq H(D(\max (u,v))),$$
and, hence, $\max (u,v)$ is a subsolution. 
\smallskip

Checking the same claim in the It\^o's formulation yields %fact for the spde, we find
$$d \max (u,v) \geqq \mathds{1}_{\{u(\cdot,t)> v(\cdot,t) \}} d u + 
\mathds{1}_{\{u(\cdot,t) \leqq v(\cdot,t) \}} d v,$$
%where $d_I$ denotes the It\^o differential.
which suggests  that $\max (u,v)$ 
is not necessarily a subsolution. % which should be the case for viscosity solutions. 
\smallskip

The final justification for considering the  Stratonovich vs It\^o's formulation when  studying, for example, the  equation 
\begin{equation*}%\label{final-just-1}
du = H (Du) \cdot dB
\end{equation*}
%instead of 
%\begin{equation}\label{final-just-2}
%du = H(Du,x)dB
%\end{equation}
comes from considering the family of problems 
\begin{equation*}
u_t^\ep = H(Du^\ep)\dot B^\ep,
\end{equation*}
where $B^\ep$ are smooth approximations of the Brownian motion   
$B$.  
If $u^\ep$ and $u$ are smooth and, as $\ep\to0$, 
$u^\ep \to u$ in $C^2(\R^d\times (0,\infty))$, 
it is not difficult to see that $u$ must solve the equation in the %\eqref{final-just-1} in the 
Stratonovich sense. 
\smallskip

Note that, under suitable assumptions on the initial datum of the regularized equation and the Hamiltonian, it is possible to show, using arguments from  the theory of viscosity solutions,   that the solutions $u^\ep$ are, uniformly in $\ep$, bounded and Lipschitz continuous in $x$, and, hence, converge uniformly along subsequences  for each $t$.  This observation is the starting point of the theory, since it provides a candidate for a possible solution of \eqref{eq0.1}.
\smallskip

%%%%%%%%%%%%%%%%%%%%%%%%%%%%%%
\section{Single versus multiple signals, the method of characteristics and nonlinear pde with linear rough dependence on time}% and  the method of characteristics.}
\subsection*{Single versus multiple signals} The next  example illustrates that there is a difference between one single  and many signals and indicates the role that rough paths may play in the theory.
\smallskip

Consider   two smooth paths $B_1$ and $B_2$ and the linear pde 
\begin{equation}\label{eq:simple-linear}
u_t = u_x \dot B_1 + f(x) \dot B_2 \  \text{in} \  Q_T  \quad  u(\cdot,0)=u_0. % \ \text{on} \ \R.
\end{equation}
% two smooth paths $B_1$ and $B_2$.
%--the argument  easily 
%extends  to the rough paths.

It is immediate that
$v(x,t) = u(x-B_1(t),t)$
solves  
$$ v_t = f(x-B_1(t)) \dot B_2 \ \text{in} \ 
Q_T \quad 
v(\cdot,0)= u_0, %\quad\text{on}\quad \R\ ,
$$
and, hence, 
$$u(x,t) = v(x+B_1(t),t) 
= u_0 (x+B_1(t)) + \int_0^t f(x+B_1(t) - B_1(s))\dot B_2(s) ds\ .$$
To extend this expression to non smooth paths, it is necessary to deal with integrals of the form 
%study for Brownian motions $W^1,W^2$ integrals of the form 
$$\int_a^b g(B_1(s))\, dB_2(s),$$
which is one of the key ingredients  of Lyons's theory of rough paths; see, for example, Qian and Lyons \cite{lyonsqian},  Lyons \cite{lyons1, lyons2}, Lejay and Lyons \cite{lejaylyons}, etc.. 
\subsection*{Nonlinear pde with linear rough dependence on time}
The calculation  above suggests, however, a possible way  to study general 
linear/nonlinear  equations with linear rough  dependence, that is, equations  of the form 
\begin{equation}\label{eq:general}
du = F (D^2 u,Du,x)dt + \langle a(x) ,Du\rangle \cdot dB_1 + c(x) u\cdot dB_2 \  \text{in} \  Q_T \quad u(\cdot,0)=u_0.% \ \text{on} \ \R^d.
\end{equation}
%which is discussed next. 
%\smallskip

Consider the system %of differential The characteristics of \eqref{eq:H-linear} are 
\begin{equation}\label{eq:H-linear-B}
\begin{cases}
dX = - a(X) \cdot dB_1\quad X(0)=x,\\
\noalign{\vskip6pt}
dP =  \langle Da(X), P\rangle \cdot dB_1 + \langle Dc(X), P\rangle  U \cdot dB_2 \quad P(0)=p,\\
\noalign{\vskip6pt}
dU = c(X) U\cdot dB_2 \quad U(0)=u,
%\noalign{\vskip6pt} 
\end{cases}\end{equation}
which, in view of the theory of rough paths, has a solution for any initial datum $(x,p,u)$. Of course, $a$ and $c$ must satisfy appropriate conditions. This, however, is not important for the ongoing discussion. 
% it is left up to the reader the the exact assumptions on the coefficients for this to be true. % are left as an exercise to the reader. 
\smallskip

It is immediate that, with initial condition  $X(0)=x, P(0)=Du_0(x), U(0)=u_0(x)$,  \eqref{eq:H-linear-B} is  the system of characteristic equations of  the linear Hamilton-Jacobi equation 
\[du=\langle a(x) ,Du \rangle\cdot dB_1 + c(x) u\cdot dB_2 \quad u(\cdot,0)=u_0.\]
%with initial condition  $X(0)=x, P(0)=Du_0(x), U(0)=u_0(x)$. 
%Note that, due to linearity of the problem, it is immediate that the map $x\to X(x,t)$ is invertible for all $t$.
%\smallskip

The next step is to make the ansatz that the solution $u$ of \eqref{eq:general}  has the form 
\begin{equation}\label{eq:unknown}
u(x,t) = v(X^{-1} (x,t),t),
\end{equation}
and to find the equation satisfied by $v$.  Note that, due to the linearity, it is immediate that the map $x\to X(x,t)$ is invertible for all $t$.

\smallskip

Substituting in \eqref{eq:general}, arguing formally (the calculation can be made rigorous using viscosity solutions when $B_1$ and $B_2$  are smooth), and rewriting \eqref{eq:unknown} as 
%the notation 
%
% The following calculation appears formal but can be easily justified using viscosity solutions arguments as long as one can make sense of the solution of the system and the resulting integration.  To simplify the notation,  it is mor e convenient  to rewrite \eqref{eq:unknown} as 
$$u(\cdot,t) = S (t) v(\cdot,t),$$
where, for any $v_0$, $S(t)v_0$ is the solution of the linear Hamilton-Jacobi equation with initial datum $v_0$, yields %and to proceed to identify the equation for $v$.
%\smallskip
%As in the previous section we denote by $(X_d,P_d,U_d)$ the characteristics
%for $\circW \equiv 1$. 
%The special form of \eqref{eq:H-linear-B} allows to ``invert'' the $X_d$ 
%characteristic globally in time.
%Indeed
%$$X_d^{-1} (x,t) = X_d (x,-t)\quad\text{and}\quad 
%X^{-1} (x,t) = X_d (x,-W(t))\ .$$
%
%Assume finally  
%that $u$ is regular enough for the following calculation to hold.
%If $u$ is only continuous, all the computations can be justified in the 
%viscosity sense. 
%
%We make the change of unknown
%\begin{equation}\label{eq:unknown}
%u(x,t) = v(X^{-1} (x,t),t)
%\end{equation}
%and find the equation satisfied by $v$. 
%This amounts to inverting the characteristics corresponding to the
%Hamilton-Jacobi part of the equation.
%At the abstract level what we do can be summarized loosely as follows:
%
%If $S_H(t)u_0$ denotes the solution of 
%$$U_t= H(DU,U,x)\circW \quad\text{in}\quad \R^d\times (0,\infty)
%\quad\text{with}\quad U= u_0\quad\text{on}\quad \R^d\times \{0\}\ , $$
%then we seek $u:\R^d \times [0,\infty)\to\R$ such that the 
%solution of \eqref{eq0.1} can be written as 
%$$u(\cdot,t) = S_H (t) v(\cdot,t)\ .$$
%As we see below, this can be done easily in the linear setting.
%For the nonlinear problem we follow the same strategy but we need to 
%address the serious difficulty that the 
%characteristics are not, in general,  invertible. 
%It follows that 
\begin{equation*}
\begin{split}
du = & d(S(t) v(\cdot,t))=dS(t)v(\cdot,t) + S(t)dv(\cdot,t)\\
& =   \langle a(x), DS (t) v(\cdot,t) \rangle \cdot dB_1 + c(x) S(t)v(\cdot,t)\cdot dB_2 
+ S(t)( v_t(\cdot,t)) \\
\noalign{\vskip6pt}
& =  \langle a(x), DS (t) v(\cdot,t) \rangle \cdot dB_1 + c(x) S(t) v(\cdot,t),x)\cdot dB_2 \\
\noalign{\vskip7pt}
&\qquad 
+ F(D^2 S (t) v(\cdot,t) , DS (t) v(\cdot,t), S(t) v(\cdot,t),x)dt,
\end{split}
\end{equation*}
and, hence, 
\begin{equation*}
S(t) dv (\cdot,t) = F(D^2 S (t) v(\cdot,t), DS (t) v(\cdot,t),
S(t) v(\cdot,t),x) dt,
\end{equation*}
and 
\begin{equation*}%\label{nonlinear-prob}
dv  = S^{-1} (t) F(D^2 S(t) v(\cdot,t) ,DS (t) v(\cdot,t),
S(t) v(\cdot,t),x) dt.
\end{equation*}
Since the last equation does not contain any singular time dependence, it is convenient to replace $dv$ by  $v_t$ and to rewrite the last equation as 
\begin{equation}\label{nonlinear-prob}
v_t  = S^{-1} (t) F(D^2 S(t) v(\cdot,t) ,DS (t) v(\cdot,t), S (t) v(\cdot,t),x).
\end{equation}
This last expression  appears to be more complicated than \eqref{eq:general}, 
but  this is only due to the notation.
\vskip.05in

The point is that \eqref{nonlinear-prob} actually 
is  simpler since the transformation eliminates the troublesome term 
$$ \langle a(x) ,Du  \rangle \cdot dB_1 + c(x) u\cdot dB_2.$$
The new equation is of the form
\[v_t=\widetilde F(D^2v, Dv,v,x,t)  \   \text{in} \ Q_T \quad  v(\cdot,0)=u_0,\] % \ \text{on} \ \R^d \]
and can be studied using the viscosity theory as long as $\tilde F$ satisfies the appropriate conditions for well-posedness.
\smallskip

The discussion above  gives an alternative way to find pathwise solutions to all the equations studied using the martingale method as well as scalar quasilinear equations of divergence form,  always with linear rough time dependence. As a matter of fact, a closer look at the existing theories for linear spde yields that the approach described above allows for the treatment of larger class of equations.  

\subsection*{Stochastic characteristics} The analysis in the previous subsection suggests  that  to handle equations with nonlinear rough dependence, it may be 
useful to look, at least when the Hamiltonians are smooth,  at the associated system of characteristics.
When the time signals are smooth this is a classical system of $2d+1$ ode. In the particular case that the rough dependence is Brownian, the stochastic characteristics were used in the work of Kunita \cite{kunita} on stochastic flows. In what follows, statements are made   without any assumptions and the details are left to the reader.
% It is left to reader For now we leave the   details to left to the  interested leave the details to the reader for now.
\smallskip

The characteristics of the Hamilton-Jacobi equation 
\begin{equation}\label{paris1}
du=\sum_{i=1}^m H^i(Du, u, x,t) \cdot dB_i   \  \text{in} \   Q_T \quad  u(\cdot,0)=u_0,% \ \text{on} \ \R^d, 
\end{equation}
are the solutions to the following system of differential equations:
%The characteristics associated with \eqref{eq:simple}    
%is the following system of ordinary differential equations (for short ode)
\begin{equation}\label{eq:odes}
\left\{
\begin{array}{l}
\ds dX = - \sum_{i=1}^m D_p H^i(P, U, X,t) \cdot dB_i \quad X(x,0) = x, \\
\noalign{\vskip6pt}
\ds dP =  \sum_{i=1}^m \left(D_x H^i(P, U, X,t) + D_uH^i (P, U, X,t) P\right) \cdot dB_i \quad P(x, 0) = Du_0(x),\\
\noalign{\vskip6pt}
\ds dU = \sum_{i=1}^m \left(H^i(P, U, x,t) - \langle D_p H^i(P, U, x,t), P\rangle \right) \cdot dB_i \quad U(x, 0) = u_0(x).\\[2mm]
%\ds X(x,0) = x\ ,\qquad P(x, 0) = Du_0(x)\ ,\qquad U(x, 0) = u_0(x).
\end{array}\right. \end{equation}
%with initial conditions
%\begin{equation}\label{eq:initial}
%X(x,0) = x\ ,\qquad P(x, 0) = Du_0(x)\ ,\qquad U(x, 0) = u_0(x)\ .
%\end{equation}
The connection between \eqref{paris1}    
and \eqref{eq:odes} 
%and \eqref{eq:initial} 
is made through the relationship 
$$U(x,t) = u(X(x,t),t)\quad\text{and}\quad P(x,t) = Du(X(x,t),t).$$
The method of characteristics works  
as long as  it is possible to invert the map    %\eqref{eq:odes},\eqref{eq:initial} and
$t \to X(x,t)$.  This  can always  be done  in some interval $(-T^*,T^*)$ for small  $T^*>0$, which depends on bounds on $H, u_0$, their derivatives and the signal, and, in general, is difficult to estimate in a sharp way.
\smallskip

It then follows that %satisfies appropriate conditions,  although
%in general it is very difficult to estimate this time, that is, there is some $T^*> 0$ such that 
\[u(x,t)= U(X^{-1}(x,t),t) \ \]
is a smooth solution to \eqref{paris1} in $ \R^d\times (-T^*,T^*).$  The latter means, for all $s,t \in  (-T^*,T^*)$ with $s < t$ and $x\in \R^d$,
\[u(x,t) =u(x,s) + \displaystyle\int_s^t \sum_{i=1}^m H^i(Du(x,r), u(x,r), x,r) \cdot dB_i(r).\]
If $m=1$, it is possible to express the solutions of \eqref{eq:odes}  using in the characteristics of the ``non rough'' equation
\[u_t=H(Du, u, x,t)  \ \  \text{in} \ \ Q_T \quad  u(\cdot,0)=u_0.\]% \ \text{on} \ \R^d.\]
Indeed if  $(X_d,P_d,U_d)$ is the solution of 
\begin{equation}\label{eq:dodes}
\left\{
\begin{array}{l}
\ds \dot X_d = - D_p H (P_d, U_d, X_d, t) \quad X_d(x,0) = x, \\
\noalign{\vskip6pt}
\ds \dot P_d =  D_x H^i(P_d, U_d, X_d, t) + D_uH^i (P_d, U_d, X_d, t) P_d  \quad P_d(x, 0) = Du_0(x),\\
\noalign{\vskip6pt}
\ds \dot U_d = H^i(P_d, U_d, X_d, t) - \langle D_p H(P_d, U_d, X_d, t), P_d\rangle \quad U_d(x, 0) = u_0(x),\\
%\noalign{\vskip6pt}
%X_d(x,0) = x\ ,\qquad P_d(x, 0) = Du_0(x), \qquad U_d(x, 0) = u_0(x), \ 
\end{array}\right. 
\end{equation}
then 
$$X(x,t)=X_d(x,B(t)), \quad P(x,t)=P_d(x,B(t)), \quad  \text{and}  \quad U(x,t)=U_d(x,B(t)),$$
and the inversion is possible as long as $|B(t)| < T^*_d$, the maximal time for which $X_d$ is invertible. 

\smallskip

This simple expression for the solution of \eqref{eq:odes} is not valid for
$m\geqq 2$ unless the Hamiltonian $H$ satisfies the involution relationship
$$\{H^i,H^j\}: =D_x H^i D_p H^j - D_x H^j D_p H^i =0\quad\text{for all}\quad i,j=1,\ldots,m.$$
%where $\{H^i,H_j\}$ denotes the usual Poisson bracket, that is 
%$\{H^i,H_j\} = 
%D_x H^i D_p H_j - D_x H_j D_p H^i\ .$
The latter yields that the solutions of the system of the characteristics commute, 
that is 
$$X(x,t) = X_d^1 (\cdot,B_1(t)) \sbullet X_d^2 (\cdot,B_2(t)) 
\sbullet \cdots\sbullet 
X_d^M (\cdot, B_m(t))(x),$$
where, for $i=1,\ldots,m$, 
$(X_d^i,P_d^i,U_d^i)$ is  the solution of \eqref{eq:odes} with 
$H\equiv H^i$ and $B_i(t) = 1$ and $\sbullet$ stands for the 
composition of maps.
\smallskip

For example, if, for all $i=1,\ldots,m$, the $H^i$'s are independent of $x,u$ and $t$, then  the involution relationship  is satisfied,  and \eqref{eq:odes} reduces to % has the form  
%\begin{equation*}
%\begin{cases}
%\ds \dot X = - \frac{\partial H^i}{\partial p} (P) \circW_1\\
%\noalign{\vskip6pt}
%\dot P = 0\\
%\noalign{\vskip6pt}
%\ds \dot U = - (P \frac{\partial H^i(P)}{\partial p} - H^i)(P) \circW_1
%\end{cases}
%\end{equation*}
$$dX = - \sum_{i=1}^m DH^i (P) \cdot dB_i, \qquad  dP=0\ ,\qquad 
dU = \sum_{i=1}^m [ H^i(P) - \langle D_p H^i(P),P\rangle  ] \cdot B_i.$$
and the $X$-characteristic is given by 
$$X(x,t) = x- \sum_{i=1}^m D_xH^i 
%\frac{\partial H^i}{\partial p} 
(Du_0(x)) B_i(t).$$
Finally, either for  $m=1$ or for space homogeneous Hamiltonians when $m\geqq 2$,  it is possible to 
find $X, P$ and  $U$ for any continuous $B$. Otherwise it is necessary to appeal to the rough path theory.

\section{Fully nonlinear equations with semilinear stochastic dependence}
%%%%%%%%%%%%%%%%%%%%%%%%%%%%%%%%%%

I describe next the work of Lions and Souganidis \cite{lionssouganidis4} about fully nonlinear equations with semilinear stochastic dependence. 
\smallskip

Consider the 
initial value problem
\begin{equation}\label{eq:initialvalue}
%\begin{cases}
\ds du = F(D^2 u,Du,u)dt + \sum_{i=1}^m H^i (u)\cdot  dB_i \ \text{in}
\  Q_T \quad  u(\cdot, 0) = u_0,\\
%\noalign{\vskip6pt}
%u(\cdot, 0) = u_0\quad\text{on}\quad \R^d,
%\end{cases}
\end{equation}
with $u_0 \in BUC (\R^d)$,  $B = (B_1,\ldots,B_m)$  is a $C^\alpha$ geometric  rough path with $\alpha\in (1/3,1/2)$, for example Brownian motion with Stratonovich, 
$F\in C({\mathcal S}^d\times \R^d) $ 
degenerate elliptic, that is, for all $(p,u) \in \R^{d+1}$ and $X,Y \in {\mathcal S}^d$, 
%the space of symmetric $N\times N$ matrices,
\begin{equation}\label{degenerateelliptic}
 \text{if} \  X\leq Y, \ \text{then} \ F(X,p,u) \leq F(Y,p,u), 
\end{equation}
and 
\begin{equation}\label{eq4.20biss} 
H= (H^1,\ldots,H^m) \in (C^{5} (\R))^m . 
\end{equation}
When $m=1$ and $B$ is continuous path, then \eqref{eq4.20biss} can be replaced by  
\begin{equation}\label{eq4.20bis} 
H \in C^{3,1} (\R). 
\end{equation}

Although the results presented here also apply to the more general equations like 
\begin{equation}\label{eq:old-0.4}		%% \leqno(0.4)
du=F(D^2u, Du,u,x,t )dt +    \sum^m_{i=1}
H^i(u,x,t)\cdot dB_i \ \text{in} \  Q_T,
\end{equation} 
%for $M=1$ (we do not know how to treat $M\geqq 2$) 
%in which case it will be necessary to assume that 	%% $M=1$ and 
%$B$ is a rough path. 
for simplicity  I concentrate on \eqref{eq:initialvalue} and assume that $m=1$. % the exact assumptions will become clear below.
%%%%%%%%%%%%%%%%%%%%%%%
%{\bf 1. Definition and main results.}
%%%%%%%%%%%%%%%%%%%%%%%
%\section{Definition and main results}
\smallskip

For $v\in\R$, consider the differential equation
%\begin{equation}\label{takis3} 
%d\Phi =\sum^m_{i=1} H^i(\Phi )\cdot dB_i  \ \text{in} \ (0,\infty) \quad 
%\Phi (v,0)=v.
%\end{equation} 	
\begin{equation}\label{takis3.1} 
d\Phi =H(\Phi )\cdot dB  \ \  \text{in} \ \ (0,\infty) \quad 
\Phi (v,0)=v.
\end{equation} 	
It is assumed that 
\begin{equation}\label{semilinearas}
\begin{cases}
%\begin{split}
 \text{there exists a unique solution  $\Phi\in C([0,T]; C^3(\R))$ % :\R \times (0,\infty)\\to\R$ 
of  \eqref{takis3.1} such that, 
for all $T >0$,}\\[1mm]
%\noalign{\vskip4pt}
%\hskip-1truein 
%\hskip-1in \Phi (\cdot,\cdot,\omega) \in C([0,T]; C^3(\R)) \ \ \ \text{and}\\[1mm]
%\text{and}\\[1mm]
%\hskip.2truein 
\hskip1in
M(T)=  \sup_{0\le t\le T} \Big[|\Phi (0,t)|+\sum^3_{i=1}\|
D^i_v\Phi (\cdot , t)\|_\infty\Big] < \infty.
%\end{array}\right.
%\end{split}
\end{cases}
\end{equation}
%When $m\rangle1$,  the $B_i$'s are taken to be geometric rough paths and \eqref{takis3} is interpreted in the rough path sense. % and the properties of their theory. 
Since $m=1$, it follows that,  %then \eqref{takis3} yields, 
for all $t >0$,
\begin{equation}\label{eq1.4} 
\Phi (v,t)=\widehat\Phi (v,B(t))\ ,%\qquad (t> 0),		
\end{equation}
where $\widehat \Phi$ solves the ode
\begin{equation}\label{eq1.5}
\dot{\widehat\Phi} =H(\widehat \Phi ) \  \text{in}  \  \R \quad 
\widehat\Phi (v,0)=v.		%% \leqno(1.5)
\end{equation}
It is then straightforward to obtain  
\eqref{semilinearas} from 
%\eqref{eq1.2} and \eqref{eq1.3} from
the analogous properties of $\hat\Phi$.
%%Moreover, it is also clear that 
%%\eqref{eq1.2},  \eqref{eq1.3} 
%\eqref{semilinearas}
%and \eqref{eq1.4} also hold, 
%if $B$ is any continuous and not only Brownian path.
\smallskip

%For the rest of the section we assume that \eqref{semilinearas} is 
%satisfied and we define the function 
Define $\widetilde F:{\mathcal S}^d\times \R\times [0,\infty )\to\R$ by
\begin{equation}\label{eq1.6} 	%% \leqno(1.6)
\widetilde F(X,p,v,t)=\frac1{\Phi'(v,t)}
F(\Phi'(v,t)X+\Phi''(v,t)(p\otimes p),
\Phi'(v,t,)p, \Phi(v,t)),
\end{equation}
where, to simplify the presentation, ``$\prime$'' denotes the
partial derivatives of $\Phi$ with respect to $v$. 
%Note that, if $H$ is linear in $u$, the transformation does not create 
%any new dependence on $v$, and the problem can be treated using the classical viscosity solution arguments.  
%% and the dependence on 
%%$x$ is omitted.
%%\smallskip
%If $H$ is nonlinear,  then $\tilde F$ depends nontrivially on $v$ and new arguments and assumptions are needed. 
%
%if $H$ is nonlinear in which case $\Phi' (\cdot,t,)$ and $\Phi''(\cdot,t)$depend on $v$, while, if $H$ is linear,  
%$\tilde F$ is independent of $v$. 
%and, as a matter of fact, 
%$$\tilde F (X,p,v,t,\omega) = \exp (-\sum \lambdA^i W_t^i) X, 
%\exp \sum \lambdA^i W_t^i p.$$
\smallskip

The following definitions are motivated by the strategy described in 
Section~3 which amounts to  inverting the 
characteristics.  For \eqref{eq:initialvalue}, the latter 
are the solutions of \eqref{takis3.1}, which, in view of the semilinear form, % of the coefficient of the rough path,  
can be inverted globally.
% and In the case at hand, of course, this can be done globally.

%For \eqref{eq:initialvalue} the latter 
%are the solutions of \eqref{takis3}. 
%%Following  %% our previous work \cite{LS1,LS2} 
%the strategy described in the previous Section, 
\smallskip
%
%For each 
%$\phi\in C^2(Q_T)$, set  %define let the smooth in $x$  function
%$$
%\Psi (x,t):=\Phi (\phi (x,t), t),
%$$
%and note that $\Psi$ is smooth in $x$.  
%\smallskip
%
The definition of weak solution of \eqref{eq:initialvalue} is introduced  next.

\begin{defn}\label{defn1.1}
%A process $u:\R^d\times [0,T]\times \Omega\rightarrow \R$ 
Fix $T>0$.  
%and $u_0\in BUC(\R^d)$ and  assume that \eqref{degenerateelliptic}, \eqref{eq4.20bis} and \eqref{semilinearas} %hold.
Then $u \in BUC(\overline Q_T)$ is a pathwise subsolution (resp. supersolution) of  
%is a pathwise  viscosity sub-solution (resp.\ super-solution) of
\eqref{eq:initialvalue},    %% \eqref{eq:old-0.1} 
if,
%\noindent
%{\rm (i)} $u(\cdot, \cdot, \omega)\in BUC(\R^d\times [0,T])$ a.s.
%
%\noindent
%{\rm (ii)} $(t,\omega )\rightarrow u(\cdot , t, \omega)\in BUC(\R^d)$ is
%$\FF_t$-measurable, and
%
%\noindent
%{\rm (iii)} 
for all $\phi\in C^2(Q_T)$ and all local
maximum (resp. minimum) points $(x_0, t_0)\in Q_T$ of
$(x,t)\rightarrow u(x,t)-\Phi (\phi (x,t), t)$, 
%if $u(x_0, t_0, \omega )=\Phi (\phi (x_0, t_0),t_0, \omega )$, then
\begin{equation}\label{eq1.7} 		%% \leqno(1.7)
\phi_t (x_0, t_0)\le
\widetilde F(D^2\phi (x_0, t_0), D\phi (x_0, t_0),  u(x_0,t_0),t_0 ),
\end{equation}
$($resp.
\begin{equation}\label{eq1.8} 		%% \leqno(1.8)
\phi_t (x_0, t_0)\ge
\widetilde F(D^2\phi (x_0, t_0), D\phi (x_0, t_0),  u(x_0,t_0), t_0 )\big)\ .
\end{equation}
A function $u\in BUC(\overline Q_T)$ is a pathwise (viscosity) solution of 
\eqref{eq:initialvalue},    %% \eqref{eq:old-0.1} 
if it is both subsolution and supersolution of  \eqref{eq:initialvalue}.
\end{defn}

Since the characteristics are globally invertible, 
%a consequence of the semilinear ``stochastic'' form of the equation, 
it is possible to introduce 
a global change of the unknown without going through test functions. This leads to the next possible definition.

%It is straightforward to check that Definition~\ref{defn1.1} is equivalent to
%Definition~\ref{defn1.2}, which is stated below.  
%Since the latter reduces directly the problem to a ``deterministic'' one 
%with random coefficients, we choose to use it in what follows.
%\smallskip
%The second definition is:

\begin{defn}\label{defn1.2}
Fix $T >0$.  
%and $u_0\in BUC(\R^d)$ and  assume that \eqref{degenerateelliptic}, \eqref{eq4.20bis} and  \eqref{semilinearas} hold.
Then $u \in BUC(\overline Q_T)$ is a pathwise subsolution (resp. supersolution) of   \eqref{eq:initialvalue}, if 
the function $v:\R^d\times [0,T]\to\R$ defined by
\begin{equation}\label{eq1.9} 		%% \leqno(1.9)
u(x,t)=\Phi (v(x,t), t)
\end{equation}
is a viscosity subsolution (resp.\ supersolution) of
\begin{equation}\label{eq1.10} 		%% \leqno(1.10)
%\begin{cases}
v_t=\widetilde F(D^2v, Dv, v, t) \ \text{in \ $Q_T$}   \quad v(\cdot,0)=u_0. % \ \text{on} \R^d.
 % \ \text{on\quad $\R^d\times \{ 0\}$.} \cr
%v=u_0&\text{on\quad $\R^d\times \{ 0\}$.}\end{cases}
\end{equation}
%$u \in BUC(\overline Q_T)$
A function $u \in BUC(\overline Q_T)$  is a pathwise solution of 
\eqref{eq:initialvalue}		%% \eqref{eq:old-0.1} 
if it is both subsolution and supersolution.
% of 
%\eqref{eq:initialvalue}.	%% \eqref{eq:old-0.1}.
\end{defn}

The two definitions are equivalent, and, 
moreover, for smooth $B$'s, the solutions 
introduced in Definitions~\ref{defn1.1} and Definition \ref{defn1.2} 
coincide with the classical viscosity solution. 
%are the same as the viscosity solutions for smooth $W$. 
\smallskip

In view of the above, the well-posedness of
solutions to  \eqref{eq:initialvalue}   %%    \eqref{eq:old-0.1} 
reduces to the study of the analogous 
questions for \eqref{eq1.10}.  
%The  address these issues.  

\smallskip
%To state the results it is necessary to consider two separate cases, namely
%$M\rangle1$ and 
%$H=(H_1,\ldots ,H_M)$  linear and, hence, $\Phi (v,t)\equiv v$, and $M=1$ and 
%$H$  nonlinear, in which case $\Phi (v,\cdot )\not\equiv v$.  
%Finally, 
After the work described above was announced, Buckdahn and Ma \cite{buckdahnma2, buckdahnma3} used  the  map  \eqref{eq1.9}, which is known as the Doss-Sussman transformation, to study equations similar to \eqref{eq:initialvalue}. The work in \cite{buckdahnma2, buckdahnma3} covers   a more restrictive class of $F$'s and  well-posedness is proved under the assumption that the transformed initial value problem admits a comparison principle. In \cite{lionssouganidis4} there is no such assumption and the comparison is proved directly. 

\smallskip

If $H$ is linear in $u$, the problem is simpler and the details are left to the reader. %can be treated using the arguments of the previous section since $F$ is independent of $v$. 
\smallskip

For the the rest of the section,  $H$ is taken to be nonlinear, and, 
to simplify the presentation,
it is also assumed  that $F$ is independent of $u$. % and $\Phi$ dependents only  of $u$. 
%Finally, note that if $H$ is linear in $u$, the problem can be treated using the arguments of the previous section since $F$ is independent of $v$. Hence
%% and $x,t$.
% and we do not 
%exhibit the dependence of $\Phi$ and $\widetilde F$ on $\omega$.
\smallskip

To deal with  $\tilde F$, it is necessary to assume that 
%it is needed  to assume, that $F$ is degenerate, that is,
%\begin{equation}\label{deg.ell}
%F \ \text{ is nondecreasing with respect to the Hessian mattrix},
%\end{equation}
%and 
\begin{equation}\label{eq1.11} %% \leqno(1.11)
F\in C^{0,1}({\mathcal S}^d\times \R^d),
\end{equation}
and 
\begin{equation}\label{eq1.12}   %%\leqno(1.12)
\begin{cases} 
\text{there exists a constant \  $C >0$ \  such that, for almost every \ $(X,p)$,}\\[1mm] 
\qquad \text{either} \quad  D_X F(X,p): X + \langle D_pF(X,p),p\rangle  - F\leqq C\\[1mm] 
\qquad   \text{or}  \quad  D_X F(X,p): X + \langle D_pF(X,p),p\rangle  - F\geqq -C.
\end{cases}
\end{equation}
It is easy to see that any linear $F$ satisfies \eqref{eq1.12}.  Moreover, 
 \eqref{eq1.11} implies that $F$ can be written as the minmax of linear functions, that is,  
$$F(X,p)=\underset{\alpha \in A}\sup \underset{\beta \in B} \inf (a_{\alpha,\beta}:X + \langle b_{\alpha,\beta},p\rangle +h_{\alpha,\beta}),$$
for $A \subset {\mathcal S}^d$ and  $B\subset \R^d$ bounded  and  $a_{\alpha,\beta}\in {\mathcal S}^d$ and $b_{\alpha,\beta} \in  \R^d$ such that
$$\underset{\alpha \in A}\sup \underset{\beta \in B} \inf  [ \|a_{\alpha,\beta}\| + |b_{\alpha,\beta}|] <\infty.$$
Since $ D_X F(X,p): X + \langle D_pF(X,P),P\rangle  - F$ is formally the derivative, at $\lambda=1$, of the map
$\lambda  \to F(\lambda X,\lambda p)-\lambda F(X,P)$, it follows that \eqref{eq1.12} is related to, a uniform in $\alpha,\beta$,  one sided bound of ${\lambda}^{-1} h_{\alpha,\beta} - h_{\alpha,\beta}$ in a neighborhood of $\lambda=1.$ 
\smallskip

%
%\left\{\begin{matrix}
%   \text{there exists a constant \  $C> 0$ \  such that either }\hfill\cr
%   \noalign{\vskip6pt}
%   \text{\quad (i)~~\ \langle(D_XF(X,p), X\rangle + \langleD_pF(X,P),P\rangle - F\leqq C$ 
%   \qquad  for almost every \ $(X,p)$,}\hfill\cr
%   \text{or}\hfill\cr
%   \text{\quad (ii)~~$\langleD_X F(X,p), X\rangle + \langleD_pF(X,P),P\rangle - F\geqq -C$ 
%   \quad\  for almost every \ $(X,p)$.}\hfill\end{matrix}\right.
%\end{equation}

%The result is:
%
%\begin{thm}\label{thm1.1} 		%%  {\bf Theorem 1.1}:  
%Fix $T> 0$, let $u_0\in BUC(\R^d)$ and  assume \eqref{degenerateelliptic}, \eqref{eq4.20bis}, \eqref{semilinearas},
%\eqref{eq1.11} and  \eqref{eq1.12}. 
%%\eqref{eq4.20bis}, %\eqref{eq:old-0.4} 
%%\eqref{semilinearas}, 
%%and, if $H$ is nonlinear,  
%%\eqref{eq1.11} and \eqref{eq1.12}.
%Then there exists at most one  pathwise 
% solution of \eqref{eq:initialvalue}.    %%  \eqref{eq:old-0.1}.
%\end{thm}

I present next two explanations for the need for an assumption like \eqref{eq1.12}. % is discussed next. Two explanations are presented. 
The first 
is based on considerations from the method of characteristics. The second relies on viscosity solution arguments.
\smallskip

Consider the following  first-order versions of \eqref{eq:initialvalue} and 
\eqref{eq1.10}, namely 
\begin{equation}\label{eq:F}
du = F(Du)dt + H(u)\cdot dB,
\end{equation}
and 
\begin{equation}\label{eq:F1}
v_t = \tilde F (Dv, v, t),
\end{equation}
%and 
%\begin{equation}\label{eq:F-tilde}
%v_t = \tilde F (Dv, v, t,\omega) \ ,
%\end{equation}
with 
\begin{equation}\label{F-condition}
\tilde F (p,v,t ) = \frac1{\Phi'(v,t)} 
F(\Phi'(v,t)p),
\end{equation}
where $d\Phi=H(\Phi)\cdot dB$, and assume that $F$, $H$, $B$ and, hence, $\tilde F$ are smooth.
\smallskip

% and that  
%there exist smooth solutions of \eqref{eq:F}% and \eqref{eq:F-tilde}, 
%which, 
%for smooth data and  up to some $T^* > 0$ 
%are given, in view of the discussion in Section~3,
%by the method of characteristics.

The characteristics of the equations in \eqref{eq:F} and \eqref{eq:F1} are respectively %and \eqref{eq:F-tilde} the latter have the form 
\begin{equation}\label{eq4x}
\left\{ \begin{array}{l}
\dot X  =  - D F(P),\\
\noalign{\vskip6pt}
\dot P  =  H' (U) P\dot B,\\
\noalign{\vskip6pt}
\dot U  =  [F(P) - \langle DF(P),P\rangle] + H(U) \dot B,
\end{array} \right.
\end{equation}
and 
\begin{equation}\label{eq5x}
\left\{ \begin{array}{l}
\dot Y  =  - D_Q \tilde F(Q,V) = -D_P F(\Phi' (V),Q),\\[1.2mm]
\noalign{\vskip6pt}
\dot Q  =  \tilde F_V Q 
= \Phi'' (V) (\Phi'(V))^{-2} Q
[D_P F(\Phi' (V) Q), \Phi'(V)Q)) - F(\Phi'(V)Q)],\\[1.2mm]
%=  \Phi'' (\Phi')^2 P [(D_p F(\Phi' Q) ,\Phi' Q)-F]\ ,\\
\noalign{\vskip6pt}
\dot V  =  \tilde F - \langle D_Q\tilde F(Q,V),Q\rangle  
= (\Phi' (V))^{-1} [F(\Phi' (V)Q) - \langle D_PF(\Phi'(V)Q),\Phi'(V)Q\rangle ].
%(\Phi')^{-1} [F - (D_p F(\Phi' q), \Phi' q)]\ .
\end{array}\right.
\end{equation}
Of course,  \eqref{eq4x} and \eqref{eq5x} are 
equivalent after a change of variables. 
It is, however, clear that some additional hypotheses are needed 
in order for \eqref{eq4x}, and, hence, \eqref{eq5x} to have unique solutions. 
For example, without any additional assumptions, the right 
hand side  of the $P$-equation in \eqref{eq4x} may not be Lipschitz continuous in $U$.
On the other hand, the right hand side of the equations for $Q$ and $V$  
in \eqref{eq5x} contain the quantity $\langle D_p F,P\rangle-F$ appearing in 
\eqref{eq1.12} and an, at least  one-sided, Lipschitz condition is  necessary to 
yield existence and uniqueness. 
\smallskip
%\end{document}

The second explanation is based on the fact that the comparison principle for the pathwise viscosity solutions 
of \eqref{eq:initialvalue} will follow from the  comparison 
in $BUC(Q_T)$ of viscosity solutions of \eqref{eq1.10}. The 
latter does not follow directly from the existing
theory unless something more is assumed; see, for example,  Barles \cite{barlesbook} and Crandall, Ishii and Lions \cite{cil}. 
\smallskip

This ``additional'' assumption is that 
for each $R >0$, there exists $C_R >0$ such that, for all $X\in {\mathcal S}^d$, 
$p\in \R^d$, $v\in [-R,R]$ and  $t\in [0,T]$, 
\begin{equation}\label{eq1star}
\frac{\partial \tilde F}{\partial v} (X,p,v,t) \leqq C_R\ .
\end{equation}

A straightforward calculation, using \eqref{eq1.12}, yields that, 
for all $X,p,v$ and $t$, 
\begin{equation}\label{eq1.13}
\frac{\partial\tilde F}{\partial v} = \frac{\Phi''}{(\Phi')^2}
[D_XF : (\Phi' X + \Phi'' p\otimes p) + \langle D_pF ,\Phi' p\rangle -F] 
+ \Phi' (\frac{\Phi''}{\Phi'})' D_XF : p\otimes p\ ;
\end{equation}
note that to keep the formula simple, the explicit dependence of $F$ and its derivatives on $\Phi' X + \Phi'' p\otimes p$ is omitted.
\smallskip

It is immediate that   $\dfrac{\partial\tilde F}{\partial v}$ cannot  
satisfy \eqref{eq1star} without an extra assumption on $F$ and control on the size of $p$. If a bound on $p$ is not available, it is necessary to know that 
%unless and assuming something like \eqref{eq1.12} on 
%$D_XF \cdot X - (D_p F,p)- F$ and having control on the size of $p$ 
%unless, for the latter, it is  known that
 $\Phi' (\Phi''(\Phi')^{-1})' \geqq 0$. 
\smallskip

%It is thus possible to compare sub- and super-solutions provided 
%one of the two is Lipschitz continuous. 
%%The former is \eqref{eq1.12}, while for the.....
%It turns out, however, that it is possible to compare general bounded 
%uniformly continuous sub- and super-solutions.
%\smallskip

The last point that needs explanation is that  \eqref{eq1.13} is nonlocal, in the sense that it depends on $v$  through $\Phi$, while  \eqref{eq1.12}
is a local one, that is $\Phi$ plays no role whatsoever. This can be taken care of  in the proof by working in uniformly small time intervals, using the local time 
behavior of $\Phi$ and then iterating in time.
\smallskip

%It is thus possible to compare sub- and super-solutions provided 
%one of the two is Lipschitz continuous. 
%The former is \eqref{eq1.12}, while for the.....
%It turns out, however, that it is possible to compare general bounded 
%uniformly continuous sub- and super-solutions.
%\smallskip

%We have: 

%The uniqueness of the stochastic viscosity solutions will follow from the a.s. 
%uniqueness in $BUC(\R^d\times [0,T])$ of viscosity solutions of \eqref{eq1.10}.
%When $H$ is nonlinear,  the uniqueness of solutions of \eqref{eq1.10} 
%does not follow directly from the existing theory.  
%Indeed, a straightforward formal calculation using \eqref{eq1.12} yields
%that, for all $X,p$ and $v$, 
%\begin{equation}\label{eq1.13}		%% \leqno(1.13)
%\text{either }\; \frac{\partial}{\partial v}\widetilde F(X,p,v)\le
%C(1+|p|^2)\;\;\text{or}\;\; \frac{\partial}{\partial v} \widetilde
%F(X,p,v)\ge -C(1+|p|^2).
%\end{equation}
%
%These gradient dependent bounds do not allow for a direct application of the 
%existing comparison results (see \cite{B}, \cite{CIL}), unless one of the 
%sub- and super-solutions to be compared is Lipschitz continuous.  
%Since the comparison is of independent interest, we state it below 
%as a separate proposition.

The comparison result is stated next.

%\begin{thm}\label{thm1.1} 		%%  {\bf Theorem 1.1}:  
%Fix $T> 0$, let $u_0\in BUC(\R^d)$ and  assume \eqref{degenerateelliptic}, \eqref{eq4.20bis}, \eqref{semilinearas},
%\eqref{eq1.11} and  \eqref{eq1.12}. 
%\end{thm}

%\begin{prop}\label{prop1.1}	%% {\bf Proposition 1.1}: 
\begin{thm}\label{thm1.11}
Assume  \eqref{degenerateelliptic}, \eqref{eq4.20bis}, \eqref{semilinearas},
 \eqref{eq1.11} and \eqref{eq1.12}.
For each $T >0$ and any geometric rough path $B$ in $C^\alpha$with $\alpha\in(1/3, 1/2]$, there exists a 
constant $C=C(F, H, B, T) >0$ such that, if  $\overline v\in BUC(\overline Q_T)
$ and $\underline v\in BUC(\overline Q_T)$ are respectively 
a subsolution and a supersolution of \eqref{eq:initialvalue}, then, for all $t\in [0,T]$, 
\begin{equation*}
\sup_{x\in \R^d} (\overline v(x,t)-\underline v(\cdot,t))_+\le C \sup_{x\in\R^d}
(\overline v(\cdot,0)-\underline v(\cdot,0))_+.
\end{equation*}
\end{thm}

%%%%%%%%%
\begin{proof}		%% [Sketch of Proof of Proposition~\ref{prop1.1}]
%I only sketch the argument when $H$ is nonlinear. 
%In the linear case  the result follows 
%from the classical theory of viscosity solution.  
%function of its arguments. 
To simplify the presentation, it is assumed that  $F$ is smooth.   The actual proof follows by writing finite differences instead of taking 
derivatives and using regularizations.
\smallskip
	   
Since $\Phi (v,0)=v$,  \eqref{semilinearas} yields that,   for 
fixed $\delta >0$, it is possible to choose $h >0$ so small that
\begin{equation}\label{eq2.1} 		%% \leqno(2.1)
\sup_{0\le t\le h}\big[ |\Phi (v,t)-v|+|\Phi'(v,t)-1|+|\Phi''(v,t)|
+|\Phi'''(v,t))\big]\le\delta.
\end{equation}
Next consider the new change of variables 
$$v=\phi (z)=z+\delta\psi (z)\ \text{ with }\ \phi' >0\ .$$

If $v$ is a subsolution (resp.  supersolution) of  \eqref{eq1.10}, then $z$ 
is a subsolution  (resp. supersolution) of
\begin{equation}\label{eq2.2}
z_t=\widetilde{\widetilde F} (D^2z, Dz,z),
\end{equation}
with  
\begin{equation}\label{eqtakis2.21}
%\begin{cases}	
\begin{split}	%% \leqno(2.2)
\widetilde{\widetilde F} (X, p, z)&=  \dfrac1{\Phi'(\phi (z),t)\phi'(z)} F\big(\Phi'(\phi
	   (z),t)[\phi'(z)X+\phi''(z)(p\otimes p)]\\[2mm]
	   &+\Phi''(\phi (z),t)(\phi'(z))^2(p\otimes p), \Phi'(\phi
	   (z),t)\phi'(z)p\big).
\end{split}
	  % \end{cases}
%\noalign{\vskip6pt}
%&\hskip1.3truein +\Phi''(\phi (z),t)(\phi'(z))^2(Dz\otimes Dz), \Phi'(\phi
%	   (z),t)\phi'(z)Dz\big)\\
%& =.
\end{equation}

The comparison result follows from the classical theory of viscosity solutions,
%(see \cite{B2} and \cite{CIL}), 
if there exists 
$C= C_R> 0$ where  $R = \max (\|\bar v\|,\|\underline v\|)$,  such that, 
for all $X$, $p$ and $z$, 
\begin{equation}\label{eq2.3}		%% \leqno(2.3)
\frac{\partial}{\partial z}
\skew3\widetilde{\widetilde F} (X,p,z)\le C.
\end{equation}

A straightforward calculation yields
\begin{equation*}
\begin{split}
\frac{\partial}{\partial z}\skew3\widetilde{\widetilde F}(X,p,z)
& =-\frac{(\Phi'\phi')'}{(\Phi'\phi')^2} F +\\
\noalign{\vskip6pt}
&\quad  \frac1{(\Phi'\phi')}\Big[ \langle D_XF, \big[(\Phi'\phi')'X+
[(\Phi'\phi'')'+(\Phi''(\phi')^2)'](p\otimes p)\rangle \big]
+ \langle D_pF, (\Phi'\phi')'p)\rangle \Big]\\
\noalign{\vskip6pt}
& = \frac{(\Phi'\phi')'}{(\Phi'\phi')^2}\big[-F + \langle D_XF, (\Phi'\phi'X
+ (\Phi'\phi''+\Phi''(\phi')^2)(p\otimes p))\rangle 
% + F_p\cdot (\Phi'\phi')p\big]\\ 
+ \langle D_p F,\Phi' \phi' p \rangle] \\
\noalign{\vskip6pt}
&\qquad + \langle D_X F,
\Big[ \frac{(\Phi'\phi''+\Phi''(\phi')^2)'}{\Phi'\phi'} -
\frac{(\Phi'\phi''+\Phi''(\phi')^2)(\Phi'\phi')''}{(\Phi'\phi')^2}
\Big](p\otimes p)\rangle,
\end{split}\end{equation*}
where, to simplify the notation, 
the arguments of $F$, $D_pF$, $D_X^2 F$, $\Phi'$,
$\Phi''$, $\phi'$  and $\phi''$ are omitted.
\smallskip

In view of \eqref{degenerateelliptic} and \eqref{eq1.11}, to obtain \eqref{eq1star} 
it suffices  to choose $\phi$ so that 
\begin{equation}\label{eq2.4} 			%% \leqno(2.4)
\frac{(\Phi'\phi''+\Phi''(\phi')^2)'}{\Phi'\phi'} -
\frac{(\Phi'\phi''+\Phi''(\phi')^2)(\Phi'\phi')''}{(\Phi'\phi')^2}\le 0
\end{equation}   
and, if the second inequality in  \eqref{eq1.12} holds, 
\begin{equation}\label{eq2.5} 			%% \leqno(2.5)
\frac{(\Phi'\phi')'}{(\Phi'\phi')^2} \le 0
\end{equation}   
or, if the first inequality in \eqref{eq1.12} holds, 
\begin{equation}\label{eq2.6} 			%% \leqno(2.6)
\frac{(\Phi'\phi')'}{(\Phi'\phi')^2} \ge 0.
\end{equation}
%\smallskip

Assumption \eqref{eq2.1} and the special choice of $\phi$ yield that 
\eqref{eq2.4} is satisfied  if $\psi'''\le -1$, and that \eqref{eq2.5} 
(resp.\ \eqref{eq2.6}) holds,  
if $\psi''\le -1$ (resp.\ $\psi''\ge 1$).
It is a simple exercise to find $\psi$ 
so that \eqref{eq2.4} and either \eqref{eq2.5} or \eqref{eq2.6} hold 
in its domain of definition.
\smallskip

The classical comparison result for viscosity solutions 
then  yields that, if $\overline v(\cdot,0)\le
\underline v(\cdot,0)$ on $\R^d$, then $\overline v\le \underline v$ 
on $\R^d\times [0,h]$.  
The same argument then yields the comparison in $[h,2h]$, etc..

\end{proof}
The existence of the pathwise solutions of 
\eqref{eq:initialvalue} 
is based on the stability properties of the ``approximating '' initial value problem %  \eqref{eq:initialvalue} by 
\begin{equation}\label{eq1.14}		%% \leqno(1.14)
%\begin{cases}
u^\ep_t =F(D^2u^\ep, Du^\ep )+\sum^m_{i=1}H^i(u^\ep)
\dot B_{i}^\ep  \ \text{ in \  $Q_T$} \quad
u^\ep(\cdot,0) =u^\ep_0, %  \  \text{on}  \  \R^d,
%\end{cases}
\end{equation}
where  $u^\ep_0\in BUC(\R^d)$, and
\begin{equation}\label{eq1.15}		%%\leqno(1.15) 
\begin{cases}
%\zeta^\ep = (\zeta^\ep_1,\ldots ,\zeta^\ep_M)
B^\ep= (B_1^{\ep},\ldots, B_m^{\ep}) \in C^1([0,\infty);\R^m),\\
\noalign{\vskip6pt}
\hbox{and, for all $T >0$, as $\ep\to 0$, }\ 
B^\ep\rightarrow B  \ \text{in the  rough path metric.}
\end{cases}
\end{equation}

Note that, if $m=1$,  the assumption in \eqref{eq1.15} can be reduced  to $B^\ep\rightarrow B$  uniformly on $[0,T]$. 
%For example, $B^\ep$ can be simply a convolution or a finite difference approximation of $B$.
 
%\marginpar{{\bf Takis:}\newline add ref \& \newline exact result 
%\newline here}
%Such approximations to $W$ arise by either simply mollifying $W$ or
%suitably scaling a given function $\tilde W$, for example 
%%$\zeta_\ep (t)=\ep^{-1}\zeta (\ep^{-2}t)$, 
%$W_t^\ep = \ep^{-1}\tilde W_{\ep^{-2}t}$,  
%which satisfies appropriate mixing conditions.
%\medskip
\smallskip

The existence result is stated next.

\begin{thm}\label{thm1.2}		%% {\bf Theorem 1.2}: 
%Assume \eqref{deg.ell}, \eqref{semilinearas}, 
%%% \eqref{eq4.20bis}, 		%\eqref{eq:old-0.4}, 
%and, if $H$ is nonlinear,  \eqref{eq1.11} and \eqref{eq1.12}. 
Assume  \eqref{degenerateelliptic}, \eqref{eq4.20bis}, \eqref{semilinearas},
 \eqref{eq1.11} and \eqref{eq1.12} and fix $T >0.$ 
Let $(\zeta^\ep )_{\ep >0}$ and $(\xi^\eta )_{\eta >0}$ 
satisfy \eqref{eq1.15} and consider the solutions $u^\ep, v^\eta\in
BUC(\overline Q_T)$ of \eqref{eq1.14} with initial datum $u^\ep_0$ and
$v^\eta_0$ respectively.  
If, as $\ep ,\eta\to 0$, 
$ u^\ep_0-v^\eta_0 \to 0$ uniformly on $\R^d$, then, as $\ep,\eta\to0$,   
$u^\ep - v^\eta \to0$ uniformly on $\overline Q_T$. 
%a.s.\ and for every $T> 0$.
In particular, each family $(u^\ep)_{\ep >0}$ is Cauchy in $\overline Q_T.$
Hence, it converges uniformly  to $u\in BUC (\overline Q_T)$, 
which is a pathwise  viscosity solution to 
\eqref{eq:initialvalue}.  	%% %\eqref{eq:old-0.1}.  
Moreover, all approximate families converge to the same limit.
\end{thm}

The proof of Theorem~\ref{thm1.2} follows from the comparison between subsolutions  and
supersolutions of \eqref{eq1.10} for different approximations 
$(\zeta^\ep)_{\ep>0}$ and $(\xi^\eta)_{\eta >0}$.
Since a similar theorem will be proved later when dealing with nonlinear 
gradient dependent $H$, the proof is omitted. %Theorem~\ref{thm1.2}. 
\smallskip

Finally the next result is  about the Lipschitz continuity of the solutions.  
Its proof is based on the comparison
estimate obtained in %Proposition~\ref{prop1.1} 
Theorem~\ref{thm1.11} and, hence, it is omitted.

\begin{prop}\label{prop1.3}		%% 
%Assume \eqref{deg.ell}, \eqref{semilinearas}, 
Fix $T$ and assume  \eqref{degenerateelliptic}, \eqref{eq4.20bis}, \eqref{semilinearas},
 \eqref{eq1.11} and \eqref{eq1.12} and let $u\in BUC (\overline Q_T)$ be the unique pathwise solution to  \eqref{eq:initialvalue} for  $u_0\in C^{0,1}(\R^d)$. 
%% \eqref{eq4.20bis},      %%  \eqref{eq:old-0.4}, 
%and, if $H$ is nonlinear,  \eqref{eq1.11}, \eqref{eq1.12}.  
%Let $u$ be the unique stochastic viscosity solution of 
%\eqref{eq:initialvalue}.  	%% %\eqref{eq:old-0.1}.  
Then $u(\cdot, t)\in C^{0,1}(\R^d)$ for all $t \in [0,T],$
and  there exists 
$C=C(F, H, B, T)>0$ such that, for all $t\in [0,T]$, 
$ \| Du(\cdot, t)\| \leqq C.
$
\end{prop}

Of course Proposition~\ref{prop1.3} is immediate if $F$ and $H$ do not depend on $x$. The point is that the clainm holds in full generality. 
%After the work described above was announced, Buckdahn and Ma \cite{buckdahnma2, buckdahnma3} used  the  map  \eqref{eq1.9}, which is known as the Sussman-Doss transformation, to study equations similar to \eqref{eq:initialvalue} but for a more restrictive class of $F$'s and proved wellposedness under the assumption that the transformed initial value problem admits a comparison principle.
%
%\end{document}

%%%%%%%%%%%%%%%%%%%%%%%%%%

\section{The extension operator for spatially homogeneous first-order problems}
%simplest nonlinear  pde with rough signals 
%as a limit of regular approximations}
%%%%%%%%%%%%%%%%%%%%%%%%%%%%%%%%%%
The object  here is the study the space homogeneous Hamilton-Jacobi equation
\begin{equation}\label{paris100}
du=\sum_{i=1}^m H^i(Du)\cdot dB_i \ \text{in} \ Q_\oo \quad u(\cdot,0)=u_0,
\end{equation}
with $B=(B_1,\ldots, B_m) \in C_0([0,\oo);\R^m)=\{B\in C([0,\oo);\R^m) : B(0)=0\}.$
\smallskip

%In this section it is assumed again that $m=1$. However, in view of the discussion in Section~3 about when the characteristics commute, the results easily extend to the multi-path setting for spatially homogenous Hamiltonians. 
%\smallskip
%
%To make some of the statements below shorter, it is convenient to introduce the space $C_0([0,\infty);\R)=\{B\in C([0,\infty);\R) :  B(0)=0\}$
%

The aim  is to show that, if $H=(H_1,\ldots,H_m)\in C^{1,1}_{\text{loc}}(\R^d; \R^m)$,  the solution operator of \eqref{paris100} with smooth paths  has a unique extension to the set of continuous paths.

%\begin{equation}\label{ShJ2}	
%du = H(Du) \cdot dB\ \text{in} \  Q_\oo \quad 
%u(\cdot, 0)= u_0, 
%\end{equation}
%for smooth paths, which is well defined even when $H$ is merely continuous,  has a unique extension 
%to the set of continuous paths.
%\smallskip
%

%Fix  $H\in C^{1,1}_{\text{loc}}(\R^d)$, $u_0\in \text{BUC}(\R^d)$  and $B\in C_0([0,\infty);\R)=\{B\in C([0,\infty);\R) :  B(0)=0\}$. Let  $(B^\ep)_{\ep >0}$ be a family of smooth paths in $C_0([0,\infty);\R)$  which approximate $B$ in $ C([0,\infty)),$ and $(u^\ep_0)_{\ep >0}$ a family in $BUC(\R^d)$ such that, as $\ep\to0$,  $u^\ep_0 \to u_0$ uniformly on $\R^d$, and 
%% that is, as $\ep\to 0$, $B^\ep \to B$ locally uniformly in $[0,\infty),$ 
% consider the family of Hamilton-Jacobi equations
%\begin{equation}\label{ShJ3}	%% \label{eq:34}
%%\begin{cases}
%u^\ep_t =H(Du^\ep)\dot B^\ep \ \text{in} \ Q_\oo \quad 
%%\noalign{\vskip6pt}
%u^\ep(\cdot,0)  = u_0^\ep, % \quad\text{on}\quad \R^d,
%%\end{cases}
%\end{equation}
%which, for each $u_0^\ep \in BUC(\R^d)$ and $H\in C(\R^d)$, have a unique solution $u^\ep \in BUC(\overline Q_\oo)$.
%%\smallskip
%\vskip.05in
\smallskip

The result  is stated next. 

\begin{thm}\label{paris1000}
Fix  $H\in C^{1,1}_{\text{loc}}(\R^d;\R^m)$, $u_0\in \text{BUC}(\R^d)$  and $B\in C_0([0,\infty);\R^m)$. 
There exists a unique $u\in \text{BUC}(\overline Q_\oo)$ such that,  for any families    $(B^\ep)_{\ep >0}$   in $C_0([0,\infty);\R^m)\cap C^1([0,\infty);\R^m)$  and $(u^\ep_0)_{\ep >0}$ in $\text{BUC}(\R^d)$ which approximate respectively  $B$ in $ C([0,\infty);\R^m)$ and $u_0$ in $\text{BUC}(\R^d)$,  if $u^\ep \in \text{BUC}(\overline Q_\oo)$ is the unique viscosity solution of $du^\ep=\sum_{i=1}^m H^i(Du^\ep)\cdot dB^\ep _i \ \text{in} \ Q_\oo \ \text{and} \ u^\ep(\cdot,0)=u^\ep_0,$ then, as $\ep\to 0$, $u^\ep \to u$ uniformly in $\overline Q_\oo$. 
\end{thm}

%$u^\ep_t =H(Du^\ep)\dot B^\ep \ \text{in} \ Q_\oo$ with $ 
%u^\ep(\cdot,0)  = u_0^\ep,$
%and 
%% that is, as $\ep\to 0$, $B^\ep \to B$ locally uniformly in $[0,\infty),$ 
% consider the family of Hamilton-Jacobi equations
%\begin{equation}\label{ShJ3}	%% \label{eq:34}
%%\begin{cases}
%u^\ep_t =H(Du^\ep)\dot B^\ep \ \text{in} \ Q_\oo \quad 
%%\noalign{\vskip6pt}
%u^\ep(\cdot,0)  = u_0^\ep, % \quad\text{on}\quad \R^d,
%\end{equation}
%then, as $\ep\to 0$, $u^\ep \to u$ uniformly in $\overline Q_\oo$. 
%\end{thm}

This unique limit will be also characterized later as the unique pathwise solution of \eqref{paris100}.

\smallskip

The claim follows from the next theorem which asserts that, if the family of smooth paths   $(B^\ep)_{\ep >0}$ and initial data $(u^\ep_0)_{\ep >0}$ are Cauchy  in $C_0([0,\infty);\R)$ and $\text{BUC}(\R^d)$ respectively, then the solutions $u^\ep \in \text{BUC}(\overline Q_\oo)$ of 
\begin{equation}\label{ShJ3}	%% \label{eq:34}
du^\ep=\sum_{i=1}^m H^i(Du^\ep)\cdot dB^\ep _i \ \text{in} \ Q_\oo \quad u^\ep(\cdot,0)=u^\ep_0, 
\end{equation}
form a Cauchy family in $\text{BUC}(\overline Q_\oo)$.  

\begin{thm}\label{thm:simplest} 
Fix  $H\in C^{1,1}_{\text{loc}}(\R^d;\R^m)$. Let $\zeta^\ep, \xi^\eta \in C_0([0,\infty);\R^m)\cap C^1([0,\infty);\R^m)$  and $u_0^\ep, v_0^\eta \in \text{BUC}(\R^d)$ be such that, as $\ep, \eta \to 0$,  $\zeta^\ep - \xi^\eta \to 0$ in  $C([0,\infty);\R^m)$ and 
$u_0^\ep - v_0^\eta \to 0$  in $\text{BUC}(\R^d)$. If $u^\ep, v^\eta \in \text{BUC}(\overline Q_\oo)$ are the viscosity solutions of \eqref{ShJ3} with respective paths and initial condition $(\zeta^\ep, u_0^\ep)$, $(\xi^\eta, v_0^\eta)$, then,  as $\ep, \eta \to 0$, $u^\ep - v^\eta \to 0$  in $\text{BUC}(\overline Q_\oo)$.
%that Fix $B\in C([0,\infty))$ with $B(0)=0$ and assume that $H\in C_{\loc}^{1,1} (\R^d)$ and $u_0 \in BUC (\R^d)$.  
%Let $(B^\ep)_{\ep >0}$ in $C^1((0,\infty))\cap C([0,\infty))$ with $B^\ep(0)=0$ and $(u_0^\ep)_{\ep > 0}$ in $BUC(\R^d)$ be such that, as $\ep \to 0$, $B^\ep \to B$ in $C([0,\infty))$ and  $u_0^\ep\to u_0$ in $BUC (\R^d)$. There exists $u\in C(\R^d \times [0,\infty))$ such that, for each $T >0$, $u\in BUC(\overline Q_T)$ and, 
%if $u^\ep \in BUC(\overline Q_T)$ is the solution to \eqref{ShJ3}, then, as $\ep \to 0$, $u_\ep\to u$ in $BUC(\overline Q_T)$.
\end{thm}

\begin{proof}
%Let $(B^\ep)_{\ep > 0}$ and 
%$(\widetilde B^\eta)_{\eta > 0}$  be two  $C^1$-approximations of $B$ in $C([0,\infty))$ such that $B^\ep(0)=\widetilde B^\eta(0)=0$, 
%% that is  two families of $C^1$-functions such that, \widetilde B^\eta
%%as $\ep,\eta \to0$, $B^\ep - \widetilde B^\eta \to0$ in $C([0,\infty))$ 
%and consider two 
%approximating families  $(u^\ep_0)_{\ep>0}$ and $(\tilde u^\eta_0)+{\eta >0}$ of $u_0$ in $BUC (\R^d).$
%%that is   $u_0^\ep,\tilde u_0^\eta \in BUC (\R^d)$ and, as $\ep,\eta\to0$, 
%%$\|u_0^\ep - \tilde u_0^\eta\|\to 0$. 
%Fix $T >0$ and let $u^\ep,\tilde u^\eta \in BUC (\overline Q_T)$ 
%be the viscosity solutions to \eqref{ShJ3} with paths $B^\ep, \widetilde B^\eta$ and initial datum  $u^\ep_0, \tilde u^\eta_0$.  
%The claim follows if it is shown that, as $\ep,\eta\to0$, $u^\ep - \tilde u^\eta\to0$ in 
%$\R^d\times [0,T]$.
%\smallskip
%
A simple density argument implies that it is enough to consider 
$u_0^\ep, v_0^\eta \in C^{0,1} (\R^d)$.
% with  $\max(\|Du_0^\ep\|,\|Dv_0^\eta\|) \leqq C$ for some $C>0$. 
Since $H$ is independent of $x$, it follows that $u^\ep, v^\eta \in C^{0,1} (\overline Q_\oo)$ and, for all $t>0$,  
 $\max (\|Du^\ep (\cdot,t) \|, \|Dv^\eta (\cdot,t)\|) 
\leqq \max (\| Du_0^\ep \|,\|D v_0^\eta\|)$. 
Hence, without any loss of generality, it may be assumed  that $H\in C^{1,1}(\R^d)$. 
%\smallskip
\vskip.075in
Notice that, for each $\ep$ and $\eta$, $u^\ep$ and $v^\eta$
are actually also Lipschitz continuous in time.  The  
Lipschitz constants in time, however,  depend on $|\dot \zeta^\ep|$ and $|\dot {\xi}^\eta |$, and, hence,  
are not bounded uniformly in $\ep, \eta.$
This is one of the main reasons behind the difficulties here. 
\smallskip

To keep the arguments simple, it is also  assumed that $u_0^\ep$ and $v_0^\ep$, and,  hence,  
$u^\ep$ and $v^\eta$ are periodic in the unit cube
$\tee^d$.
This simplification allows  not be concerned about infinity, and, more 
precisely, the possibility that the suprema  below are not achieved.
The periodicity can be eliminated as an assumption by introducing
appropriate penalizations at infinity that force the sup's to be actually maxima.
\smallskip

Finally, from now on I assume that $m=1$. This is only done to keep the notation simpler. Since the equation does not depend on the space variable, the extension to $m>1$
is immediate
\smallskip

The general strategy in the theory of viscosity solutions 
to show that,  
as $\ep,\eta\to0$, $u^\ep - v^\eta\to0$ in $\overline Q_\oo$, 
is to  double the 
variables and consider the function 
$$z(x,y,t) = u^\ep (x,t) - v^\eta (y,t),$$
which satisfies the so-called   ``doubled'' initial value problem 
\begin{equation}\label{doubled1}	%% \label{eq:35} 
%\begin{cases} 
z_t = H(D_x z)\dot \zeta^\ep - H(-D_y z) \dot {\xi}^\eta \ \text{in} \  
\R^d\times\R^d\times (0,\infty) \quad 
%\noalign{\vskip6pt}
z(x, y,0)= u_0^\ep (x) - v_0^\eta (y). % \ \text{in} \  \R^d\times\R^d,
%\end{cases}
\end{equation}
%where, for $\lambda > 0$, 
%with  
%$$z_0 (x,y) = u_0^\ep (x) - u_0^\eta (y).$$%+ \sup(u_0^\ep -\tilde u_0^\eta) \ .$$

The assumptions on $u_0^\ep$ and  $v_0^\eta $ imply that, 
for $\lambda >0$, there exists $\theta (\lambda )>0$ such that $\theta (\lambda )\to0 \ \text{as} \  \lambda\to\infty$ and 
\begin{equation}\label{takis1} z_0 (x,y)   %% = u_0^\ep (x) - u_0^\eta (y) 
\leqq \lambda |x-y|^2 +\theta(\lambda) + \sup(u_0^\ep -v_0^\eta).
\end{equation}

To conclude, it suffices  to show  that there exists 
%will follow, if we can establish a bound like 
%$$z(x,y,t) \leqq U^{\ep,\eta,\lambda} (x,y,t)$$
%% $$z \leqq U^{\ep,\eta,\lambda} \ \text{ in }\ \R^d\times \R^d 
%% \times [0,T]\ , $$
$U^{\ep,\eta,\lambda}: \R^d\times\R^d\times [0,T]\to\R$ 
such that, as $\ep,\eta\to 0$ and $\lambda \to \infty$,
$$U^{\ep,\eta,\lambda}(x,x,t) \to 0 \ 
\text{uniformly in} \ \overline Q_\oo \   \text{and} \ 
z \leqq U^{\ep,\eta,\lambda} \ \text{ in }\ \R^d\times \R^d
\times [0,\oo).$$

It would then follow that 
$$\lim_{\ep,\eta \to0}  \sup_{(x,t) \in \overline Q_\oo} z(x,x,t) =0,$$
which is one part of the claim. 
The other direction is proved similarly.
%we consider $u^\eta (x,t) - u^\ep (y,t)$ and argue similarly.
\smallskip

Again, as in  the general ``non rough'' theory,    it is  natural to try to show that  there exists, for   some $C >0$ and $a(\lambda) >0$ such that  $a(\lambda)\to 0$ as $\lambda\to\infty$, a supersolution of \eqref{doubled1} of  the form  
$$ U^{\ep,\eta,\lambda} (x,y,t) = 
 C \lambda |x-y|^2 + a(\lambda).$$
%with $a(\lambda)\to 0$ as $\lambda\to\infty$.
This is, however, the main difficulty, since both the $C$ and  $a(\lambda)$ will depend on  
$|\dot \zeta^\ep|$ and $|\dot {\xi}^\eta|,$ 
which  are not bounded uniformly in $\ep$ and $\eta$.
\smallskip

The first new idea to circumvent this difficulty is to find sharper 
upper bounds by considering the solution 
$\phi^{\lambda,\ep,\eta}(x,y,t) $ 
of
\begin{equation}\label{paris200} \phi^{\lambda,\ep,\eta}_t =H(D_x \phi^{\lambda,\ep,\eta})\dot \zeta^\ep - H(-D_y \phi^{\lambda,\ep,\eta}) \dot {\xi}^\eta \ \text{in} \  
\R^d\times\R^d\times (0,\infty) \quad 
%\noalign{\vskip6pt}
\phi^{\lambda,\ep,\eta}(x, y,0)=\lambda |x-y|^2,
\end{equation}

which, in view of the spatial homogeneity of $H$ and the fact that $\phi^{\lambda,\ep,\eta}(\cdot,\cdot,0)$ depends on $x-y$, is given by 
%
%%\eqref{doubled1} 
%%\begin{equation}\label{doubled2}	%% \label{eq:36}
%%\begin{cases}
%%\phi_t = H(D_x\phi) \dot W^\ep - H(-D_y\phi)\dot W^\eta\quad\text{in}\quad
%%\R^d\times\R^d\times  (0,\infty)\ ,\\
%%\noalign{\vskip6pt}
%%\phi = \phi_0 \quad\text{on}\quad \R^d\times\R^d\times  \{0\} \ ,
%%\end{cases}
%%\end{equation}
%with 
%$\phi_0 (x,y) = \lambda |x-y|^2\ .$
%\smallskip
%
%Note that, since $H$ depends only on the gradient and $\phi_0$ 
%depends only on $x-y$, %it follows that 
$$\phi^{\lambda,\ep,\eta} (x,y,t) = \Phi^{\lambda,\ep,\eta} (x-y,t),$$
with  $\Phi^{\lambda,\ep,\eta}$ solving the initial value problem 
\begin{equation}\label{doubled3}	
\Phi_t = H(D\Phi) (\dot \zeta^\ep - \dot {\xi}^\eta) \ \text{in} \ Q_\oo
 \quad 
\Phi (\cdot, 0) =\lambda |z|^2, %\quad\text{on} \  \R^d.
\end{equation}
%with  $\Phi_0 (z): = \lambda |z|^2.$
%\smallskip
which is well-posed for each $\ep, \eta$.
\smallskip

The classical  comparison principle for viscosity solutions yields, that 
for all $x,y\in\R^d$, 
$t\geqq 0$, $\lambda >0$ and $\ep, \eta$,
$$z(x,y,t) \leqq \phi^{\lambda,\ep,\eta} (x,y,t) +  
\max_{x,y\in \R^d} (z(x,y,0) - \lambda |x-y|^2), $$
and, hence,  
for all $x\in \R^d$ and $t\geqq0$,
$$u^\ep (x,t) - u^\eta (x,t) 
\leqq \phi^{\lambda,\ep,\eta}(x,x,t) 
+ \theta (\lambda) + \sup(u^\ep_0- v^\eta_0) .$$
To conclude, it is necessary to show that there exists $\Theta(\lambda) >0$ such that $\lim_{\lambda \to \infty}\Theta(\lambda)=0$ and  
%such that, as $\lambda\to\infty$, $\Theta (\lambda)\to0$ and 
$$\varlimsup_{\ep,\eta\to0} \sup_{x\in\R^d} \phi^{\lambda,\ep,\eta}(x,x,t) 
\leqq \Theta (\lambda),$$
a fact that apriori may present 
a problem since 
the ``usual'' viscosity theory yields the existence of $\phi^{\lambda,\ep,\eta}$ but not the 
desired uniform estimate.
\smallskip

Here comes the second new idea, namely, to use the characteristics to construct a smooth solution 
$\phi^{\lambda,\ep,\eta}$, at least for a small time, which, of course, depends on $\ep$ and $\eta$. 
The aim then will be to show that, as $\ep,\eta\to0$, 
the interval of existence becomes of order one.
\smallskip

The characteristics of  \eqref{paris200} 
are 
\begin{equation}\label{sch}	
\begin{cases}
\dot X =  -DH (P) \dot \zeta^\ep \quad X(0)= x \qquad \dot Y = -DH (Q)\dot {\xi}^\eta \quad Y(0) =y,\\
\noalign{\vskip6pt}
\dot P = 0\ \quad \dot Q=0 \quad P(0) = Q(0) = 2\lambda (x-y), \\
\noalign{\vskip6pt}
\dot U = \left(H(P)-\langle D_p H(P),P\rangle\right) \dot \zeta^\ep - (H(Q) - \langle DH (Q), Q\rangle) \dot {\xi}^\eta \quad U(0) = \lambda |x-y|^2.
%\noalign{\vskip6pt}
%X(0)= x\ ,\quad Y(0) =y\ ,\qquad 
%P(0) = Q(0) = 2\lambda (x-y)\ ,\qquad 
%U(0) = \lambda |x-y|^2.\\
%\noalign{\vskip6pt}
\end{cases}
\end{equation}
%Note that to simplify the notation and keep the equations in \eqref{sch} 
%``symmetric'' we write the characteristic for $Q$ instead of $-Q$.
%\smallskip
Note that to keep the equations simpler  the system is written for $Q(x,t)=-Dv^\eta (Y(t),t)$ instead of $Dv^\eta (Y(t),t).$ Similarly I ignore the dependence on $\lambda, \ep$ and $\eta$. 
\smallskip

The method of characteristics provides a classical 
solution of \eqref{paris200}  for some short time $T_{\ep,\eta,\lambda}^*$ 
as long the map $(x,y) \to (X(t), Y(t))$ 
is invertible. 
%%One way to guarantee  this is to show that the Jacobian of this map, 
%which at $t=0$ is $1$, does not vanish in $[0,T_{\ep,\eta,\lambda}^*)$.
%Finding a sharp estimate for the first time the Jacobian vanishes is, in general, a 
%difficult problem. In the case at hand  this is not a problem, since the simplest estimate does the job.
%%while with some work, it is possible to obtain a rough one. 
%It is, however,  one of the major technical difficulties when trying 
%to study equations which depend on $x$ and have multiple paths. 
\smallskip

The special structure of \eqref{sch} yields that,  for all $t\geq 0$,  
%, which follows from the fact that $H$ depends only on the gradient, yields that, for all $t\geq 0$,  
%$P$ and $Q$ are consthomogeneous in space yields that allows for many simplifications. 
%Indeed, it follows from \eqref{sch} that $P$ and $Q$ are constant in 
%time and, in view of the choice of the initial datum, 
$$P(t) = Q(t) =  2\lambda (x-y), $$
%Then
%$$%% \mathop{X\!-\!Y}\limits^{\textstyle\dot\frown}  %\frac{d}{dt} (X-Y) 
%\dot X - \dot Y 
%= - DH(2\lambda (x-y)) (\dot B^\ep - \dot {\widetilde B}^\eta)\ ,$$
%and, hence, for $t\geqq0$,
and 
$$(X-Y) (t) = (x-y) - DH(2\lambda (x-y)) (\zeta^\ep(t) -\xi ^\eta(t)).$$
To simplify the notation, let $z= x-y$ and $Z(t) = X(t)-Y(t)$, in 
which case the last  equation can be rewritten as 
$$Z(x,t) = z-DH (2\lambda  z) (\zeta^\ep(t) -\xi ^\eta(t)).$$
Note that 
$z \to Z(z,t)$ is the position characteristic associated with the simplified 
initial value problem \eqref{doubled3}, and, in the problem at hand, 
is  the only map  that needs to be inverted. 
%Hence the aim is to estimate the time during which $z\mapsto Z$ is a
%one-to-one map. 
Since 
$$D_zZ(z,t)=I + 2\lambda D^2 H ( 2\lambda z) (\zeta^\ep(t) -\xi^\eta(t)),$$
it follows that the map 
%and clearly its Jacobian does not vanish 
$z\mapsto Z$ is invertible
as long as 
\begin{equation}\label{paris201}
\sup_{t\in[0,T]} |(\zeta^\ep(t) -\xi ^\eta(t))|\, \|D^2 H\|_\infty 
< (2\lambda)^{-1}.
\end{equation}

This is, of course, possible for any  $T$ and $\lambda$ provided  
$\ep$ and $\eta$ are small, since, as $\ep,\eta\to0$, $\zeta^\ep -\xi^\eta \to 0$ 
in $C([0,\oo))$. 
\smallskip

The above estimates depend on having $H\in C^2$. 
Since  the interval of existence depends only on the $C^{1,1}$ bounds of $H$, it can be 
assumed that $H$ has this regularity and then conclude 
introducing yet another level of approximations.
%If $H$ is less regular than $C^{1,1}$, 
%then there is a nontrivial interaction between 
%the regularities of $H$ and $W$.
%This is described and analyzed in detail in the next section.
\smallskip

It now follows that %Next  observe that   the characteristics yield  
\begin{equation*}
\begin{split}
&\phi^{\lambda,\ep,\eta} (X(t),Y(t),t) 
 = \lambda |x-y - DH (2\lambda (x -y)) (\zeta^\ep(t) -\xi ^\eta(t))|^2 \\
\noalign{\vskip6pt}
&\qquad + 
\left[ H(2\lambda (x - y)) - \langle D H(2\lambda (x- y)), 
2\lambda (x - y)\rangle \right]  (\zeta^\ep(t) -\xi ^\eta(t))\ .
\end{split}
\end{equation*}

Moreover, it follows from \eqref{paris200}, that   % Since $|\zeta^\ep(t) -\xi ^\eta(t)| \, \|D^2 H\|_\infty <  
there exists  $C >0$ depending  only on  
$\|H\|_{C^{1,1}}$ such that 
$$|\phi^{\lambda,\ep,\eta} (X(t),Y(t),t) - \lambda |x - y|^2|
\le \lambda C \sup_{0\leqq t\leqq T} 
|\zeta^\ep(t) -\xi ^\eta(t)|\ .$$

Returning to the $x,y$ variables, the above estimate gives that, for each fixed 
$\lambda >0$ and $T >0$ and as $\ep,\eta\to0$,
%$$\phi_{\ep,\eta}^\lambda \to \frac{\lambda}2 |x-y|^2$$
$$\sup_{\substack{x,y\in\R^d\\ t\in [0,T]}} 
(\phi^{\lambda,\ep,\eta}  (x,y,t) - \lambda |x-y|^2) \to 0.$$
%which concludes the proof.
\end{proof}

\subsection*{A summary of the general startegy} Since the approach and  the arguments of the  proof above are used several times in the theory and the notes, it is helpful 
to  present a brief summary of the main points. 
% and the strategy presented 
%above so that it is possible to refer to it later in the notes.
\smallskip

The conclusion of the theorem is that it is possible to construct, 
using the classical theory of viscosity solutions, a (unique) 
$u\in BUC (\overline Q_\oo)$, which is the candidate for the solution 
of  \eqref{paris100} for any $B$ continuous as long as $H\in C_{\text{loc}}^{1,1}$.
\smallskip

The key technical step in the proof was the fact that, if, 
as $\ep,\eta\to0$, 
$\zeta^\ep - \xi^\eta \to0$ in $C([0,\infty))$, then, for each $\lambda > 0$ 
and $T>0$, as $\ep,\eta\to0$ 
%$$v(z,t) = v_{\ep,\eta}^\lambda(z,t) \to \lambda |z|\ ,$$
$$\sup_{\substack{z\in \R^d\\  t\in [0,T]}} 
\big(v_{\ep,\eta}^\lambda (z,t) - \lambda |z|\, \big)\to 0\ ,$$
where  $v= v_{\ep,\eta}^\lambda$ is the solution of the initial value problem 
\begin{equation}\label{paris202}
v_t = H(Dv) (\dot \zeta^\ep - \dot {\xi}^\eta) \ \text{in} \  Q_\oo
 \quad 
v(z,0) = \lambda |z|. \end{equation}

The proof presented earlier used $\lambda |z|^2$ as initial condition in  \eqref{paris202}.
It is not hard to see, however, that the same argument will work for 
initial datum $\lambda |z|$. 
Indeed it is enough to consider regularizations like $(\delta + |z|^2)^{1/2}$ 
and to observe that the estimate on $u^\ep (\cdot,t) -{v}^\eta(\cdot,t)$ 
is uniform 
on $\delta$ in view of the assumption that $H\in C^{1,1}$. 
The conclusion for $\lambda |z|$ then follows from the stability properties 
of viscosity solutions. 
\smallskip

The result about the extension can be summarized as follows.
\smallskip

Given sufficiently regular paths  $B_n=(B_{1,n},\ldots, B_{m,n}):[0,\infty)\to \R^m$, $H=(H^1,\ldots, H^m) \in C(\R^d;\R^m)$   and $\lambda >0$,   let $v_{n,t}^\lambda \in \text{BUC}(\overline Q_\oo)$ be the solution of 
\begin{equation}\label{paris210}
v_{n,t}^\lambda = \sum_{i=1}^mH^i(Dv_n^\lambda) \dot B_{n,i}\ \text{ in } \  Q_\oo \quad v_n^\lambda(z,0)= \lambda |z|.
\end{equation}

The following theorem gives a sufficient condition for the existence of the extension.

\begin{thm}\label{paris210.5}
If for every $B_n\in C_0([0,\oo);\R^m)\cap C^1 ([0,\infty);\R^m)$ such that, 
as $n\to\infty$,  $B_n \to0$  in $C([0,\oo);\R^m)$,  and $T>0$, the solution $v_{n,t}^\lambda$ of \eqref{paris210}
has the property 
\begin{equation}\label{paris211} 
\lim_{n\to \oo, \lambda \to \oo} \sup_{(z,t)\in\R^d\times [0,T]} 
\big (v_n^\lambda (z,t) - \lambda |z| \big ) =0
\end{equation}
then there is an extension.
\end{thm}
%Assuming that  $\text{C}_{\loc}^{1,1}$ is, however, rather 
%restrictive.
%For example, the typical Hamiltonian $H(p) = |p|$ arising in front 
%propagation does not have this regularity.
%On the other hand, the only assumption made on $B$ is continuity.
%It turns out %(see Section~8) 
%that it is actually  possible to relax the regularity of $H$, if more is  
%assumed about $B$. 

%As far as the regularizations $W^\ep$'s are concerned, the canonical 
%example is
%$$W^\ep = W * \rho_\ep\ ,$$
%where $\rho_\ep$ is a standard mollifier in $\R$.
%
%In the cases of the Brownian motion, there are other more sophisticated 
%approximations --- see, for example, the book by Ikeda and Watanabe
%\cite{IK}. 
%
%A particular example is given by 
%$$\dot W^\ep_t (\omega) = \frac1{\ep} \dot B (\frac{t}{\ep^2},\omega)\ ,$$
%where the random field $B:\R\times\Omega\to \R$, which is 
%differentiable in $t$, 
%satisfies some appropriate mixing conditions explained below. 
%
%If $\A$ is a collection of subsets of $\Omega$, we denote by 
%$\sum (\A)$ the smallest $\sum$-algebra containing all sets in $\A$.
%For $s,t\in [0,\infty]$ with $s\leqq t$ let   
%\begin{gather*}
%\M_s^t = 
%\sum (\{ \dot B (u,\cdot )\}_{s\leqq u\leqq t} )\ \text{ and}\\
%\noalign{\vskip6pt}
%\beta (t) = \sup_{s\geqq 0} \{ |P(A\cap B) - P(A) P(B)| : 
%A \in \M_0^s,\ B\in \M_{t+s}^\infty \}\ .
%\end{gather*}
%Assume that 
%$$\int_0^\infty \beta^{1/4} (t)\,dt \langle\infty\ .$$
%
%It then follows that, as $\ep\to0$ and a.s. in $\omega$, $W^\ep\to W$ in 
%$C([0,\infty))$.
%
%%%%%%%%%%%%%%%%%%%%%%%%%%%%%%%%%%
%% old sec.9 -- now sec.8
%\end{document}
%\newpage

\section{Pathwise solutions for  equations with non-smooth Hamiltonians}

It is important to extend the class of Hamiltonians for which the solution operator of \eqref{paris200} with smooth paths has an extension.  The assumption that   $H\in \text{C}_{\loc}^{1,1}$ is rather 
restrictive. For example, the typical Hamiltonian $H(p) = |p|$ arising in front 
propagation does not have this regularity. 
\smallskip

The aim of this section is to provide a necessary and sufficient condition on $H$ to have an extension as well as to investigate if it is possible to assume less in $H$ by ``increasing'' the regularity of the paths, while still covering many cases of interest.
% without, of course, going all the way up to $C^1$. 
\smallskip

An important question and tool in this direction is to understand/control  the cancellations arising from the oscillations of the paths. And for this, it is useful to investigate if there are some formulae for the solutions in the presence of sign changing driving signals.

%On the other hand, the only assumption made on $B$ is continuity.
%It turns out %(see Section~8) 
%that it is actually  possible to relax the regularity of $H$, if more is  
%assumed about $B$. 
%
%It is mathematically interesting but also important, in view of the applications discussed earlier, to extend the theory 
%to nonsmooth Hamiltonians. Since in this case, it will not be possible to construct solutions of the system of characteristics, it is necessary to develop another method to construct solutions. This will be a density argument, which shows that the solution operator for smooth paths has a unique extension to the set of merely continuous paths. To prove this fact, we need to investigate the possibility of having some general formulae for solutions as well as the class of Hamiltonians for which it is possible to have the extension. An important step in this direction is to develop a good control of the cancellations  of the solutions arising from the oscillations of the paths. 

\subsection*{Formulae for solutions} The simplest possible formulae for the solutions of 
\begin{equation}\label{paris230}
u_t=H(Du) \ \text{in} \ Q_T \quad u(\cdot,0)=u_0, 
\end{equation}
are the well known Lax-Oleinik and Hopf formula which require convexity for $H$ and $u_0$ respectively.  In the Appendix the reader can find an extensive discussion about these formulae, their relationship and possible extensions.
\smallskip

When $H$ is convex,  the Lax-Oleinik formula is 
\begin{equation}\label{LO1} 
u(x,t) = \sup_{y\in \R^d} \left[u_0 (y) - tH^* (\frac{y-x}t)\right],  %\quad \text{(resp. } \ u(x,t) = \sup_{y\in \R^d} [u_0 (y) - t H^* (\frac{y-x}t)]).
\end{equation}
%and if is convex, 
%\begin{equation}\label{LO2} 
%u(x,t) = \sup_{y\in \R^d} [u_0 (y) - \frac1t H^* (\frac{x-y}t)].
%\end{equation}
where, given a convex function $w:\R^d\to \R$, $w^*(q) = \sup [\langle q,p \rangle -w (p)]$ is its Legendre transform. 
% (resp.  $H^*(p) = \inf [(p,q) +H (q)]$).  
%while in  \eqref{LO1} is the convex dual.
\smallskip

%Similarly if  $H$  convex, then $H(p)=\sup [(p,q) -H^*(q)]$ and   
%\begin{equation}\label{LO2} 
%u(x,t) = \sup_{y\in \R^d} [u_0 (y) - \frac1t H^* (\frac{y-x}t)]\ .
%\end{equation}
The Hopf formula, which  is the ``dual'' of the Lax-Oleinik one, says  that, if  $u_0$ is convex, then 
\begin{equation}\label{Hf1}%% \label{eq:25} 
u(x,t) = \sup_{p\in\R^d} \big[ \langle p,x \rangle +tH(p) - u_0^* (p)\big].
\end{equation}

In general, neither formula extends to the solutions of % to have either formula for the solution of 
\begin{equation} \label{ivp3}	%% \label{eq:Hopf}
%\begin{cases}
u_t = H(Du)\dot \xi\  \text{in} \    Q_T  \quad 
%\noalign{\vskip6pt}
u(\cdot, 0)= u_0, %\ \ \text{on} \ \ \R^d,
%\end{cases} 
\end{equation}
except in time  intervals where the path is either increasing or decreasing in which case it is possible to change time. 
%A unless $B$ is either increasing or decreasing, in which case it is possible to change the time. 
\smallskip

Indeed,  if $H\in C(\R^d)$, $\xi \in C^1$ and $u_0$ convex, the natural extension of \eqref{Hf1} should be 
%natural extension of the classical Hopf formula 
$$\sup_{p\in\R^d} \big[\langle p, x\rangle +\xi(t) H(p) - u_0^* (p)].$$
The formula above is a subsolution, as the ``sup'' of solutions $\langle p, x\rangle + B(t) H(p) - u_0^* (p)$,
but, in general, is not   a solution of  \eqref{ivp3}. The heuristic reason is that shocks are not reversible.
\smallskip

For example, if  $H(p) = |p|$ and $u_0 (x) = |x|$, then 
$$\sup_{p\in\R} (p x +\xi(t)|p| - |\cdot|^* (p)) = (|x| +\xi(t))_+.$$
On the other hand, the following is true.

\begin{prop}\label{prop6.1}
The unique viscosity solution of \eqref{ivp3} %$u_t = |Du|\circW$ 
with $\xi \in C^1$, $\xi(0)=0$, $H(p)= |p|$
 and $u_0 (x) = |x|$ is \begin{equation}\label{Hf2}	%% \label{eq:26}
u(x,t) =\max \left[( |x| + \xi(t))_+,  (\max_{0\leqq s\leqq t} \xi(s))_+\right].
\end{equation}
\end{prop}
Although the regularity of $\xi$ is used in the 
proof of \eqref{Hf2}, 
the actual formula extends by density to arbitrary continuous $\xi$'s.
\smallskip

%I refer to the Appendix for a more extensive discussion about extensions of the Lax-Oleinik and Hopf formulae.
% to Hamiltaonians

It is possible to give two different  proofs for \eqref{Hf2}. 
One is based on dividing $[0,T]$ into intervals where $\dot \xi $ is positive 
or negative and iterating  the Hopf formula. The second is a direct justification that  
\eqref{Hf2} is the viscosity solution to the problem. The details can be found in \cite{lionssouganidisbook}.
\smallskip

From the analysis point of view, the difficulty is related to the fact that when the signal changes sign,  the convexity properties of Hamiltonian also change. This leads to the possibility of using the formulae provided by the interpretation of the solution as the value function of a  two-player, zero-sum differential games, which I briefly recall next. %yield formulae without requiring any convexity/concavity of $H$, 
\smallskip

Assume that
\[ H(p)= \sup_{\alpha \in A} \inf_{\beta \in B} (\langle f(\alpha,\beta),  p\rangle +h(\alpha, \beta)), \]
where, for simplicity, the sets $A$ and $B$ are assumed  to be compact subsets of $\R^p$ and $\R^q$ and 
$f :  A\times B \to \R^d$ and $h:  A\times B \to \R$ are bounded;
note that any Lipschitz continuous Hamiltonian $H$ can be written  as a max/min of linear maps. 
\smallskip

It was shown in Evans and Souganidis \cite{evanssouganidis} that the unique viscosity solution of the initial value problem
\[ u_t=H(Du) \  \text{in} \ Q_T \quad u(\cdot,0)=u_0, \]
admits the representation
\[ u(x,t)= \underset{\alpha \in \Gamma(T-t)}\sup \underset{z \in N(T-t)} \inf \left\{ \int_{T-t}^T h(\alpha[z] (s),z(s))ds  \ + \ u (x(T))\right\}, \]
where $N(T-t)$ is the set of controls $z:[T-t,T] \to B$,  $\Gamma(T-t)$ is the set of nonanticipating strategies which map $B$-valued controls to $A$-valued ones, and $(x(s))_{s\in[T-t,T]}$ is the solution of the ode
% in $[T-t,T]$  the ordinary differential equation
\[ \dot x= f(\alpha[z] (s),z(s)) \quad  x(T-t)=x.\] 

%It is reasonable to expect that such a formula will extend to the problem with rough driving paths. There are, however, serious difficulties, which in, general, cannot be circumvented as it is explained below.
%%\smallskip
%
An attempt to extend this formula to \eqref{Hf1} meets immediately difficulties. Assume, for example, that $\xi\in C^1$. Then  
\[H(p)\dot \xi(t)= H(p)\dot \xi(t)_+ - H(p)  \dot \xi(t)_- ,\] 
and it easy to check that, in  general, it is not possible to find compact sets $C$ and $D$, vectors $f_\pm :C\times D \to \R^d$ and scalars $h_\pm: C\times D \to \R$ such that 

\[H(p)\dot \xi(t)=\sup_{c\in C} \inf_{d\in D} \langle((f_+(c,d)\dot \xi(t)_+ + f_-(c,d)\dot \xi(t)_-), p\rangle + (h_+(c,d)\dot \xi(t)_+ + h_-(c,d)\dot \xi(t)_-)).\]

Of course,  if the above holds, then the solution of \eqref{Hf1} is given  
by the formula 
\[
\begin{split} 
&u(x,t)=\\
&\underset{\alpha \in \Gamma(T-t)}\sup \underset{z \in N(T-t)} \inf \left [ \int_{T-t}^T\big( h_+( \alpha[z] (s),z(s)) \dot \xi(T-s)_+)) 
+ h_-(  \alpha[z] (s),z(s) \dot \xi(T-s)_-) \big) ds  \ + \ u(x(T)) \right],
\end{split} 
\]
where for the  control $z\in N(T-t)$ and strategy $a\in \Gamma(T-t)$, $x(s)_{s\in [T-t,T]}$ solves 
\[  \dot x(s)= f_+(\alpha[z] (s),z(s))\dot \xi(T-s)_+ + f_-(\alpha[z] (s),z(s))\dot \xi(T-s)_- \ \ \  x(T-t)=x.\]

%%%%%%%%%%%%%%%%%%%%%%%%%%%%%%%%%%
%%%%%%%%%%%%%%%%%%%%%%%%%%%%%%%%%%
\subsection*{Pathwise  solutions for   nonsmooth Hamiltonians.} When $H$ is less regular than in Theorem~\ref{thm:simplest}, 
%%%%%%%%%%%%%%%%%%%%%%%%%%%%%%%%%%
it is also possible to prove the unique extension property for the solution operator for smooth paths, but the argument is different and does not rely on inverting the characteristics. It is, however,  possible to use  the general  strategy summarized in Theorem~\ref{paris210.5} to identify the conditions on $H$ that will allow for the extension to exist for the initial value problem
\begin{equation}\label{paris260}
du=\sum_{i=1}^m H^i(Du)\cdot dB_i \ \text{in} \ Q_T \quad u(\cdot,0)=u_0.
\end{equation}
%after the proof of Theorem~\ref{thm:simplest}. % can be used here too. 
%\smallskip

The main result is stated next.

\begin{thm}\label{paris261}
The solution operator of \eqref{paris260} with smooth paths  has a unique extension to continuous paths if and only if $H^i$ is the difference of two convex functions for $i=1,\ldots, m$.  
\end{thm}

Identifying  the class of Hamiltonians  $H$ which can be as written  as the difference of two 
convex functions is a difficult  question. 
\smallskip

When $d=1$, a  necessary and sufficient condition for $H$ to be the difference of two convex functions 
is that $H' \in \text{BV}$. 
Indeed in this case, in the sense of distributions, 
$H'' = H''_1 -H''_2$
with $H''_1$ and $H''_2$ nonnegative distributions and, hence, locally 
bounded measures. Conversely, if $H' \in \text{BV}$, then $H'' = (H'')_+ - (H'')_-$.
\smallskip

When $d\geqq 2$, if $H = H_1 - H_2$ with $H_1,H_2$ convex, then, as above, 
$DH\in \text{BV}$. The converse is, however, false.
Functions with gradients in $\text{BV}$ may not have directional derivatives at 
every point, while differences of convex functions do.

\smallskip
Finally, if $H\in C^{1,1}$, then $H$ is clearly the difference of convex functions.
Indeed since, for some $c >0$,  
$D^2 H \geqq -2cI $,  then $H = H_1 - H_2$ with 
$H_1 (p) = H(p) +  c|p|^2$ and $H_2(p) =  c|p|^2$.
\smallskip

The proof of Theorem~\ref{paris261} is divided in several perts and requires a number of ingredients which are developed next.
\smallskip

\begin{prop}\label{DIFCONVEX}
Assume that the extension operator exists for all continuous paths. Then $H$ must be the difference of two convex functions.
\end{prop} 
\begin{proof}
In what follows, I assume for simplicity that $m=1$ and the problem is set in $Q_1.$
\smallskip

The necessity follows from the criterion summarized in Theorem~\ref{paris210.5}.  Since the extension must hold for any continuous path, it is possible to construct a sequence of paths satisfying the assumptions of Theorem~\ref{paris210.5} such that \eqref{paris211} implies that $H$ must be the difference of two convex functions. 
\smallskip

Consider a partition of $[0,1]$ of $2n$ intervals of length $1/2n$ and define the piecewise linear paths $B_n:[0,1]\to \R$ with slope $\dot B_n=\pm\mu$ and, for definiteness, assume that $\dot B_n=\mu$ in the first interval. 
\smallskip
It follows that 
\[\sup_{t\in [0,1]} |B_n(t)| \leqq \frac{\mu}{2n} \quad \text{and},  \text{ if } \  \mu/2n \to 0, \quad B_n\to 0 \ \text{in} \  C([0,\infty)).\]
%and, hence, 
%$$B_n\to 0\quad \text{ if }\quad \mu/2n \to 0\ .$$
%%\smallskip
Fix $\lambda >0$ and let $v_n^\lambda:\overline Q_1\to \R$ be the solution of 
\begin{equation}\label{paris260}
v_{n,t}^\lambda = H(Dv_{n}^\lambda) \dot B_n \ \text{in} \ Q_1 \quad v_{n}^\lambda(z,0)=\lambda |z|.
\end{equation}

Assume that, for some $\delta >0$, $\mu =  2n\delta$.  The claim   follows if it is shown that the $v_{n}^\lambda$'s  blow up, as $n\to\infty$,  if 
$H$ is not the difference  of two convex functions in a ball of radius $\lambda$. 
%A The conclusion then follows using a diagonalization argument as 
%$\delta\to0$. 
\smallskip

Recall % that %To investigate the behavior of the $v_n$'s, as $n\to\infty$,  observe
that, in each time interval of length $1/2n$, the equation in \eqref{paris260}  are   
$$\text{either }\quad v_{n,t} = 2n\delta H(Dv_n)\quad\text{or}\quad 
v_{n,t} = - 2n\delta H(Dv_n),$$
or, after rescaling, 
$$\text{either }\quad U_{n,t} = \delta H(DU_n)\ \text{ in }\ \overline Q_1\quad 
\text{or}\quad  V_{n,t}= - \delta H(DV^n)\ \text{ in }\ \overline Q_1;$$
%\smallskip
here, for notational simplicity, I omit the explicit dependence on $\lambda$.
\smallskip
 
The $V_n$'s  are constructed by a repeated iteration of Hopf's formula. This procedure yields sequences 
$(V_{2k+1}^*)_{k=0}^\infty$ and 
$(V_{2k}^*)_{k=0}^\infty$ which, as $k\to \infty$, either  
blow up   or converge, uniformly in $\overline B_\lambda$,
to $\overline V_1^*$ and $\overline V_2^*$ respectively. 
\smallskip

In the latter case, it follows that 
$$\overline V_2^* = (\overline V_1^* - \delta H)^{**}\quad\text{and}\quad 
\overline V_1^* = (\overline V_2^* + \delta H)^{**},$$
and, therefore, 
%$$\bar V_1^* = \bar V_2^* + \delta H\ ,$$
%and
$$\delta H = \overline V_1^* - \overline V_2^*,$$
which yields that $H$ is the difference of two convex functions.
\smallskip

If the sequences $(V_{2k+1}^*)_{k=0}^\infty$ and $(V_{2k}^*)_{k=0}^*$ blow up, then 
a diagonal argument, in the limit $\delta\to0$, shows that 
\eqref{paris211} cannot hold.
\smallskip

Indeed, since, for each $\delta >0$ and  as $k\to\infty$, 
$V_{2k+1}^* \to -\infty\ \text{ and }\  
V_{2k}^* \to -\infty\ \text{ in }\ \overline B_\lambda,$
choosing  $\delta = 1/m$  along a sequence $k_m \to \infty$  yields %we must have 
$V_{2 k_m}^* \leqq -1.$
\smallskip

Going back to the original scaled problem, it follows that 
$v_{k_m} \leqq -1,$
while $B_{k_m} \to0$ in $C([0,\infty))$ and $v_{k_m} (0,0)=0$.
%\end{proof}

\end{proof}
The next step is to show that,  if $H$ is the difference of two convex functions, then there exists a unique extension of the solution operator with $B$ smooth to the class of merely continuous $B$.
\smallskip

The main difficulties are  the  lack of differentiability of $B$ and  how to  control the oscillation of the solutions with respect to time. This was actually already exploited  in the proof of Proposition~\ref{DIFCONVEX}. For the sufficiency,  it is important to obtain a more explicit estimate. 
\smallskip

Controlling the cancellations due to the oscillations of the paths is very much related to the irreversibility of the equations due to the formation of shocks. ``Some memory'', however,  remains resulting in cancellations taking place as it can be seen in the next result. 
\smallskip

%The other direction is: 
%
%\begin{prop}\label{DIFCONEX2}
%Assume \eqref{difconvex}. Then \eqref{eq:36} holds for any sequence $B_n \in C([0,\infty))$ 
%such that $B_n(0)=0$ and, as $n\to 0$, $B_n \to 0$ in 
%$C([0,\infty))$. %Then  $H$ must be the difference of two convex functions.
%\end{prop}
%
%%To conclude the discussion about the necessary and sufficient condition, it needs to be shown 
%%that, if \eqref{difconvex} holds, then 
%%there exists a unique extension of the solution 
%%operators for $B$ smooth to the class of merely continuous $B$. 
%%\smallskip
%%
%The proof of Proposition~\ref{DIFCONEX2} is more complicated, since it is not possible to use the characteristics 
%as it was done in the proof of Theorem~\ref{thm:simplest}. Instead, it is necessary to find a way to control the oscillations of the solutions in time to deal with the  
%difficulty that $\dot B$ may not be defined. 
% it is 
%also important to control the oscillation of the solutions with respect 
%to time. 
%This is what we was done in the argument outlined above.
%The issue is very much connected to the fact that, due to 
%\smallskip
%This  is very much related to the fact that, due to the formation of shocks, the equations are not reversible. 
%In other words, solving the problem backwards does not give the same function. 
%On the other hand, ``some memory'' remains, resulting in cancellations taking place as it can be seen in the next result, 
%which is about initial value problems of the form 

Consider the initial value problems 
\begin{equation} \label{mshJ}	%% \label{eq:67bis}
%\begin{cases} 
u_t =\sum_{i=1}^m H^i(Du)\dot B^i  \ \text{in}
\   Q_T\quad u(\cdot,0)= u_0,\\
%\noalign{\vskip6pt}
\Eq
and
\begin{equation} \label{mshJ1}
v_t^i = H^i (Dv^i) \dot B^i \  \text{in}  \ Q_T \quad 
%\noalign{\vskip6pt}
v^i(\cdot, 0) = v_0^i,
\end{equation}
%and
%\begin{equation}\label{mshJ1}	%% \label{eq:67bis2}
%%\begin{cases}
%v_t^i = H^i (Dv^i) \dot B^i \ \text{in} \ Q_T \quad 
%%\noalign{\vskip6pt}
%v^i(\cdot, 0) = v_0^i \ \text{on} \  \R^d\,
%%\end{cases}
%\end{equation}
where, for each $i=1,\ldots,m$, 
\begin{equation}\label{reg}	%% \label{eq:68bis}
H^i \in C (\R^d)\ ,\quad B^i \in C^1 ([0,\infty))
\quad \text{and}\quad 
u_0,v_0^i \in BUC (\R^d).
\end{equation}

It is known that  both initial value problems in \eqref{mshJ}  and \eqref{mshJ1} have  unique viscosity 
solutions. 
In the  statement below,  $S_{H^i}(t) v$ is the 
solution of \eqref{mshJ1} with $\dot B^i\equiv 1$ at time $t > 0$.

\begin{thm}\label{prop:IVP}
Assume, in addition to \eqref{reg}, that, for each $i=1,\ldots,m$,
$H^i$ is convex and $DH^i(p_i)$ exists for some  $p_i\in\R^d$,  % be such that $DH^i(p_i)$ exists. 
and let $u\in BUC (\overline Q_T)$ be the viscosity solution 
of \eqref{mshJ}. 
Then, for all $(x,t) \in \overline Q_T$, 
\begin{equation}\label{canc1}
\begin{split}
&\prod_{i=1}^m S_{H^i} (-\min_{0\leqq s\leqq t} B^{i}(s)) u_0\left(x 
+ \sum_{i=1}^m DH^i (p_i)  (\min_{0\leqq s\leqq t}( B^{i}(s) - B^i(t))\right) \\
\noalign{\vskip6pt}
&\qquad \qquad 
+ \sum_{i=1}^m H^i (p_i) (\min_{0\leqq s\leqq t} (B^{i}(s) - B^i(t))) 
\leqq u (x,t)\leqq \\
\noalign{\vskip6pt}
&\qquad \prod_{i=1}^m S_{H^i} 
(\max_{0\leqq s\leqq t} B^{i}(s)) u_0 \left(x
+ \sum_{i=1}^m DH^i (p_i) (\min_{0\leqq s\leqq t} (B^{i}(s) -B^i(t)))\right)\\
\noalign{\vskip6pt} 
&\qquad \qquad 
- \sum_{i=1}^m H^i (p_i) (\max_{0\leqq s\leqq t} (B^{i}((s) -B^i(t))).
\end{split}
\end{equation}
%\begin{equation}\label{canc1}
%\begin{split}
%&\prod_{i=1}^m S_{H^i} (\max_{0\leqq s\leqq t} B^{i,-}(s)) u_0(x 
%- \sum_{i=1}^m DH^i (p_i)  (\max_{0\leqq s\leqq t} B^{i,-}(s) - B^i(t))) \\
%\noalign{\vskip6pt}
%&\qquad \qquad 
%- \sum_{i=1}^m H^i (p_i) (\max_{0\leqq s\leqq t} B^{i,-}(s) - B^i(t)) 
%\leqq u (x,t) \\
%\noalign{\vskip6pt}
%&\qquad \leqq \prod_{i=1}^m S_{H^i} 
%(\max_{0\leqq s\leqq t} B^{i,+}(s)) u_0 (x
%- \sum_{i=1}^m DH^i (p_i) (\max_{0\leqq s\leqq t} B^{i,+}(s) -B^i(t)))\\
%\noalign{\vskip6pt} 
%&\qquad \qquad 
%- \sum_{i=1}^m H^i (p_i) (\max_{0\leqq s\leqq t} B^{i,+}((s) -B^i(t))\ .
%\end{split}
%\end{equation}
\end{thm}

The proof of \eqref{canc1}, which is complicated,  is  based on repeated use of the 
Lax-Oleinik and Hopf formulae. The details can be found in \cite{lionssouganidisbook}. 
\smallskip

The following remark is  useful for what follows. 
\smallskip

If, in addition,  
%the $H^i$'s are nonnegative with minimum $0$
$$\min H^i =0\quad\text{ for }\quad i=1,\ldots,m\ ,$$
then \eqref{canc1} can be simplified considerably to read
\begin{equation}\label{canc1.1}
\prod_{i=1}^m S_{H^i} (\max_{0\leqq s\leqq t}  B^{i,-}_s) u_0 (x) 
\leqq u (x,t) \leqq \prod_{i=1}^m S_{H^i} (\max_{0\leqq s\leqq t} 
B^{i,+}_s) u_0(x)\ .
\end{equation}

The bounds in  \eqref{canc1} are  sharp. Indeed, recall that in the particular case 
%
%Before we present the proof of Proposition~\ref{prop:IVP}, which is long 
%and involved, we discuss few remarks and state and prove the 
%extension result.
%
%%The proof of \eqref{prop:IVP} 
%If in addition 
%%the $H^i$'s are nonnegative with minimum $0$
%$$\min H^i =0\quad\text{ for }\quad i=1,\ldots,M\ ,$$
%then \eqref{canc1} can be simplified considerably to read
%\begin{equation}\label{canc1.1}
%\prod_{i=1}^M S_{H^i} (\max_{0\leqq s\leqq t}  W^{i,-}_s) u_0 (x) 
%\leqq u (x,t) \leqq \prod_{i=1}^M S_{H^i} (\max_{0\leqq s\leqq t} 
%W^{i,+}_s) u_0(x)\ .
%\end{equation}
$$H(p) = |p|\ ,\quad u_0(x) = |x|\quad\text{and}\quad m=1\ ,$$
it was already  claimed  that the solution of \eqref{mshJ} is given by 
$$u(x,t) =\max \left [(|x| + B(t))_+ , \max_{0\leqq s\leqq t} B(s))_+\right ].$$
Evaluating the formula at $x=0$ yields  that  the upper bound in 
Proposition~\ref{prop:IVP} is sharp, since, in this case, 
$$S_H (\max_{0\leqq s\leqq t} B(s)) u_0 (0) 
= \max_{0\leqq s\leqq t}  B(s) \quad \text{and}  \quad  \max( B_+(t) , \max_{0\leqq s\leqq t} B(s)) = \max_{0\leqq s\leqq t} 
B(s).$$
%and 
%$$\max( B^+(t) , \max_{0\leqq s\leqq t} B^+(s)) = \max_{0\leqq s\leqq t} 
%B^+(s).$$
%
%Instead of the lower bound in \eqref{canc1} and \eqref{canc1.1}, 
%it is also possible to use the fact 
%$$u(x,t) \geqq \sup_{p\in\R^d}
% ((p,x) - v_0^* (p) + \sum_{i=1}^m H^i(p) W^i_t)\ ,$$
%which follows from the observation that, for time dependent problems, 
%the Hopf formula yields only a sub-solution. 
%
Using Theorem~\ref{prop:IVP} it is now possible to prove the sufficient part of Theorem~\ref{DIFCONVEX}. 

\begin{prop}\label{takis51}
Assume that, for each other $i=1,
\ldots,m,$  $H^i\in C(\R^d)$ is the difference of two convex functions. % satisfies  \eqref{difconvex}. 
Then the solution operator of  \eqref{paris260} on the class of smooth paths has a unique extension to the space of continuous paths. 
\end{prop} 

% consequence is the following theorem.
%Proposition~\ref{DIFCONEX2}. 

%\begin{proof}[The proof of Proposition~\ref{DIFCONEX2}]
 %Then the solution operator $S$ of \eqref{mshJ} has a unique extension.
%\end{thm}
%\begin{proof}
\begin{proof} 
The proof is based again on Theorem~\ref{paris210.5}. Fix $\lambda >0$, let $B_n=(B^n_1,\ldots B^n_m) \in C_0([0,\oo),\R^m)\cap C^1((0,\oo);\R^m)$ be sequence of signals such that, as $n\to \oo$,  $B^n \to 0$ in $C([0,\oo),\R^m)$, and consider the solution $\phi^n \in \text{BUC}(\overline Q_\oo)$ of 
%$\phi^{\lambda,\ep,\eta}$  solves 
\begin{equation}\label{paris265}
\phi^{n}_t = \sum_{i=1}^m H^i (D\phi^{n}) 
\dot B^n_i \  \text{in}  \ Q_\oo \quad   \phi^{n} (z,0) = 
\phi_0(z)  
= \lambda |z|.\end{equation}
It shown here that the assumption on $H$ yields,  for each $\lambda >0$ and $T>0$,
$$  %% \phi_{\ep,\eta}^\lambda (x) - \lambda |x| \to0\ , \qquad 
\lim_{n} \max_{(z,t) \in\overline Q_T} 
(\phi^{n} (z,t) -\lambda |z|)=0.$$
%\,|\to 0\ , $$
%where $(\xi^{i,\ep})_{\ep> 0}$ and $( \zeta^{i,\eta})_{\eta> 0}$ are 
%smooth approximations in $C([0,\infty))$  of the given 
%${B\in C([0,\infty))}$.
For each $i=1, \ldots, m$, $H^i=H^i_1 - H^i_2$ with $H^{i}_{1} , H^{i}_{2}$ convex. To simplify the presentation, it is assumed that 
each  $H^i_1$ and $H^i_2$ has minimum $0$ which is attained at $p=0$. Then, it is possible to use  \eqref{canc1}. %follows from 
%%with 
%%$H^i_1(0)=H^i_2(0)=\min H^i_1=\min H^i=0$ and $DH^i_1(0)=DH^i_2(0)$,  
% \eqref{canc1}, that   %\eqref{canc1.1} holds. 
%\begin{equation}\label{canc1.1}
%\prod_{i=1}^m S_{H^i} (-\min_{0\leqq s\leqq t}  B^i(s)) u_0 (x) 
%\leqq u (x,t) \leqq \prod_{i=1}^m S_{H^i} (\max_{0\leqq s\leqq t} 
%B^i(s)) u_0(x).
%\end{equation}
%
%In view of the previous discussion assume next that, for all $i=1,\ldots,m$, 
%both $H^i^1$ and $H^i^2$ are nonnegative. 
%We leave it up to the reader to carry out the calculations in the 
%general case. 
\smallskip

Rewriting  \eqref{paris265} as %the equation satisfied by $\phi^{\lambda,\ep,\eta}$ as 
$$\phi^{n}_t
= \sum_{i=1}^m H^{i,1} (D\phi^{n}) 
B^n_i
+\sum_{i=1}^m H^{i,2} (D\phi^{n}) (-\dot B^n_i),
$$
and using  \eqref{canc1.1} yields,   
for all $x\in \R^d$, 
\begin{equation*}
\begin{split}
&\prod_{i=1}^m S_{H^{i,1}} (-\min_{0\leqq s\leqq t} B^n_i(s)) 
\prod_{i=1}^m S_{H^{i,2}} (\max_{0\leqq s\leqq t}  B^n_i(s))
\phi_0(x) \\
\noalign{\vskip6pt}
&\qquad  \leqq 
\phi^{\lambda,\ep,\eta} (x,t) 
\leqq\prod_{i=1}^M S_{H^{i,1}}(\max_{0\leqq s\leqq t} B^n_i(s))
\prod_{i=1}^m S_{H^{i,2}}  (-\min_{0\leqq s\leqq t}  B^n_i(s))\phi_0(x),
\end{split}
\end{equation*}
and the claim now follows since, $\lim_{\ep,\eta\to0}(\max_{i=1,\dot,m}\max_{s\in [0,T]}|B^n_i|)=0.$

\end{proof}
%Proposition~\ref{DIFCONEX} and  Proposition~\ref{DIFCONEX2} are then combined in the  following theorem.
%\begin{thm}\label{takis51}
%Assume that $H\in C(\R^d)$ satisfies  \eqref{difconvex}. Then the solution operator for \eqref{ivp3} on the class of smooth paths has a unique extension to the space of continuous paths. 
%\end{thm} 
%
%It is immediate that the left and right hand sides of the above inequality 
%converge, and  as $\ep,\eta\to0$, to $\phi_0$  for $T> 0$ and $\lambda$ fixed,  
%uniformly in $\R^d\times [0,T]$. % as $\ep,\eta\to0$, to $\phi_0$.
%\end{proof}

%Theorem 6.1. Assume that, for each i = 1, . . . , m, Hi ? C(RN ) is the difference of two convex functions Hi1,Hi2 ? C(RN) and Bi ? C([0,°)). Then the solution operator S of (6.10) has a unique extension.

Another consequence of the ``cancellation'' estimates of 
Theorem~\ref{prop:IVP} is an explicit error estimate between two solutions with different signals.
% in terms of the 
%difference of the signals and initial data. 
\smallskip

In what follows, for $k=1,2$, $u^k\in \text{BUC}(\overline Q_\oo)$,  $B^k \in C_0([0,\oo);\R^m)$ and $u^k_0 \in C^{0,1}(\R^d)$ is the solution of the initial value problem
%
% 
%where $u^\ep$ and $v^\eta$ are the solutions of 
%%\begin{equation}\label{cep1}
%%\begin{split}
%%&\begin{cases}
%%u_t^\ep = \sum_{i=1}^M  H^i (Du^\ep)\dot W^{i,\ep}\ 
%%\text{ in }\ \R^d\times (0,\infty)\ ,\\
%%\noalign{\vskip6pt}
%%u^\ep = u_0\ \text{ on }\ \R^d\times \{0\}\ ,
%%\end{cases}\\
%%\text{ and }& \\
%%&\begin{cases}
%%u_t^\eta =\sum_{i=1}^M  H^i(Du^\eta) \dot W^{i,\eta}\ 
%%\text{ in }\ \R^d\times (0,\infty)\ , \\ 
%%\noalign{\vskip6pt}
%%u^\eta = u_0\ \text{ on }\ \R^d\times\{0\}\ .
%%\end{cases}\end{split}
%%\end{equation}
%%%%% tried on one line  %%%%%%%%%%
\begin{equation}\label{cep1}
u_t^k= \sum_{i=1}^m H^i (Du^k)\dot B^k_i \ 
\text{ in } \ Q_\oo \quad u^k(\cdot,0) = u_0^k. % \ \text{ on } \ \R^d,\\[1mm]
%u^\ep = u_0\ \text{ on }\ \R^d\times \{0\}\ ,
%\end{cases}
%\text{ and }\\[1mm]
%v_t^\eta =\sum_{i=1}^m H^i(Dv^\eta) \dot {\zeta}^{i,\eta} 
%\text{ in } \ Q_T \quad 
%%\noalign{\vskip6pt}
%v^\eta(\cdot,0) = u_0 \ \text{ on } \ \R^d.
%\end{cases}
\end{equation}
In \eqref{cep1}, the solution is either  a classical viscosity solution if the signal is smooth, or  the function obtained by the extension operator if $B^k$ is continuous. 
%%%%%%%%%%%%%%

\begin{thm}\label{thm:diff}
Assume that, for each $1,\ldots,m$, $H^{i}\in C(\R^d)$ is the difference of 
two nonnegative convex  functions $H^{i,1}, H^{i,2}.$  For $k=1,2$, $B^k \in C_0([0,\oo);\R^m)$ and $u^k_0 \in C^{0,1}(\R^d)$.  
%$(\xi^{i,\ep})_{\ep> 0}$ and $(\zeta^{i,\eta})_{\eta> 0}$ are two approximations 
%in $C([0,+\infty))$ of 
%$B^i \in C([0,\infty))$ with $\xi^{i,\ep}(0)=\zeta^{i,\eta}(0)=0$, and $u_0 \in C^{0,1} (\R^d)$. 
Let $u^k \in \text{BUC}(\overline Q_\oo)$ be the solution of \eqref{cep1}.  
There exists $C>0$ depending on $\|u_0^k\|$ and $\|Du_0^k\|$ 
and the growth of $H^i$'s  such that, for all $t >0$,  
  \begin{equation*}
   \sup_{x\in\R^d} |u^1 (x,t) - v^1 (x,t)| 
    \leqq C \max_{i=1,\ldots,m}  \max_{0\leqq s\leqq t} 
    |B^1_{i}(s) - B^2_{i}(s)| + \sup_{x\in\R^d} |u^1_0 (x)- u^2_0 (x)|.
\end{equation*}
%Moreover, if, for $u_{1,0}, u_{2,0} \in C^{0,1} (\R^d)$ 
%and $B^{i,1}, B^{i,2} \in C([0,+\infty))$, 
%$u_1,u_2 \in BUC(\overline Q_T)$ 
%are the extensions of the ``non-rough'' solution operators,  then 
%%the solutions of \eqref{mshJ} for $W_i^1, W_i^2 \in C([0,\infty))$ 
%%obtained by Theorem~\ref{takis51}, 
%% there exists $M> 0$  such that, for all $t> 0$,  
%\begin{equation*}
%\sup_{x\in \R^d} |u_1 (x,t) - u_2 (x,t)| \leqq 
%C \max_{i=1,\ldots,m} \max_{0\leqq s\leqq t} |B_i^{1}(s) - B_i^{2}(s)| + \sup_{x\in \R^d} |u_{1,0} (x) - u_{2,0} (x)| \ .
%\end{equation*}
\end{thm}

\begin{proof} 
Only  the estimate for $\sup (u^1 - u^2)$ is shown here. The one for $\sup (u_2-u_1)$ follows similarly.
Moreover,  the claim is proven under the additional assumption that the signals are smooth. The general case follows by density. 
%direction follows similarly. 
\smallskip

Let $L=\max_{k=1,2}\|Du_0^k\|$.  
Since the Hamiltonians are $x$-independent, it is immediate from the 
contraction property that, for all $t\geqq 0$, 
$ u^1 (\cdot,t) , u^2(\cdot,t) \in C^{0,1}(\R^d)$ 
and 
$\max_{k=1,2} \|Du^k (\cdot,t)\| \leqq L.$
%Without loss of generality then it  may be assumed that 
%$$|H^i| \leqq \max_i \max_{|p|\leqq L} |H(p)|\ .$$
The standard comparison estimate for viscosity solutions %argument in the proof of Theorem~\ref{thm:sharp} 
implies that, for all $(x,t)\in \overline Q_T$,  
\begin{equation*}
u^1 (x,t) - u^2 (x,t) - \phi^{L} (x,x,t) 
\leqq \sup_{x,y\in \R^d} \left[u_0^1 (x) - u_0^2(y) - L |x-y|\right ] 
\leqq 0,%\frac{L^2}{4\lambda} \ .
\end{equation*}
where $\phi^L$ is the solution of the usual doubled equation with $\phi^L(x,y,0)=L|x-y|.$ 
%Standard estimates from the theory of viscosity solutions yield 
%that, if $w\in BUC (\R^d\times (0,\infty))$ solves, for $\zeta$ smooth,  
%$$w_t = H(Dw) \dot\zeta\ \text{ in }\ \R^d\times (0,\infty)\ ,$$
%and $w(\cdot,0)$ is Lipschitz with constant $\lambda$, then 
%$$|w(z,t) - w(z,0)| \leqq 
%\Big(\sup_{p\in B(0,\lambda)} |H(p)|\Big) |\zeta (t)|\ .$$
\smallskip

Basic estimates from 
the theory of viscosity solutions yields that, for any $\tau >0$ and $w\in C^{0,1}(\R^d)$ with $\|Dw\|\leq L$,
$$\max_i \|S_{H^i} (\tau) w-w\| 
\leqq \left(\max_i \max_{|p|\leqq L} |H^i (p)|\right)\tau\ .$$

It follows that 
\begin{equation*}
\begin{split}
\phi^{L} (x,x,t) 
& \leqq L |x-x| + m \max_{1\leqq i\leqq m} 
[\max_{|p|\leqq L} |H^i (p)| 
\max_{0\leqq s\leqq t} |B^1_i(s)-B^2_i(s)]\\
\noalign{\vskip6pt}
& = m \max_{1\leqq i\leqq m} 
[\max_{|p| \leqq L} |H^i (p)| 
\max_{0\leq s\leq t} |\xi^{i,\ep}(s) - \zeta^{i,\eta}(s)|].
\end{split}
\end{equation*}

Combining the upper bounds for $u^1 (x,t) -u^2 (x,t)$ and 
$\phi^{L} (x,x,t)$ %and optimizing over $\lambda$ 
gives the claim.

\end{proof}

 \subsection*{Control of cancellations for spatially dependent  Hamiltonians} It is both interesting and important for the study of qualitative properties of the pathwise solutions, see, for example, section 7 of the notes, to extend the results about the cancellations  
to  spatially dependent Hamiltonians $H=H(p,x)$ and the initial value problem 
 \beq\label{gal10}
 du=H(Du,x)\cdot d\xi \  \text{in} \  Q_\oo.
\end{equation}

The basic cancellation estimate reduces to whether if, for any $u\in \text{BUC}(\R^d)$ and any $a >0$,
% 
%
%Another important open question is whether it is possible to control the cancellations for spatially dependent Hamiltonians. 
%An important step towards such a conclusion is to show that, if 
% $H$ is convex and nonnegative, then 
%%%Then, for all $a> 0$ and $u\in BUC(\R^d)$, 
\beq\label{gal11}
S_H (a) S_{-H}(a)u \leq u \leq S_{-H} (a) S_H (a) u, 
\end{equation}
where $S_{\pm H}$ is the solution operator of \eqref{gal10} with Hamiltonians $\pm H$.

A consequence of a counterexample of Gassiat \cite{g} presented in the next section is that such a result cannot be expected for nonconvex Hamiltonians, since it would imply a domain of dependence property which is shown in \cite{g} not to hold for a very simple nonconvex problem. 
\smallskip

A first step towards an affirmative result was shown some time ago by Lions and the author. This was extended lately by Gassiat, Gess, Lions and Souganidis \cite{ggls} who established the following.
% and It is shown next that if, however, 
\begin{thm}\label{gal12}
Fix $\xi \in C_0([0,\oo);\R)$ and assume that 
\beq\label{gal13}
H=H(p,x):\R^d \times \R^d \to \R \ \ \text{is convex and Lipschitz continuous in $p$ uniformly in $x$.}.
\end{equation}
Then \eqref{gal11} holds. 
\end{thm}

\begin{proof}  The result is shown for  $\xi\in C^1([0,\oo))$. The general conclusion follows by density. Moreover, since the arguments are identical, I only work with the inequality on the left.
\smallskip

For notational simplicity, I assume that $\|D_pH\|=1$. If $L$ is the Legendre transform of $H$, it follows that 
% It then follows that 
\[H(p,x)= \underset{B_1(0)} \sup \{\langle p, v\rangle -L(v,x)\}.\]
%where $L(v,x)= \sup_{\R^d}\{\langle v, p\rangle-H(p,x)\}.$
\smallskip

Let $\mathcal A=L^{\infty}(\R_{+}; \overline B_1(0))$. The  control representation of the solution $u$ of \eqref{gal10} (see, for example, Lions~\cite{lbook}) with $\xi_{t}\equiv t$ and $u_0 \in \text{BUC}(\R^d)$ gives  %{[}my calculations gave the roles of $H$, $-H$ reversed? but kept expressions from Lions-Souganidis{]}
\[
S_{H}(t)u_{0}(x)=\sup_{q\in\mathcal A}\left\{ u_{0}(X(t))-\int_{0}^{t}L(q(s), X(s))ds:\,X(0)=x,\,\dot{X}(s)=q(s) \ \text{for $s\in [0,t]$} \right\} ,
\]
and 
\[
S_{-H}(t)u_{0}(y)=\inf_{r\in\mathcal A}\left\{ u_{0}(Y(t))+\int_{0}^{t}L(r(s),Y(s))ds:\,Y(0)=y,\,\dot{Y}(s)=-r(s) \ \text{for $s\in [0,t]$}\right\}.
\]
%where $\mathcal A=L^{\infty}(\R_{+}; \overline B_1(0))$ is the set of controls. 
%\smallskip

It follows that 
\begin{align*}
S_{H}(t)\circ S_{H}(-t)u_{0}(x)=\sup_{q\in\mathcal A}\inf_{r\in\mathcal A}\Big\{ & u_{0}(Y(t))+\int_{0}^{t}L(r(s),Y(s))ds-\int_{0}^{t}L(q(s),X(s))ds:\\
 & Y(0)=X(t),\,\dot{Y}(s)=-r(s),\,X(0)=x,\,\dot{X}(s)=q(s)\ \text{for $s\in [0,t]$}\Big\}.
\end{align*}
Given $q\in\mathcal A$  choose $r(s)=q(t-s)$ in the infimum above. Since  $Y(s)=X(t-s)$, it follows that 
\begin{align*}
S_{H}(t)\circ S_{H}(-t)u_{0}(x) & \le\sup_{q\in\mathcal A}\Big\{ u_{0}(X(0))+\int_{0}^{t}L(q(t-s),X(t-s))ds-\int_{0}^{t}L(q(s),X(s))ds:\\
 & \quad\quad\quad X(0)=x,\,\dot{X}(s)=q(s)  \ \text{for $s\in [0,t]$} \Big\}\\
 & =u_{0}(x).
\end{align*}

\end{proof} 

As a matter of fact, Lions and Souganidis came up recently with  a more refined form of \eqref{gal11}, which is stated below without proof. % was shown very recently by L
\begin{thm}
Fix $\xi \in C_0([0,\oo);\R)$ and assume  \eqref{gal13}. For every $a, b, c>0$ such that $b\leq \min(a,c)$,
\[
S_H(c)S_H(b)S_H(a)=S_H(a+c-b).\]
\end{thm}

\subsection*{The interplay between the regularity of the Hamiltonians and the paths}
%A summary follows of what is  known so far 
%about the interplay between the regularity 
%of $H$ and the paths in order to have a unique extension.
%\smallskip
%
The classical theory of viscosity solutions  applies when 
$H\in C$ and $B\in C^1$;  
actually  it is possible to consider $B\in C^{1,1}$ or even 
discontinuous $B$ as long as $\dot B\in L^1.$
%(see Perthame and Lions \cite{PL}, Ishii \cite{I1}, Barles \cite{B1}, etc.).
It was also shown here that, when $H\in C^{1,1}$ or, more generally, if $H$ is the difference of two 
convex (or half-convex) functions, 
there exists a unique extension for any $B\in C ([0,\infty))$. 
%Recall that $H$ is called half-convex, if for some $C> 0$, 
%$D^2H \geqq -cI$, a fact which is equivalent to saying that $p\to H(p) 
%+ (c/2)|p|^2$ is convex. 
%Arguments similar to the one's we will develop next allow us to show that 
%there exists a unique extension to $H\in C^{2(1-\alpha)+\ep}$ and 
%$W\in C^{0,\alpha}$  for $\alpha \in (0,1)$ and $\ep> 0$. 
%Whether it is possible to avoid the additional 
%$\ep$-regularity is an open problem
\smallskip

Arguments similar to the ones presented  next yield a unique extension for 
$B\in C^{0,\alpha} ([0,\infty))$ with $\alpha \in(0,1)$ and 
$H\in C^{2(1-\alpha)+\ep} (\R^d)$ for $\ep>0$; recall that, 
for any $\beta \in [0,\infty)$,   $C^\beta (\R^d)$ is the space 
$C^{[\beta],\beta-[\beta]} (\R^d)$. It is not clear, however,  if the additional $\ep$-regularity is necessary. % is an 
\smallskip

The conclusion  resembles nonlinear interpolation.  Indeed, 
consider the solution mapping $T(B,H) = u$, which is 
a bounded map from $C^1 \times C^0$ into $C$ and 
$C^0\times C^2$ into $C$. 
Typically, if $T$ is bilinear, abstract interpolation results 
would imply that $T$ must be a bounded map from $C^{0,\alpha} \times 
C^{[2(-1-\alpha)], 2(1-\alpha)- [2(1-\alpha)]}$ into $C$. 
But $T$ is far from being bilinear. 
%Therefore we cannot 

\smallskip

Next it is stated without proof (see \cite{lionssouganidisbook} for the details) that, in the particular case  $\alpha = 1/2$, it is 
possible to have a unique extension if $H\in C^{1,\delta}$ for $\delta > 0$. 
Of course, the goal is to show that is enough to have $H\in C^1$ or 
even $H\in C^{0,1}$. 
Questions related to the issues described above are studied  in an ongoing work  by Lions, Seeger and Souganidis \cite{lss}. 
%apply such abstract results to the problem we are considering here.
%\smallskip
%
%Next it is shown that, in the particular case  $\alpha = 1/2$, it is 
%possible to have a unique extension if $H\in C^{1,\delta}$ for $\delta> 0$. 
%Of course, the goal is to show that is enough to have $H\in C^1$ or 
%even $H\in C^{0,1}$. This is question studied in \cite{lss}. 
\smallskip

%To state the result we fix $W\in C^0$ and 
A sequence $(B_n)_{n\in\N}$ in $ C^1 ([0,\infty))$ is said to 
approximate $B\in C([0,\infty))$ in $C^{0,\,1/2}$ if, as $n\to\infty$,  
$$B_n \to  B\ \text{ in $C([0,\infty))$ 
\quad and\quad} 
\sup_n \|\dot B_n \|_\oo \, \| B_n-B\|_\oo < \infty\ .$$
Given $B\in C^{0,1/2}([0,\infty))$, it is possible to 
find  at least two classes of such approximations. 
The first uses convolution with a suitable smooth 
kernel, while the second relies on finite differences. 
\smallskip

Let $\rho_n (t) = n\rho (nt)$ with 
$\rho$ a smooth nonnegative kernel with compact support in $[-1,1]$
such that 
$\text{$\int z\rho (z) dz=0$  and $\int \rho (z)dz=1$,}$
and consider the smooth function $B_n = B*\rho_n$.   
If 
$C = (\|\rho'\| +\|\rho\| +1) [B]_{0,\, 1/2}$, then 
$$\| \dot B_n \| \leqq C\sqrt{n}\quad\text{and}\quad \|B_n -B\|
\leqq C/\sqrt{n}\ .$$
%
%Indeed 
%$$|W_t - W^n_t | = | \int (W_t -W_{t-n^{-1} s} ) \rho (s)ds| 
%\leqq C/\sqrt{n}$$
%and 
%\begin{equation*}
%\begin{split}
%|\dot W^n_t |  = |\int W_{t-s}  n^2 \rho' (ns)ds| 
%& = |\int (W_{t-s} - W_t ) n^2\rho' (ns)ds|\\
%\noalign{\vskip6pt}
%& = |\int (W_{t-s/n}  - W_t ) n\rho' (s)ds| 
%\leqq C\sqrt{n}\ .
%\end{split}
%\end{equation*}
%\smallskip
For the second approximation, subdivide $[0,T]$ into intervals of length
$\Delta = T/n$ and construct $B_n$ by a linear interpolation of 
$(B_{k\Delta})_{k=1,\ldots,n}$. 
Then 
$$|\dot B_n|  = \frac{|B_{(k+1)\Delta} - B_{k\Delta}|}{\Delta} 
\leqq \frac{[B]_{0,\,1/2}}{\sqrt{\Delta}} = C\sqrt{n} \quad \text{and} \quad \|B-B_n\| \leqq  [B]_{0,\,1/2} \sqrt{\Delta} = \frac{C}{\sqrt{n}}. $$
%and 
%$$\|B-B_n\| \leqq  [B]_{0,\,1/2} \sqrt{\Delta} = \frac{C}{\sqrt{n}}\ .$$
The next result says that 
$C^{0,\frac12}$-approximations of $C^{0,\frac12}$ paths yield 
a unique extension for $H\in C^{1,\delta} (\R^d)$ with $\delta >0$. 
%there exists a unique extension provided we use such kind of regularizations. 
As a matter of fact the result not only gives an extension but also an 
estimate. For the proof I refer to \cite{lionssouganidisbook}. 

%We have:

\begin{thm}\label{thm:extension}
Assume that $B\in C^{0,1/2}([0,\infty))$  and  
$H\in C^{1,\delta}(\R^d)$  for some $\delta >0$, and fix $T > 0$ and $u_0 \in BUC(\R^d)$. 
For any $(\xi_n)_{n\in \N}, (\zeta^m)_{m\in\N} \in C_0([0,\infty))$  and $u_{0,n}, v_{0,m} \in BUC(\R^d)$, which are  respectively $C^{0,1/2}$-approximations 
 of $B$  and $u_0$ in $BUC(\R^d),$ let 
 $u_n, v_m \in BUC(\overline Q_T)$ be the solutions of the corresponding initial value problems.
Then there exists  $u\in BUC (\overline Q_T)$ such that, as $n, m\to\infty$, $u_n, v_m \to u$ in $BUC(\overline Q_T).$ Moreover,  if  $\|u_{0,n} - u_0\| \leqq Cn^{-\beta}$ for some $\alpha, \beta >0$, then there exist
$\gamma > 0$ and $C>$ such that 
$|u_n-u|\leqq Cn^{-\gamma}$ in $\overline Q_T$.
%
%If, as $n,m \to\infty$, $u_{0,n} - \tilde u_{0,m} \to 0$ uniformly in $\R^d$,  
%%%  \in BUC (\R^d)$, 
%there exists  $u\in BUC (\overline Q_T)$ such that, as $n\to\infty$
%and $m\to\infty$,
%$u_n,\tilde u_m\to u$ uniformly in  $\overline Q_T$.
%If, in addition, $|u_{0,n} - u_0| \leqq Cn^{-\beta}$ for some $C, \beta> 0$, then there exists
%$\gamma > 0$, such that 
%$|u_n-u|\leqq Cn^{-\gamma}$ in $\overline Q_T$.
%% for some $\gamma > 0$.
\end{thm}

A discussion follows about  the need to have conditions on $H$. The key step in the proof of Theorems~\ref{thm:extension} % and \ref{thm:9.3}
can be reformulated as follows.
Let $(B_n)_{n\in\N}$ be a sequence of $C^1$-functions such  that, 
as $n\to\infty$, 
%\begin{equation}\label{adreg}
$B_n \to0\quad\text{and}\quad \sup_n \|B_n\|_\oo\, \|\dot B_n\|_\oo < \infty,$ 
%\end{equation}
and consider the solution $v_n$ of 
%$$%\begin{cases}
$v_{n,t} = H(Dv_n) \dot B_n \ \text{in}\ Q_T$ and 
%\noalign{\vskip6pt}
$v_n(x,0) = \lambda |x|.$
As before, it suffices to show that, for each fixed $T >0$ and for all $(x,t)\in \overline Q_T$,% as $n\to\infty$,     
$$\lim_{n\to \oo}\sup_{(x,t)\in \overline Q_T} [v_n (x,t) - \lambda |x|]\to 0\ .$$
Next let $\dot B$ be piecewise constant such that,  
for $t_i = \frac{T}{k} i$,  
%$$\circW = \begin{cases} 
%\Delta_1 &\text{ on $[0,t_1]$}\\
%-\Delta_2&\text{ on $[t_1,t_2]$}\end{cases}\quad 
%\text{where}\quad t_i = \frac{T}k i\ .$$
$$\dot B= \Delta_1\ \text{ in }\ [t_{2k},t_{2k+1}] \quad \text{ and }\quad 
\dot B = - \Delta_2 \ \text{ in }\ [t_{2k+1},t_{2k}],
$$
and, for simplicity, take $\lambda =1$.
Arguments similar to the ones earlier in this section and the fact that $v_k$ is 
convex, since $v_k(\cdot,0)$ is, yield a sequence $w_k = v_k^*$ 
such that 
$$w_0 = 0 {\mathds 1}_{\{|p| \leqq 1\}} + \oo {\mathds 1}_{\{|p| > 1\}}
%\begin{cases}0&\text{ if $|p| \leqq 1$}\ ,\\
%+\infty &\text{ if $|p| \geqq 1$\ ,}\end{cases}
$$
and
$$w_{2k+1}= (w_{2k} + \Delta_{2k} H)^{**} 
\ \text{and}\  
w_{2k}= (w_{2k-1} - \Delta_{2k-1} H)^{**}\ \text{where} \ \Delta_i = k \big[B(\frac{(i+1)T}k) - B(\frac{iT}k)\big].$$
%where 
%$$\Delta_i = k \bigg[B(\frac{(i+1)T}k) - B(\frac{iT}k)\bigg]\ .$$
The convexity of the $v_k$'s and Hopf's formula implies that that the sequence $w_k$ is decreasing. 
Then convergence will follow if there is a lower bound for the $w_k$'s.. %To be able to prove convergence, it is necessary to obtain a lower bound.
\smallskip

Consider next the particular case 
$H(p) = |p|^\theta$
and assume that, for all $i$,  
$\Delta^i = \sqrt{k}$. %\quad\text{for all}\quad i\ .$
\smallskip

If  $\tilde w_k$ is constructed similarly to $w_k$ but with 
$\Delta_i\equiv 1$, 
it is immediate that 
$w_k = k^{-1/2} \tilde w_k,$
%If  $\tilde w_k$ is constructed similarly to $w_k$ but with 
%$\DeltA^i\equiv 1$.
and,  since  
$\tilde w_{k+1} = ((\tilde w_k \pm |p|^\theta )^{**} \mp |p|^\theta)^{**},$
it follows that $\tilde w_{k+1} \leqq \tilde w_k$ and 
$\tilde w_k = +\infty$ if $|p| >1$. 
\smallskip

Let $m_k = - \inf_{|p| <1} w_k(p)$.  Since $H$ is not the difference of two 
convex functions if  $\theta \in(0,1/2)$, it must be that 
%In view of Theorem we must have that 
$m_k\to\infty\quad\text{as}\quad k\to\infty. $
%since, for $\theta \in(0,1/2)$, $H$ is not the difference of two 
%convex functions.
\smallskip

It turns out, and this is tedious computation, 
that there exists $c >0$ such that 
$\tilde w_k \leqq -ck^{1-\theta}\ .$
\smallskip

It follows that,  if $\theta <1/2$,
$$w_k = k^{-1/2} \tilde w_k \leqq -ck^{1/2 -\theta} \to-\infty\quad\text{as}
\quad k\to\infty\ .$$

The above calculations show that, if $H\in C^{0,\alpha} (\R^d)$ 
with 
$\alpha \in (0,\frac12)$ and $\sup_n \|B_n\|_{C^{0,\frac12}}< \infty$, then 
there is blow up, and, hence, not a good solution.
%\smallskip
On the other hand,  if $H\in C^{0,\frac12}(\R^d)$,  there is no 
blow up. 
\smallskip

\section{qualitative properties}

Recently there has been great interest  in the study and understanding of various qualitative properties of the solutions. In this section, I focus manly on the initial problem
\begin{equation}\label{london1}
du=H(Du,x)\cdot dB \ \text{in} \ Q_\oo \quad u(\cdot,0)=u_0,
\end{equation}
and I discuss the following three qualitative behaviors: domain of dependence and finite speed of propagation ,  intermittent regularizing effect and regularity,  and long time behavior of the pathwise solutions. 
\smallskip

\subsection*{Domain of dependence and finite speed of propagation} Given that the pathwise solutions are obtained as uniform limits of solutions of hyperbolic equations with  domain of dependence and finite speed of propagation property, it is natural to ask if this property remains true in the limit. 
 
\smallskip

In the context of the ``non-rough''  viscosity solutions, it is  known 
that, if $H$ is  Lipschitz continuous with constant $L$, and 
$u^1, u^2 \in BUC(\overline Q_T)$ solve the initial value problems 
$$u^1_t=H(Du^1) \ \text{in} \ Q_T \quad u^1(\cdot,0)=u^1_0 \quad \text{and} \quad u^2_t=H(Du^2) \ \text{in} \ Q_T \quad u^2(\cdot,0)=u^2_0,$$
%and 
%$$u^2_t=H(Du^2) \ \text{in} \ Q_T \quad u^2(x,0)=u^2_0(x),$$
then
%$$w_t=H(Dw) \ \text{in} \ Q_T$$
$$\text{if} \ u_0^1 = u_0^2 \  \text{ in } \   B(0,R), \quad \text{then} \  u^1(\cdot,t)=u^2(\cdot,t) \  \text{in}  \  B(0,R-Lt).$$
%at $(x_0,t_0)$ depends only on $u_0$ in the ball 
%$B(x_0,L(t_0-t))$ for all $t\in[0,t_0]$. 
%

The first positive but partial result in this direction for pathwise solutions was proved \cite{lionssouganidisbook}. The claim is  the following. %that, if $u$ solves \eqref{london1} and 
%which yields something weaker than finite speed of propagation, is about 
%the initial value problem 
%\begin{equation}\label{ShJ4}   %% \label{eq:70}
%u_t = (H_1 (Du) - H_2 (Du))\dot B
%\  \text{in}\  Q_T, 
%\end{equation}
%with 
%\begin{equation}\label{reg111}	%% \label{eq:71}
%H_1,H_2\ \text{ convex and bounded from below.}
%\end{equation}
\begin{prop}\label{cor:Lmax}
Assume that $H=H_1-H_2$ with $H_1$ and $H_2$ convex and bounded from below, and $u_0\in C^{0,1}(\R^d)$.  Let $L$ be the Lipschitz constant of $H_1$ and $H_2$ in $B(0, \|Du\|)$ and consider the solution  $u\in BUC (\overline Q_T)$ of \eqref{london1}. If, for some $A\in \R$ and $R > 0$, 
$$u(\cdot,0)\equiv A \ \text{ in} \  B(0,R),$$
%and assume that \eqref{reg1}   holds. 
%Let $L> 0$ be the Lipschitz constant of $H_1$ and $H_2$ in $B(0,2R_0)$.
then 
$$u(\cdot,t) \equiv A \ \text{in } \  B(0,R- L( \max_{0\leqq s\leq t} B(s) -\min_{0\leqq s\leq t}B(s)).$$ 
\end{prop}

\begin{proof}
Without loss of generality, the problem may be  reduced to Hamiltonians with the additional property
\begin{equation}\label{reg2}	%% \label{eq:71bis}
H_1,H_2 \ \text{ nonnegative and }\ H_1(0) = H_2(0) =0.
\end{equation}
As long as $R > L( \max_{0\leqq s\leq t} B(s) -\min_{0\leqq s\leq t}B(s))$, and, since 
$H_1(0) = H_2(0)=0$, 
%Proposition~\ref{prop:constantL} gives
the finite speed of propagation of the initial value problem with $B(t)=t$ yields 
\begin{equation*}
S_{H_1} (\max_{0\leqq s\leqq t} B^\pm(s)) u_0 
= S_{H_2} (\max_{0\leqq s\leqq t} B^\pm(s)) u_0  = A,
\end{equation*}
and  the claim then follows using the estimate in Theorem~\ref{prop:IVP}.
\end{proof} 

The following example in  \cite{g}  shows  that, when the Hamiltonian is neither convex nor concave, the initial value problem does not have the finite speed of propagation property.

\smallskip

Fix $T> 0$ and $\xi \in C_0([0,\oo);\R)$. The total variation $V_{0,T} (\xi)$ of $\xi$  in $[0,T]$ is 
\[
V_{0,T} (\xi):= \underset{ (t_0,\ldots,t_n)\in \mcl P} \sup \sum _{i=0}^{n-1} |\xi(t_{i+1}) - \xi(t_i)|,
\]
where $\mcl P= \{0 = t_0 < t_1 < \cdots, t_n = T\}$ is  a partition of $[0,T]$. 
% that is,
%and $\abs{ \mcl P}$ its mesh size, that is,
%\[ \mcl P := \{0 = t_0 < t_1 < \cdots, t_n = T\}.\]
%Given a path $\xi \in C([0,T])$, its total variation $V_{0,T} (\xi)$ over $[0,T]$ is denoted by 
%\[
%V_{0,T} (\xi):= \underset{ (t_0,\ldots,t_n)\in \mcl P} \sup \sum _{i=0}^{n-1} |\xi(t_{i+1}) - \xi(t_i)|.
%\]
\smallskip

The result is stated next. 

\begin{prop}\label{gassiat0}
Given $\xi \in C_0([0,T];\R)$, let $u\in \text{BUC}(\R^2 \times [0,T])$ be the solution of
\beq\label{gassiat1}
du= (|u_x|-|u_y|) \cdot d\xi \ \text{in} \   \R^2 \times [0,T] \quad u(x,y,0)=|x-y| +\Theta(x,y),
\end{equation}

with  $\Theta \in \text{BUC}(\R^2)$  nonnegative and such that, for some $R> 0$,  $\Theta(x,y) \geq R$ if $\min(x,y) \geq R.$
Then
\beq\label{gassiat2}
%u(0,0,T)\geq \left( \underset{ (t_0,\ldots,t_n)\in \mcl P} \sup \frac{\sum _{i=0}^{n-1} |\xi(t_{i+1}) - \xi(t_i)|- \R}{n} \right )_+ \wedge 1.
u(0,0,T)\geq \left( \underset{ (t_0,\ldots,t_n)\in \mcl P} \sup \dfrac{\sum _{i=0}^{n-1} |\xi(t_{i+1}) - \xi(t_i)| -R}{n} \right)_+ \wedge 1. %\frac{V_{0,T} (\xi) -R}{n} \right )_+ \wedge 1.
\end{equation}
In particular, $u(0,0,T)> 0$ as soon as $V_{0,T} (\xi) > R. $
\end{prop}

If  $\xi$ is a Brownian motion, then $V_{0,T} (\xi)=+\oo$ for all $T>0$.  Then \eqref{gassiat1} implies there is no finite speed of propagation property for any $R>0$. 
\smallskip

\begin{proof}[The proof of Proposition~\ref{gassiat0}]
The argument is based on the differential games representation formula discussed earlier in the notes, which is possible to have for the very special Hamiltonian considered here.%in Theorem~\ref{gassiat0}. 
\smallskip

Arguing by density, I assume that $\xi\in C^1$. A simple calculation shows that, for all $p,q \in \R$,
\[(|p|-|q|)\dot\xi(t)=\underset{|a|\le 1}\max \underset{|b|\le 1}\min \left\{(a \dot\xi(t)_+ +  b \dot\xi(t)_- ) p+ (b \dot \xi(t)_+  + a \dot\xi(t)_- )q\right\}.\]
It follows that, for any $T>0$,
\[u(0,0,T)=\sup_{\alpha \in \Gamma (T)} \inf_{z\in N(T)} J(\alpha[z],z),\]
where, for each pair $(w,z)$ of controls in $[0,T]$, 
\[J(w,z)= |x^{w,z}(T)-y^{w,z}(T)| + \Theta (x^{w,z}(T),y^{w,z}(T)), \]
and 
\[
\begin{split} 
& \dot x^{w,z} (s)= w(s) \dot\xi(T-s)_+ +  z(s) \dot\xi(T-s)_- \quad x^{w,z}(0)=0,\\[1.2mm]
& \dot y^{w,z} (s)= z(s) \dot \xi(T-s)_+  + w(s) \dot\xi(T-s)_- \quad  y^{w,z}(0)=0.
 \end{split} \]

I refer to \cite{gassiat} for the rest of the argument, which is based on the choice, for each partition of $[0,T]$, of a suitable pair of strategy and control, and the assumption on $\Theta$.

\end{proof}

Motivated by the general question and  the partial result and counterexample discussed above, \cite{ggls} considered the case of convex, spatially dependent Hamiltonians. Using the cancellation property discussed in the previous section, it is  proven in \cite{ggls}  that, in this setting,  there is a finite speed of propagation.  This required the use of what is known as ``skeleton'' of the path. The details are presented next.
\smallskip

%
%The notion of domain of dependence is fist reformulated as speed of propagation.
%\smallskip
%\end{document}

%Given $H :\R^d \times \R^d \to \R$ and $\xi\in C([0,T])$, consider the equation 
%\begin{equation}\label{gal1}
%du=H(Du,x)\cdot d\xi \  \text{in}  \ \R^d\times [0,\oo), % \quad u(\cdot,0)=u_0.
%\end{equation}
%%\end{document}
%where throughout the discussion in this subsection it is assumed that
%\begin{equation}\label{gal0}
%\begin{split}
%&H:\R^d\times\R^d \to \R \ \text{ is convex and }\\
%&\text{Lipschitz continuous with constant $L$  in the first argument.}
%\end{split}
%\end{equation}
%%\end{document}
%
%The speed of propagation of $H$ at time $T$ is defined by
%\begin{align}
%\rho_{H}(\xi,T):=\sup\Big\{ R\ge0:\  \text{there exist solutions} & \ \ u^{1},u^{2}\mbox{ of }\eqref{gal1} \text{ and } x \in \R^d, \\  \text{ such that  }  u^{1}(\cdot, 0)=u^{2}(\cdot, 0) &\mbox{ in }B_R(x)\ %\label{eq:defR-2}
%\text{and }u^{1}(x, T)\ne u^{2}(x, T)\Big\}. \nonumber 
%\end{align}
%where $B_R(x)$ is the ball in $\R^d$ centered at $x$ with radius $R$.
%\smallskip

Given $\xi \in C([0,T])$, if $\arg\min_{[a,b]}$ (resp.  $\arg\max_{[a,b]}$) denotes the set of minima (resp. maxima) points of $\xi$ on the interval $[a,b]\subseteq[0,T]$,   the sequence $(\tau_{i})_{i\in\mathbb{Z}}$ of successive extrema of $\xi$  is defined by 
\begin{equation}
\tau_{0}=\sup\left\{ t\in[0,T]: \xi(t)=\max_{0\leq s\leq T}\xi(s)\mbox{ or }\xi(t)=\min_{0\leq s\leq T}\xi(s)\right\} ,\label{eq:extrema1}
\end{equation}
where, for all  $ i\geq0$, %inductively set 
\begin{equation}
%\forall i\geq0,
\tau_{i+1}=\left\{ \begin{array}{ll}
\sup\arg\max_{[\tau_{i},T]}\xi & \mbox{ if }\ \ \xi(\tau_{i})<0,\\[1.5mm]
\sup\arg\min_{[\tau_{i},T]}\xi & \mbox{ if } \ \ \xi(\tau_{i})> 0,
\end{array}\right.\label{eq:extrema2}
\end{equation}
and, for all  $ i\leq0$,   
\begin{equation}
%\forall i\leq0,
\tau_{i-1}=\left\{ \begin{array}{ll}
\inf\arg\max_{[0,\tau_{i}]}\xi & \mbox{ if } \ \  \xi(\tau_{i}) < 0,\\[1.5mm]
\inf\arg\min_{[0,\tau_{i}]}\xi & \mbox{ if }\ \  \xi(\tau_{i})> 0.
\end{array}\right.\label{eq:extrema3}
\end{equation}
\smallskip 
The skeleton (resp. full skeleton) or  reduced (resp. fully reduced) path ${R}_{0,T}(\xi)$ (resp. $\tilde{R}_{0,T}(\xi)$) of $\xi \in C_0([0,T])$ are defined as follows. %(see Figure \ref{fig:reduced_path}). 
\begin{defn}
Let $\xi \in C([0,T])$. % be a continuous path with $\xi(0)=0$. 

%\begin{enumerate}
(i)~The reduced path $R_{0,T}(\xi)$ is a piecewise linear function which agrees  with $\xi$ on $(\tau_{i})_{ i\in\mathbb{Z}}$.

% for $t\in\{0,T,\tau_{i}:\ i\in\mathbb{Z}\}$ and is piecewise linear in between. 
(ii)~The fully reduced path  $\tilde{R}_{0,T}(\xi)$ is a piecewise linear function  agreeing with % the fully reduced path and  which coincides with
 $\xi$ on $(\tau_{-i})_{i\in\mathbb{N}} \cup \{T\}$.
 
  %$for $t\in\{0,T,\tau_{-i}:\ i\in\mathbb{N}\}$ and is piecewise linear in between. 
(iii)~ A path $\xi \in C_0([0,T])$ is reduced (resp. fully reduced) if $\xi=R_{0,T}(\xi)$ (resp. $\xi=\tilde{R}_{0,T}(\xi)$). %See Figure  for a pictorial representation. 
%\end{enumerate}
\end{defn}
%\begin{figure}
%\input{reduced_path.pdf_tex}
%\caption{The (fully) reduced path}
%\label{fig:reduced_path}
%\end{figure}
Note that the reduced and the fully reduced paths coincide prior to the global extremum $\tau_0$. While the reduced path captures the max-min fluctuations also after $\tau_0$, the fully reduced path is affine linear on $[\tau_0,T]$ and, in this sense, is more ``reduced''.
\smallskip

%Next I consider \eqref{london1}. It is assumed that 
%Given $H :\R^d \times \R^d \to \R$ and $\xi\in C([0,T])$, consider the equation 
%\begin{equation}\label{gal1}
%du=H(Du,x)\cdot d\xi \  \text{in}  \ \R^d\times [0,\oo), % \quad u(\cdot,0)=u_0.
%\end{equation}
%%\end{document}
%where throughout the discussion in this subsection it is assumed that

Throughout the discussion, it is assumed that 
\begin{equation}\label{gal0}
\begin{cases}
H:\R^d\times\R^d \to \R \ \text{ is convex and }\\[1.2mm]
\text{Lipschitz continuous with constant $L$  in the first argument.}
\end{cases}
\end{equation}
%\end{document}

The speed of propagation of \eqref{london1} at  time $T$ is defined by
\begin{align}
\rho_{H}(\xi,T):=\sup\Big\{ R\ge0:\  \text{there exist solutions} & \ \ u^{1},u^{2}\mbox{ of }\eqref{london1} \text{ and } x \in \R^d, \\  \text{ such that  }  u^{1}(\cdot, 0)=u^{2}(\cdot, 0) &\mbox{ in }B_R(x)\ %\label{eq:defR-2}
  \text{and }u^{1}(x, T)\ne u^{2}(x, T)\Big\}. \nonumber 
\end{align}
%where $B_R(x)$ is the ball in $\R^d$ centered at $x$ with radius $R$.
%\smallskip

To keep track of the dependence of the solution on the path, in what follows I use the notation $u^\xi$ for the solution of \eqref{london1} with path $\xi$. The main observation is that 
%Let $u^{\xi}$ denote a   solution of \eqref{gal1}. It follows from arguments introduced earlier that, for any $T > 0$,
%  The next claim expresses the solution at time of $T$ in terms  is to show that Theorem \ref{thm:reduction} in the next section  that %(cf.~Theorem \ref{thm:reduction} below)	
\begin{equation}
u^{\xi}(\cdot, T)=u^{R_{0,T}(\xi)}(\cdot, T),\label{eq:intro-reduced}
\end{equation}
which immediately implies the following result about the speed of propagation.
%in the statement $L$ is the Lipschitz constant of $H$ from \eqref{gal0}.
\begin{thm}\label{takis3}
Assume \eqref{gal0}. %that $H:\R^{d}\to\R$ is convex. 
Then, %\begin{enumerate}
for all $\xi \in C([0,T])$,  % we have 
\begin{equation} \label{eq:main-speed}
\rho_{H}(\xi,T)\leq L \ \|R_{0,T}(\xi)\|_{TV([0,T])}.
\end{equation}
%where $R_{0,T}(\xi)$ is the reduced path associated to $\xi$. 
%\item 
\end{thm}
The second main result of \cite{ggls} concerns the total variation of the reduced path of a Brownian motion. To state it, it is necessary to introduce the random variable $\theta: [0,\infty) \to [0,\infty)$ given by 
%recall that the range $\theta (t)$ of the Brownian motion is given by 
\begin{equation}\label{takis10005}
\theta(a)=\inf \{ t\geq0:\;\;\max_{[0,t]}B-\min_{[0,t]}B=a\},
\end{equation}
which is the first  time that the range, that is $\max-\min$ of a Brownian motion equals $a$.
% introduce the following notation for the range of the Brownian motion, that is the difference 
\smallskip

It is proved in \cite{ggls}, where I refer for the details,  that the length of the reduced path is a random variable with almost Gaussian tails. It is also shown that if the range, that is,  the maximum minus the minimum of $B$,  is fixed instead of the time horizon $T$,  then the length has Poissonian tails.

\begin{thm}\label{takis4}
%\begin{thm}
\label{thm:BM} Let $B$ be a Brownian motion and fix $T> 0$. Then, for each $\gamma \in (0,2)$, there exists $C=C(\gamma, T)> 0$ such that, for any $x\geq2$, 
\begin{equation}
\P\left(\left\Vert R_{0,T}(B)\right\Vert _{TV([0,T])}\geq x\right)\leq C \exp\left(-Cx^{\gamma}\right),\label{eq:RGaussian}
\end{equation}
%Let $\theta(1):=\inf\left\{ t\geq0,\;\;\max_{[0,t]}B-\min_{[0,t]}B=1\right\} $. Then, 
and
\begin{equation}\label{eq:RPoisson}
\displaystyle{\underset{x\to\infty} \lim \frac{\ln \P\left(\left\Vert R_{0,\theta(1)}(B)\right\Vert _{TV([0,\theta(1)])}\geq x\right)}{x\ln(x) }= -1.} %\sim_{x\to\infty}-x\ln(x).\label{eq:RPoisson}
\end{equation}
\end{thm}
\smallskip
A related result, proving that the expectation of the total variation of the so-called piecewise linear oscillating running max/min function of Brownian motion is finite, has been obtained independently by Hoel, Karlsen, Risebro, and Storr{\o}sten in \cite{HKRS}.
\smallskip

The following remark shows that  that upper bound in Theorem~\ref{takis3} is actually sharp. %  for simplicity I  only treat the case $H(p)=|p|$. 
% is the Euclidean norm.%(??Maybe try for more general ??). For a given , we call  the half-reduced path which is obtained similarly to  but only considering the  for . 
\begin{prop}\label{taki4}
%\label{thm:opt} 
Let $H(p)=|p|$ on $\R^{d}$  with $d\ge1$. 
Then, for all $T >0 $ and  $\xi \in C_0([0,T];\R)$, 
\begin{equation}
\rho_{H}(\xi, T)\geq\|\tilde{R}_{0,T}(\xi)\|_{TV([0,T])}.
\end{equation}
When $d=1$, then  %higher dim
\[
\rho_{H}(\xi, T)=\|\tilde{R}_{0,T}(\xi)\|_{TV([0,T])}.
\]
\end{prop}

Here  I only sketch  the proof of the first result. %Theorem~\ref{takis3}.

\begin{proof}[A sketch of the proof of Theorem~\ref{takis3}] The first step is \eqref{gal11}. 
\smallskip

The second is a monotonicity property for piecewise linear paths. Let $\xi_{t}={\mathds 1}_{t\in[0,t_{1}]}(a_{0}t)+\\{\mathds 1}_{t\in[t_{1},T]}(a_{1}(t-t_{1})+a_{0}t_{1})$ and, for $s <t$, set  $\xi_{s,t}=\xi_t -\xi_s$.
\smallskip

If  $a_{0}\ge0$ and $a_{1}\le0$ (resp. $a_{0}\le0$ and  $a_{1}\ge0$), then 
\beq\label{gal20}
S_{H}^{\xi}(0,T)\ge S_{H}(\xi_{0,T}) \qquad (resp. \ \ S_{H}^{\xi}(0,T)\le S_{H}(\xi_{0,T}).)
\end{equation}

Since the claim is immediate if $a_{0}=0$ or $a_{1}=0$, next it is assumed that 
% then there is nothing to show. Assume next that 
%Since both claims are proved similarly, we only consider the case 
$a_{0}>0$ and $a_{1} <0$ %(see Figure \ref{fig:reduction1}). 

If $\xi_{0,T}\le0$, then
\begin{align*}
S_{H}(a_{1}(T-t_{1})) & =S_{-H}(-a_{1}(T-t_{1}))=S_{-H}(-a_{1}(T-t_{1})-a_{0}t_{1})\circ S_{-H}(a_{0}t_{1})\\
 & =S_{-H}(-\xi_{0,T})\circ S_{-H}(a_{0}t_{1})=S_{H}(\xi_{0,T})\circ S_{H}(-a_{0}t_{1}),
\end{align*}
and, hence, in view of \eqref{gal13},
\begin{align*}
S_{H}^{\xi}(0,T) & =S_{H}(\xi_{0,T})\circ S_{H}(-a_{0}t_{1})\circ S_{H}(a_{0}t_{1}) 
  \ge S_{H}(\xi_{0,T}).
\end{align*}
If $\xi_{0,T}\ge0$ %(see Figure \ref{fig:reduction1}), 
then, again, \eqref{gal13} yields 
\begin{align*}
S_{H}^{\xi}(0,T) & =S_{H}(a_{1}(T-t_{1}))\circ S_{H}(-a_{1}(T-t_{1})+a_{0}t_{1}+a_{1}(T-t_{1}))\\
 & =S_{H}(a_{1}(T-t_{1}))\circ S_{H}(-a_{1}(T-t_{1}))\circ S_{H}(a_{0}t_{1}+a_{1}(T-t_{1}))
  \le S_{H}(\xi_{0,T}).
\end{align*}
For the second inequality,  note that $S_{-H}^{-\xi}(0,T)=S_{H}^{\xi}(0,T)$, $S_{-H}(-t)=S_{H}(t)$. It then follows from  the first part that 
\[
S_{H}^{\xi}(0,T)=S_{-H}^{-\xi}(0,T)\ge S_{-H}(-\xi_{0,T})=S_{H}(\xi_{0,T}).
\]
\smallskip

The next observation provides the first indication of the possible reduction encountered when using the max or min of a given path. For the statement, given a piecewise linear path $\xi$,  set 
\[ \tau_{max}=\sup \left\{t\in[0,T]:\ \xi_{t}=\max_{s\in[0,T]}\xi_{s}\right \}  \ \text{ and} \  \tau_{min}={\inf}\left\{t\in[0,T]:\ \xi_{t}=\min_{s\in[0,T]}\xi_{s}\right\}.\]

\begin{lem}
\label{lem:reduction3} Fix  a piecewise linear path $\xi$.
% and set $\tau_{max}:=\sup\{t\in[0,T]:\ \xi_{t}=\max_{s\in[0,T]}\xi_{s}\}$ and  $\tau_{min}:=\sup\{t\in[0,T]:\ \xi_{t}=\min_{s\in[0,T]}\xi_{s}\}$. 
 Then 
\[
S_{H}^{\xi}(\tau_{max},T)\circ S_{H}(\xi_{0,\tau_{max}})\leq S_{H}^{\xi}(0,T)\leq S_{H}(\xi_{\tau_{min},T})\circ S_{H}^{\xi}(0,\tau_{min}).
\]
\end{lem}

\begin{proof}
Since the proofs of both inequalities are similar, I only show the details for the first. %one on the left. %We prove the first inequality, the second inequality can be shown analogously. 
\smallskip

Without loss of generality, it is assumed that  $\text{sign}(\xi_{t_{i-1},t_{i}})=-
\text{sign}(\xi_{t_{i},t_{i+1}})$ for all $[t_{i-1},t_{i+1}]\subseteq[0,\tau_{max}]$. 
It follows that, if  $\xi_{|[0,\tau_{max}]}$ is linear, then $S_{H}^{\xi}(0,\tau_{max})=S_{H}(\xi_{0,\tau_{max}})$.
\smallskip

If not, since $\xi_{0,\tau_{max}}\ge0$, there is an index $j$ such that $\xi_{t_{j-1},t_{j+1}}\ge0$ and $\xi_{t_{j-1},t_{j}}\le0$. It then follows from \eqref{gal20} that 
\[
S_{H}^{\xi}(0,\tau_{max})\le S_{H}^{\tilde{\xi}}(0,\tau_{max}), 
\]
 where $\tilde{\xi}$ is piecewise linear and coincides with $\xi$ for all $t\in\{t_{i}:\ i\ne j\}$. 
 %and is piecewise linear otherwise. 
\smallskip
 
A simple iteration yields $S_{H}^{\xi}(0,\tau_{max})\le S_{H}(\xi_{0,\tau_{max}})$, and, since $S_{H}^{\xi}(0,T)=S_{H}^{\xi}(\tau_{max},T)\circ S_{H}^{\xi}(0,\tau_{max})$, this concludes the proof.

\end{proof}
The previous conclusions and  lemmata are combined to establish the following monotonicity result. 
\begin{cor} \label{cor:Monotone}
Let $\xi,\zeta$ be piecewise linear, $\xi(0)=\zeta(0)$, $\xi(T)=\zeta(T)$ and $\xi\leq\zeta$ on $[0,T]$. Then
\begin{equation} \label{eq:Monotone}
S_{H}^{\xi}(0,T)\leq S_{H}^{\zeta}(0,T).
\end{equation}
\end{cor}

\begin{proof}
Assume that $\xi$ and $\zeta$ are piecewise linear on each interval $[t_i,t_{i+1}]$ on a common partition  $0=t_0 \leq \ldots \leq t_N=T$ of $[0,T]$.
\smallskip

If  $N=2$, then, for all  $\gamma \geq 0$ and all $a, b \in \R$,
\begin{equation} \label{eq:Monotone2}
%\forall \gamma \geq 0, \forall a, b \in \R, \;\; 
S_H(a+\gamma) \circ S_H(b-\gamma) \leq S_H(a) \circ S_H(b).
\end{equation}
If $a\geq 0$, this follows from the fact that, in view of  \eqref{gal20},  %Lemma \ref{lem:bump}, 
\[S_H(\gamma) \circ S_H(b-\gamma) \leq S_H(b).\] % by Lemma \ref{lem:bump}. 
If $a+\gamma \leq 0$, then again \eqref{lem:reduction3} yields  \[S_H(a) \circ S_H(b) = S_H(a+ \gamma) \circ S_H(-\gamma) \circ S_H(b) \geq S_H(a+ \gamma) \circ S_H(b-\gamma).\]
% again by Lemma \ref{lem:bump}. 
 Finally, if $a \leq 0 \leq a+\gamma$ we have \[S_H(a) \circ S_H(b) \geq S_H(a+b) \geq S_H(a+\gamma) \circ S_H(b-\gamma).\]
%\smallskip

The proof for $N> 2$ follows  by induction on $N$. Let $\rho$ be piecewise linear on the same partition and coincide with $\zeta$ on $t_0, t_1$, and with $\xi$ on $t_2, \ldots, t_N$. The induction hypothesis then yields
\[
S_{H}^{\xi}(0,t_2)\leq S_{H}^{\rho}(0,t_2) \ \  \text{and} \ \ S_{H}^{\rho}(t_1,T)\leq S_{H}^{\zeta}(t_1,T)
\]
from which we deduce
\[
S_{H}^{\xi}(0,T)\leq S_{H}^{\rho}(0,T) \leq S_{H}^{\zeta}(0,T).
\]
 \end{proof}
 
 To complete the study of the cancellations, it is necessary to use a density argument, which, itself,  requires a result 
 about the uniform continuity of the solutions with respect to the paths. Such a result was shown earlier in the notes  for spatially-independent Hamiltonians which are the difference of two convex functions and for spatially dependent under some additional conditions on the joint dependence but not convexity. The most general result available without additional assumptions other than convexity was obtained in \cite{ls8}. Here it is stated without a proof.
 
 \begin{thm}\label{gal30}
 Assume \eqref{gal13}. Then, for each $u_0$ $\in$ $BUC(\R^d)$ and $T \geq 0$, the family
\[ \left\{ S_H^\xi (0,T)u_0: \;\;\; \xi \mbox{ piecewise linear} \right\}
\]
has a uniform  modulus of continuity.
\end{thm} 
 
 An immediate consequence is the following extension result which is stated as a corollary without proof; see \cite{ggls} for the details.
 \begin{cor}
\label{cor:continuity}The map $\xi\mapsto S_{H}(\xi)$ is uniformly continuous in the sup-norm in the sense that, if $(\xi^n)_{n\in \N}$ is a sequence of  piecewise-linear functions on $[0,T]$ with $\lim_{n,m\to\infty}\|\xi^n-\xi^m\|_{\infty,[0,T]}= 0$, then, for all $u \in BUC(\R^d)\times (0,\infty)$, 
\begin{equation}\label{takis40}
\underset{n,m\to\infty}\lim %_{n,m\to\infty}
 \|S^{\xi^n}_H(0,T)u-S^{\xi^m}_H(0,T)u\|_{\infty} = 0.
\end{equation}
\end{cor}

Combining all the results above completes the proof.

\end{proof}

%\end{document}

\subsection*{Stochastic intermittent regularization} A very interesting question is whether there is  some kind of stochastic regularization-type property for the pathwise  solutions of 
\begin{equation}\label{gal40}
du=H(Du)\cdot d\zeta \  \text{in} \ Q_\oo.
\end{equation}

It is  assumed that 
\begin{equation}\label{hconvex2}
H\in C^2(\R^d) \ \ \text{is uniformly convex},
\end{equation}
which implies that there exist $\Theta \geq \theta >0$ such that, for all $p \in \R^d$ and in the sense of symmetric matrices,
\begin{equation}\label{hconvex1}
\theta I \leq D^2 H(p) \leq \Theta I.
%\Theta |\xi|^2 \geq \sum_{i,j=1}^d H_{p_i p_j} \xi_i\xi_j \geq \theta |\xi|^2.
\end{equation}
 The upper bound in \eqref{hconvex1} can be relaxed when dealing  with Lipschitz  solutions of \eqref{gal40}.  
\smallskip

Motivated by a recent 
observation of Gassiat and Gess \cite{gg} for the very special case that $H(p)=(1/2) |p|^2$, recently Lions and the author \cite{lsregularity} investigated this question. A summary of these results is presented next without proofs. The
details can be found in  \cite{lsregularity}. 

%The results for the stochastic equation follow from two new and sharp analogous results for the deterministic initial value problems
%\begin{equation}\label{takis41}
%u_t=\pm H(Du) \ \ \text{in} \ \  \R^d\times (0,\infty), 
%\end{equation}
%with
%\begin{equation}\label{hconvex2}
%H\in C^2(\R^d) \ \ \text{is uniformly convex}
%\end{equation}
%satisfying \eqref{gal44}.
%
\smallskip

The possible intermittent regularizing results  follow from iterating regularizing and propagation of regularity-type results for the ``non rough'' problem
\beq\label{gal65}
u_t=\pm H(Du) \ \text{in}  \ Q_\oo.
\Eq

It turns out that the quantity to measure the regularizing effects is the symmetric matrix
$$F(p)=\sqrt{D^2H(p)},$$
the reason being that,  if, for example, $u $ is a smooth solution of \eqref{gal65}, then a simple calculation yields that the matrix $W(x,t)=F(Du(x,t))$ satisfies the matrix valued problem
\[W_t=DH(Du)DW \pm |DW|^2.\]
%\end{document}

The first claim is about the regularizing effect of \eqref{gal65}.  In what follows  all the inequalities and solutions below should be understood in the viscosity sense.
\begin{thm}\label{gal60} Assume \eqref{hconvex2}. If $u\in \text{BUC}(\R^d \times (0,\infty))$
is a solution of $u_t=H(Du) $  
(resp. \\ $u_t=-H(Du)$) \ $\text{in \ $\R^d \times (0,\infty)$}$ and, for some $C \in (0,\infty]$,  
\begin{equation}\label{manos1}
-F(Du(\cdot, 0))D^2u(\cdot, 0) F(Du(\cdot, 0)) \leq CI \ \text{in} \ \R^d, 
\end{equation}
(resp.
 \begin{equation}\label{manos2}
-F(Du(\cdot, 0))D^2u(\cdot, 0) F(Du(\cdot, 0)) \geq - CI  \ \text{in} \ \R^d),
\end{equation}
%(resp. $-F(Du(\cdot, 0))D^2u(\cdot, 0) F(Du(\cdot, 0)) \geq - CI  \ \text{in} \ \R^d$,) 
then, for all $t > 0$, 
\begin{equation}\label{gal66}
 -F(Du(\cdot, t))D^2u(\cdot, t) F(Du(\cdot, t)) \leq \dfrac{C}{1+Ct} I  \ \text{in} \ \R^d,
 \end{equation}
 (resp.
\begin{equation}\label{gal67}
 -F(Du(\cdot, t))D^2u(\cdot, t) F(Du(\cdot, t)) \geq -\dfrac{C}{1+Ct} I  \ \text{in} \ \R^d).
 \end{equation}
 \end{thm}
 
% \end{document}
Estimates  \eqref{gal66} and \eqref{gal67}  are sharper versions of the classical regularizing effect-type results for viscosity solutions (see Lions \cite{lbook}, Lasry and Lions~\cite{LL}), which say that, if $u_t=H(Du)$  
(resp. $u_t=-H(Du)$) in $Q_\infty$, and, for some $C \in (0,\infty]$, $-D^2u(\cdot, 0) \leq C I$ (resp. $-D^2u(\cdot, 0) \geq -C I) \ \text{in} \ \R^d$, then, for all $t > 0$, 
\begin{equation}\label{gal61}
 -D^2 u(\cdot, t) \leq \dfrac{C}{1+\theta Ct} I  \ \text{in} \ \R^d
 \end{equation}
 (resp.
\begin{equation}\label{gal62}
 -D^2 u(\cdot, t)  \geq -\dfrac{C}{1+\theta Ct} I  \ \text{in} \ \R^d.)
 \end{equation}
 Note that, when $C=\infty$, that is,  no assumption is made on $u(\cdot,0)$, then \eqref{manos1} and \eqref{manos2} reduce to 
\beq\label{manos22.1}  -F(Du(\cdot, t))D^2u(\cdot, t) F(Du(\cdot, t)) \leq \dfrac{1}{t} \quad \left(\text{resp.} -F(Du(\cdot, t))D^2u(\cdot, t) F(Du(\cdot, t)) \geq -\dfrac{1}{t}\right),
\end{equation}
%(resp. 
%\beq\label{takis23.1}-F(Du(\cdot, t))D^2u(\cdot, t) F(Du(\cdot, t)) \geq -\dfrac{1}{t}), \Eq
which are sharper versions of \eqref{gal61} and \eqref{gal62}, in the sense that they do not depend on $\theta$, of the classical 
estimates 
$$-D^2 u(\cdot, t) \leq \dfrac{1}{\theta t} \quad (\text{resp.} \  D^2 u(\cdot, t) \geq -\dfrac{1}{\theta t}).$$
%and the sharpness of the new regularizing effect is clear since the estimate is independent of $\theta$. 
%\smallskip

To continue with  the propagation of regularity result,  I first recall that it was shown in \cite{LL} that, if $u$ solves $u_t=H(Du)$ (resp. $u_t=-H(Du)$)  in $\R^d \times [0,\infty)$, with $H$ satisfying \eqref{hconvex1}, then,   
\begin{equation}\label{takis2221}\text{ if \ $-D^2 u(\cdot, 0) \geq - C I$, \ then \  $-D^2u(\cdot,t) \geq - \dfrac{C}{(1-\Theta Ct)_+}$},\end{equation}
%\end{document}
(resp.
\beq\label{takis2222} \text{if  \ $-D^2u(\cdot, 0) \leq C I$, \ then \ $-D^2u(\cdot,t) \leq  \dfrac{C}{(1-\Theta Ct)_+}$.})\Eq

The new propagation of regularity result depends on the dimension. In what follows, it is said that  $H:\R^d\to \R$ is quadratic, if there exists a symmetric matrix $A$ which satisfies \eqref{hconvex1}  such that 
$$H(p)=(Ap,p).$$

\begin{thm}\label{takis26.1} Assume \eqref{hconvex2} and let   $u\in \text{BUC}(\overline Q_\oo)$ solve $u_t=H(Du) $   
(resp.   $u_t=-H(Du)$) \ $\text{in $Q_\oo$}$.  Suppose that  either  $d=1$ or $H$ is quadratic.  If, for some  $C> 0$, 
\begin{equation}\label{manos4}
 -F(Du(\cdot, 0))D^2u(\cdot, 0) F(Du(\cdot, 0)) \geq -CI \ \text{in} \ \R^d, 
 \end{equation}
 (resp. 
 \begin{equation}\label{manos5}
 -F(Du(\cdot, 0))D^2u(\cdot, 0) F(Du(\cdot, 0)) \leq CI  \ \text{in} \ \R^d), 
 \end{equation}
  then, for all $t > 0$, 
\begin{equation}\label{takis27}
 -F(Du(\cdot, t))D^2u(\cdot, t) F(Du(\cdot, t)) \geq - \dfrac{C}{(1-Ct)_+} I  \ \text{in} \ \R^d,
 \end{equation}
 (resp.
\begin{equation}\label{takis28}
 -F(Du(\cdot, t))D^2u(\cdot, t) F(Du(\cdot, t)) \leq \dfrac{C}{(1-Ct)_+} I  \ \text{in} \ \R^d.)
 \end{equation}
% If $d=1$, then \eqref{manos4} (resp. \eqref{manos5}) imply \eqref{takis27}  (resp. \eqref{takis28}),
% for any $u(\cdot,0) \in \text{BUC}(\R^d).$
\end{thm}

The result for $d\geq 2$ and general $H$ requirer more regularity for the initial condition.

%Since the dimension plays a role we state them as separate theorems, the first for $d\rangle1$ and the second for $d=1$.
\begin{thm}\label{takis26} Assume that  $d > 1$ and that $H$ satisfies \eqref{hconvex2} but is  not quadratic.
 Let   $u\in \text{BUC}(\R^d \times [0,\infty))$ solve $u_t=H(Du) $   
(resp.  $u_t=-H(Du)$) \ $\text{in $\R^d \times (0,\infty)$}$ and assume that $u(\cdot,0) \in C^{1,1}(\R^d)$. If, for some  $C> 0$, 
\begin{equation}\label{manos4}
 -F(Du(\cdot, 0))D^2u(\cdot, 0) F(Du(\cdot, 0)) \geq -CI \ \text{in} \ \R^d, 
 \end{equation}
 (resp. 
 \begin{equation}\label{manos5}
 -F(Du(\cdot, 0))D^2u(\cdot, 0) F(Du(\cdot, 0)) \leq CI  \ \text{in} \ \R^d), 
 \end{equation}
  then, for all $t > 0$, 
\begin{equation}\label{takis27}
 -F(Du(\cdot, t))D^2u(\cdot, t) F(Du(\cdot, t)) \geq - \dfrac{C}{(1-Ct)_+} I  \ \text{in} \ \R^d,
 \end{equation}
 (resp.
\begin{equation}\label{takis28}
 -F(Du(\cdot, t))D^2u(\cdot, t) F(Du(\cdot, t)) \leq \dfrac{C}{(1-Ct)_+} I  \ \text{in} \ \R^d.)
 \end{equation}
% If $d=1$, then \eqref{manos4} (resp. \eqref{manos5}) imply \eqref{takis27}  (resp. \eqref{takis28}),
% for any $u(\cdot,0) \in \text{BUC}(\R^d).$
\end{thm}
%\smallskip

It turns out that  the assumption that $u(\cdot,0) \in C^{1,1}(\R^d)$ if $d>1$ and $H$ is not quadratic   is necessary to have estimates like \eqref{takis27} and \eqref{takis28}. This is the claim  of  the next result.
\begin{thm}\label{takis29} Assume \eqref{hconvex2} and  $d>1$.  If  \eqref{takis27} holds for all solutions $u \in \text{BUC}(\overline Q_\oo)$ of $u_t=H(Du) $   
(resp.  $u_t=-H(Du)$) in  $Q_\oo$ with $u \in C^{0,1}(\R^d)$ satisfying 
\eqref{manos4} (resp. \eqref{manos5}), then the map  $\lambda \to (D^2H(p +\lambda \xi) \xi^\perp, \xi^\perp)$ must be concave (resp. convex).  In particular, both estimates hold  without any restrictions on the data if and only if 
$H$ is quadratic.
\end{thm}

The motivation behind Theorem~\ref{gal60} and Theorem~\ref{takis26.1} and Theorem~\ref{takis26} is twofold. The first is to obtain as sharp as possible regularity results for solutions of \eqref{gal65}. The second is to  obtain  intermittent regularity results for  \eqref{gal40}, like the ones obtained in \cite{gg} in the specific case that $H(p)=\frac{1}{2}|p|^2$, where, of course, $\theta=\Theta=1$, 
$F(Du)D^2uF(Du)=D^2u$  and the ``new'' estimates are the same as the old ones, that is, \eqref{takis2221} and \eqref{takis2222},  which hold without any regularity conditions.  
\smallskip

The regularity results of \cite{gg} follow from an iteration of \eqref{gal66}, \eqref{gal67}, \eqref{takis2221} and \eqref{takis2222}. As  shown next,  the iteration scheme cannot work when $H$ is not quadratic unless $d=1$. 
\smallskip

To explain the problem, I consider the first two steps of the possible iteration for  $u\in \text{BUC}(Q_\oo)$ solving
\[\begin{split}& u_t=H(Du) \ \text{in} \ \R^d\times (0,a], \quad  u_t=-H(Du) \ \text{in} \ \R^d\times (a,a+b] \quad\\[1.2mm]
 & \text{and} \quad u_t=H(Du) \ \text{in} \ \R^d\times (a+b,a+b+c].\end{split}\]
If the only estimates available were \eqref{gal66},  \eqref{gal67},  \eqref{takis2221} and  \eqref{takis2222}, we find, after some simple algebra, that 
$$D^2u(\cdot, a) \geq -\dfrac{1}{\theta a} I, \quad D^2 u(\cdot, a+b) \geq - \dfrac{1}{(\theta a - \Theta b)_+} I  \quad \text{and} \quad  D^2 u(\cdot, a+b+c) \geq - \dfrac{1}{(\theta a - \Theta b)_+ + \theta c} I.$$
It is immediate that the above estimates cannot be iterated unless there is a special relationship  between the time intervals and the convexity constants, something which will not be possible for arbitrary continuous paths $\xi$. 
\smallskip

If it were possible, as is the casewhen $d=1$, to use the estimates of Theorem~\ref{takis26} without any regularity restrictions, then    Theorem~\ref{gal60}, Theorem~\ref{takis26.1} and Theorem~\ref{takis26}  would imply
$${\mathcal W} (a) \geq -\dfrac{1}{a} I, \quad {\mathcal W} (a+b) \geq -\dfrac{1}{(a-b)_+} I \quad \text{and} \quad {\mathcal W} (a+b+c) \geq -\dfrac{1}{(a-b)_+ +c} I,$$
which can be further iterated, since the estimates are expressed only in terms of  increments $\zeta$.
\smallskip

Before turning to the intermittent regularity results, it is necessary to make some additional remarks. For the sake of definiteness, I continue the discussion in the context of the example above. Although $u(\cdot, a)$ may not be in $C^{1,1}$, it follows from \eqref{takis2221} and \eqref{manos5} that, for some $h\in (0,b]$ and $t\in (a, a+h)$, $u(\cdot, t)\in C^{1,1}$. There is no way, however, to guarantee that $h=b$. Moreover, as was shown in \cite{lsregularity}, in general,  it is possible to have $u$ and $h> 0$ such that $u_t=-H(Du)$ in $\R^d\times (-h,0]$, $u_t=H(Du)$  in $\R^d\times (0,h]$, $u(\cdot,t) \in C^{1,1}$ for $t\in (-h,0) \cup (0,h)$ and $u(\cdot,0) \notin C^{1,1}$. The implication is that when $d >1$ and $H$ is not  quadratic, there is no hope to obtain after iteration smooth solutions. 
\smallskip

To state the  results about intermittent regularity, 
 %about the long time behavior of the solutions of \eqref{takis100}, 
it is convenient to   introduce  the running maximum and minimum functions $M:[0,\infty)\to \R$ and $m:[0,\infty)\to \R$ of a path $\zeta \in 
C_0([0,\infty);\R)$ defined respectively by 
\begin{equation}\label{running}
M(t)=\underset{0\leq s \leq  t}\max \zeta(t) \quad \text{and} \quad m(t)=\underset{0\leq s \leq  t}\min \zeta(t).
\end{equation}
\begin{thm}\label{takis310}  Assume   \eqref{hconvex2} and  either  $d=1$ or  $H$ is quadratic  when  $d>1$, fix  $\zeta \in C_0([0,T);\R)$ and let $u \in \text{BUC}(\overline Q_\oo$ be a solution of  \eqref{gal40}. 
% and, for $\zeta \in C_0([0,T);\R)$,  let $u \in \text{BUC}(\R \times [0,\infty)$ be a solution of  \eqref{takis100}
Then, for all $t> 0$, 
\begin{equation}\label{takis311}
- \dfrac{1}{M(t)-\zeta(t)}\leq -F(Du(\cdot,t))D^2u(\cdot,t)F(Du(\cdot,t)) \leq \dfrac{1}{\zeta(t)-m(t)}.  %\leq F(Du(\cdot,t))D^2u(\cdot,t)F(Du(\cdot,t)) \leq \dfrac{1}{M(t)-\zeta(t)}.
\end{equation}
\end{thm}
\smallskip

Note that when \eqref{takis311} holds, then, at times $t$ such that $m(t) < \zeta (t) < M(t)$, $u(\cdot, t)\in C^{1,1}(\R^d)$ and  
% of  \eqref{takis712} when  $d=1$ is in $C_x^{1,1}$ and the bound is  
\eqref{takis311} implies that, for all $t> 0$, 
\begin{equation}\label{IRR}
|F(Du(\cdot, t)) D^2u (\cdot, t))F(Du(\cdot, t))| \leq  \max\left[\frac{1}{\zeta(t) -m(t)}, \frac{1}{M(t) -\zeta(t)}\right].
\end{equation}
\smallskip

When, however, \eqref{takis311} is not available, the best regularity estimate available, which is also new, is a decay on the  Lipschitz constant $\|Du\|$. % in what follows $\|\cdot\|$ stands for the usual $L^\infty$-norm.
\begin{thm}\label{takis32} Assume \eqref{hconvex2},  fix  $\zeta \in C_0([0,T);\R)$, and  let $u \in \text{BUC}(\overline Q_\oo$ be a solution of  \eqref{gal40}. Then, for all $t> 0$,
\begin{equation}\label{takis33}
\|Du(\cdot,t)\| \leq \sqrt{\dfrac{2\|u(\cdot,t)\|}{\theta (M(t)-m(t))}}.
\end{equation}
\end{thm}
\smallskip

It follows from \eqref{takis33} that, for any $t >0$ such that $m(t) < M(t)$, any solution of \eqref{gal40} is actually Lipschitz continuous.

\smallskip

An immediate consequence of the estimates in Theorem~\ref{takis32} and Theorem~\ref{takis310}, which is based on well known properties of the Brownian motion (see, for, example, Peres~\cite{p})  is the following observation.
\begin{thm}\label{takis1000}
Assume that $\zeta$ is a Brownian motion and $H$ satisfies \eqref{hconvex2}. There exists a random  uncountable subset of $(0,\infty)$ with no isolated points and of Hausdorff measure $1/2$, which depends on $\zeta$, off of which,  any stochastic viscosity solution of \eqref{gal40} is in  $C^{0,1}(\R^d)$ with a bound satisfying \eqref{takis33}. If  $d=1$ or $H$ is quadratic, for  the same set of times,  the solution is in  $C^{1,1}(\R^d)$ and satisfies \eqref{takis311}. % in space.  % {\bf I did not understand PL comment}
\end{thm}

\subsection*{Long time behavior of the ``rough''  viscosity solutions}  I begin with a short introduction about the  long time behavior of solutions of Hamilton-Jacobi equations. 
In order to avoid technicalities due to the behavior of the solutions  at infinity,  throughout this subsection, it is assumed that solutions are periodic functions in ${\mathbb {T}}^d$. 
\smallskip

To explain the difficulties,  I first look at two very simple cases. In the first case, fix some  $p\in \R^d$ and consider the linear initial value
problem
$$du=(p, Du) \cdot d\zeta \ \text{in} \   Q\oo \quad u(\cdot,0)=u_0.$$
%for some fixed $p \in \R^d$.
%\smallskip

Its solution is
$u(x,t)=u_0( x + p \zeta(t)),$
and clearly it is not true that $u(\cdot,t)$ has, as $t\to \infty$, a uniform limit.
\smallskip

The second example is about  \eqref{gal40} with $H$ satisfying \eqref{hconvex2},  and $\dot\xi > 0$ and $\underset{t \to \infty}\lim \xi(t)=\infty$.
Since\\  %It is immediate from the formula of the solution
$u(x,t)=\underset{y\in \R^d} \sup \left[u_0(y) -tH^\star(\dfrac{x-y}{\xi(t)})\right],$
%where $L$ is the convex dual of $H$, 
it is immediate that, as $t\to \infty$ and uniformly in $x$, $u(x,t) \to \sup u$.
\smallskip

The intermittent regularizing results yield information about the long time behavior of the solutions of \eqref{gal40} under the rather weak assumption that 
%\beq\label{gal70}
%du=H(Du)\cdot dB \ \ \text{in} \ \ \R^d\times(0,\infty),
%\Eq
%where $B$ is a standard Brownian motion and 
\beq\label{gal71}
H\in C(\R^d) \ \ \text{is convex \ \ and } \ \ H(p)>H(0)=0 \ \ \text{for all} \ \ p\in \R^d\setminus\{0\}.
\Eq
%The result about the long time behavior of solutions of \eqref{gal}, which is an immediate consequence of Theorem~\ref{takis32} is stated next.
\begin{thm}\label{takis34}
Assume \eqref{gal71}, fix  $\zeta \in C_0([0,T);\R)$, and let $u \in \text{BUC}(\overline Q_\oo$ be a space periodic solution of  \eqref{gal40}.  If there exists  $t_n\to \infty$ such that  $M(t_n)-m(t_n) \to \infty$, then there exists  $u_\infty \in \R$ such that, as $t\to\infty$ and uniformly in space, $u(\cdot,t) \to u_\infty$.
\end{thm}

In the particular case that $\xi$ is a standard Brownian motion the long time result is stated next.

\begin{thm}\label{gal72}
Assume \eqref{gal71}. For almost every Brownian path $B$, if $u \in \text{BUC}(\overline Q_\oo)$ is a periodic solution of $du=H(Du)\cdot dB \  \text{in} \   Q_\oo$,  
% For almost every  Brownian path $\zeta$, %and $u_0$  periodic and Lipschitz continuous 
there exists a constant $u_\infty =u_\infty (B, u(\cdot, 0))$ such that, 
%if $u$ is the stochastic viscosity solution of \eqref{takis1}, then, 
as $t\to \infty$ and uniformly in $\R^d$, $u(\cdot, t) \to u_\infty$.  Moreover, the random variable is, in general, not constant. 
\end{thm}
\begin{proof}  The contraction property and the fact that $H(0)=0$ yield that  the family  $(u(\cdot,t))_{t \geq 0}$ is uniformly bounded. 
\smallskip

It is  assumed next that the Hamiltonian satisfies \eqref{hconvex2}. It follows from the intermittent regularizing  property, the a.s. properties of the running max and min of the Brownian motion, and the fact that  the Lipschitz constant of the solutions decreases in time that, as $t\to \infty$,  $\|Du(\cdot,t)\| \to 0$. 
\smallskip

In view of the periodicity, it follows  that, along subsequences $s_n\to \infty$, the $u(\cdot,s_n)$'s converge uniformly to constants. 

\smallskip

It remains to show that the whole family converges to the same constant. This is again a consequence of the intermittent regularizing result and  the fact that  the periodicity, the contraction property of the solutions of \eqref{gal40} and $H(0)=0$ yield that 
\begin{equation}\label{gal80}
%\begin{cases}
t \to \max_{x \in \R^d} u(x,t) \quad \text{is nonincreasing, and}  \quad t \to \min_{x \in \R^d} u(x,t) \quad \text{is nondecreasing}.
\Eq

It remains to remove the assumption that the Hamiltonians satisfy \eqref{hconvex2}.  Indeed, if   \eqref{gal71} holds, $H$ can be  approximated  uniformly by a sequence $(H_m)_{m\in \N}$ of Hamiltonians satisfying \eqref{hconvex2}. Let $u_m$ be the solution of the \eqref{takis1} with Hamiltonian $H_m$ and same initial datum. Since, as $m\to \infty$,  $u_m \to u$ uniformly in $Q_T$ for all $T >0$, it follows that, for all $t> 0$, 
$$\int_{\mathbb{T}} H(Du(x,t))  dx \leq \underset{m \to \infty} \liminf \int_{\mathbb{T}} H(D u_m (x,t)) dx.$$

Choose  the sequence $t_n$ and $s_n$  as before to conclude.

\end{proof}

I conclude with an example that shows that, in the stochastic setting,  the limit constant $u_\infty$ must  be random.  
\smallskip

Consider the initial value problem
\begin{equation}\label{1}
du=|u_x| \cdot dB \  \text{in} \ Q_\oo \quad u(\cdot,0)=u_0,
\end{equation}
with  $u_0$  a $2$-periodic extension on $\R$ of    $u_0(x)=1-|x-1| \ \text{on} \  [0,2]$. 
%and $\zeta$ is a Brownian motion; to emphasize that we are working with stochastic viscosity solutions below we write $du=|u_x|\circ dB.$ 
%\smallskip
%
Let $c$ be the limit as $t \to \infty$ of $u$. Since $1-u_0(x)=u_0(x+1)$ and $-B$ is also a Brownian motion with the same law as $B$, if $\mathcal L(f)$ denotes the law of the random variable $f$,  it follows that 
\begin{equation}\label{2}
\mathcal L(c)=\mathcal L(1-c).
\end{equation}

If the limit $c$ of the solution of \eqref{1} is deterministic, then  \eqref{2} implies that $c= 1/2$.  It is shown  next that this is not the case.

\smallskip

Recall that the pathwise solutions are Lipschitz with  respect to paths. Indeed,  if $u, v$ are two pathwise solutions of \eqref{1} with paths respectively $B, \xi$ and $u(\cdot, 0)\equiv v(\cdot,0)$, then there exists $L> 0$, which depends on $\|u_{x}(\cdot,0)\|$ such that, for any $T >0$,
\begin{equation}\label{3}
\underset{x\in \R, t \in  [0,T]}\max |u(x,t)-v(x,t)| \leq L \underset{t \in  [0,T]} \max|\zeta(t)-\xi(t)|.
\end{equation}

Next fix $T=2$ and  use \eqref{3}  to compare the solutions of \eqref{1} with $\zeta \equiv B$ and $\xi(t)\equiv t$ and  $\zeta \equiv B$ and $\xi (t)\equiv -t$.
\smallskip

When $\xi(t)\equiv t$ (resp. $\xi \equiv -t$) the solution $v$  of \eqref{1} is given by 
$$v(x,t)=\underset{|y|\leq t}\max \ u_0(x+y) \quad (\text{resp.} \ \  
v(x,t)=\underset{|y|\leq t}\min  \ u_0(x+y)).$$   

It  is then simple to check that, if $\xi(t)\equiv t$, then $v(\cdot, 2)\equiv 1,$ while, when $\xi(t)\equiv -t$, $v(x,2)=0$. 
\smallskip

Fix $\ep=1/4L$ and consider the events
\begin{equation}\label{11}
A_+:=\left\{\underset{t\in [0,2]} \max|B(t) -t| < \ep \right\} \quad \text{and} \quad A_-:=\left\{\underset{t\in [0,2]} \max|B(t) + t| > \ep\right\}.
\end{equation}
Of course, 
\begin{equation}\label{12}
\mathbb{P}(A_+) > 0 \quad \text{and} \quad \mathbb{P}(A_-) > 0.
\end{equation}

Then  \eqref{3} implies 
\begin{equation}\label{takis12}
u(x,2) \geq 1-L\ep = 3/4  \ \ \text{on} \ \ A_+ \quad \text{and} \quad u(x,2) \leq L\ep = 1/4   \ \ \text{on} \ \ A_-.
\end{equation}

It follows that the random variable $c$  cannot be constant since in $A_+$ it must be bigger than $3/4$ and in $A_-$ smaller than $1/4$.

\smallskip

In an upcoming publication (Gassiat, Lions and Souganidis~\cite{gls}) we are visiting this problem and obtain in a special case more information about $u_\infty$.

\section[Stochastic viscosity solutions]
{Pathwise solutions for fully nonlinear,
second-order  PDE with rough signals and 
 smooth, spatially homogeneous Hamiltonians}
%%%%%%%%%%%%%%%%%%

Consider the initial value problem 
\begin{equation}\label{FNL1}
%\begin{cases}
du = F(D^2 u,Du,u,x,t)\,dt 
+ \sum_{i=1}^m H^i (Du)\cdot dB_i \  \text{ in }\ 
Q_\oo\  \quad u(\cdot,0)= u_0,\\
%\noalign{\vskip6pt}
%u(\cdot,0)= g\ \text{ in }\ \R^d,
%\end{cases}
\end{equation}
with 
\begin{equation}\label{RegH}
H=(H^1,\ldots, H^m) \in C^2 (\R^d;\R^m),
\end{equation}
\begin{equation}\label{paris272}
B=(B_1,\dots, B_m)  \in C_0([0,\infty);\R^m)
\end{equation}
and
\begin{equation}\label{paris270}\
\text{$F$ is degenerate elliptic.}
\end{equation}
% that is, for all $(p,r,x,t)\in \R^d\times
%\R\times\R^d\times [0,\infty)$ and all $X,Y\in S^N$ such that $X\leqq Y$,
%\begin{equation}\label{de2}
%F(X,p,r,x,t) \geqq F(Y,p,r,x,t)\ .
%\end{equation}
%\begin{equation}\label{de2}
%\left\{\begin{array}{l}
%\text{$F$ degenerate elliptic, i.e., for all $(p,r,x,t)\in \R^d\times
%\R\times\R^d\times [0,\infty)$}\\
%\text{and all $X,Y\in S^N$ such that $X\leqq Y$,}\\
%\hskip.5in F(X,p,r,x,t) \geqq F(Y,p,r,x,t)\ .
%\end{array}\right.
%\end{equation}

The case of ``irregular'' Hamiltonians requires different arguments. Spatially dependent regular Hamiltonians are discussed later. 
\smallskip

An important question is if the Hamiltonian's can depend on $u$ and $Du$ at the same time. 
The theory for Hamiltonians depending only on $u$ was developed  in Section~3.
The case where $H$ depends both on $u$ and $Du$ is an open problem 
with the exception of a few special cases, like, for example, linear 
dependence on $u$ and $p$, which are basically an exercise.  %We also remark that to develop such a the theory it is necessary to make 
%additional hypotheses on $F$ similar to the ones needed for the analogous 
%deterministic theory. 
%We introduce these conditions when they are needed. 
\smallskip

The  theory of viscosity solutions for equations like 
\eqref{FNL1} with $H \equiv 0$ is based on using smooth test functions to 
test the equation at appropriate points. % (see Appendix~A).
As already discussed earlier this can not be applied directly to 
\eqref{FNL1}. 
\smallskip
 
Recall that, when $H$ is sufficiently regular,  it is possible to construct, using the characteristics, local in time smooth solutions to \eqref{paris100}. % using  characteristics were constructed in  Section~XXX.
These solutions, for special initial data,  play the role of the smooth test functions for 
\eqref{FNL1}.

\begin{defn}\label{defn8.2}
Fix $B\in C([0,\infty); \R^m)$ and $T >0$. A function  
$u \in BUC (\overline Q_T )$ 
is a pathwise subsolution (resp. supersolution)  of  \eqref{FNL1} 
if, for any maximum (resp. mimum)  $(x_0,t_0)\in Q_\oo$  of $u-\Phi -\psi$, 
where $\psi \in C^1((0,\infty))$ and $\Phi$ is a smooth solution of 
$d\Phi = \sum_{i=1}^m H^i (D\Phi) \circ dB^i$ 
in $\R^d\times (t_0-h, t_0+h)$ for some $h> 0$, then 
\begin{equation}\label{sub}
\psi'(t_0)\leqq F(D^2 \Phi (x_0,t_0), D\Phi (x_0,t_0), u(x_0,t_0), x_0,t_0)
\end{equation}
\begin{equation}\label{super}
(\text{resp.} \quad  \psi'(t_0)\geqq F(D^2 \Phi (x_0,t_0), D\Phi (x_0,t_0), u(x_0,t_0), x_0,t_0).)
\end{equation}
%A function $u \in BUC (\overline Q_T )$ is 
%a super-solution of \eqref{FNL1}
%if, for any minimum $(x_0,t_0)\in Q_T$ of 
%$u-\Phi-\psi$, where $\psi\in C^1 ((0,\infty))$ and  $\Phi$ is a smooth 
%solution of $\Phi_t = \sum_{i=1}^m H^i (D\Phi)\circ dW_i$ 
%in $\R^d\times (t_0-h,t_0+h)$ for some $h> 0$, then 
%\begin{equation}\label{super1}
%\psi' (t_0) \geqq F(D^2\Phi (x_0,t_0), D\Phi (x_0,t_0), u(x_0,t_0),x_0,t_0)\ .
%\end{equation}
Finally, $u \in BUC (\overline Q_T )$ 
is a solution of \eqref{FNL1} if
it is both a subsolution  and supersolution.
\end{defn} 

As in the. classical ``non rough'' theory, it is possible to have  upper-semicontinuous subsolutions, lower-semicontinuous supersolutions and discontinuous solutions.  For simplicity, this is avoided here. Such weaker ``solutions'' are used to carry out the Perron construction in Section~10.
\smallskip

Although somewhat natural,  
the definition introduces several difficulties at 
the technical level. 
One of the advantages of the theory of viscosity solutions is the flexibility associated with the choice of the test functions.
This is not, however, the case here.  
As a result,  it is necessary to work very hard to obtain facts which were 
almost trivial in the deterministic setting. 
For example, in the definition, it is often useful to assume that the max/min is 
strict. Even this fact, which is trivial for classical viscosity solutions, in the current setting 
requires a more work. 
\smallskip

It is also useful to point out the relationship 
between the approach used  for equations with linear dependence on $Du$
and the above definition. 
Heuristically, in Definition~\ref{defn8.2},  one  inverts locally 
the characteristics in an attempt to ``eliminate'' the bad term 
involving $dB$.
Since the problem is nonlinear and $u$ is not regular, it is, of course, 
not possible to do this globally.
In a way consistent with the spirit of the theory of viscosity solutions,  
 this difficulty is overcome by working at the level of the test functions, 
where, of course, it is possible to invert locally the characteristics.
The price to pay for this  is that the test functions  used
here are very robust and not as flexible as the ones used in the classical 
deterministic theory.
%Recall that instead of quadratics here it is needed to use short time smooth 
%evolution of quadratics.
This leads to several technical difficulties, since all the theory has 
to be revisited.

The fact that Definition~\ref{defn8.2} is good in the sense that 
it agrees with the classical (deterministic) one if $B \in C^1$,  is left as an exercise. There are also several other preliminary facts about short time behavior, etc., which are omitted.
\smallskip

The emphasis here is on establishing a comparison principle and some stability properties. The existence
follows  either by a density argument or by Perron's method. The latter was established lately in a very general setting by Seeger \cite{seeger} for $m\geq1$.

\smallskip

The next result is about the  stability properties of the pathwise viscosity solutions.
Although it can be stated in a much more general form using 
``half relaxed limits'' and lower- and upper-semicontinuous envelopes, here 
it is  presented in a simplified form.

\begin{prop}\label{prop:stability}
Let $F_n, F$ be degenerate elliptic, $H_{n}, H \in C^2 (\R^d;\R^m)$, $B_n, B\in C([0,\infty);\R^m)$ 
%and $u_{0,n}, u_0\in BUC (\R^d)$ 
be such that 
$\sup_{i,n} \|D^2 H_{i,n} \| < \infty$  and, as $n\to\infty$  and locally uninformly, $F_n\to F$,  
$H_n \to H$ in $C^2(\R^d;\R^m)$,  and 
$B_n\to B$ in $C([0,\infty);\R^m)$. % and $u_{0,n}\to u_0$ in $C(\R^d)$. 
If $u_n$ is a pathwise  solution  of
\eqref{FNL1} with  nonlinearity $F_n$, Hamiltonian $H_n$ and path $B_n$ and 
 $u_n\to u$ in $C(\overline Q_T)$, then $u$ is a pathwise   solution of \eqref{FNL1}.
\end{prop}

The assumptions that $H_n\to H$ in $C^2(\R^d;\R^m)$ instead of just in 
$C(\R^d)$ and 
$\sup_n \|D^2 H_n\| <\infty$ are not needed  for the ``deterministic'' 
theory. 
Here they are dictated by the nature of the test functions. 
%We will revisit this issue later when we consider nonsmooth $H$'s.

\begin{proof}[Proof of Proposition~\ref{prop:stability}]
Let $(x_0,t_0)\in\R^d\times (0, T]$ be a strict maximum of 
$u-\Phi-\psi$ where $\psi\in C^1 ((0,\infty))$ and,  for some $h> 0$, $\Phi$ is a smooth 
solution of \eqref{paris100} in $(t_0-h,t_0+h)$.
% for some $h> 0$, 
%given by the characteristics.  
\smallskip

Let $\Phi_n$ be the  smooth solution of 
\begin{equation*}
%\begin{cases}
\Phi_{nt} = H_n (D\Phi_n) \dot B_n \ \text{in}\  \R^d \times 
(t_0-h_n,t_0+h_n)\quad 
%\noalign{\vskip6pt}
\Phi_n(\cdot, t_0) =\Phi(\cdot, t_0).
%\end{cases}
\end{equation*}
%given by the method of characteristics.
%\smallskip

The assumptions on the $H_n$ and $B_n$ imply that,  as $n\to \infty$, 
$\Phi_n\to \Phi$,
$D\Phi_n \to D\Phi$ and $D^2\Phi_n\to D^2\Phi$  in 
$C(\R^d\times (t_0-h', t_0+h'))$, for some, uniform in $n$, $h'\in (0,h)$; note that this is  the place 
where $H_n\to H$ in $C^2(\R^d)$ and 
$\sup_n \|D^2 H_n\| < \infty$ are used.
\smallskip

Let $(x_n,t_n)$ be a maximum point of $u_n - \Phi_n-\psi$ in 
$\R^d\times [t_0-h', t_0+h']$. 
Since $(x_0,t_0)$ is a strict maximum of $u-\Phi -\psi$, there exists
a subsequence such that $(x_n,t_n)\to (x_0,t_0)$.
The definition of viscosity solution then gives 
$$\psi' (t_n) \leqq 
F(D^2\Phi_n (x_n,t_n) ,D\Phi_n (x_n,t_n),u_n(x_n,t_n),x_n,t_n)
\ .$$

Letting $n\to\infty$ yields the claim.

\end{proof}

The next result is  the comparison principle for  pathwise viscosity solutions of  the first-order initial value  problem, that is,
\begin{equation}\label{paris273}
du =  \sum_{i=1}^m H^i (Du)\cdot dB_i \  \text{ in }\ 
Q_\oo\  \quad u(\cdot,0)= u_0.\\
\end{equation}

\begin{thm}\label{prop:uniqueness} 
Assume that \eqref{RegH}, \eqref{paris272} and  $u_0 \in BUC(\R^d)$. Then   \eqref{paris273} has a unique pathwise solution
$u\in BUC (\overline Q_\oo)$ which agrees with the ``solution'' obtained from the extension operator.  % for all $T> 0$.
\end{thm}

The proof follows from the arguments used to prove the next result about the extension operator for \eqref{FNL1} which is stated next, hence it is omitted. 
\smallskip

The next result is about the extension operator for \eqref{FNL1}. As before, it is shown that the solutions to  initial value problems \eqref{FNL1} with smooth time signal approximating the given rough one form a Cauchy family in 
$BUC(\overline Q_T)$ and, hence, all converge to the same function which is a pathwise viscosity solution to \eqref{FNL1}.
\smallskip

The next result  provides an extension from  smooth  to 
arbitrary continuous paths $B$.  For simplicity the the dependence of $F$ on $u,x$ and $t$ is omitted. 

\begin{thm}\label{exist}
Assume   \eqref{RegH}, \eqref{paris272} and \eqref{paris270}  and fix $u_0\in \text{BUC}(\R^d)$ and $B\in C_0([0,\infty);\R^m)$.
 Consider   
 two families $(\zeta_\ep)_{\ep> 0}$, $(\xi_\eta)_{\eta> 0}$ in $C_0([0,\infty);\R^m) \cap C^1([0,\infty);\R^m)$
and $(u_{0,\ep})_{\ep> 0}$, $( v_{0,\eta})_{\eta> 0}\in BUC(\R^d)$ such that, as $\ep,\eta\to0$,  $\zeta_\ep$ and  $\xi_\eta$ converge to $B$ in  $C([0,\infty);\R^m)$  and  
$u_{0,\ep}$ and $ v_{0,\eta}$ converge to $u_0$ uniformly in $\R^d$. 
Let $(u_\ep)_{\ep> 0}, (v_\eta)_{\eta> 0}\in BUC(\overline Q_\oo)$ 
be the unique viscosity 
solutions of \eqref{FNL1} with signal and initial datum $\zeta_\ep, u_{0,\ep}$ and $\xi_\eta, v_{0,\eta}$ respectively. %and \eqref{FNL4}. 
Then, for all $T >0$,  
as $\ep,\eta\to0$, $u_\ep- v^\eta \to0$ uniformly in 
$\overline Q_T$. 
In particular, the family $(u_\ep)_{\ep> 0}$ is Cauchy in 
$BUC(\overline Q_T)$ and all  approximations converge to the same limit.
\end{thm}

\begin{proof}
Fix $T>0$ and consider the doubled initial value problem 
\begin{equation}\label{eq100}
\begin{cases} 
dZ^{\lambda,\ep,\eta} = \sum_{i=1}^m H(D_x Z^{\lambda,\ep,\eta} ) \dot \xi_{i,\ep} 
- \sum_{i=1}^m H(-D_y Z^{\lambda,\ep,\eta}) \dot{\zeta}_{i,\eta}, \ \text{ in }\ 
\R^d\times\R^d \times (0,T)\\
\noalign{\vskip6pt}
Z^{\lambda,\ep,\eta} (x,y,0) = \lambda |x-y|^2 .
\end{cases}
\end{equation}
It is immediate that 
%\begin{equation}\label{eq102} 
$Z^{\lambda,\ep,\eta} (x,y,t) = \Phi^{\lambda,\ep,\eta} (x-y,t),$ 
%\end{equation} 
where 
\begin{equation}\label{eq103} 
%\begin{cases}
\Phi_t^{\lambda,\ep,\eta} = \sum_{i=1}^mH(D_z\Phi^{\lambda,\ep,\eta}) (\dot \xi_{i,\ep} 
- \dot \zeta_{i,\eta}) \ \text{ in }\  Q_T\quad 
%\noalign{\vskip6pt}
\Phi^{\lambda,\ep,\eta} (z,0) = \lambda |z|^2.
%\end{cases}
\end{equation}

As discussed earlier,  there exists $T^{\lambda,\ep,\eta} > 0$ such that 
$\Phi^{\lambda,\ep,\eta}$ is given by the method of characteristics in 
$\R^d\times [0,T^{\lambda,\ep,\eta})$ and
% $\lim_{\ep,\eta\to0}T^{\lambda,\ep,\eta}=\infty$
%as $\ep,\eta\to0$, $T^{\lambda,\ep,\eta}\to +\infty$ and, for all $T> 0$,   
\begin{equation}\label{eq101} 
\lim_{\ep,\eta\to0}T^{\lambda,\ep,\eta}=\infty \ \text{and} \  \lim_{\ep,\eta\to0}\sup_{(z,t) \in \R^d\times [0,T] }\left(\Phi^{\lambda,\ep,\eta}(z)  - \lambda |z|^2\right)=0.
%\ \text{ uniformly in }\ \R^d\times [0,T]\ .
\end{equation}

The conclusion will follow as soon as it established that 
\begin{equation}\label{e1}
\lim_{\lambda\to\infty} \varlimsup_{\ep,\eta\to0} 
\sup_{(x,y)\in \R^{2N},t\in [0,T]} 
(u^\ep (x,t) - v^\eta (y,t) - \lambda |x-y|^2) = 0\ .
\end{equation}
Consider next the function  
$$\Psi^{\lambda,\ep,\eta} (x,y,t) = u^\ep (x,t) - v^\eta (y,t) 
- \Phi^{\lambda,\ep,\eta} (x-y,t)\ .$$
The classical theory of viscosity solutions (see \cite{cil}) 
yields that the map  
%\marginpar{{\bf PL:} \em should\\ we give details here?}
$$t\longmapsto M^{\lambda,\ep,\eta}(t) = 
\sup_{x,y\in\R^d} 
[u^\ep (x,t) -  v^\eta (y,t) - \Phi^{\lambda,\ep,\eta} (x-y,t)]$$
is nonincreasing  in $[0,T^{\lambda,\ep,\eta})$. 
\smallskip

Hence, for $x,y\in \R^d$ and $t\in [0,T^{\lambda,\ep,\eta})$, 
$$u^\ep (x,t) -  v^\eta (y,t) - \Phi^{\lambda,\ep,\eta}(x-y,t) 
\leqq \sup_{x,y\in\R^d} (u_0^\ep (x) -  v_0^\eta (y) 
-\lambda |x-y|^2)\ .$$

The claim now follows from the assumptions on $u_0^\ep$ and $v_0^\eta$. 

\end{proof}

The uniqueness of the pathwise  viscosity solutions of \eqref{FNL1}
is considerably more complicated than the one for \eqref{paris273}.
This is consistent with the deterministic theory, where the uniqueness theory 
of viscosity solutions for second-order degenerate, elliptic equations is 
by far more complex than the one for Hamilton-Jacobi equations.
For the same reasons as for the existence, I will present the argument 
omitting the dependence on $u,x$ and $t$.
\smallskip

The proof follows the general strategy outlined in the ``User's Guide''. 
The actual arguments  are, however, different and more complicated.
\smallskip

Recall that in the background of the ``deterministic'' proof are 
the so called sup- and inf-convolutions. 
These are particular regularizations that yield approximations which 
have parabolic expansions almost everywhere and are also subsolutions and 
supersolutions of the nonlinear pde.
\smallskip

This is exactly where the pathwise case becomes different. 
The ``classical'' sup-convolutions and inf-convolutions of pathwise viscosity solutions 
do not have parabolic expansions. To deal with this serious difficulty, it is necessary 
%In accordance with the general strategy about pathwise viscosity solutions,
to  change the sup-convolutions and inf-convolutions by replacing the 
quadratic weights by short time smooth solutions of the first-order part 
of \eqref{FNL1}. 
The new regularizations have now parabolic expansions---the reader should 
think that the new weights remove the ``singularities'' due to the roughness
of $B$.
%The rest of the argument follows then the proof of the deterministic 
%uniqueness, although the details are more intricate. 

\begin{thm}\label{uniq}
Assume \eqref{RegH}, \eqref{paris272} and \eqref{paris270}. Let $u,v\in BUC(\overline Q_\infty)$ be respectively a viscosity subsolution and 
supersolution of \eqref{FNL1}. 
Then, for all $t\geq 0$,
\begin{equation}\label{comp}
\sup_{x\in\R^d} (u-v) (x,t) \leqq 
\sup_{x\in\R^d} (u-v) (x,0)\ .
\end{equation}
\end{thm}
\begin{proof}
To simplify the presentation below it is assumed that $m=1$. 
Recall that, for any $\phi \in C^3(\R^d\times \R^d)\cap C^{0,1} (\R^d\times\R^d) $, there exists some $a> 0$ such that the doubled initial value problem 
\begin{equation*}
%\begin{cases}
dU=[H(D_xU)- H(-D_y U)]\circ dB\ \text{ in }\ \R^d\times (t_0-a, t_0+a)\quad U(x, y,t_0)=\phi(x, y),%\ \text{ in }\ \R^d\times \R^d. \\ 
%\noalign{\vskip6pt}
%U=\phi\ \text{ on }\ \R^d\times \{t_0\}\ ,
%\end{cases}
\end{equation*}
has a  smooth solution which, for future use, is denoted by  $S_H ^d (t-t_0,t_0)\phi$.
\smallskip
%which exists by the method of characteristics, if  
%$\phi \in C^3 (\R^d\times\R^d)\cap C^{0,1} (\R^d\times\R^d)$ 
%for some $a> 0$ depending on $\|\phi\|_{C^2}$.

If $\phi$ is of separated form,  that is,  $\phi (x,y) = \phi_1(x) +\phi_2(y)$, making if necessary, the interval of existence smaller,
it is immediate that 
\begin{equation*}
S^d (t-t_0,t_0) \phi (x,y) 
= S_H^+ (t-t_0,t_0) \phi_1 (x) + S_H^- (t-t_0,t_0)\phi_2 (y)\ ,
\end{equation*}
where, as before, $S^{\pm}_{H^{\pm}}$ denote the smooth short time 
solution operators to  $du = \pm H(Du) \cdot dB$.% and $du = - H(-Du)\circ dW$ 
\smallskip

Moreover, for any $\lambda > 0$ and $t,t_0\in \R$, it is obvious that 
\begin{equation*} 
S^d (t-t_0,t_0) (\lambda |\cdot  - \cdot|^2) (x,y) = \lambda|x-y|^2\ .
\end{equation*}

Finally, again for smooth solutions, 
\begin{equation*}
S_H^- (t-t_0,t_0)\phi_2 (y)  = - S_H^+ (t-t_0,t_0) (-\phi_2) (y).
\end{equation*}

Fix $\mu >0$. The claim is that,  for large enough $\lambda$,
\begin{equation*} 
\Phi (x,y,t) = u (x,t) - u(y,t) - \lambda |x-y|^2 - \mu t
\end{equation*}
cannot have  a maximum %$(x_\lambda,y_\lambda,t_\lambda)
in $ \R^d\times \R^d \times (0,T]$. 
This leads to the desired conclusion as in the classical proof of the maximum principle. 
\smallskip

Arguing by contradiction, it is assumed that there exists 
$(x_\lambda,y_\lambda,t_\lambda) \in \R^d\times\R^d\times (0,T]$ 
such that, for all $(x,y,t) \in \R^d\times\R^d \times [0,T]$,
\begin{equation}\label{eq:uniq} 
\Phi (x,y,t) = u(x,t) - v(y,t) - S^d(t-t_\lambda,t_\lambda) 
(\lambda |\cdot - \cdot|^2) (x,y) - \mu t
\leqq \Phi (x_\lambda,y_\lambda,t_\lambda)\ .
\end{equation} 

%\begin{proof}[Proof of Theorem~\ref{uniq}]
%As in the proof of the comparison principle for the deterministic version 
%of \eqref{FNL2}, \eqref{comp} follows if we show that the assumption that, 
%for some $\mu> 0$, the function 
%\begin{equation}\label{eq:uniq}
%\left\{
%\begin{array}{rcl}
%\Phi (x,y,t) & = & u(x,t) - v(y,t) - \lambda |x-y|^2 -\mu t\\
%&=& u(x,t) - v(y,t) - \widetilde S (t-t_\lambda,t_\lambda)
%[\lambda |\cdot -\cdot |^2] (x,y)-\mu t
%\end{array}\right.
%\end{equation}
%attains a maximum $(x_\lambda,y_\lambda,t_\lambda)\in\R^d\times\R^d
%\times (0,T)$ as $\lambda\to\infty$, leads to contradiction.

To handle the behavior at infinity and assert the existence of a maximum, it is 
necessary to consider 
$ S^d(t-t_\lambda,t_\lambda)[\lambda|\cdot-\cdot|^2 +\beta \nu(\cdot)]
(x,y)$ instead of 
$S^d(t-t_\lambda)[\lambda|\cdot-\cdot|^2]$ in \eqref{eq:uniq}, for 
$t-t_\lambda$ small, $\beta\to0$, and a smooth approximation $\nu(x)$ of $|x|$.
Since this adds some tedious details which may obscure the main ideas 
of the proof, below it is assumed that a maximum exists. % and obtain a contradiction.
\smallskip

Elementary computations and a straightforward application of the 
Cauchy-Schwarz inequality yield, for all $\ep > 0$ and $\xi,\eta\in\R^d$, 
\begin{equation}\label{eq.cs} 
\begin{split}
|x-y|^2 - |x_\lambda - y_\lambda|^2 
&\leqq 2\langle x_\lambda - y_\lambda, x-x_\lambda - \xi\rangle + 
- 2\langle x_\lambda - y_\lambda,y-y_\lambda-\eta\rangle\\
&\qquad
  + 2\langle x_\lambda-y_\lambda, \xi-\eta\rangle\\
&\qquad 
+ (2+\ep^{-1}) (|x-x_\lambda-\xi|^2 + |y-y_\lambda -\eta|^2) 
+ (1+2\ep) |\xi-\eta|^2.
\end{split}
\end{equation} 
Let 
\begin{equation*} 
p_\lambda = \lambda (x_\lambda - y_\lambda) \ ,\quad 
\lambda_\ep = \lambda (2+\ep^{-1})\quad\text{and}\quad 
\beta_\ep = \lambda (1+2\ep)\ .
\end{equation*} 

%We follow the proof of the uniqueness of ``deterministic' viscosity 
%solutions presented in \cite{CIL}. 
%A straightforward application of the Cauchy-Schwarz inequality yields, 
%for $\ep> 0$ and $\xi,\eta\in\R^d$,
%\begin{equation*}
%\lambda| x-y|^2 \leq \lambda_\ep (|x-\xi|^2 + |y-\eta|^2) 
%+\beta_\ep |\xi-\eta|^2\ ,
%\end{equation*}
%where 
%\begin{equation*}
%\lambda_\ep = \ep^{-1} + 2\lambda\quad\text{ and }\quad 
%\beta_\ep =  \lambda\lambda_\ep\ .
%\end{equation*}
The comparison of local in time smooth solutions of stochastic 
Hamilton-Jacobi equations, which are easily obtained by the method of 
characteristics, and the facts explained  before the beginning 
of the proof yield that the function 
\begin{equation*}
\begin{split} 
\Psi (x,y,\xi,\eta,t) &= u(x,t) - v (y,t) - S_H^+ (t-t_\lambda,t_\lambda) 
(2\langle p_\lambda,\cdot- x_\lambda - \xi\rangle 
+ \lambda_\ep |\cdot - x_\lambda - \xi|^2)(x)\\
& \qquad 
-  S_H^- (t-t_\lambda,t_\lambda)(-2\langle p_\lambda, \cdot -y_\lambda-\eta \rangle
+ \lambda_\ep |\cdot - y_\lambda - \eta|^2) (y)\\ 
& \qquad
-2\langle p_\lambda, \xi-\eta\rangle - \beta_\ep |\xi-\eta|^2 - \mu t
\end{split}
\end{equation*}
achieves,  for $h\le h_0 = h_0 (\lambda,\ep^{-1})$,
its maximum in $\R^d\times\R^d\times\R^d\times\R^d 
\times (t_\lambda - h,t_\lambda +h)$ at 
$(x_\lambda, y_\lambda,0,0,t_\lambda)$. 
\smallskip

Note that here it is necessary to  take $t-t_\lambda$ sufficiently small  
to have local in time smooth solutions for the doubled as well as the
$H$ and $-H(-)$ equations given by the characteristics. 
\smallskip

For $t\in (t_\lambda -h,t_\lambda +h)$ define the modified sup- and 
inf-convolutions
\begin{equation*}
\bar u(\xi,t) = \sup_{x\in\R^d} [u(x,t) - S_H^+ (t-t_\lambda,t_\lambda) 
(2\langle p_\lambda,\cdot-x_\lambda-\xi\rangle
+ \lambda_\ep|\cdot - x_\lambda - \xi|^2)(x)]
\end{equation*}
and
\begin{equation*}
\un{v}(\eta,t) = \inf_{y\in\R^d} [v(y,t) +S_H^- (t-t_\lambda, t_\lambda)
(-2\langle p_\lambda,\cdot -y_\lambda -\eta\rangle
+ \lambda_\ep |\cdot - y_\lambda -\eta|^2)(y)]\ .
\end{equation*}
It follows that, for $\delta >0$,  
\begin{equation*}
G (\xi,\eta,t) = \bar u(\xi,t) - \un{v}(\eta,t) 
- (\beta_\ep+\delta) |\xi-\eta|^2 -2\langle p_\lambda, \xi-\eta\rangle- \mu t
\end{equation*}
attains its maximum in $\R^d\times\R^d\times (t_\lambda - h,t_\lambda +h)$
at $(0,0,t_\lambda)$.
\smallskip

Observe next that 
% the key properties of $\hat u$ and $\hat v$.  Its proof is outlined: 
there exists a constant $K_{\ep,\lambda}> 0$ such that, in 
$\R^d\times (t_\lambda-h,t_\lambda+h)$,
\begin{equation}\label{lsided}
D_\xi^2 \bar u \ge - K_{\ep,\lambda},\ \ 
D_\eta^ 2\un{v}\le K_{\ep,\lambda}, \quad 
\bar u_t \le K_{\ep,\lambda},\ \text{ and }\ \un{v}_t \ge -K_{\ep,\lambda}\ .
\end{equation}
with the inequalities understood both in the viscosity and 
distributional sense.
\smallskip

The one sided bounds of $D_\xi^2\bar u$ and $D_\eta^2 \un{v}$ are an immediate 
consequence of the definition of $\bar u$ and $\un{v}$ and the regularity 
of the kernels, which imply that, for some $K_{\ep,\lambda} > 0$ and 
in $\R^d \times (t_\lambda - h, t_n +h)$, 
\begin{equation*}
\begin{split}
&|D_\xi^2 S_H^+ (\cdot - t_\lambda,t_\lambda)
(2\langle p_\lambda, \cdot - x_\lambda -\xi\rangle + \lambda_\ep 
|\cdot - x_\lambda - \xi|^2)| \\
&\qquad 
+ |D_\eta^2 S_H^- (\cdot - t_\lambda,t_\lambda) 
(-2\langle p_\lambda,\cdot - y_\lambda - \eta\rangle +\lambda_\ep 
|\cdot - y_\lambda - \eta|^2)| 
\leqq K_{\ep,\lambda}\ .
\end{split}
\end{equation*}

The bound for $\bar u_t$ is shown next;
the argument for $\un{v}_t$ is similar. 
Note that, in view of the behavior of $B$, such a bound cannot be expected to hold 
for $u_t$. 
Indeed take $F\equiv 0$ and $H\equiv 1$, in which case $u(x,t) = B_t$. 
\smallskip

Assume that, for some smooth function $g$ and for $\xi$ fixed,  
the map $(\xi,t)\to  \bar u(\xi,t) - g(t)$ has a max at $\hat t$. 
It follows that 
$$(x,t)\mapsto 
u(x,t) - S_H^+ (t-t_\lambda,t_\lambda) 
\left(2\langle p_\lambda, \cdot - x_\lambda-\xi\rangle + \lambda_\ep| \cdot - \xi|^2\right)(x) 
- g(t)$$
has a max at $(\hat x,\hat t)$, where $\hat x$ is a point where that 
supremum in the definition of $\bar u(\xi,t)$ is achieved, that is, 
$$\bar u(\xi,\hat t) = u(\hat x,\hat t) - S_H^+ (\hat t- t_\lambda,t_\lambda)
\left(2\langle p_\lambda,\cdot - x_\lambda -\xi\rangle + 
\lambda_\ep| \cdot - \xi|^2 \right)(\hat x)\ .$$

In view of the definition of the pathwise  viscosity sub-solution, 
it follows that there exists some $\hat K_{\ep,\lambda}$ depending on $K_{\ep,\lambda}$ 
and $H$, such that 
$g(\hat t)\le \hat K_{\ep,\lambda}$, and, hence, the claim follows.
\smallskip

The one-sided bounds \eqref{lsided} 
 yield the existence of $p_n,q_n,\xi_n,\eta_n \in\R^d$
and $t_n> 0$ such that, as $n\to\infty$, 

(i) \ $(\xi_n,\eta_n,t_n)\to (0,0,t_\lambda)$, $p_n,q_n\to0$, 
\smallskip

(ii) the map $(\xi,\eta,t) \to \bar u(\xi,t) - \un{v}(\eta,t) 
- \beta_\ep |\xi-\eta|^2 - \langle p_n,\xi\rangle - \langle q_n,\eta\rangle - 2\langle p_\lambda,\xi-\eta\rangle - \mu t$\\ 
has a maximum at $(\xi_n,\eta_n,t_n)$, 
\smallskip

(iii) $\bar u$ and $\un{v}$ have parabolic second-order expansions from above and below at 
$(\xi_n,t_n)$ and $(\eta_n,t_n)$ respectively, that is, there exist $a_n, b_n \in \R$ such that 
\begin{equation*}
\begin{split}
\bar u(\xi,t) & \leq \bar u(\xi_n,t_n) +  
%\bar u_t(\xi_n,t_n)(t-t_n) + (D_\xi \bar u(\xi_n,t_n), \xi-\xi_n)\\
a_n(t-t_n) + (D_\xi \bar u(\xi_n,t_n), \xi-\xi_n)\\
&\qquad + \frac12 (D_\xi^2 \bar u (\xi_n,t_n)(\xi-\xi_n),\xi-\xi_n) 
+ o(|\xi-\xi_n|^2 + |t-t_n|)\ ,
\end{split}
\end{equation*}
and
\begin{equation*}
\begin{split} 
\un{v}(\eta,t) & \geq \un{v}(\eta_n,t_n) + 
b_n
%\un{v}_t (\eta_n,t_n)
(t-t_n) + (D_\eta \un{v}(\eta_n,t_n), \eta-\eta_n)\\
&\qquad + \frac12 ( D_\eta^2 \un{v}(\eta_n,t_n)(\eta-\eta_n),\eta-\eta_n)
+ o(|\eta-\eta_n|^2 + |t-t_n|)\ ,
\end{split}
\end{equation*}
and, finally,
\smallskip

(iv) $a_n=b_n +\mu$, \enspace
%$\bar u_t (\xi_n,t_n) = \un{v}_t (\eta_n,t_n)+\mu$,\enspace
$D_\xi \bar u(\xi_n,t_n) = p_n + 2p_\lambda + 2\beta_\ep (\xi_n-\mu_n)$, \enspace
$D_\eta \un{v}(\xi_n,t_n) = -q_n + 2p_\lambda 2\beta_\ep (\xi_n -\eta_n)$ and \enspace
$D_\xi^2 \bar u (\xi_n,t_n) \le D^2_\eta \un{v}(\eta_n,t_n)$.
\smallskip

It follows that, for some $\theta >0$ fixed, $t <t_n$, $(\xi,t)$ near 
$(\xi_n,t_n)$ and $(\eta,t)$ near $(\eta_n,t_n)$, the maps 
\begin{equation*}
(x,\xi,t)\to 
u(x,t) - S_H^+(t-t_\lambda,t_\lambda) 
( 2\langle p_\lambda,\cdot-x_\lambda -\xi\rangle 
+ \lambda_\ep |\cdot -x_\lambda -\xi|^2)(x) - \Phi (\xi,t),
\end{equation*}
and 
\begin{equation*}
(y,\eta ,t) \to v(y,t) + S_{H}^- 
(t-t_\lambda,t_\lambda)
(-2\langle p_\lambda,\cdot - y_\lambda -\eta\rangle
+ \lambda_\ep|\cdot-\eta|^2)(y) -\Psi (\eta,t),
\end{equation*}

attain respectively  a maximum at $(x_n,\xi_n,t_n)$ and 
a minimum at $(y_n,\eta_n,t_n)$, where 
\begin{equation*}
\begin{split}
\Phi(\xi,t) & = \bar u(\xi_n,t_n) + (\bar u_t (\xi_n,t_n)-\theta) 
(t-t_n) + (D_\xi \bar u(\xi_n,t_n), \xi-\xi_n)\\
&\qquad + 
\frac12 ( (D_\xi^2 \bar u(\xi_n,t_n) + \theta I)(\xi-\xi_n), \xi-\xi_n) 
\end{split}
\end{equation*}
and 
\begin{equation*}
\begin{split}
\Psi (\eta,t) &= \un{v}(\eta_n,t_n) 
+(\un{v}_t (\eta_n,t_n) +\theta) 
(t-t_n) + (D_\eta \un{v}(\eta_n,t_n), \eta-\eta_n)\\
&\qquad + \frac12 ((D_\eta^2 \un{v}(\eta_n,t_n) +\theta I) 
(\eta - \eta_n),\eta-\eta_n)\ .
\end{split}
\end{equation*}
\smallskip

Next, for sufficiently small $r> 0$,  let ${\mathcal B}(\xi_n,t_n,r_n)= B(\xi_n,r) \times (t_n-r,t_n]$ and define
\begin{equation*}
\overline \Phi (x,t)=
 \inf [\Phi (\xi,t) + S^+_H(t-t_\lambda,t_\lambda)
(2 \langle p_\lambda, \cdot - x_\lambda -\xi\rangle 
+ \lambda_\ep |\cdot -x_\lambda- \xi|^2) (x) : (\xi,t) \in {\mathcal B}(\xi_n,t_n,r_n)],
% : (\xi,t) \in B(\xi_{n,r}) \times (t_n-r,t_n)],
%\end{split}
\end{equation*} 
and 
\begin{equation*}
\un{\Psi} (y,t) =
\sup [\Psi (\eta,t) - S^-_H (t-t_\lambda, t_\lambda)
(-2\langle p_\lambda,\cdot - y_\lambda -\eta\rangle
+ \lambda_\ep |\cdot - y_\lambda -\eta|^2) (x) : (\xi,t) \in {\mathcal B}(\xi_n,t_n,r_n)].
\end{equation*} 
%\end{gather*}
\smallskip

It follows that $u-\bar \Phi$ and $v-\un{\Psi}$ 
attain a local max at $(x_n,t_n)$ and a local min at $(y_n,t_n)$.
Moreover, $\bar\Phi$ and $\un{\Psi}$ are  smooth solutions
of $du = H(Du)\cdot dB$ for $(x,t)$ near $(x_n,t_n)$ and 
$dv = -H(-D_yv)\cdot dB$ for $(y,t)$ near $(y_n,t_n)$. 
This last assertion for $\bar\Phi$ and $\un{\Psi}$ follows,  
using  the inverse function theorem, from the 
fact that, at $(x_n,t_n)$ and $(y_n,t_n),$ 
there exists a unique minimum in the definition of $\bar\Phi$ and $\un{\Psi}$.  
This in turn comes from the observation that for $\lambda >\lambda_0$, 
at $(\xi_n,x_n,t_n)$ and $(\eta_n,y_n,t_n),$ 
$$D^2 \overline \Phi (\xi_n,t_n) + (\lambda_\ep +\theta) I> 0
\quad\text{and}\quad 
D^2 \un{\Psi} (\eta_n,t_n) - (\lambda_\ep+ \theta) I <0\ .$$
%\smallskip
Finally, elementary calculations also yield that 
$$D_\xi^2 \Phi (\xi_n,t_n) \geqq D_x^2 \bar\Phi (x_n,t_n)
\quad\text{and}\quad 
D_\eta^2 \Psi (\eta_n,t_n) \leqq D_y^2 \un{\Psi} (y_n,t_n)\ .$$

Applying now the definitions of the pathwise  subsolution  and 
supersolution to $u$ and $v$ respectively,  yields
\begin{equation*}
\begin{split}
\bar u_t (\xi_n,t_n) - \theta 
& \le F(D_x^2 \bar\Phi (x_n,t_n) ,D_x \bar\Phi (x_n,t_n))\\
& \le F(D_\xi^2 \Phi (\xi_n,t_n), D_\xi\Phi (\xi_n,t_n)) 
= F(D_\xi^2 \bar u(\xi_n,t_n) + \theta I, D_\xi\bar u(\xi_n,t_n))
\end{split}
\end{equation*}
and 
\begin{equation*}
\un{v}_t (\eta_n,t_n) + \theta 
\ge F(-D_\xi^2 \un{v}(\xi_n,t_n) - \theta I, D_\eta \un{v}(\eta_n,t_n)).
\end{equation*}
Hence
\begin{equation*}
\begin{split}
\mu - 2\theta & \leq a_n - b_n - 2\theta\\
& \le \sup
[F(A+\theta I,p+p_n) - F(A-\theta I,p+q_n) :  |p_n|, |q_n| 
\le n^{-1}, |A| \le K_{\ep,\lambda}].
\end{split}
\end{equation*}
The conclusion now follows choosing $\ep = (2\lambda)^{-1}$ and 
letting $\lambda \to \infty$ and $\theta \to0$.
%\end{proof}
%\smallskip

\end{proof}
It is worth remarking that, in the course of the previous proof,  
it was  shown that,  
for $0 < h\le \hat h_0$, with $\hat h_0=\hat h_0(\lambda,\ep) \le h_0$,
$\bar u$ $($resp. $\un{v})$ is a viscosity subsolution 
(resp. supersolution) of 
$$\bar u_t \le F(D_\xi^2 \bar u,D_\xi \bar u)\quad (\text{resp. } 
\un{v}_t \ge F(D_\eta^2 \un{v},D_\eta\un{v})) \ \text{ in }\ 
\R^d\times (t_\lambda - h,t_\lambda +h).$$
%\end{proof}
%%%%%%%%%%%%%%%%%%%%%%%%%%%%%%%%%%%%%%%%%%%%%%
\section{Pathwise solutions to fully nonlinear first and second order pde with spatially dependent smooth Hamiltonians}
\subsection*{The general problem, strategy and difficulties} The next step in the development of the theory is to consider spatially dependent Hamiltonians and, possibly, multiple paths. 
\smallskip

Most of this section is about pathwise solutions of initial value problems of the form 
%the discussion about  stochastic viscosity solutions of 
%$x$-dependent fully nonlinear second-order pde of the form 
\begin{equation}\label{xdep1}
%\begin{cases}
du = F(D^2 u, Du,u,x) +H(Du,x) \cdot dB 
\ \text{ in }\ Q_\oo \quad 
u(\cdot, 0)=u_0,
%\end{cases}
\end{equation}
with only one path and, as always,  $F$   degenerate elliptic.
\smallskip

%It was already discussed that it is not yet understood how to define solutions when the Hamiltonian depends nonlinearly both on $Du$ and $u$. When $H$ depends linearly on the gradient, then a change of unknowns similar to the one described in Section~3 leads to problems with semilinear rough path dependence, which can be studied using the methodology of  Section~4. When $H$ depends linearly on $u$, then again a change of variables to remove the $u$-dependence from the Hamiltonian leads to an equation that can be studied following the results of this section.
%\smallskip

Extending the theory to equations with multiple rough time dependence had been an open problem until very recently, when  Lions and Souganidis \cite{lionssouganidis7} came up with a way to resolve the difficulty. A brief discussion about this appears at the end of this section. % and work out the details. The strategy is different. 
\smallskip

Finally, to study equations for nonsmooth Hamiltonians, it is necessary to modify the definition of the solution using now as test functions solutions of the doubled equations constructed for non smooth Hamiltonians.   The details  appear in \cite{lionssouganidisbook}.
\smallskip

%the approach of Section~XXXX and prove the exis

%Our knowledge about such equations is more limited than the 
%$x$-independent case.
%When $H$ is nonlinear, our theory only applies to ``regular'' Hamiltonians---the
%degree of regularity is specified later---and to paths 
%$W$ which are indeed Brownian---the stochastic nature of $W$ plays a 
%role here---and, very importantly, 
%to only one Brownian motion---we explain this difficulty in the course of 
%discussion.
%As was the case in the $x$-independent case, we do not know at the moment 
%how to handle Hamiltonians depending both on $Du$ and $u$.
%Finally, under the appropriate hypotheses the theory also applies to 
%$x$-dependent Brownian motions as long as the spacial dependence is 
%sufficiently smooth.

The strategy of  the proof of the comparison is similar to the one followed for spatially homogeneous Hamiltonians.
The pathwise  solutions are defined using as test functions smooth solutions of 
\begin{equation}\label{xdep2}
%\begin{cases}
du = H(Du,x)\cdot dB\ \text{ in }\ \R^d\times (t_0-h,t_0+h), \quad 
u(\cdot,t_0) = \phi,
%\end{cases}
\end{equation}
which under the appropriate assumptions on $H$ exist for each $t_0 >0$ and smooth 
$\phi$ in $(t_0-h,t_0+h)$ for some small $h$.
\smallskip

%
%Stochastic viscosity solutions of \eqref{xdep1} 
%are defined   exactly as in Definition~\ref{defn8.2},
%hence we do not repeat it here.
The aim in this section is to prove that  pathwise  solutions are well posed.
% that is,   exist, and satisfy a comparison principle. 
To avoid many technicalities, the discussion here  is restricted to  Hamilton-Jacobi initial value problems 
\begin{equation}\label{takis61}
du = H(Du,x) \cdot dB 
\  \text{ in } \ Q_T \quad 
u(\cdot, 0)=u_0.
\end{equation}
The general problem \eqref{xdep1} is studied using he arguments of this and the previous sections; some details can be found in  \cite{seeger}. 
%We succeed in showing this under some restrictive assumptions on $H$ 
%in terms of regularity.
\smallskip

Similarly to the spatially homogeneous case, the main technical issue is to control the length of the interval of existence of  smooth solutions of  the  doubled equation with quadratic initial datum and smooth approximations to $\zeta^\ep$ and $\xi^\eta$ of the path $B$, that is
\begin{equation}\label{xdepds}
\begin{cases}
dz = H(D_xz,x)\cdot d\zeta^\ep - H (-D_y z,y) \cdot d\xi^\eta 
\ \text{ in }\ \R^d\times\R^d\times (t_0-h, t_0+h),\\[1.2mm]
z(x,y,0) = \lambda |x-y|^2.
\end{cases}
\end{equation}
%where $B^\ep$ and $\widetilde B^\eta$ are regular approximations of $B$. 
%\smallskip

As already discussed earlier, the most basic estimate  is that   $h=\text{O}(\lambda^{-1})$, which, as  is explained below, is too small   to carry out the comparison proof. The challenge, therefore,  is to take advantage of the cancellations, due to  the special form of the initial datum as well as of  the doubled Hamiltonian, to obtain smooth solutions in a longer time interval. 
\smallskip

Since the smooth solutions to \eqref{xdepds} are constructed by the method of characteristics, the technical issue is to control the length of the interval of invertibility of the characteristics.  This can be done by estimating the interval of time in which the Jacobian does not vanish. It is here  that using a single path helps, because, after a change of time, the problem reduces to studying the analogous question for homogeneous in time odes.  % without rough ti
\smallskip

To further simplify the presentation, the ``rough'' problem discussed in the sequel  is not \eqref{xdepds} but rather the doubled equation with the rough path,  that is 
\begin{equation}\label{xdepd3}
%\begin{cases}
dw =H(D_xw,x)\cdot dB -H(-D_y w,y)\cdot dB\ 
\text{ in }\ \R^d\times\R^d\times (t_0-h,t_0+h) \quad 
w(x,y,0) = \lambda' |x-y|^2. %\frac{\lambda}{2} |x-y|^2\ .
%\end{cases}
\end{equation}
In what follows, to avoid cumbersome expressions, %$ \frac{\lambda}{2} |x-y|^2$ is written as  $\lambda' |x-y|^2$
$\lambda =2 \lambda'$ 
\smallskip

The short time smooth solutions of  \eqref{xdepd3} are given by  $w(x,y,t)=U(x,y,B(t)-B(t_0))$, where $U$ is the short time smooth solutions 
to the ``non-rough'' doubled initial value problem
\begin{equation}\label{takis63}
U_t = H(D_xU,x) - H (-D_y U,y) 
\ \text{ in }\ \R^d\times\R^d\times (-T^*,T^*) \ \quad 
U(x,y,0) =  \lambda' |x-y|^2,
\end{equation}
and  $T^*>0$ and $h$ are such that $\sup_{s \in (t_0-h,t_0+h)} |B(s)-B(t_0)| \leq T^*.$
%note that this is the place where assuming that there is only one rough path enters. The argument for multiple paths is rather different. 
\smallskip

The smooth solutions of  \eqref{takis63} are constructed  by inverting the map  $(x,y)\to (X(x,y,t), Y(x,y,t))$ of  the  corresponding system of characteristics, that is 
\begin{equation}\label{takis64}
\begin{cases}
\dot X=- D_p H(P,X) \quad  X(x,y,0)=x \quad \dot Y=-D_q H(Q, Y) \quad Y(x,y,0)=y,\\[1mm]
\dot P = D_x H (P,X)  \quad   \dot Q = D_y H (Q, Y) \quad P(x,y,0)=Q(x,y, 0)=\lambda (x-y),\\[1mm]
\dot U= H(P,X) - \langle D_p H(P,X),P\rangle -H(Q,Y) + \langle D_q H(Q,Y), Q\rangle \ \  U(x,y,0)=\lambda'|x-y|^2.
\end{cases}
\end{equation}
A crude estimate, which does not take into account the special form of the system and the initial data, gives that  the map $(x,y)\mapsto (X(x,y,t),Y(x,y,t))$
is invertible at least in a time interval of length $\text{O}(\lambda^{-1})$ % $[-c\lambda^{-1},c\lambda^{-1}]$,
with the constant depending  on $\|H\|_{C^2}$.
\smallskip

This implies that the characteristics of \eqref{xdepd3} are invertible as 
% It follows that the map $(x,y)\to (X(x,y,t),Y(x,y,t))$ where 
%\begin{equation*}
%X(x,t) = X(x,y,W(t)0\quad\text{and}\quad Y=Y(x,y,W(t))
%\end{equation*}
%is invertible as 
long as 
$$\sup_{s \in (t_0-h,t_0+h)} |B(s)-B(t_0)|\leqq \text{O}(\lambda^{-1}).$$
It turns out, as it is shown below,  that this interval is not long enough to yield a comparison for the pathwise solutions. Taking, however,  advantage of the special structure of \eqref{takis63} and  \eqref{takis64} and under suitable assumptions on $H$ and its derivatives, it is possible to improve the estimate of  the time interval.
% to $O(\lambda^{-2})$ which will 
%again not be enough to prove comparison.
%With additional regularity and some decay assumptions on $H$ it is also 
%possible to come up with interval of existence of order~1.
\smallskip

The discussion next aims to explain the need of  intervals of invertibility that are longer  than $\text{O}(\lambda^{-1})$, and serves as a blueprint for the strategy of the actual proof. 
%invertibility for the characteristics, next we go through the arguments needed
%to prove uniqueness of stochastic viscosity solutions of \eqref{xdep2}, which
%are similar to the ones used in the $x$-independent case.
%\medskip
%%%%%%%%%%%%%%%%%%
%\subsection*{Sketch of a (potential) proof of comparison for 
%\eqref{xdep2} for separated Hamiltonians}$\quad$
%%%%%%%%%%%%%%%%%%
\smallskip

Assume  that $u$ and $v$ are respectively a subsolution and a supersolution 
of \eqref{xdep2}.  As in the $x$-independent case, it is  assumed that,  
for some $\alpha > 0$ and $\lambda > 0$,
$(x_0,y_0,t_0)$ with $t_0> 0$ is a maximum point of 
$$u(x,t) - v(y,t) - \lambda'| x-y|^2 - \alpha t.$$
Then, for $h> 0$ and all $x,y$,
$$u(x,t_0-h) - v(y,t_0-h) \leqq \lambda' |x-y|^2 - \alpha h 
+ u (x_0,t_0) - v(y_0,t_0) -\lambda' |x_0-y_0|^2\ .$$
Since 
%As in the $x$-independent case it is immediate that 
$w(x,y,t) = u(x,t) - v(y,t)$ solves the doubled equation  \eqref{xdep2}, to obtain the comparison 
it is enough to compare $w$ with the small time smooth solution $z$ to \eqref{xdepd3} starting at $t_0-h$.
\smallskip

It follows that 
%Since, as explained earlier, it is always possible to compare viscosity 
%solutions with smooth solutions, the previous inequality for some 
%It follows that 
%$h= h_\lambda (\omega)$ yields 
$$u(x_0,t_0) - v(y_0,t_0) \leqq w(x_0,y_0,t_0) + u(x_0,t_0) 
- v(y_0,t_0) - w(x_0,y_0,t_0-h)-\alpha h,$$
%\lambda' |x_0-t_0|^2 - \alpha h\ ,$$
%where $U_\lambda$ is the local time smooth solution of \eqref{xdepd3} which 
%exists by the method of characteristics.  
and, hence, %  we must have 
%$$\alpha h \leqq z (x_0,y_0,t_0) -z(x_0,y_0,t_0-h),$$ 
%%\lambda' |x_0-y_0|^2,$$
%and, finally,
$$\alpha \leq \frac{w (x_0,y_0,t_0) -w(x_0,y_0,t_0-h)}{h}.$$ 
Recall that $h$ depends on $\lambda$ and, to conclude, this dependence must be such that 
$$\limsup_{\lambda\to \infty}\frac{w (x_0,y_0,t_0) -w(x_0,y_0,t_0-h)}{h}\leq 0.$$
On the other hand, it will be shown that, if $z$ is a smooth solution to \eqref{xdepd3}, then
$$
w (x_0,y_0,t_0) -w(x_0,y_0,t_0-h)\lesssim  \sup_{s \in (t_0-h,t_0+h)} |B(s)-B(t_0)| h^{-1}\lambda^{-\frac 12}.
$$
Combining the last two statements implies  that, to get a contradiction,  $h=h(\lambda)$ must be such that 
\begin{equation}\label{takis66}
\limsup _{\lambda \to \infty} \sup_{s \in (t_0-h,t_0+h)} |B(s)-B(t_0)| h^{-1}\lambda^{-\frac 12}=0.
\end{equation}
%
%
%
%
%
% by obtaining a contradiction it suffice to show that 
%
%
%Recall The previous discussion
%On the other hand, as long as $z$ is smooth, it follows that 
%$$U_\lambda (x,y,t) = \widetilde S(W(t_0) - W(t_0-h)) (\lambda 
%|\cdot - \cdot|^2) (x,t)\ ,$$
%where $\widetilde S$ is the solution of \eqref{xdepd3} with $W \equiv 1$.
%
%Of course the above estimate will hold as long as \eqref{xdepd3} has a 
%local time smooth solution.
%In view of our previous discussion, is some $\omega (\lambda^{-1})$?
%In general we must have $W(t_0) - W(t_0,t) = O(\lambda^{-1})$. 
%
%If the interval of existence of smooth solutions to \eqref{takis63} is $\text{O}(\lambda^{-1})$ and, for example, $B \in C^{0,1/2}([0,\infty)),$ then $h=\text{O}(\lambda^{-2})$ and, as a result, 
%$$\sup_{s \in (t_0-h,t_0+h)} |B(s)-B(t_0)| h^{-1}\lambda^{-\frac 12}\lesssim \lambda^{\frac 12}$$
%which, of course, does not do the job. 
%\smallskip
%
%It is clear from the above discussion that, to have \eqref{takis66} for $B \in C^{0,1/2}([0,\infty)),$ it is necessary to show that   $\lim_{\lambda \to \infty}h\lambda^{-1}=0$, and  this requires considerable effort.
%some more care.
%is not possible to get \eqref{takis66}. 
%\smallskip
The next argument indicates  that there is indeed a problem 
if the smooth solutions of the ``deterministic'' doubled problem exist only for times of order $O(\lambda^{-1})$. 
%$e only have an interval of existence  of smooth $U_\lambda$'s given by 
%the characteristic in a interval of length $O(\lambda^{-1})$.
\smallskip

Indeed in this case, the proof of the comparison argument outlined above, yields  
$$\alpha h\leqq o(1) |B (t_0) - B(t_0-h)|\lambda^{-1/2},$$
and, if $B\in C^{0,\beta}([0,\infty))$, it follows that 
%As before, if $W$ is H\"older continuous with exponent $\beta$, 
 $h^\beta \approx \lambda^{-1},$
and the above inequality yields 
$\alpha \leqq o(1) h^{\frac{3\beta}2 -1}$
%% \lambda^{\frac1{\beta} - \frac32}\ ,$$
in which case it is not possible to obtain a contradiction, if 
$\beta < 2/3$, which, of course, is the case for Brownian paths. 
%unless there is an explicit knowledge of the $O(1)$. 
\smallskip

It appears, at least for the moment formally,  
that for this case the Brownian case ``optimal'' interval of 
existence is $O(\lambda^{-1/2})$. 
Indeed if this is the case then 
we must have $|B(t_0) - B(t_0-h)| \approx \lambda^{-\frac12}$,
and, hence, $h\approx \lambda^{-1}$. 
This leads to 
%% $$\alpha \frac1{\lambda} \leqq \frac{o(1)}{\lambda}$$
$\alpha \leqq o(1)$ and, hence, a contradiction.
\subsection*{Improvement of  the interval of existence of smooth solutions} The problem is to find longer than $\text{O}(\lambda^{-1})$ intervals of existence of smooth solution of the doubled deterministic Hamilton-Jacobi equation \eqref{takis63}. 
\smallskip

%\begin{equation}\label{detdHJ}
%\begin{cases}
%U_t = H (D_x U,x) - H(-D_y U,y) \ \text{ in }\ \R^d\times (-T,T),\\[1mm]
%U (x,y,0) = \lambda' |x-y|^2.
%\end{cases}
%\end{equation}
Two general sets of conditions will be modeled by two particular classes of Hamiltonians, namely separated and linear $H$'s.
\smallskip
 
%We will introduce several different sets of assumptions. 
To give the reader a flavor of the type of arguments that will be 
involved, it is convenient to begin with ``separated'' Hamiltonians 
of the form
\beq\label{paris50}
H(p,x) = H(p) + F(x),
\end{equation}
in which case the doubled equation and its characteristics are 
\begin{equation}\label{detds}
%\begin{cases} 
U_t = H (D_xU) - H(-D_y U) + F(x) - F(y)\ \text{ in }\ \R^d\times \R^d \times (-T,T) \quad U(x,y,0) = \lambda' |x-y|^2,
%,
%\\\[1mm]
%U(x,y,0) = \lambda' |x-y|^2,
%\end{cases}
\end{equation}
and %the characteristics are 
\begin{equation}\label{detchds}
\begin{cases} \dot X = -D H(P) \ \ X(0) =x   \quad \dot Y = - D H(Q) \ \ Y(0) = y ,\\[1mm]
\dot P = D F(X) \quad \dot Q = D F(Y) \ \  P(0) = Q(0) =\lambda (x-y),\\[1mm]
\dot U = H(P) - \langle D H(P),P\rangle - H(Q) + \langle D H(Q),Q\rangle + F(X) - F(Y) \ \ U(0) = \lambda' |x-y|^2.%\\[1mm]
%X(0) =x \quad Y(0) = y \quad P(0) = Q(0) =\lambda (x-y) \quad  U(0) = \lambda' |x-y|^2.
%P(0) = \lambda (x-y) \end{cases} 
%\quad\text{and}\qquad
%\begin{cases} \dot Y = - D_q H(Q)\\
%\dot Q = D_y F(Y)\\
%Y(0) = y\\
%Q(0) = \lambda (x-y) 
\end{cases}
\end{equation}
%and 
%\begin{equation}\label{detchds2}
%\begin{cases}
%\dot U = H(P) - (D_p H(P),P) - H(Q) + (D_q H(Q),Q) 
%+ F(X) - F(Y)\\
%\noalign{\vskip6pt}
%U(0) = \lambda' |x-y|^2\ .
%\end{cases}
%\end{equation}
%As before, to simplify the notation since everything becomes more symmetric, 
%$Q= -D_y u(Y(t),t)$.
%%\smallskip
%\vskip.05in
\vskip.075in

Let $J(t)$ denote the Jacobian of the map 
$(x,y)\mapsto (X(x,y,t), Y(x,y,t))$ at time t.  In what follows,  
to avoid the rather cumbersome notation involving determinants, all the calculations below are presented 
for $d=1$, that is $x,y\in \R$. % the reader, of course, can easily see how to extend everything to higher dimensions. 
\smallskip

It follows that 
$$J = \frac{\partial X}{\partial x} \ \frac{\partial Y}{\partial y}
- \frac{\partial X}{\partial y}\ \frac{\partial Y}{\partial x} \quad \text{and} \quad J(0)=1.$$
The most direct way to find an estimate for the time of existence of smooth solutions is, for example, to obtain a bound for the first time $t_\lambda$ such that $J(t_\lambda)=\frac{1}{2}$, and, for this, it is convenient to calculate and estimate the derivatives of $J$ with respect to time at $t=0$. 
\smallskip

Hence, it is necessary to derive the odes  satisfied  
by $\dfrac{\partial X}{\partial x},\dfrac{\partial X}{\partial y},   
\dfrac{\partial Y}{\partial x}$ and $\dfrac{\partial Y}{\partial y}$. 
Writing $\dfrac{\partial X}{\partial\alpha},\dfrac{\partial Y}{\partial\alpha},
\dfrac{\partial P}{\partial\alpha}$, and $\dfrac{\partial Q}{\partial\alpha}$ 
with $\alpha = x$ or $y$, differentiating \eqref{detchds} and omitting the subscripts for the derivatives of $H$ and $F$ yields the systems % we find the ode 
\begin{equation*} 
\begin{cases}
\ds \frac{\Circfrown{\partial X}}{\partial\alpha} 
= - D^2H(P)\dfrac{\partial P}{\partial\alpha}, \quad \begin{cases} \dfrac{\Circfrown{\partial X}}{\partial x} (x,y,0)=1, \\ \dfrac{\Circfrown{\partial X}}{\partial y} (x,y,0)=0,\end{cases} \\
\noalign{\vskip6pt}
\ds \dfrac{\Circfrown{\partial P}}{\partial \alpha}
= D^2F(X)\dfrac{\partial X}{\partial \alpha}, \quad \begin{cases} \dfrac{\Circfrown{\partial P}}{\partial x} (x,y,0)= \lambda,\\  \dfrac{\Circfrown{\partial P}}{\partial y} (x,y,0) = - \lambda\ , \end{cases}
%\noalign{\vskip6pt}
%\ds \frac{\Circfrown{\partial X}}{\partial x} (x,y,0)=1,\ \ 
%\frac{\Circfrown{\partial P}}{\partial x} (x,y,0)= \lambda,\\
%\noalign{\vskip6pt}
%\ds\frac{\Circfrown{\partial X}}{\partial y} (x,y,0)=0,\ \ 
%\frac{\Circfrown{\partial P}}{\partial y} (x,y,0) = - \lambda\ ,
\end{cases}
\hskip-.2in  \text{ and } \  
\begin{cases}
\ds\frac{\Circfrown{\partial Y}}{\partial\alpha} 
= - D^2H(Q)\dfrac{\partial Q}{\partial \alpha}, \quad \begin{cases}\dfrac{\Circfrown{\partial Y}}{\partial x}(x,y,0) =0,\\
\dfrac{\Circfrown{\partial Y}}{\partial Y}(x,y,0) =1,\ \end{cases}\\
\noalign{\vskip6pt}
\ds\frac{\Circfrown{\partial Q}}{\partial \alpha} 
= D^2F(Y) \dfrac{\partial Y}{\partial\alpha}, \begin{cases} \dfrac{\Circfrown{\partial Q}}{\partial x}(x,y,0) = \lambda,\\ \dfrac{\Circfrown{\partial Q}}{\partial y}(x,y,0) = - \lambda.\end{cases}
%\noalign{\vskip6pt}
%\ds\frac{\Circfrown{\partial Y}}{\partial x}(x,y,0) =0,\ \ 
%\frac{\Circfrown{\partial Q}}{\partial x}(x,y,0) = \lambda,\\
%\noalign{\vskip6pt}
%\ds\frac{\Circfrown{\partial Y}}{\partial Y}(x,y,0) =1,\ \ 
%\frac{\Circfrown{\partial Q}}{\partial y}(x,y,0) = - \lambda.
\end{cases}
\end{equation*}
%to keep the notation simple the dependence of $H$ and $F$ and their derivatives on $P,Q,X$ and 

\begin{prop}\label{intex1}
Assume that $DH$, $DF$, $D^2F$, $D^2H$, $|D^3H|(1+|p|)$ and $D^4H$ are bounded. 
If $t_\lambda$ is the first time  that $J(t_\lambda)=1/2$, 
then, for some uniform constant $c> 0$ which depends on the bounds on  $H,F$ and
their derivatives, and for all $x,y \in \R^d$, 
%%$$t_\lambda \geqq \min [(2c)^{-1/2}, (2c)^{-1/3} \lambda^{-1/3}]\ .$$
$$t_\lambda \ge c \min (1, \lambda^{-1/3})\ .$$
\end{prop}
\begin{proof} 
%To keep the notation simpler throughout the proof, we write $J$ instead of 
%$J_\lambda$.
%
Straightforward calculations that take advantage of the separated form of the Hamiltonian  yield
%calculating the first and second derivatives of $J$. 
%This is the place where the separated form of the Hamiltonian leads to simpler 
%expressions.
%We have: 
$$\dot J = 
- D^2H (P) \left(
\frac{\partial P}{\partial x} \frac{\partial Y}{\partial y} 
- \frac{\partial P}{\partial y} \frac{\partial Y}{\partial x}\right) 
- D^2H(Q) \left(
\frac{\partial X}{\partial x} \frac{\partial Q}{\partial y} 
- \frac{\partial X}{\partial y} \frac{\partial Q}{\partial x}\right) 
$$
and 
\begin{equation*}
\begin{split}
\ddot J &= - \left(D^2H(P) D^2F(X) + D^2H(Q) D^2F(Y)\right) J 
+ 2D^2H(P) D^2H(Q) \left(\frac{\partial P}{\partial x}\frac{\partial Q}{\partial y} 
- \frac{\partial P}{\partial y} \frac{\partial Q}{\partial x}\right)\\
\noalign{\vskip6pt}
&\qquad 
- D^3H(P) DF(X) \left(\frac{\partial P}{\partial x} \frac{\partial Y}{\partial y}
- \frac{\partial P}{\partial y} \frac{\partial Y}{\partial x} \right) 
- D^3H(Q) DF(Y)\left(\frac{\partial X}{\partial x} \frac{\partial Q}{\partial y} 
- \frac{\partial X}{\partial y}\frac{\partial Q}{\partial x}\right)\ .
\end{split}
\end{equation*}
To simplify the expressions for $\dot J$ and $\ddot J$, it is convenient to write $\dfrac{\partial X}{\partial\alpha},\dfrac{\partial Y}{\partial\alpha},
\dfrac{\partial P}{\partial\alpha}$and $\dfrac{\partial Q}{\partial\alpha}$ in terms of the solutions $(\eta_1,\psi_1,\eta_2,\psi_2)$ of the linearized system %Consider next the solutions $(\xi_1,\phi_1,\xi_2,\phi_2)$ and 
%$(\eta_1,\psi_1,\eta_2,\psi_2)$ of the simplified system of ode 
\begin{equation*}
\begin{cases}
\dot{\xi}_1 = - D^2H(P) \phi_1,\ \ \xi_1(0) =1,\\
\noalign{\vskip6pt}
\dot{\phi}_1 = D^2F(X) \xi_1,\ \ \phi_1 (0) =0,\\
\noalign{\vskip6pt}
\dot{\xi}_2 = - D^2H(P) \phi_2,\ \ \xi_2 (0)=0,\\
\noalign{\vskip6pt}
\dot{\phi}_2 = D^2F(X)\xi_2,\ \ \phi_2 (0)=1,
\end{cases}\quad \text{and}\quad 
\begin{cases}
\dot{\eta}_1 = -D^2H(Q) \psi_1,\ \ \eta_1 (0)=1,\\
\noalign{\vskip6pt}
\dot{\psi}_1 = D^2F(Y) \eta_1,\ \ \psi_1 (0)=0,\\
\noalign{\vskip6pt}
\dot{\eta}_2 = -D^2H (Q)\psi_2,\ \ \eta_2 (0)=0,\\
\noalign{\vskip6pt}
\dot{\psi}_2 = D^2F(Y)\eta_2,\ \  \psi_2(0)=1,
\end{cases}
\end{equation*}
which are bounded in $[0,1]$ and satisfy 
\begin{equation*}
\begin{cases} 
\xi_1(t) =1 + O(1)t,\quad \xi_2(t)=O(1)t,\\
\noalign{\vskip6pt}
\phi_1(t)=O(1) t,\quad \phi_2(t)=1+O(1)t,\end{cases}
\quad\text{and}\quad 
\begin{cases}
\eta_1(t) = 1+O(1)t,\quad \eta_2 =O(1)t,\\
\noalign{\vskip6pt}
\psi_1(t) = O(1)t,\quad \psi_2 = 1+ O(1)t,
\end{cases}
\end{equation*}
where $O(1)$ denotes different quantities for each functions which 
are uniformly bounded in $[0,1]$; % and  are diffrom equation to equation. 
note  that the assumption that $D^2 H$ and $D^2 F$ are bounded is used here. 
\smallskip

A direct substitution yields 
\begin{equation*}
\begin{cases}
\frac{\partial X}{\partial x} = \xi_1 + \lambda\xi_2  \quad   \frac{\partial X}{\partial y} = - \lambda \xi_2,\\[3.5mm]
%\noalign{\vskip6pt}
\frac{\partial P}{\partial x} = \phi_1 + \lambda\phi_2 \quad   \frac{\partial P}{\partial y} = - \lambda \phi_2,
\end{cases}
%\end{equation*}
\quad \text{and} \quad 
%\quad\text{and}\quad
%\begin{cases}
%\ds \frac{\partial X}{\partial y} = - \lambda \xi_2,\\
%\noalign{\vskip6pt}
%\ds \frac{\partial P}{\partial y} = - \lambda \phi_2,\end{cases}
%\end{equation*}
%and 
%\begin{equation*}
\begin{cases}
\ \frac{\partial Y}{\partial x} = \lambda\eta_2 \quad \frac{\partial Y}{\partial y} = \eta_1 - \lambda\eta_2,\\[3.5mm]
 \frac{\partial Q}{\partial x} = \lambda\psi_2  \quad \frac{\partial Q}{\partial y} = \psi_1 - \lambda\psi_2.
\end{cases}
%\quad\text{and}\quad
%\begin{cases}
%\ds\frac{\partial Y}{\partial y} = \eta_1 - \lambda\eta_2,\\
%\noalign{\vskip6pt}
%\ds\frac{\partial Q}{\partial y} = \psi_1 - \lambda\psi_2.
%\end{cases}
\end{equation*}
%\end{proof}
Using the observations above gives
\begin{equation*}
\begin{split}
\frac{\partial P}{\partial x} \frac{\partial P}{\partial y} 
- \frac{\partial P}{\partial y} \frac{\partial Q}{\partial x} 
& = (\phi_1 + \lambda\phi_2) (\psi_1 - \lambda\psi_2) 
- (-\lambda\phi_2)\lambda\psi_2\\
%\noalign{\vskip6pt}
& = \phi_1 \psi_1 + 2\lambda (\phi_2 \psi_1 - \phi_1 \psi_2) 
= O (1) (1+2\lambda t),
\end{split}
\end{equation*}
since 
$$\phi_1\psi_1 =O(1)\quad\text{and}\quad  
\phi_2 \psi_1 - \phi_1\psi_2 = (1+O(1)t) O(1) t 
- O(1) t(1+ O(1)t) = O(1)t\ .$$

Similarly, since 
\begin{gather*}
\phi_2 \eta_1 - \phi_1\eta_2 = (1+O(1)t) (1+O(1)t) - O(1) tO(1)t = 1+ O(1)t
\quad\text{and}\\
\noalign{\vskip6pt}
\xi_2\psi_1 - \xi_1 \psi_2 = O(1) tO(1)t - (1+O(1)t)(1+O(1)t) = O(1)t-1, 
\end{gather*}
it follows that 
\begin{equation*}
\begin{split}
\Partial{P}{x} \Partial{Y}{y} - \Partial{P}{y}\Partial{Y}{x}
& = (\phi_1 + \lambda\phi_2)(\eta_1 - \lambda\eta_2) - (-\lambda\phi_2)
\lambda\eta_2\\
& = \phi_1 \xi_1 + \lambda (\phi_2 \eta_1 - \phi_1\eta_2 ) 
= O(1) (1+\lambda t) + \lambda ,
\end{split}
\end{equation*} 
and 
\begin{equation*}
\begin{split}
\Partial{X}x \Partial{Q}y - \Partial{X}y \Partial{Q}x 
& = (\xi_1 + \lambda\xi_2)(\psi_1 - \psi_2) - (-\lambda\xi_2)\lambda \psi_2\\
& = \xi_1 \psi_1 + \lambda (\xi_2 \psi_1 - \xi_1\psi_2) 
= O(1) (1+\lambda t) -\lambda\ .
\end{split}
\end{equation*} 
Inserting all the above in the expression for $\ddot J$ yields 
$$\ddot J = O(1) J + O(1) (1+\lambda t) 
+ \lambda (D^3H (Q) DF(Y) - D^3H(P)DF(X)).$$
Set 
$$A:= (D^3H(Q) - D^3H(P)) DF(Y)
\quad\text{and}\quad 
D:= D^3H(P) (DF(Y) - DF(X)).$$
It is immediate that 
$$\lambda |A| \leqq \lambda O (\|DF\|_\infty |Q-P|) = 
\lambda O(1)t,$$
with the last estimate following from the observation that 
\begin{equation*}
\begin{split}
(P-Q) (t) & = \lambda (x-y) + \int_0^t DF(X(s))ds 
- (\lambda  (x-y)t \int_0^t DF(Y(s))ds\\
\noalign{\vskip6pt}
& = \int_0^t (DF(X(s)) - DF(Y(s)) ds = O(1) t.
\end{split}
\end{equation*}
As far as $D$ is concerned, observe that 
$$|D| \leqq \|D^3H \|_\infty |X-Y| \lesssim \frac{|X-Y|}{1+|P|},$$
and recall that 
$$|X-Y| \leqq |x-y| + O(1)t \ \ \text{and} \ \ |P| = |\lambda (x-y) + O(1)t|.$$
%and
%$$|P| = |\lambda (x-y) + O(1)t|.$$
%
Hence, 
$$\lambda |D| \lesssim \left[\frac{\lambda |x-y|}{1+|\lambda (x-y)+O(1)t|}
+ \frac{\lambda O(1)t}{1+|P|}\right]
\lesssim \left[\frac{|P(0)|}{1+|P(0) + O(1)t|} 
+ \frac{\lambda O(1)t}{1+|P|}\right];$$
the second term in the bound above comes  from $\lambda O(1)t$, while an additional argument is needed for the first. 
%it is necessary to argue a bit more.
\smallskip

Choose $t\leq t_1$ so that the $O(1)t$ term in $P$ is such that 
$|O(1)t| \leqq \frac12.$ 
%
%%Then
%%$$\frac{\lambda |x-y|}{1+ |\lambda (x-y) + O(1)t|}
%%= \frac{|P(0)|}{1+|P(0) + O(1)t|}\ .$$
%
If $|P(0)| \leqq 1$, then 
$$\frac{|P(0)|}{1+|P(0) + O(1)t|} \leqq 1$$
while, if $|P(0)| >1$, 
$$1+|P(0) + O(1)t| 
\geqq 1+ |P(0)| - |O(1)t| 
\geqq |P(0)| + \frac12$$
and 
$$\frac{|P(0)|}{1+|P(t)|} \leqq 
\frac{|P(0)|}{\frac12 + |P(0)|} \leqq 1\ .$$
Combining the estimates on $\lambda A$ and $\lambda D$ gives  
$$\ddot J = O(1) J + O(1) \lambda t + O(1).$$
Iit is also immediate that
$$J (0) =1\quad\text{ and }\quad \dot J (0)=0;$$
this is another place where the separated form of the Hamiltonian and the symmetric form of the test function play a  role. 
%The first equality follows from the definition while the 
%second is a consequence of the separated form of the Hamiltonian.  
\smallskip

It follows  there exists 
%Let $t_\lambda$ be the first time that 
%$$J_\lambda (t_\lambda) = \frac12.$$
%
%The previous calculations In view of the previous calculations we have, for some 
$s_\lambda \in  (0,t_\lambda)$  such that 
$$\frac12 = 1+\frac12 t_\lambda^2 \ddot J(s_\lambda)$$
and, hence, 
$$|t_\lambda^2 \ddot J (s_\lambda)| = 1,$$
which implies 
$$1\lesssim t_\lambda^2 (1+ \lambda t_\lambda).$$
%The previous Therefore, for some uniform $C> 0$,  
%$$Ct_\lambda^2 (1+ \lambda t_\lambda) \geqq 1$$
It follows that 
$$1\lesssim \lambda t_\lambda^3 + t_\lambda^2,$$
and the claim is proved. 
% hence, for some other uniform $c> 0$,  
%$$\lambda t_\lambda^3 + t_\lambda^2 \geqq c,$$
%which yields the claim.

\end{proof}

Having established a longer  than $\text{O}(\lambda^{-1})$ interval of existence for the 
solution $U_\lambda$ of \eqref{detds}, it is now possible to obtain the following 
comparison result for pathwise solutions to Hamilton-Jacobi equations with 
separated Hamiltonians.
\begin{thm}\label{comp1}
Let $u\in BUC(\overline Q_T)$ and $v \in BUC(\overline Q_T)$ be respectively a  subsolution and a  supersolution of \eqref{xdep2}
in $Q_T$ with $H$ as in \eqref{paris50}, that is, of separated form,  satisfying  
the assumptions of Proposition~\ref{intex1}. 
Moreover, assume that $B\in C^{0,\beta} ([0,\infty])$ with $\beta\geq 2/5$.
Then, for all $t\in [0,T]$, 
$$\sup_{x\in\R^d} (u(x,t) - v(x,t)) 
\leqq \sup_{x\in \R^d} (u(\cdot,0) - v(\cdot,0))\ .$$
\end{thm}
The following lemma, which is stated without a proof since it is rather classical,  will be used in the proof of Theorem~\ref{comp1}. % Its proof is presented after th
\begin{lem}\label{sharp}
Assume that  $H\in C(\R^d)$ and $F\in C^{0,1} (\R^d)$ and let $U_\lambda$ be the viscosity solution of  the doubled equation $w_t=H(D_xw)+F(x)-H(-D_yw,y) -F(y) \ \text{in} \ Q_T$  
%for the Hamiltonian $H+F$ 
with initial datum $\lambda |x-y|^2$.  %\eqref{takis63} with $H$ as above. 
Then, for all $x,y \in \R^d$ and $t\in [0,T]$, 
$$|U_\lambda (x,y,t) - \lambda |x-y|^2 | \leqq t \|DF\|\, |x-y|\ .$$
\end{lem}
%
%\begin{proof}
%A straightforward calculation  shows  that $\lambda' |x-y|^2 \pm t \|DF\|\, |x-y|$ are, 
%respectively, super- and sub-solutions of 
%$$w_t=H(D_xw)+F(x)-H(-D_yw,y) -F(y) \ \text{in} \ Q_T \ \  \text{and} \ \ w(x,y,0)=\lambda' |x-y|^2,$$
%%\eqref{detds}.
%and then the  comparison for viscosity solutions implies 

\begin{proof}[The proof of Theorem~\ref{comp1}]
Assume that $(x_0,y_0,t_0)$ with $t_0 > 0$ is a maximum point of 
$u(x,t) - v(y,t) - \lambda |x-y|^2 - \alpha t$.  Repeating the arguments at the end of the previous subsection and using Lemma~\ref{sharp} yields 
%leading to Proposition~\ref{intex1} we 
%assume taht $(x_0,y_0,t_0)$ with $t_0> 0$ is a maximum point of 
%$$u(x,t) - v(y,t) - \lambda |x-y|^2 - \alpha t\ .$$
%
%Then, as in previous discussion, we are led to the inequality 
\begin{equation}\label{en}
\alpha h \leqq \|DF\|_\infty |x_0-y_0|\, |B(t_0) - B(t_0 -h)|\ . 
\end{equation}
Recall that, in view of Proposition~\ref{intex1}, the above inequality 
 holds as long as 
$$|B(t_0) - B(t_0-h) | \lesssim \lambda^{-1/3}\ .$$
Since $B\in C^{0,\beta} ([0,\infty])$, $h=h(\lambda)$can be chosen so that 
%leads to the following relationship between $\lambda$ and $h$:
$$\lambda^{-1} \approx h^{-3\beta}.$$
%between $\lambda$ and $h$.
Moreover,  $(x_0,y_0) \in \R^d\times\R^d$ being  a maximum of 
$u(x,t_0) - v(y,t_0) - \lambda |x-y|^2$
yields
$\lambda |x_0-y_0|^2 \leqq \max (\|u\|, \|v\|)$
and, if  $\omega$ is the modulus of continuity of $u$, 
$\lambda |x_0-y_0|^2 \leqq \omega (\lambda^{-1/2} \max (\|u\|,\|v\|)^{1/2})
= \text{O}(1),$
%\smallskip
and, hence,
%where $\omega$ is the modulus of continuity of $u$, and, 
$|x_0-y_0| \lesssim \lambda^{-1/2}.$
\smallskip

Inserting all the observations  above in \eqref{en}, gives
%we must have, for some $C> 0$,
$\alpha h  \lesssim  % h^{\alpha +\frac{3\alpha}2} 
%= C\ O  (1)
\text{o}(1) h^{\frac{5\beta}2}\ ,$
%\smallskip
and, thus, %which yields 
$\alpha \lesssim  \text{o}(1) h^{\frac{5\beta}2 -1},$
which leads to a contradiction as $\lambda\to\infty$.
% and, therefore, 
%$h\to0$ provided $\alpha \rangle2/5$.

\end{proof}
Note  that it is possible  to assume less on $B$  in Theorem~\ref{comp1}, if more information is available about  the modulus of continuity of either $u$ or $v$.
%regularity is an explicit modulus assuming less on $W$ if we knew more explicitly 
%the modulus of continuity of either $u$ or $v$.
\smallskip

For Hamiltonians that are  not of separated form, the situation is more complicated. 
Indeed the ``canonical'' assumption on $H$ for the deterministic theory is 
that, for some modulus $\omega_H$ and all $x,y,p\in\R^d$, 
\begin{equation}\label{Hass1} 
|H(p,x) - H(p,y)| 
\leqq \omega_H (|x-y| (1+ |p|)).
\end{equation}
On the other hand, the proof of the  comparison yields %, it was shown that at the maximum point is alwasy known that  we found that we always have 
$$\lambda |x_0-y_0|^2 \leqq 2 \max (\|u\|,\|v\|),\quad \text{and} \quad 
\lambda |x_0-y_0|^2\leqq \max(\omega_u (|x_0-y_0|), \omega_v (|x_0-y_0|))\ ,$$
and, hence, 
$$|x_0-y_0|^2 \leqq \lambda^{-1} \max (\omega_u,\omega_v) 
(2(\lambda^{-1} \max (\|u\|, \|v\|))^{1/2} )\ .$$

If either $u$  or $v$ is Lipschitz continuous, then the above
estimate can be improved to 
$$|x_0-y_0| \leqq \min (\|Du\|,\|Dv\|) \lambda^{-1}\ .$$
The next technical result replaces Lemma~\ref{sharp}. Its proof is again classical and it is omitted.
\begin{lem}\label{sharp2} 
Assume that $H$ satisfies \eqref{Hass1} with $\omega_H(r) = Lr.$ 
Let $U_\lambda$ be the viscosity solution of  the doubled initial value problem \eqref{takis63}. 
%equation $w_t=H(D_xw)+F(x)-H(-D_yw) -F(y) \ \text{in} \ Q_T$ with initial datum $\lambda' |x-y|^2$.  %\eqref{takis63} with $H$ as above. 
Then there exists  $C>0$ depending on $L$ such that, for all $x,y \in \R^d$ and $t\in [0,T]$, 
%$$|U_\lambda (x,y,t) - \lambda |x-y|^2 | \leqq t \|DF\|\, |x-y|\ .$$
%There exists $C> 0$ depending on $L$, such that 
\begin{equation}\label{e2} 
|U_\lambda (x,y,t) - \lambda' e^{Ct} |x-y|^2| 
\leqq (e^{Ct} -1)|x-y|\ . 
\end{equation}
If, in addition  $|x-y| \lesssim \lambda^{-1}$, then %\eqref{e2} reduces to 
$$|U_\lambda (x,y,t) - \lambda' e^{Ct} |x-y|^2 | 
\lesssim  t\lambda^{-1}\ .$$
\end{lem}
A discussion follows about how to  
``increase'' the length of the interval of existence of solutions given by the method of characteristics for $H$'s which are not  separated.  
%\smallskip
To keep the notation simple, it is again convenient to argue for $d=1$.
% and  leave it to the reader to check the details in  the general case. 
\smallskip

The characteristic odes for the deterministic doubled pde \eqref{detds} are 
\begin{equation*}
\begin{cases} 
\dot {X} = - D_p H(P,X) \ \ X(0) =x  \quad \dot{Y} = - D_q H(Q,Y) \ \ Y(0) = y,\\[1mm]
%\noalign{\vskip6pt}
\dot{P} = D_x H(P,X)  \quad \dot{Q} = D_y H(Q,Y) \quad P(0)=Q(0) = \lambda (x-y)\\[1mm]
%\noalign{\vskip6pt}
\dot{U} = H(P,X) - \langle D_p H(P,X,P)\rangle - H(Q,Y) + \langle D_Q H(Q,Y) Q)\rangle \ \  U(0) = \lambda' |x-y|^2.
%\\[1mm]
%X(0) =x \ \ Y(0) = y \ \  P(0)=Q(0) = \lambda (x-y) \ \  U(0) = \lambda' |x-y|^2.
\end{cases}
%\quad\text{and}\quad
%\begin{cases} 
%\dot{Y} = - D_q H(Q,Y),\\
%\noalign{\vskip6pt}
%\dot{Q} = D_y H(Q,Y),\\
%\noalign{\vskip6pt}
%Y(0) = y,\ Q(0) = \lambda (x-y)\ ,
%\end{cases}
\end{equation*}
%and 
%$$\begin{cases}
%\dot{U} = H(P,X) - D_p H(P,X,P) - H(Q,Y) + (D_Q H(Q,Y) Q)\ ,\\
%\noalign{\vskip6pt}
%U(0) = \lambda' |x-y|^2.
%\end{cases}
%$$
%again to keep the formul{\ae} symmetric we write the 
%equation $Q$ instead of $-Q$.
Recall that the Jacobian is given by 
\begin{equation*}
J = \Partial{X}x \Partial{Y}y - \Partial{X}{y} \Partial{Y}x,
\end{equation*}
and, for $\alpha =x$ or $y$,
\begin{equation*} 
\begin{cases}
\ds \Partial{X}{\alpha} = - D_p^2 H(P,X) \Partial{P}{\alpha} - D_{px}^2 H(P,X)
\Partial{X}{\alpha}\ , \quad  \Partial{X}{\alpha}(0) = \begin{cases} 1&\text{if } \alpha=x\\
0&\text{if } \alpha \ne x\end{cases}\ ,           \\
\noalign{\vskip9pt}
\ds \Partial{P}{\alpha} = D_{px}^2 H(P,X)\Partial{P}{\alpha} 
+ D_x^2 H(P,X) \Partial{X}{\alpha}\ ,  \quad  \Partial{P}{\alpha}(0) = \begin{cases} \lambda&\text{if } \alpha=x\\
-\lambda &\text{if } \alpha=y\end{cases}\ ,
\end{cases}
  \\
%\noalign{\vskip9pt}
%\ds \Partial{X}{\alpha}(0) = \begin{cases} 1&\text{if } \alpha=x\\
%0&\text{if } \alpha \ne x\end{cases}\ , 
%\quad  \Partial{P}{\alpha}(0) = \begin{cases} \lambda&\text{if } \alpha=x\\
%-\lambda &\text{if } \alpha=y\end{cases}\ ,
%\end{cases}
\end{equation*}
and 
\begin{equation*}
\begin{cases}
\ds\Partial{Y}{\alpha} = - D_q^2 H(Q,Y)\Partial{Q}{\alpha} 
- D_{qy}^2 H(Q,Y) \Partial{Y}{\alpha}\ , \quad \Partial{Y}{\alpha}(0) = \begin{cases} 0&\text{if } \alpha=x\\
1&\text{if } \alpha=y\end{cases}\ ,\\
\noalign{\vskip9pt}
\ds\Partial{Q}{\alpha} = D_{yq}^2 H(Q,Y)\Partial{Q}{\alpha}
+ D_y^2 H(Q,Y)\Partial{Y}{\alpha}\ , \quad  \Partial{Q}{\alpha}(0) = \begin{cases} \lambda&\text{if } \alpha=x\\
-\lambda &\text{if } \alpha=y\end{cases}.\\
%\noalign{\vskip9pt}
%\ds\Partial{Y}{\alpha}(0) = \begin{cases} 0&\text{if } \alpha=x\\
%1&\text{if } \alpha=y\end{cases}\ ,
%\quad 
%\Partial{Q}{\alpha}(0) = \begin{cases} \lambda&\text{if } \alpha=x\\
%-\lambda &\text{if } \alpha=y\end{cases}.
\end{cases}
\end{equation*}
It is also convenient to consider, for $i = 1,2$ and $z=x-y$, the linearized auxiliary systems
\begin{equation}\label{aux}
\begin{cases}
\dot\xi_i \!=\! -D_{pp}^2 H(P,X) 
(1+\lambda |z|) \phi_i \! -\!  D^2_{xp} H(P,X)\xi_i\ , \quad \xi_1 (0) =1,\ \xi_2 (0) = 0\ \\
\noalign{\vskip6pt}
\dot\phi_i = D_{xp}^2 H(P,X)\phi_i + 
\ds \frac{D_{xx}^2 H(P,X)\xi_i}{1+\lambda |z|}\ , \quad \phi_1 (0)=0,\ \phi_1 (0) = \ds \frac1{1+\lambda |z|}\ , \\
%\noalign{\vskip6pt}
%\xi_1 (0) =1,\ \xi_2 (0) = 0\ ,\\
%\noalign{\vskip6pt}
%\phi_1 (0)=0,\ \phi_1 (0) = \ds \frac1{1+\lambda |z|}\ ,
\end{cases}
\end{equation}
%\text{ and } 
and
\begin{equation}\label{takis67}
\begin{cases} 
\dot\eta^i \! =\!  - D_q^2 H(Q,Y)(1+\lambda |z|) \psi_i \! -\!  
D_{qy}^2 H(Q,Y)\eta^i\ , \quad \eta_1 (0)=1,\ \eta_2 (0)=0\ ,\\
\noalign{\vskip6pt}
\dot\psi_i = D_{qy}^2 H(Q,Y) \psi_i + D_{yy}^2 H (Q,Y)
\ds\frac{\eta^i}{1+\lambda |z|}\ , \quad \psi_1 (0)=0,\ \psi_2 (0) = \ds\frac1{1+\lambda|z|}.\\
%\noalign{\vskip6pt}
%\eta_1 (0)=1,\ \eta_2 (0)=0\ ,\\
%\noalign{\vskip6pt}
%\psi_1 (0)=0,\ \psi_2 (0) = \ds\frac1{1+\lambda|z|}.
\end{cases}
\end{equation}
%where, for the $(\xi_i,\phi_i)$-system, $H$ is evaluated at $(X,P)$, while 
%for the $(\eta,\psi)$-system $H$ is evaluated at $(Y,Q)$.
It is immediate that 
\begin{equation*}
\begin{cases}
\ds\Partial{X}{x} = \xi_1 +\lambda\xi_2,\ \ \ \Partial{X}{y} = -\lambda\xi_2, \ \  \ \Partial{Y}{x} = - \lambda\eta_2, \ \  \ \Partial{Y}{x} = - \lambda\eta_2,\\
\noalign{\vskip6pt}
\ds\Partial{P}{x} = (\phi_1+\lambda\phi_2) (1+\lambda |z|) \qquad  \Partial{Q}{x} = -\lambda\psi_2 (1+\lambda |z|), \\
\noalign{\vskip6pt}
\ds\Partial{P}{y} = - \lambda\phi_2(1+\lambda|z|) \qquad  \Partial{Q}{y} = (\psi_1+\lambda\psi_2)(1+\lambda|z|).
\end{cases}
%\ \text{ and } \ \ 
%\begin{cases}
%\ds\Partial{Y}{x} = - \lambda\eta_2,\\ 
%\ds\Partial{Y}{y} = \eta_1 +\lambda\eta_2,\\
%\noalign{\vskip6pt}
%\ds\Partial{Q}{x} = -\lambda\psi_2 (1+\lambda |z|),\\
%\noalign{\vskip6pt}
%\ds\Partial{Q}{y} = (\psi_1+\lambda\psi_2)(1+\lambda|z|).
%\end{cases}
\end{equation*}
Assume next that, for all $p,x\in\R^d$, 
\begin{equation}\label{Hass3}
\begin{split}
&|D^2_{xp}H(p,x)| \lesssim  1,\quad 
|D_p^2H(p,x| \lesssim  1,\quad 
(1+|p|)|D_x^2H(p,x| \lesssim  1,\quad 
|D_{xxp}^3 H(p,x| \lesssim  1,\\[1.5mm] 
& (1+|p|)|D_{xpp}^3 H(p,x)| 					%% \frac{C}{(1+|p|)}
 \lesssim  1\ \text{ and }\ 
(1 +|p|)^2|D_p^3 H(p,x)|  \lesssim  1.
\end{split}
\end{equation}
It follows that there exists $C= C(T)> 0$ such that, for all $t\in [-T,T]$,  
\begin{equation}\label{Hass4} 
|\xi_1(t) | \leqq C,\quad |\eta_1(t)|\leqq C,\quad 
|\xi_2(t)|\leqq \frac{Ct}{1+\lambda |z|}\ \ \text{ and }\ \ 
|\eta_2(t)| \leqq \frac{Ct}{1+\lambda |z|}.
\end{equation}
Consider the matrices
\begin{equation*}
A^x = \begin{pmatrix}
-D_{xp}^2 H(P,X) & - D_{pp}^2 H(P,X) (1+\lambda |z|)\\
\noalign{\vskip6pt}
\ds \frac{D_x^2 H(P,X)}{1+\lambda |z|} & D_{xp}^2 H(P,X)
\end{pmatrix}% \quad\text{and}\quad 
\end{equation*}
and 
\begin{equation*}
A^y = \begin{pmatrix} 
-D_{yq}^2 H(Q,Y) & - D_{qq}^2 H(Q,Y)(1+\lambda |z|) \\
\noalign{\vskip6pt}
\ds \frac{D_{yy}^2 H(Q,Y)}{1+\lambda |z|}&D_{yq}^2 H(Q,Y).
\end{pmatrix}
\end{equation*}
%where again, $H$ is evaluated along $(X,P)$ in $A^x$ and along 
%$(Y,Q)$ in $A^y$.
The next lemma, which is stated without proof,  is important for the development of the rest of the theory here as well as for the theory of pathwise
conservation laws. 
\begin{lem}\label{mat}
Assume that,  in addition to \eqref{Hass3}, for all $p,x \in \R^d$, 
$|D_pH(p,x)|$ 
and $(1+|p|)^{-1}|D_x H(p,x)|$  are bounded. 
Then there exist  $\ep_0> 0$ and $C> 0$  such that, for all $t\in (0,\ep_0)$,
\begin{equation*}
\|A^x - A^y\| \leqq C|z|\ .
\end{equation*}
\end{lem}
Lemma~\ref{mat} implies   that, for all 
$t\in (0,\ep_0)$, 
\begin{equation}\label{Hass5} 
|\xi_1 - \eta_1 | \leqq C|z|t\quad\text{ and }\quad 
|\xi_2 - \eta_2 | \leqq \frac{C|z|t}{1+\lambda |z|},
\end{equation}
and, since
\begin{equation*}
\lambda (\xi_2 \eta_1 - \xi_1 \eta_2) = \lambda (\xi_2 -\eta_2) \eta_1 
+ \lambda \eta_2 (\eta_1 - \xi_1),
\end{equation*}
it follows from \eqref{Hass4} and \eqref{Hass5} that 
\begin{equation}\label{Hass6}
|\lambda (\xi_2 \eta_1 - \xi_1 \eta_2)| \leqq Ct.
\end{equation}

Similar arguments  allow to 
obtain an interval of invertibility of the 
characteristics that is uniform in $\lambda$, and, hence, a $O(1)$-interval of existence of smooth solutions of the doubled equation if  either one of 
the following three groups of possible assumptions  hold 
for all $(x,p) \in \R^d\times\R^d$:% and some $c> 0$:
%They are 
\begin{align}
&\begin{cases} 
|D_x^2 H| \lesssim 1 ,\quad |D_{xxx}^3 H| \lesssim 1 ,\quad 
|D_{xp}^2 H| \lesssim 1 ,\quad |D_{xxp}^3 H| (1+|p|) \lesssim 1 ,\\
\noalign{\vskip6pt}
|D_p^2 H| \lesssim 1 ,\quad |D_{xpp}^3 H| (1+|p|) \lesssim 1 ,\quad 
|D_p^3 H| (1+|p|) \lesssim 1 .
\end{cases} \label{Hcase1}\\
\noalign{\vskip9pt}
&\begin{cases}
|D_x^2H|\lesssim 1 ,\quad |D_x^3H| \lesssim (1+|p|),\quad |D_{xp}^2H| \lesssim 1,
\quad |D_{xxp}^3 H| \lesssim 1,\\
\noalign{\vskip6pt}
|D_p^2 H| (1+|p|) \lesssim 1,\quad |D_{ppx}^3 H| (1+|p|)\lesssim 1,\quad 
|D_p^3 H| (1+|p|^2) \lesssim 1.
\end{cases}\label{Hcase2}\\
\noalign{\vskip9pt}
&\begin{cases}
|D_p^2 H| \lesssim 1,\quad |D_{xp}^2H| \lesssim 1,\quad |D_x^2 H| (1+|p|) \lesssim 1,\\
\noalign{\vskip6pt}
|D_p^3 H| (1\!+\!|p|) \lesssim 1,\ \ |D_{xxp}^3H| (1\!+\!|p|^2) \lesssim 1,\ \ 
|D_{xpp}^3 H| \lesssim 1,\ \ |D_x^3 H| (1\!+\!|p|) \lesssim 1.
\end{cases} \label{Hcase3}
\end{align}
Note that  \eqref{Hcase1} contains the split variable case, and 
linear-type Hamiltonians are a special case of \eqref{Hcase2}.
%\smallskip
\vskip.05in

Calculations similar to the ones used in  the split variable case  yield
$$|\xi_2 \eta_1 -\eta_2 \xi_1| \lesssim t^2,$$
and, as was already seen,  $t_\lambda =\lambda^{-1/3}$.
%\smallskip
Note that, if $|DH|$, $|D^2H|$ and $|D^3H|$ are all bounded, then 
$t_\lambda = \lambda^{-1/2}$.
%\smallskip

\subsection*{The necessity of the assumptions} An important question is whether  conditions like the ones stated above are actually necessary to have well posed problems for Hamiltonians that depend on $p,x$. That some conditions are needed is natural since the argument is based on inverting characteristics and, hence, staying away from shocks. In view of this, assumptions that control the behavior of $H$ and its derivatives for large $|p|$ are to be expected.
\smallskip

On the other hand, some of the restrictions imposed are due to the specific choice of the initial datum of the doubled equation, which, in principle, does not ``interact well'' with the cancellation properties of the given $H$.   
\smallskip

Consider, for example,  the Hamiltonian 
\begin{equation}\label{Hspecial}
H(p,x)= F(a(x)p),
\end{equation}
with 
\begin{equation}\label{Hspecial1}
a, F\in C^2 (\R)\cap C^{0,1}(\R) \ \ \text{ and } \ \  a >0.
\end{equation} 
The characteristics are 
\begin{equation*}
\dot X = - F' (a(X)P) a(X)\quad\text{ and }\quad 
\dot P = F'(a(X)P) a'(X)P.
\end{equation*}
Let $\phi \in C^2 (\R)$ be such that 
$\phi' = \dfrac1a,$
%consider the new characteristics
$\hat X = \phi (X)$ and 
$\hat P= a(X) P.$
Then 
$$\dot{\hat X} = \phi' (X)\dot X = - a^{-1} (X) F' (a(X)P) a(X) 
= - F'(\hat P)$$
and
$$\dot{\hat P} = a' (X)\dot X P + a(X) \dot P 
= - a' (X) F'(\hat P) a(X)P + a(X) F'(a(X)P) a'(X)P = 0\ .$$
The observations above yield that it is better to use $\lambda |\phi^{-1}(x)-\phi^{-1}(y)|^2$ instead of $\lambda |x-y|^2$ in the comparison proof. % instead of $\lambda |x-y|^2$, it is  can be 
\smallskip

At the level of the pde
$$du = F(a(x) u_x)\cdot dB,$$
the above transformation yields that, if 
$u(x,t) = U(\phi (x),t),$
then 
$$dU = F(U_x) \cdot dB\ ,$$
a problem which is, of course, homogeneous in space, and, hence, as 
already seen, there is a $O(1)$-interval of existence for the doubled pde.
\smallskip

This leads to the question if it is possible to find, instead of $\lambda |x-y|^2$,  an initial datum for the doubled pde, which is still coercive, and, in the mean time, better adjusted to the structure of the doubled equation.  This is the topic of the next subsection.

%This question was answered affirmatively  for equations with quadratic Hamiltonians corresponding to Riemannian metrics in Friz, Gassiat, Lions and Souganidis \cite{fgls}.  Positively homogeneous and convex in $p$ Hamiltonians were considered by  Seeger \cite{seeger2}. The final result, due to Lions and the author \cite{ls8}, applies to problems with  general and convex in $p$. In what follows, I sketch the results \cite{ls8}. For simplicity, I am only considering the well-posedness of 
%\beq\label{paris10}
%du=H(Du,x)\cdot dB \ \ \text{in} \ \  Q_T \quad u(\cdot,0)=u_0.
%\Eq
% Hamiltonians in Lions and Souganidis \cite{lionssouganidis8}.  Finally, in an about to be completed investigation, Lions and Souganidis extend  the well-posedness theory   of 
%\beq\label{takis2000}
%du = 
%F (D^2 u, D u, u,x,t) dt + \sum_{i=1}^m H^i(Du, x,t) \cdot dB_i \ \ 
%\text{in} \  \ Q_T, 
%\Eq
%where the $B_i$'s are Brownian paths and the $H^i$'s satisfy suitable regularity assumptions and are nothing necessarily convex.

\subsection*{Convex Hamiltonians and a single path} The example discussed was the motivation behind several works which eventually led to a new class of well-posedness results in the case of a single path and convex Hamiltonians. 
\smallskip

The first result in this direction which applied to quadratic Hamiltonians corresponding to Riemannian metrics is due to Friz, Gassiat, Lions and Souganidis \cite{fgls}. A more general version of the problem (positively homogeneous and convex in $p$ Hamiltonians)  was studied in  \cite{seegerh}.  The final and definitive results, which apply to general convex in $p$ Hamiltonians with minimal regularity conditions,  were obtained by  Lions and Souganidis \cite{ls8}.  These  results  are sketched next. 
\smallskip

To keep the ideas simple, here I only discuss the first-order problem 
\beq\label{paris10}
du=H(Du,x)\cdot dB \ \ \text{in} \ \  Q_T \quad u(\cdot,0)=u_0.
\Eq
%Finally, in an about to be completed investigation, Lions and Souganidis extend  the well-posedness theory   of 
%
%expand here about the well-posedness results for smooth, convex and  spatially dependent Hamiltonians, which do not satisfy the conditions of the previous subsection.
%
%% The key idea here is to replace $\lambda |x-y|^2$ in the doubled equation by something that is more adapted to the Hamiltonians. This is related to the ``distance'' function associated with $H$. 
%\smallskip

To motivate the question,  I recall that the basic step of any comparison proof for viscosity solutions is to maximize functions  like $u(x,t)-v(y,t) -\lambda |x-y|^2$. The properties of $\lambda |x-y|^2$ used in the proofs are that 
\beq\label{aa1}
D_x L_\lambda =-D_y L_\lambda, \ \  L_\lambda \geq 0, \ \  L_\lambda(x,x)=0 \ \ \text{and} \ \   L_\lambda(x,y) \to \infty \ \text{if}\ |x-y|>0.
\Eq

The difficulty is that in the spatially dependent problems this choice of $ L_\lambda$ leads to expressions like
$H(\lambda (x-y),x)-H(\lambda (x-y),y)$ and, hence, error terms that are difficult to estimate when dealing with rough signals.

\smallskip

To circumvent this problem it seems to be natural to ask if it is possible to replace $ \lambda |x-y|^2$ by some $L_\lambda (x,y)$ that has similar continuity and coercivity properties and is better suited to measure the ``distance'' between  $H(\cdot,x)$ and $H(\cdot, y)$.
\smallskip

In particular, it is necessary  to find $L_\lambda :\R^d\times \R^d \to \R$ such that  
\beq\label{aa2}
\begin{split}
&H(D_x L_\lambda,x)=H(-D_y L_\lambda,y), \ \  L_\lambda \gtrsim -\lambda^{-1}, L_\lambda(x,y) \underset{\lambda \to \infty}\to \oo \ \text{if} \ x\neq y, \\[1mm]
 &L_\lambda (x,x)\underset{\lambda \to \infty}\to 0, \ \ \text{and} \ \  L_\lambda\in C^1_{x,y} \ \text{in a neighborhood of} \ \{x=y\}.
\end{split}
\Eq

It turns out (see \cite{ls8}) that this is possible if $H$ is convex or, more generally, if there exists $H_0$ convex such that the pair $H, H_0$ is an involution, that is, $\{H,H_0\}=0$. Here I concentrate on the convex case. 
\smallskip

Given $H$ convex with Legendre transform $L$, define
\beq\label{aa3}
L_\lambda (x,y)=\inf \left \{\int_0^{\lambda^{-1}} L(-\dot x(s), x(s)) ds: \  x(0)=x, \ x(\lambda^{-1})=y,  \ x(\cdot)\in C^{0,1}([0,\lambda^{-1}])\  \right\}.
\Eq
It follows, see, for example, Crandall, Lions and Souganidis \cite{cls} and \cite{lbook},   that $L_\lambda (x,y)=\overline L(x,y, \lambda^{-1})$, where $\overline L$  is the unique solution of 
\begin{equation*}
\begin{split}
&\overline L_t + H(D_x \overline L, x)=0 \ \ \text{in} \ \ \R^d\times(0,\infty) \quad  \overline L(x,y,0)={\bf \delta}_{\{y\}}(x)\\[1.2mm]
&\overline L_t + H(-D_y \overline L, y)=0 \ \ \text{in} \ \ \R^d\times(0,\infty) \quad  \overline L(x,y,0)={\bf \delta}_{\{x\}}(y),
\end{split}
\end{equation*}
where ${\bf \delta}_A(x)=0 \ \text{if} \ x\in A$ and ${\bf \delta}_A(x)= \oo \ \text{otherwise}.$ % \end{cases}$
\smallskip

Note that, at least formally, the above imply that $H(D_x L_\lambda,x)=H(-D_y L_\lambda,y)$.  From the remaining properties in \eqref{aa2} the most challenging one is the regularity.
% since the limiting behavior as $\lambda \to \oo$ is already encoded in the initial value problems.
\smallskip

%The special case that $H(p,x)=(\sum_{i,j=1}^{d} a_{ij}(x)p_ip_j)^{1/2}$ with $(a_{ij})_{i,j\in \N}$ positive definite, in which case $\overline L$ is related to the Riemmanian distance associated to $H$, was considered in \cite{gfls}, while positively homogeneous and convex $H$, in which case $\overline L$ is related with the associated Finsler metric were studied in \cite{seeger}. 
%\smallskip

I summarize next without proofs the main result of \cite{ls8}. In what follows $\nu$ and $\mu$ denote respectively constants for lower and upper bounds.
\smallskip

The assumption on $L:\R^d\times \R^d \to \R$ is that there exist positive constants $q,\nu, \mu$ and $C\geq 0$ such that, for all $\xi \in \R^d$, 
\beq\label{aa4}
\begin{split}
&\nu |p|^q-C\leq L\leq \mu |p|^q +C,  \ \ |D_x L| \leq  \mu |p|^{q} +C,\ \ |D_{px} L| \leq  \mu |p|^{q-1} +C,  \\[1mm]
& \nu |p|^{q-2}|\xi|^2\leq  \langle D^2L_{p}\xi, \xi\rangle \leq  (\mu |p|^{q-2} +C)|\xi|^2 \ \ \text{and} \ \ |D^2L_{x}| \leq \mu |p|^{q} +C;
\end{split}
\Eq
%
%
%\beq\label{aa6}
%\nu |p|^{q-2}\leq  |L_{pp}| \leq  \mu |p|^{q-2} +C,
%\Eq
%\beq\label{aa7}
%|D_{px} L| \leq  \mu |p|^{q-1} +C,
%\Eq
%and
%\beq\label{aa8}
%|L_{xx} \leq \mu |p|^{q} +C;
%\Eq
notice that it is important that $D^2_pL$ is positive definite. 
%  there is not a constant in the lower bound for $|L_{pp}|$ and this is important. 
\smallskip

The result is stated next.

\begin{thm} Assume \eqref{aa4}. Then:

(i). If $q \leq  2$, then there exists $\lambda_0$ such that, if $\lambda >  \lambda_0$, $L_\lambda \in C^1_{x,y}(\{|x-y| < \lambda^{-1}\}).$ 

(ii). If $q>2$ and $C> 0$, then, in general, (i) above is false, and, in fact, $\overline L(x,x,\lambda^{-1})$ may not be differentiable for any $\lambda$. 

(iii).  If $q >2$ and $C=0$, then there exists $\lambda_0$ such that, if $\lambda >\lambda_0$, $L_\lambda \in C^1_{x,y}(\{|x-y| < \lambda^{-1}\}).$ 

(iv).  In all cases, $\overline L$ is semiconcave in both $x$ and $y$.
\end{thm} 

It follows that, when $q\leq 2$ or $q> 2$ and $C=0$,  the pathwise solutions of the stochastic Hamilton-Jacobi initial value problem are well posed. The result extends to the full second order problem, because the semiconcavity  is enough to carry out the details. 

\subsection*{Multiple paths} I sketch here  briefly the strategy that Lions and the author developed in \cite{lionssouganidis7} to establish the well-posdeness of the pathwise solutions in the multi-path spatially dependent setting with Brownian signals. The argument is rather technical and to keep the ideas as simple as possible I only discuss the first-order problem
\beq\label{aa1}
du=\sum_{i=1}^m H^i(Du,x)\cdot dB_i \ \ \text{in} \ \ Q_T \quad u(\cdot,0)=u_0,
\Eq
and provide some hints about  the difficulties and the methodology.
\smallskip

As in the single-path case, the main step is to obtain a sufficiently long interval of existence of smooth solutions of the doubled  initial value problem 
\beq\label{a2}
dU=\sum_{i=1}^m \left [H^i(D_xU,x) - H^i(-D_yU,y)\right ]\cdot dB_i \ \ \text{in} \ \ \R^d\times \R^d \times (0,T] \quad U(x,y,0)=\lambda |x-y|^2.
\Eq

The semiformal argument presented earlier suggests that it is necessary to have an interval of existence of order $\lambda^{-\alpha}$ for an appropriately chosen small  $\alpha > 0$ which depends on the properties of the path. This was accomplished by reverting to the ``non rough'' time homogeneous doubled equation, 
% and the fact that for short time the solution of the stochastic problem were obtained by rescaling the one of the deterministic one, 
 something that is not possible for \eqref{a2}.
\smallskip

The new methodology developed in \cite{lionssouganidis7} consists of several steps. The first is to  provide a large deviations-type estimate about the error, in terms of powers of $\lambda^{-1}$, between the stochastic characteristics and their linearizations and the Jacobian,
% that is, between the paths $(X, P, Y, Q, \dfrac {\partial X}{\partial x}, \dfrac {\partial X}{\partial x}, \dfrac {\partial X}{\partial y}, \dfrac {\partial Y}{\partial x}, \dfrac {\partial Y}{\partial y})$, 
 and  their second-order expansion in terms of $B$ and its Levy areas. This would be straightforward, if it were not for the fact that the error must be uniform in $(x,y)$ such that $|x-y|=\text{O}(\lambda^{-1/2})$. 
\smallskip

Next I describe this problem for  the solution $S$ of  a stochastic differential equation $dS=\sigma (S) dB$ with $S(0)=s$. The aim  is to obtain an exponentially small estimate for the probability of the event that $\sup_{s\in K} |S(t)-(s+\sigma(s) B(s)+(1/2)\sigma \sigma'(s) B^2(t))| > \lambda^{-\beta}$, where is $A$ is a subset of $\R$ which may depend on $\lambda$.  In other words we need an estimate for  the probability of the $\sup$ instead of the $\sup$ of the probability. Obtaining such a result requires a new approach based on estimating $L^p$-norms of events for large $p$. 
\smallskip

Having  such estimates allows for a local in time comparison result off a set of exponentially small probability in terms of $\lambda^{-1}$. An ``algebraic'' iteration of this local comparison provides the required result at the limit $\lambda \to \infty$.

\section{Perron's method}

Perron's method is a general way to obtain solutions of equations which satisfy a comparison principle. The general argument is that the maximal subsolution is  actually a solution. The  idea is that, at places where it fails to be a solution, a subsolution can be strictly increased and maintain the subsolution property. This is a local argument  which has been carried out successfully for ``deterministic'' viscosity solutions. This locality creates, however, serious technical difficulties in the rough path setting due to the rigidity of the test functions.  
\smallskip

In this section I discuss this method in the context of 
%The existence of solutions as uniform limits of solutions of equations with smooth paths was discussed earlier in the notes.
the simplified initial value problem
\begin{equation}\label{E:Perrons}
	du = F(D^2 u, Du)\;dt + \sum_{i=1}^m H^i(Du, x) \cdot dB_i \ \  \text{in } \ \ Q_T \quad  u(\cdot,0) = u_0 ,
\end{equation}
where $u_0 \in BUC(\R^d)$, $T >0$, and $B = (B^1, \ldots, B^m)$ is a Brownian path.  The method can be a extended to problems with $F$ depending also on $(x,t)$ and $B$ a geometric rough path that is $\alpha$-H\"older continuous for some $\alpha \in (1/3,1/2]$. For details I refer to \cite{seeger}. 
%on the rough-path theory, see [reference]. 
%Without loss of generality, it may be assumed in this section that $B$ is a Brownian path and \eqref{E:Perrons} is interpreted in the Stratonovich sense.
\smallskip

Throughout the discussion it is  assumed that
\beq\label{paris20}
\begin{cases}
 F: {\mathcal S}^d \times \R^d \to \R \ \text{is continuous, bounded for bounded $(X,p) \in {\mathcal S}^d \times \R^d$} \\[1.2mm]
\text{ and  degenerate elliptic,}
\end{cases}
\Eq
 and the Hamiltonians are sufficiently regular, for example, 
\begin{equation}\label{A:Hregularity}
	H \in C^4_b(\R^d \times \R^d; \R^m),
\end{equation}
to allow for the construction of local-in-time, $C^2$ in space solutions of $dw=\sum_{i=1}^m H^i(Du, x) \cdot dB_i.$
\smallskip

As mentioned in Section~3, if the Poisson brackets of the $\{H^i\}$ vanish, for example, if $m =1$ or there is no spatial dependence, then it suffices to have $H \in C^2_b(B_R \times \R^d;\R^m)$ for all $R > 0$.
\smallskip

%These assumptions are sufficient to establish the most important fact for Perron's method, that is, the comparison principle for \eqref{E:Perrons}. The result
%
%This has already been proved earlier in the notes, but for the convenience of the reader it is recalled next. 

The result is stated next. 

\begin{thm} \label{T:perrons}
	Assume \eqref{paris20} and \eqref{A:Hregularity}. Then \eqref{E:Perrons} has a unique solution $u\in BUC(Q_T)$, which is given by
	\begin{equation} \label{E:solution}
		u(x,t) = \sup\left\{ v(x,t) : v(\cdot,0) \le u_0 \text{ and $v$ is a subsolution of \eqref{E:Perrons}} \right\}.
	\end{equation}
\end{thm}

As has been discussed earlier, more assumptions are generally required for $F$, $H$, and $B$ in order for the comparison principle to hold. This is especially the case when $H$ has nontrivial spatial dependence even when $m = 1$. Apart from the assumptions that yield the comparison, the only hypotheses used for the Perron construction are \eqref{A:Hregularity}. %the basic boundedness and continuity assumptions for $F$.
\smallskip

%Next I  outline the steps of the proof of Theorem \ref{T:perrons}. For the sake of easing the presentation full generality is avoided.  For the details, I  refer the reader to \cite{seeger}.
%\smallskip

As before,  $S(t,t_0): C^2_b(\R^d) \to C^2_b(\R^d)$ be the solution operator for local in time, spatially smooth solutions of 
\begin{equation}\label{E:HJpart}
	d\Phi = \sum_{i=1}^m H^i(D\Phi,x) \cdot dB_i \ \  \text{in } \ \ \R^d \times (t_0 - h, t_0 + h) \quad \Phi(\cdot,t_0) = \phi.
\end{equation}
%which exists where, for some sufficiently small $h > 0$ depending on $\nor{D^2\phi}{\oo}$, $\nor{H}{C^4}$, and $B$, the solution $\Phi$ is given, for $t \in (t_0 - h, t_0 + h)$, by $\Phi(\cdot,t) = S(t,t_0)\phi$, and $\Phi$ belongs to $C((t_0 - h, t_0 + h), C^2_b(\R^d))$.

It is clear from the definition of stochastic viscosity subsolutions that the maximum of a finite number of subsolutions is also a subsolution, with a corresponding statement holding true for the minimum of a finite number of supersolutions. This observation to can be generalized to infinite families.

\begin{lem} \label{L:stability}
	Given a family  $\mcl F$ of subsolutions (resp. supersolutions) of \eqref{E:Perrons}, let  
		$U(x,t) = \sup_{ v \in \mcl F} v(x,t) \quad \pars{\text{ resp. } \inf_{v \in \mcl F} v(x,t) }.$
	If  $U^*<\oo$ (resp. $U_*> -\oo$),  then $U^*$ (resp. $U_*$) is a subsolution  (resp. supersolution) of \eqref{E:Perrons}.
\end{lem}

\begin{proof}
I only a sketch of the proof of the subsolution property. %, since it is almost identical for super-solutions.
\smallskip
	
Let $\phi \in C^2_b(\R^d)$, $\psi \in C^1([0,T])$, $t_0 > 0$, and $h > 0$ be such that $S(\cdot,t_0)\phi \in C((t_0-h,  t_0 + h),C^2_b(\R^d))$, assume that
	$U^*(x,t) - S(t,t_0)\phi(x) - \psi(t)$
attains a strict local maximum at $(x_0,t_0) \in \R^d \times (t_0 - h, t_0 + h)$, and set $p = D\phi(x_0)$, $X = D^2\phi(x_0)$, and $a = \psi'(t_0)$.  The goal is to show that
\[
	a \le F(X,p).
\]

%Upon modifying the  test functions $\phi$ and $\psi$ appropriately, it may be assumed that the maximum is strict. Note carefully that the comparison property of the solution operator $S(t,t_0)$ is used to justify this.
%
The definition of upper-semicontinuous envelopes and arguments from the classical viscosity solution theory imply that  there exist sequences $(x_n,t_n) \in \R^d \times (t_0-h,t_0+h)$ and $v_n \in \mcl F$ such that $\lim_{n \to \oo} (x_n,t_n) = (x_0,t_0)$, $\lim_{n \to \oo} v_n(x_n,t_n) = U^*(x_0,t_0)$, and
\[
	v_n(x,t) - S(t,t_0)\phi(x) - \psi(t)
\]
attains a local maximum at $(x_n,t_0)$. Applying the definition of stochastic viscosity subsolutions and letting $n \to \oo$ completes the proof.

\end{proof}

The second main step of the Perron construction is discussed next.
%to show that, if a ``strict'' subsolution has its values increased in a certain way in a sufficiently small open cylinder, then the resulting function is another subsolution. This ``bump'' construction is less straightforward than in the classical viscosity solution setting, due to the limited flexibility in the choice of test functions. 

\begin{lem} \label{L:bump}
Suppose that $w$ is a subsolution of \eqref{E:Perrons}, and that $w_*$ fails to be a supersolution. Then there exists $(x_0,t_0) \in \R^d \times (0,T]$ such that, for all $\kappa > 0$, \eqref{E:Perrons} has  a subsolution $w_\kappa$ such that 
	\[
		w_\kappa \ge w, \quad
		\sup(w_\kappa - w) > 0, \quad \text{and} \quad
		w_\kappa = w \ \ \text{in} \ \  Q_T \backslash \pars{B_{\kappa}(x_0) \times (t_0 -\kappa, t_0 + \kappa)}.
	\]
\end{lem}

\begin{proof}
By assumption, there exist $\phi \in C^2_b(\R^d)$, $\psi \in C^1([0,T])$, $(x_0,t_0) \in \R^d \times (0,T]$, and $h \in (0,\kappa)$ such that $S(\cdot,t_0)\phi \in C((t_0 - h, t_0 + h),C^2_b(\R^d))$, 
\[
	w_*(x,t) - S(t,t_0)\phi(x) - \psi(t)
\]
attains a local minimum at $(x_0,t_0)$, and
\begin{equation} \label{E:inequality}
	\psi'(t_0) - F(D^2\phi(x_0),D\phi(x_0)) < 0.
\end{equation}
Assume,  without loss of generality, that $x_0 = 0$, $\phi(0) = 0$, and $\psi(t_0) = 0$,  set $X = D^2\phi(0)$, $p = D\phi(0)$, and $a = \psi'(t_0)$, fix $\gamma \in (0,1)$, $r \in (0,\kappa)$, and $s \in (0,h)$, and choose $\hat \eta \in C_b^2(\R^d)$ and $h>0$ so that 
\[\hat \eta(x) = p\cdot x + \frac{1}{2} \langle Xx, x\rangle  - \gamma|x|^2 \text{ in } B_r(x_0), \ \  
	\hat \eta \le \phi \ \text{ in } \ \R^d, \ \  \text{and} \ \ S(\cdot,t_0)\hat \eta \in C((t_0 - h, t_0 + h);C^2_b(\R^d)).
\]
%\[
%	\left\{
%	\begin{split}
%	&\hat \eta(x) = p\cdot x + \frac{1}{2} Xx\cdot x - \gamma|x|^2 \text{ in } B_r(x_0) \quad \text{and} \\
%	&\hat \eta \le \phi \text{ in } \R^d.
%	\end{split}
%	\right.
%\]
%If necessary, shrink $h > 0$ so that $S(\cdot,t_0)\hat \eta \in C((t_0 - h, t_0 + h), C^2_b(\R^d))$. 

For $(x,t) \in \R^d \times (t_0 - h, t_0 + h)$ and $\delta > 0$, define
\[
	\widehat w(x,t) = w_*(0,t_0) + \delta + S(t,t_0)\hat \eta(x) + a(t-t_0) - \gamma(|t - t_0|^2 + \delta^2)^{1/2}.
\]
In view of the strict inequality in \eqref{E:inequality}, the continuity of the solution map $S(t,t_0)$ on $C^2_b(\R^d)$, and the continuity of $F$, if $\gamma$, $r$, $s$, and $\delta$ are sufficiently small, then $\widehat w$ satisfies the subsolution property in $B_r(0) \times (t_0-s,t_0+s)$.
\smallskip

The most important step in the proof is to show that, with all parameters sufficiently small, there exist $0 < r' < r$ and $0 < s' < s$ such that
\begin{equation}\label{E:bump}
	w > \widehat w  \quad \text{in } \left( B_{r}(0) \times (t_0 - s, t_0 +s) \right) \backslash \overline{ B_{r'}(0) \times (t_0 - s',t_0 + s')}. 
\end{equation}

Achieving the inequality in \eqref{E:bump} for points of the form $(x,t_0)$ can be done using classical arguments. However, this is much more  difficult for arbitrary $t \ne t_0$, because, in view of the definition of $\widehat w$, it is necessary to study the local in time, spatially smooth solution operator $S(t,t_0)$. 
\smallskip

This difficulty is overcome by establishing a finite speed of propagation for such local in time, spatially smooth solutions. As has been discussed earlier in the notes, such a result cannot be true in general. Here it relies on access to the system of rough characteristics. Indeed, the domain of dependence result is proved by estimating the deviation of characteristics from their starting points, using tools from the theory of rough or stochastic differential equations.
\smallskip

Once \eqref{E:bump} is established,  define
\[
	w_{\kappa}(x,t) = 
	\begin{cases}
		\max( \widehat{w}(x,t), w(x,t) ) & \text{for } (x,t) \in B_{r}(0) \times (t_0-s, t_0+s),  \\
		w(x,t) & \text{for } (x,t) \notin B_{r}(0) \times (t_0-s,t_0+s).
	\end{cases}
\]
Then $w_\kappa \ge w$, and $w_\kappa = w$ outside of $B_{\kappa}(0) \times (t_0 -\kappa,t_0+\kappa)$. If $(x_n,t_n)$ is such that $\lim_{n \to \oo} (x_n,t_n) = (0,t_0)$ and $\lim_{n \to \oo} w(x_n,t_n) = w_*(0,t_0)$, then 
\[
	\lim_{n \to \oo} \pars{ w(x_n,t_n) - \hat w(x_n,t_n) }= -(1 - \gamma)\delta < 0,
\]
so that
\[
	\sup_{B_{\kappa}(0) \times (t_0-\kappa,t_0+\kappa)} \pars{ w_\kappa - w} > 0.
\]
Finally, $w_\kappa$ is a subsolution. This is evident outside of $B_{r}(0) \times (t_0 -s, t_0 +s)$, as well as in the interior of $B_{r}(0) \times (t_0-s,t_0+s)$, because there, $w_\kappa$ is equal to the pointwise maximum of two subsolutions. It remains to verify the subsolution property on the boundary of $B_{r}(0) \times (t_0 -s,t_0+s)$, and this follows because, in view of \eqref{E:bump}, $w_\kappa = w$ in a neighborhood of the boundary of $B_{r}(0) \times (t_0-s,t_0+s)$.

\end{proof}

\begin{proof}[Proof of Theorem \ref{T:perrons}]
	The first step is to verify that $u$ is well defined and bounded. This follows from the comparison principle, and the fact that, in view of the assumptions,  it is possible to construct a subsolution  and a supersolution with respectively  initial datum $- \nor{u_0}{\oo}$ and $\nor{u_0}{\oo}$. 
\smallskip
	
	Fix $\eps > 0$ and let $\phi^\eps \in C^2_b(\R^d)$ be such that 
	\[
		\phi^\eps - \eps \le u_0 \le \phi^\eps + \eps \quad \text{on } \R^d.
	\]
	It is possible to construct a subsolution  and a supersolution $\uline{u}^\eps$ and $\oline{u}^\eps$ which are continuous in a neighborhood of $\R^d \times \{0\}$ and achieve respectively the initial datum $\phi^\eps - \eps$ and $\phi^\eps + \eps$. This can be done by using the solution operator $S(t_{k+1},t_k)$ on successive, small intervals $[t_k,t_{k+1}]$ and  the boundedness properties of $F$. Once again, see \cite{seegerp} for the details.
\smallskip
		
	The comparison principle yields
	\[
		\uline{u}^\eps \le u_* \le u \le u^* \le \oline{u}^\eps \quad \text{in } Q_T,
	\]
	and, in view of the continuity of $\uline{u}^\eps$ and $\oline{u}^\eps$ near $\R^d \times \{0\}$,
	\[
		\phi^\eps - \eps \le u_*(\cdot, 0) \le u^*(\cdot, 0) \le \phi^\eps + \eps.
	\]
	Since $\eps$ is arbitrary, it follows that $u(\cdot, 0) = u_0$ and $\lim_{(x,t) \to (x_0,0)} u(x,t) = u_0(x_0)$ for all $x_0 \in \R^d$.
\smallskip

	Lemma \ref{L:stability} now implies that $u^*$ is a subsolution of \eqref{E:Perrons} with $u^*(x,0) \le u_0(x)$. The formula \eqref{E:solution} for $u$ then yields $u^* \le u$, and, therefore, $u^* = u$. That is, $u$ is itself upper-semicontinuous and a subsolution.
\smallskip
	
	On the other hand, $u_*$ is a supersolution. If this were not the case, then Lemma \ref{L:bump} would imply the existence of a subsolution $\tilde u \ge u$ and a neighborhood $N \subset \R^d \times (0,T]$ such that $\tilde u = u$ in $(\R^d \times [0,T]) \backslash N$ and $\sup_{N} (\tilde u - u) > 0$, contradicting the maximality of $u$.
\smallskip
	
	The comparison principle gives $u^* \le u_*$, and, as a consequence of the definition of semicontinuous envelopes, $u_* \le u^*$. Therefore, $u = u_* = u^*$ is a solution of \eqref{E:Perrons} with $u = u_0$ on $\R^d \times \{0\}$. The uniqueness of $u$ follows from yet another application of the comparison principle.

\end{proof}

\section{Approximation schemes, convergence and error estimates}

Here I discuss a general program for constructing convergent (numerical)  approximation schemes for the pathwise viscosity solutions and obtain, for first-order equations,    explicit error estimates.
\smallskip

The  presentation focuses on the initial value problem 
\begin{equation}\label{E:schemes}
	du = F(D^2 u, Du)\; dt + \sum_{i=1}^m H^i(Du)\cdot dB_i \ \  \text{in } \ \ Q_T \quad 
	u(\cdot,0) = u_0, 
\end{equation}
where  $T > 0$ is a fixed finite horizon, $F \in C^{0,1}(\mathcal{S}^d \times \R^d)$ is degenerate elliptic, $H \in C^2(\R^d)$, $B = (B_1,  \ldots, B_m) \in C([0,T]; \R^m)$, and $u_0 \in BUC(\R^d)$. 
\smallskip

\subsection*{The scheme operator} Following the general methodology for constructing convergent  schemes for ``non-rough'' viscosity solutions put forward by Barles and Souganidis \cite{barlessouganidis2},  the approximations are constructed using  a ``scheme'' operator, which, for $h > 0$, $0 \le s \le t \le T$, and $\zeta \in C([0,T]; \R^m)$, is a  map $S_h(t,s; \zeta): BUC(\R^d) \to BUC(\R^d)$.
\smallskip

Given a partition $\mcl P= \{0 = t_0 < t_1 < \cdots< t_N = T\} $ of $[0,T]$ with mesh size  $\abs{ \mcl P}$ and
%=\max_{n = 0, 1, \ldots, N-1} \pars{ t_{n+1} - t_n}$, and 
% its mesh size, that is,
%\[
%	%\left\{
%	%\begin{split}
%	\mcl P := \{0 = t_0 < t_1 < \cdots< t_N = T\} \quad \text{and} \quad \abs{ \mcl P} := \max_{n = 0, 1, \ldots, N-1} \pars{ t_{n+1} - t_n}.
%	%&\abs{ \mcl P} := \max_{n = 0, 1, \ldots, N-1} \pars{ t_{n+1} - t_n}.
%	%\end{split}
%	%\right.
%\]
%Given a partition $\mcl P$ and 
a path $\zeta \in C_0([0,T]; \R^m)$, usually a piecewise linear approximation of $B$,  the (approximating) function $\tilde u_h(\cdot; \zeta, \mcl P)$ is defined by
\begin{equation} \label{E:introapproxsolution}
	\left\{
	\begin{split}
		&\tilde u_h(\cdot,0;\zeta, \mcl P ) := u_0 \quad \text{and}  \\[1.2mm]
		&\tilde u_h(\cdot,t; \zeta, \mcl P) := S_h(t,t_n; \zeta)\tilde u_h(\cdot,t_n; \zeta, \mcl P) & \text{for } n = 0, 1, \ldots, N-1,\; t \in (t_n, t_{n+1}].
	\end{split}
	\right.
\end{equation}
The strategy is to choose families of approximating paths $\{B_h\}_{h > 0}$ and partitions $\{\mcl P_h\}_{h > 0}$ satisfying
\begin{equation} \label{A:introapproximators1}
	\lim_{h \to 0^+} \nor{B_h - B}{\oo} = 0 = \lim_{h \to 0^+} \abs{ \mcl P_h},
\end{equation}
in such a way that the function
\begin{equation} \label{E:introschemedef1}
	u_h(x,t) := \tilde u_h(x,t; B_h, \mcl P_h)
\end{equation}
is an efficient approximation of the solution of \eqref{E:schemes}. 
\smallskip

The main restriction on  the scheme operator is that it has to be monotone, that is,
\begin{equation}\label{A:monotone}
	\left\{
	\begin{split}
		\text{if } t_n \le t \le t_{n+1}, \;& t_n, t_{n+1} \in \mcl P_h, \text{ and } u,v \in BUC(\R^d) \ 
		\text{such that} \     u\leq v \ \text{in} \ \R^d, \text{ then} \\
		%&u(x) \le v(x) \quad \Rightarrow \quad 
		&S_h(t,t_n;B_h)u \le S_h(t,t_n;B_h) \ \text{in} \ \R^d.
	\end{split}
	\right.
\end{equation}
It will also be necessary for the scheme operator to commute with constants,  that is, for all $u \in BUC(\R^d)$, $h > 0$, $0 \le s \le t < \oo$, $\zeta \in C_0([0,T],\R^m)$, and $k \in \R$,
\begin{equation}\label{A:constantcommute}
	S_h(t,s;\zeta)\pars{ u + k} = S_h(t,s;\zeta)u + k.
\end{equation}
Finally, the scheme operator must be ``consistent'' with the equation in some sense. This point, as well as the motivation for the above assumptions, are explained below.

\subsection*{The method of proof} I give here  a brief sketch  of the proof.
% to illustrate the utility of the definition of pathwise viscosity solutions and its amenability to the approximation theory outlined here. 
All the  details and concrete examples can be found in Seeger \cite{seegerschemes}.
\smallskip

Assume for the moment that $\lim_{h \to 0} u_h = u$ locally uniformly for some $u \in BUC(Q_T)$. In fact, a rigorous proof involves studying the so-called half-relaxed limits of $u_h$, but I omit these cumbersome details. % here.
\smallskip

The goal is to show that $u$ is the unique pathwise  solution of \eqref{E:schemes}. To that end, suppose that
\[
	u(x,t) - \Phi(x,t) - \psi(t)
\]
attains a strict maximum at $(y,s) \in \R^d \times I$, where $\psi \in C^1([0,T])$ and, for some small open interval $I \subset [0,T]$, $\Phi \in C(I;C^2(\R^d))$ is a local in time, smooth in space solution of
\begin{equation}\label{E:schemetestfn}
	d\Phi = \sum_{i=1}^m H^i(D\Phi) \cdot dB_i \ \  \text{in } \ \ \R^d \times I.
\end{equation}
I will show that % arrive at the inequality
\[
	\psi'(s) \le F(D^2\Phi(y,s), D\Phi(y,s)),
\]
which implies that $u$ is a subsolution. The argument to show it is a supersolution is similar.
\smallskip

For $h > 0$, let $\Phi_h$ be the local in time, smooth in space solution of
\begin{equation}\label{E:approximateHJ}
	\Phi_{h,t} = \sum_{i=1}^m H^i(D\Phi_h) \dot B_{i,h} \ \
	 \text{in } \ \ \R^d \times I \quad \Phi_h(\cdot,s) = \Phi(\cdot,t_0).
\end{equation}

Recall that such a solution can be shown to exist using  the method of characteristics. The interval $I$ may need to be shrunk, if necessary, but its length  is uniform in $h$. Since $\lim_{h\to 0} B_h=B$  uniformly on $[0,T]$, it follows that,  as $h \to 0$, $\Phi_h$ converges to $\Phi$ in $C(I;C^2(\R^d))$. % as $h \to 0$. 
\smallskip

As a result, there exists $\{(y_h,s_h)\}_{h > 0} \subset \R^d \times I$ such that $\lim_{h \to 0} (y_h,s_h) = (y,s)$ and
\[
	u_h(x,t) - \Phi_h(x,t) - \psi(t)
\]
attains a local maximum at $(y_h,s_h)$. 
\smallskip

That  $\underset{h\to 0}\lim \abs{ \mcl P_h}=0$ yields that,  for $h$ sufficiently small, there exist $n \in \N$ depending on $h$ such that
\[
	t_n < s_h \le t_{n+1} \quad \text{and} \quad t_n,t_{n+1} \in I.
\]

It then follows that 
\[
	u_h(\cdot,t_n) - \Phi_h(\cdot,t_n) - \psi(t_n) \le u_h(y_h,s_h) - \Phi_h(y_h,s_h) - \psi(s_h),
\]
or, after  rearranging terms,
\begin{equation}\label{applyscheme}
	u_h(\cdot,t_n) \le u_h(y_h,s_h) + \Phi_h(\cdot,t_n) - \Phi_h(y_h,s_h) + \psi(t_n) - \psi(s_h).
\end{equation}

This is the place where  the monotonicity \eqref{A:monotone} and the commutation with constants \eqref{A:constantcommute} of the scheme come into play. Applying $S_h(s_h,t_n; W_h)$ to both sides of \eqref{applyscheme}, %using the fact that the scheme commutes with constants, 
and evaluating the resulting expression at  $x = y_h$ give 
\[
	u_h(y_h,s_h) \le u_h(y_h,s_h) + S_h(s_h,t_n; B_h)\Phi_h(\cdot,t_n)(y_h) - \Phi_h(y_h,s_h) + \psi(t_n) - \psi(s_h),
\]
whence
\[
	\frac{\psi(s_h) - \psi(t_n)}{s_h - t_n} \le \frac{ S_h(s_h,t_n; B_h)\Phi_h(\cdot,t_n)(y_h) - \Phi_h(y_h,s_h)}{s_h - t_n}.
\]

As $h \to 0$, the left-hand side converges to $\psi'(s)$. The construction of a convergent scheme then reduces to creating a scheme operator, partitions $\mcl P_h$, and paths $W_h$ satisfying \eqref{A:monotone} and \eqref{A:constantcommute}, as well as the consistency requirement
\begin{equation}\label{A:consistent}
	\lim_{s,t \in I, \; t-s \to 0} \frac{ S_h(t,s; B_h)\Phi_h(\cdot,s) - \Phi_h(\cdot,s)}{t-s} = F(D^2\Phi,D\Phi)
\end{equation}
whenever $\Phi$ and $\Phi_h$ are as in respectively \eqref{E:schemetestfn} and \eqref{E:approximateHJ}.

\subsection*{The main examples} Presenting a full list of the types of schemes that may be constructed is beyond the scope of these notes. Here, I give a few specific examples that are representative of the general theory. More schemes  and details can be found in \cite{seegerschemes}.
\smallskip

Here I focus mainly on finite difference schemes. To simplify the presentation, assume $d = m = 1$, $F$ and $H$ are both smooth, and $F$ depends only on $u_{xx}$, so that \eqref{E:schemes} becomes
\begin{equation}\label{E:simpleeqintro}
	du = F(u_{xx})\;dt + H(u_x)\cdot dB \  \text{in }  \ Q_T \quad  u(\cdot,0) = u_0,
\end{equation}
and, in the first-order case when $F \equiv 0$,
\begin{equation} \label{E:simpleeqintrofirstorder}
	du = H(u_x)\cdot dB  \ \text{in } \   Q_T \quad u(\cdot,0) = u_0.
\end{equation}

I present next a number of different partitions $\mcl P_h$ and approximating paths  $B_h$  for which the program in the preceding subsection may be carried out. While technical, these are all made with the same idea in mind, namely, to ensure that the approximation $B_h$ is ``mild'' enough with respect to the partition. In particular, for any consecutive points $t_n$ and $t_{n+1}$ of the partition $\mcl P_h$ and for sufficiently small $h$, the ratio
\[
	\frac{\abs{ B_h(t_{n+1}) - B_h(t_n)}}{h}
\]
should be less than some fixed constant. This is a special case of the well-known  Courant-Lewy-Friedrichs (CFL) conditions required for the monotonicity of schemes in the ``non-rough'' setting. 
%schemes discussed in more detail in the forthcoming sections.
\smallskip 

For some $\eps_h > 0$ to be determined, define
\begin{equation} \label{E:introLF}
	\begin{split}
	&S_h(t,s; \zeta)u(x) = u(x) + H \pars{ \frac{u(x+h) - u(x-h)}{2h} } (\zeta(t) - \zeta(s)) \\[1.2mm]
	&+ \left[ F\pars{ \frac{u(x+h) + u(x-h) - 2u(x)}{h^2} } 
	+ \eps_h \pars{ \frac{u(x+h) + u(x-h) - 2u(x)}{h^2} } \right](t-s).
	\end{split}
\end{equation}
\smallskip

The first result, which is qualitative in nature, applies to the simple setting above as follows.

\begin{thm} \label{T:introLFresult}
	Assume that, in addition to \eqref{A:introapproximators1}, $B_h$ and $\mcl P_h$ satisfy
	\[
		\abs{ \mcl P_h} \le \frac{h^2}{\nor{F'}{\oo}} \quad \text{and} \quad \eps_h = h \|\dot B_h\|
		 \xrightarrow{h \to 0} 0.
	\]
	Then, as $h \to 0$, the function $u_h$ defined by \eqref{E:introschemedef1} using the scheme operator \eqref{E:introLF} converges locally uniformly to the solution $u$ of \eqref{E:simpleeqintro}.
\end{thm}

The condition in Theorem \ref{T:introLFresult} on the approximating path $B_h$ can be satisfied in  several different ways. For example, $B_h$ could be a piecewise linear approximation of $B$ of step-size $\eta_h > 0$, with $\lim_{h \to 0} \eta_h = 0$ in such a way that $\lim_{h \to 0} h \|\dot B_h\|= 0$.
\smallskip

By quantifying the method of proof in the previous subsection, it is possible to obtain explicit error estimates for finite difference approximations of the pathwise Hamilton-Jacobi equation \eqref{E:simpleeqintrofirstorder}. The results below are stated for the following scheme, which is defined, for some $\theta \in (0,1]$, by

\begin{equation} \label{E:introLFfirstorder}
\begin{split}
	S_h(t,s; \zeta)u(x) &= u(x) + H \pars{ \frac{u(x+h) - u(x-h)}{2h} } (\zeta(t) - \zeta(s)) \\
	&+ \frac{\theta}{2} \pars{ u(x+h) + u(x-h) - 2u(x) };
	\end{split}
\end{equation} 
note that this corresponds to choosing $\eps_h = \dfrac{\theta h^2}{2(t-s)}$ in \eqref{E:introLF}. 
\smallskip

Assume that $\omega: [0,\oo) \to [0,\oo)$ is the modulus of continuity of the fixed continuous path $B$ on $[0,T]$, define, for $h > 0$, $\rho_h$ implicitly by
\begin{equation} \label{A:ctspathCFL}
	\lambda = \frac{ (\rho_h)^{1/2} \omega( (\rho_h)^{1/2}) }{h} < \frac{\theta}{\nor{H'}{\oo}},
\end{equation}
and choose  the partition $\mcl P_h$ and path $B_h$ so that 
\begin{equation} \label{A:introregularpathspartition} 
	\left\{
	\begin{split}
		&\mcl P_h = \{ n \rho_h \wedge T \}_{n \in \N_0}, \; M_h := \lfloor (\rho_h)^{-1/2} \rfloor, \\
		&\text{and, for } k \in \N_0 \text{ and } t \in [kM_h\rho_h, (k+1)M_h\rho_h), \\
		&B_h(t) = B(kM_h\rho_h) + \pars{ \frac{ B((k+1)M_h\rho_h) - B(kM_h\rho_h)}{M_h\rho_h} } \pars{ t - kM_h\rho_h}.
	\end{split}
	\right.
\end{equation}
	
\begin{thm}\label{T:introLFresultfirstorder}
	There exists $C > 0$ depending only on $L$ such that, if $u_h$ is constructed using \eqref{E:introschemedef1} and \eqref{E:introLFfirstorder} with $\mcl P_h$ and $B_h$ as in \eqref{A:ctspathCFL} and \eqref{A:introregularpathspartition}, and $u$ is the pathwise viscosity solution of \eqref{E:simpleeqintrofirstorder}, then
	\[
		\sup_{(x,t) \in \R^d \times [0,T]} \abs{ u_h(x,t) - u(x,t)} \le C(1+T) \omega( (\rho_h)^{1/2}).
	\]
\end{thm}

If, for example, $B \in C^{0,\alpha}([0,T])$, then \eqref{A:ctspathCFL} means that $\rho_h = \text{O}(h^{2/(1+\alpha)})$, and the rate of convergence in Theorem \ref{T:introLFresultfirstorder} is $\text{O}(h^{\alpha/(1+\alpha)})$.
\smallskip

I  describe next some examples  in the case that $B$ is a Brownian motion.
% for which \eqref{E:simpleeqintrofirstorder}. % becomes the stochastic Hamilton-Jacobi equation
%\begin{equation}\label{E:introstochHJ}
%	du = H(u_x)\circ dB \ \  \text{in } \ \ \R \times (0,T] \quad u(\cdot,0) = u_0.
%\end{equation}
%
\smallskip

As a special case of Theorem \ref{T:introLFresultfirstorder}, the approximating paths and partitions may be taken to satisfy \eqref{A:introregularpathspartition} with $\rho_h$ given by
\begin{equation} \label{A:Brownianregularpartition}
	\lambda = \frac{ (\rho_h)^{3/4} \abs{ \log \rho_h}^{1/2}}{h} < \frac{\theta}{\nor{H'}{\oo}},
\end{equation}
in which case the scheme operator will be monotone almost surely for all $h$ smaller than some (random) threshold $h_0 > 0$. 
\smallskip

It is also possible to define the partitions and approximating paths using  certain stopping times that ensure that the scheme is monotone almost surely for all $h > 0$. More details can be found in   \cite{seegerschemes}.

\begin{thm} \label{T:introLFpathwiseBM}
	Suppose that $B$ is a Brownian motion, and assume that $\mcl P_h$ and $B_h$ are as in \eqref{A:introregularpathspartition} with $\rho_h$ defined by \eqref{A:Brownianregularpartition}. If $u_h$ is constructed using \eqref{E:introschemedef1} and \eqref{E:introLFfirstorder}, and $u$ is the solution of \eqref{E:simpleeqintrofirstorder}, then there exists a deterministic constant $C > 0$ depending only on $L$ and $\lambda$ such that, with probability one,
	\[
		\limsup_{h \to 0} \sup_{(x,t) \in \R^d \times [0,T]} \frac{ \abs{ u_h(x,t) - u(x,t)}}{h^{1/3} \abs{ \log h}^{1/3}} \le C(1+T).
	\]
\end{thm}

The final  result presented here is about a scheme that  converges in distribution in the space $BUC(\R^d \times [0,T])$ equipped  with the topology of local uniform convergence.
\smallskip

Recall that,  given random variables $(X_\delta)_{\delta > 0}$ and $X$ taking values in some topological space $\mcl X$, it is said that  $X_\delta$  converges, as $\delta \to 0$ in distribution (or in law) to $X$,  if the law $\nu_\delta$ of $X_\delta$ on $\mcl X$ converges weakly to the law $\nu$ of $X$. That is, for any bounded continuous function $\phi: \mcl X \to \R$,
\[
	\lim_{\delta \to 0} \int_{\mcl X} \phi \; d\nu_\delta = \int_{\mcl X} \phi\;d\nu.
\]

Below, the paths $B_h$ are taken to be appropriately scaled simple random walks, and, as a consequence, $B_h$ converges in distribution to a Brownian motion $B$ (see for instance Billingsley \cite{Bill}). This corresponds above to $\mcl X = C([0,T];\R^m)$ and $\nu$ the Wiener measure on $\mcl X$.
\smallskip

Let $\lambda$, $\rho_h$, $B_h$, and $\mcl P_h$ be given, for some probability space $(\mcl A, \mcl G, \mbf P)$, by
\begin{equation} \label{A:simplerandomwalks}
\left\{
\begin{split}
	&\lambda = \frac{(\rho_h)^{3/4}}{h} \le \frac{\theta}{\nor{H'}{\oo}}, \quad M_h := \lfloor (\rho_h)^{-1/2} \rfloor,  \qquad \mcl P_h = \{t_n\}_{n =0}^N = \left\{ n \rho_h \wedge T \right\}_{n \in \N_0},\\
	%&\mcl P_h := \{t_n\}_{n =0}^N = \left\{ n \rho_h \wedge T \right\}_{n \in \N_0}, \\
	&\{\xi_n\}_{n=1}^\oo: \mcl A \to \{-1, 1\} \text{ are independent,} \qquad \mbf P(\xi_n = 1) = \mbf P(\xi_n = -1) = \frac{1}{2},  \\[1.5mm]
	&B(0) = 0, \quad\text{and} \quad \text{for } k \in \N_0, \; t \in [kM_h \rho_h , (k+1)M_h \rho_h), \\
	&B_h(t) = B_h(kM_h \rho_h) + \frac{ \xi_k }{\sqrt{M_h \rho_h} }(t - kM_h \rho_h).% \\
	%&\text{for } k \in \N_0, \; t \in [kM_h \rho_h , (k+1)M_h \rho_h).
\end{split}
\right.
\end{equation}

\begin{thm} \label{T:introLFdistributionBM}
	If $u_h$ is constructed using \eqref{E:introschemedef1} and \eqref{E:introLFfirstorder} with $B_h$ and $\mcl P_h$ as in \eqref{A:simplerandomwalks}, and $u$ is the solution of \eqref{E:simpleeqintrofirstorder}, then, as $h \to 0$, $u_h$ converges to $u$ in distribution.
\end{thm}

\subsection*{The need to regularize  the paths}  A short discussion follows about the necessity to consider regularizations $B_h$ of the continuous path $B$ in all of the results above. To keep the presentation simple, I  concentrate on the one-dimensional, pathwise Hamilton-Jacobi equation \eqref{E:simpleeqintrofirstorder}.
\smallskip

Consider the following naive attempt at constructing a scheme operator by setting
\begin{equation}
	\begin{split}
	S_h(t, s)u(x) &= u(x) + H \pars{ \frac{u(x+h) - u(x-h)}{2h} }(B(t) - B(s))\\
	& + \eps_h \pars{ \frac{ u(x+h) + u(x-h) - 2u(x)}{h^2} }(t-s).
	\end{split} \label{E:pathwiseguess}
\end{equation}

A simple calculation reveals that $S_h(t,s)$ is monotone for $0 \le t-s \le \rho_h$, if $\rho_h$ and $\eps_h$ are such that, for some $\theta \le 1$,
%\[
%	\eps_h = \frac{\theta h^2}{2(t-s)}
%\] 
%and
\begin{equation}\label{E:intropathCFL}
	\eps_h = \frac{\theta h^2}{2(t-s)} \ \ \text{and} \ \  \lambda = \max_{|t-s| \le \rho_h} \frac{ \mathrm{osc}(B, s,t) }{h} \le \lambda_0 = \frac{\theta}{\nor{H'}{\oo}}.
\end{equation}
\smallskip

On the other hand, for any $s,t \in [0,T]$ with $|s-t|$ sufficiently small, spatially smooth solutions $\Phi$ of \eqref{E:simpleeqintrofirstorder} have the expansion 
\begin{equation}\label{E:introsmoothexpansion}
	\begin{split}
	\Phi(x,t) &= \Phi(x,s) + H(\Phi_x(x,s))(B(t)-B(s)) \\
	&+ H'(\Phi_x(x,s))^2 \Phi_{xx}(x,s) (B(t)-B(s))^2 + O(\abs{ B(t) - B(s)}^3).
	\end{split}
\end{equation}

It follows that, if $0 \le t-s \le \rho_h$, there exists  $C > 0$ depending only on $H$ such that 
\begin{equation}\label{E:pathconsistencyconsequence}
	\begin{split}
	\sup_{\R} \abs{ S_h(t , s)\Phi(\cdot,s) -\Phi(\cdot,t)} &\le C \sup_{r \in [s,t]} \nor{D^2\Phi(\cdot,r)}{\oo} \pars{ \abs{B(t) - B(s)}^2 + h^2} \\
	&\le C \sup_{r \in [s,t]} \nor{D^2\Phi(\cdot,r)}{\oo} (1 + \lambda_0^2)h^2.
	\end{split}
\end{equation}
Therefore, in order for the scheme to have a chance of converging, $\rho_h$ should satisfy
\begin{equation}\label{E:introlimits}
	\lim_{h \to 0} \frac{h^2}{\rho_h} = 0.
\end{equation}
Both \eqref{E:intropathCFL} and \eqref{E:introlimits} can be achieved when $B$ is continuously differentiable or merely Lipschitz continuous by setting
\[
	\rho_h ={\lambda h}{ \nor{\dot B}{\oo}}^{-1}.
\]
More generally, if $B\in C^{0,\alpha}([0,T])$ with $\alpha > \frac{1}{2}$ and 
\begin{equation} \label{E:alphaCFL}
	(\rho_h)^\alpha = \frac{\lambda h}{ [W]_{\alpha,T}},
\end{equation}
then both \eqref{E:intropathCFL} and \eqref{E:introlimits} are satisfied, since 
\[
	\frac{h^2}{\rho_h} = \pars{ \frac{[B]_{\alpha,T} h^{2\alpha-1}}{\lambda} }^{1/\alpha} \xrightarrow{h\to 0} 0.
\]
However, this approach fails as soon as the quadratic variation path
\[
	\ip{B}_T := \lim_{\abs{ \mcl P} \to 0}\sum_{n=0}^{N-1} \abs{ B(t_{n+1}) - B(t_n)}^2
\]
is non-zero, as \eqref{E:intropathCFL} and \eqref{E:introlimits} together imply that $\ip{B}_T = 0$. This rules out, for instance, the case where $B$ is the sample path of a Brownian motion, or, more generally, any nontrivial semimartingale.
\smallskip

Motivated by the theory of rough differential equations, it is natural to explore whether the scheme operator \eqref{E:pathwiseguess} can be somehow altered to refine the estimate in \eqref{E:pathconsistencyconsequence}, potentially allowing \eqref{E:introlimits} to be relaxed and $\rho_h$ to converge more quickly to zero as $h \to 0^+$. 
\smallskip

More precisely, the next term in the expansion \eqref{E:introsmoothexpansion} suggests taking $B \in C^{0,\alpha}([0,T]; \R^m)$ with $\alpha > \frac13$, or,  more generally, $B$ with $p$-variation with  $p < 3$,   and defining
\begin{equation} \label{E:higherorderguess}
	\begin{split}
	&S_h(t, s)u(x) = u(x) + H \pars{ \frac{u(x+h) - u(x-h)}{2h} }(B(t) - B(s))\\
	& + \frac{1}{2} H'\pars{ \frac{u(x+h) - u(x-h)}{2h}}^2 \pars{ \frac{u(x+h) + u(x-h) - 2u(x) }{h^2}} \pars{ B(t) - B(s)}^2 \\
	& + \frac{\theta}{2} \pars{  u(x+h) + u(x-h) - 2u(x) }.
	\end{split}
\end{equation}
As can easily be checked, \eqref{E:higherorderguess} is monotone as long as \eqref{E:intropathCFL} holds, and 
\[
	\Lip(u) \le L, \quad \theta + \nor{H'}{\oo} \lambda^2 \le 1, \quad \text{and} \quad \lambda \le \frac{\theta}{\nor{H'}{\oo} \pars{ 1 + 2 L \nor{H''}{\oo}}}.
\]
The error in \eqref{E:pathconsistencyconsequence} would then be of order $h^2 + \abs{ B(t) - B(s)}^3$, which again leads to a requirement like \eqref{E:introlimits}. This seems to indicate that it is necessary to   incorporate higher order corrections in \eqref{E:higherorderguess} to deal with the second-order spatial derivatives of $u$. However, this will disrupt, in general, the monotonicity of the scheme, since it will no longer be possible to use discrete maximum principle techniques.
\smallskip

For this reason, it is more  convenient to concentrate on the more effective strategy of regularizing the path $B$. If $\{B_h\}_{h > 0}$ is a family of smooth paths converging uniformly, as $h \to 0$, to $B$, then $\ip{B_h}_T = 0$ for each fixed $h > 0$, and therefore, $B_h$ and $\rho_h$ can be chosen so that \eqref{E:intropathCFL} and \eqref{E:introlimits} hold for $B_h$ rather than $W$. 

\section{Homogenization}

I present a variety of results regarding the asymptotic properties, for small $\ep > 0$, of equations of the form
\begin{equation}\label{E:generalhomog}
		u^\eps_t + \sum_{i=1}^m H^i(Du^\eps, x/\eps) \dot \zeta^{\eps}_i = 0   \  \text{in } \   Q_\oo \quad  u^\eps(\cdot,0) = u_0. % \quad \text{in } \R^d.
\end{equation}
Many proofs and details are omitted here, and can be found in Seeger \cite{seegerh}.
\smallskip

Each Hamiltonian $H^i$ in \eqref{E:generalhomog} is assumed to have some averaging properties in the  variable $y = x/\eps$. The paths $\zeta^\eps = (\zeta^{\eps}_1,  \cdots, \zeta^{\eps}_m)$, which converge locally uniformly to some limiting path $\zeta \in C_0([0,\oo); \R^m)$, will be assumed to be piecewise $C^1$, although I  present some results where they are only continuous.
\smallskip

One motivation for considering such problems is to study general equations of the form
\begin{equation}\label{E:introscaling}
	%\left\{
	%\begin{split}
		u^\eps_t + \frac{1}{\eps^{\gamma}} H\pars{ Du^\eps, \frac{x}{\eps}, \frac{t}{\eps^{2\gamma}}} = 0 \ \ \text{in } \ \  Q_\oo \quad
		u^\eps(\cdot,0) = u_0.% \quad \text{in } \R^d.
%	\end{split}
%	\right.
\end{equation}
In addition to the averaging dependence on space, the Hamiltonian $H$ is assumed to have zero expectation, so that, on average, $u^\eps$ is close to its initial value $u_0$. The dependence on time, meanwhile, is assumed to be ``mixing'' with a certain rate, so that, with the scaling of the central limit theorem, $\eps^{-\gamma} H(\cdot,\cdot, t \eps^{-2\gamma})$ will resemble,  as $\eps \to 0$, to white noise in time. % as $\eps \to 0$. 
\smallskip

When $\gamma = 1$, \eqref{E:introscaling} arises naturally as a scaled version of 
\begin{equation}\label{E:unscaled}
	%\begin{dcases}
	u_t + H(Du, x, t) = 0 \   \text{in } \   Q_\oo \quad 
	u(\cdot,0) = \eps^{-1} u_0(\eps \cdot),  % & \text{in } \R^d.
	%\end{dcases}
\end{equation}
with $u$ and $u^\eps$  related by $u^\eps(x,t) = \eps u(x/\eps, t/\eps^2)$. 
\smallskip

Studying the $\eps \to 0$ limit of $u^\eps$ then amounts to understanding the averaged large space, long time behavior of solutions of \eqref{E:unscaled} with large, slowly-varying initial data. 
\smallskip

Although it is of interest to examine \eqref{E:introscaling} for different values of $\gamma$,  it turns out that  the nature of the limiting behavior does not change for different values of $\gamma$.  Hence, from a practical point of view, $\eps$ and $\delta = \eps^\gamma$ can be viewed as small, independent parameters. It should be, however,  noted that for technical reasons, some results can only be proved under a mildness assumption on the approximate white noise dependence, which translates to a smallness condition on $\gamma$.
\smallskip

The Hamiltonians considered  in  \eqref{E:introscaling}  have the form
\begin{equation}\label{A:linearnoise}
	H(p,y,t) = \sum_{i=1}^m H^i(p,y) \xi_i(t),
\end{equation}
where the random fields $\xi_i: [0,\oo) \to \R$ are defined on a probability space $(\Omega, \mcl F, \P)$ and are assumed to be mixing with rate $\rho$ as explained below.
\smallskip

For $0 \le s \le t \le \oo$, consider the sigma algebras $\mcl F^i_{s,t} \subset \mcl F$ generated by $\{\xi_i(r)\}_{r \in [s,t]}$.  The mixing rate is then defined by 
\begin{equation}\label{mixingrate}
	\rho(t) = \max_{i=1,2,\ldots,m} \sup_{s \ge 0} \sup_{A \in \mcl F^i_{s+t, \oo}} \sup_{ B \in \mcl F^i_{0,s} } \abs{ \P(A \mid B) - \P(A)}.
\end{equation}
The quantitative mixing assumptions for the $\xi^i$ are that
\begin{equation} \label{A:xi}
	\left\{
	\begin{split}
	&t \mapsto \xi_i(t) \text{ is stationary,} 
	\quad \rho(t) \xrightarrow{t \to \oo} 0, \quad \int_0^\oo \rho(t)^{1/2}\;dt < \oo,\\
	&\E[\xi_i(0)] = 0, \text{ and } \E[\xi_i(0)^2] = 1.
	\end{split}
	\right.
\end{equation}

Above stationarity means that
\[
	(\xi(s_1), \xi(s_2), \ldots, \xi(s_M)) \quad \text{and} \quad (\xi(s_1 + t), \xi(s_2 + t), \ldots, \xi(s_M + t))
\]
have the same joint distribution for any choice of $s_1, s_2, \ldots, s_M \in [0,\oo)$ and $t \ge - \min_j s_j$. %However, the conclusions reached above are unchanged if $\xi$ is stationary only with respect to integer shifts. This will be the case if, for instance, $\zeta$ is a linearly-interpolated random walk, in which case the convergence of $\zeta^\delta$ to a Brownian motion is a consequence of Donsker's invariance principle.

It follows from  the ergodic theorem, the stationarity  and the centering assumptions  that 
%\[
%	\zeta^i(t)= \int_0^t \xi^i(s)\;ds 
%\]
%then 
\[\lim_{\delta \to 0} \delta \int_0^{\frac{t}{\delta}} \xi_i(s)\;ds = 0.\] 

The properties of  the long time fluctuations of $\zeta=\int_0^t \xi(s)ds$ around $0$ can be studied  using the central limit theorem scaling. Indeed setting $\zeta^{\delta}_i(t) = \delta \zeta_i(t/\delta^2)$, it is well-known that,  as $\delta \to 0$, $\zeta^{\delta}_i$ converges in distribution and locally uniformly to a standard Brownian motion. Indeed, with $\delta = \eps^\gamma$,  \eqref{E:introscaling} is then a specific form of \eqref{E:generalhomog}.
%\smallskip

%The assumption of stationarity in \eqref{A:xi} means that
%\[
%	(\xi(s_1), \xi(s_2), \ldots, \xi(s_M)) \quad \text{and} \quad (\xi(s_1 + t), \xi(s_2 + t), \ldots, \xi(s_M + t))
%\]
%have the same joint distribution for any choice of $s_1, s_2, \ldots, s_M \in [0,\oo)$ and $t \ge - \min_j s_j$. However, the conclusions reached above are unchanged if $\xi$ is stationary only with respect to integer shifts. This will be the case if, for instance, $\zeta$ is a linearly-interpolated random walk, in which case the convergence of $\zeta^\delta$ to a Brownian motion is a consequence of Donsker's invariance principle.
%
\subsection*{The difficulties and general strategy} Here I  discuss some of the difficulties in the study of the $\ep\to 0$ behavior of \eqref{E:generalhomog} and the strategies that can be used  to overcome them. To keep things simple, I  only consider Hamiltonians  that are periodic in space.
\smallskip

The starting (formal) assumption  is  that the noise is ``mild'' enough to allow for averaging behavior in space, and therefore, $u^\eps$ is closely approximated by a solution $\oline{u}^\eps$ of an equation of the form $\oline{u}^\eps_t + \oline{H}^\eps(D\oline{u}^\eps,t) = 0$. 
\smallskip

More precisely, following  the standard strategy of  the homogenization theory, it is assumed that there exists some auxiliary function $v: \T^d \times [0,\oo) \to \R$, so that $u^\eps$ has the formal expansion
\[
	u^\eps(x,t) \approx \oline{u}^\eps(x,t) + \eps v(x/\eps,t).
\]
An asymptotic analysis yields that, for fixed $p \in \R^d$ (here, $p = D\oline{u}^\eps(x,t)$ and $y = \frac{x}{\eps}$), $v$ solves the so called ``cell problem'' 
\begin{equation}\label{E:generalcellproblem}
	\sum_{i=1}^m H^i(D_y v + p,y)\xi_i = \oline{H}(p, \xi) \   \text{in }. \ \R^d, 
\end{equation}
where the fixed parameter $\xi \in \R^m$ stands in place of the mild white noise $\eps^{-\gamma}\xi(t/\eps^{2\gamma})$. %The equation \eqref{E:generalcellproblem} is known as the ``cell problem,'' and,
\smallskip

It is standard  the theory of periodic homogenization of Hamilton-Jacobi equations that, under the right conditions, there is a unique constant $\oline{H}(p, \xi)$ for which \eqref{E:generalcellproblem} has periodic solutions, which are called ``correctors.''
\smallskip

Taking this fact for granted for now and always arguing formally yields that $u^\eps$ will be closely approximated by $\oline{u}^\eps$ which solves 
\begin{equation} \label{E:generaletaeq}
	\oline{u}^\eps_t + \frac{1}{\eps^\gamma} \oline{H}\pars{ D \oline{u}^\eps,\xi \pars{ \frac{t}{\eps^{2\gamma} } } }= 0 \   \text{in } \  Q_\oo \quad  \oline{u}^\eps(\cdot,0) = u_0.% \quad \text{in } \R^d.
\end{equation}
Note that, in deriving \eqref{E:generaletaeq}, it was used  that $\xi \mapsto \oline{H}(\cdot,\xi)$ is positively homogenous, which follows  from multiplying \eqref{E:generalcellproblem} by a positive constant and using the uniqueness of the right-hand side. 
\smallskip

If 
\begin{equation}\label{centeredeffective}
	\E \left[ \oline{H}(p,\xi(t)) \right]= 0 \ \  \text{for all } \ \ p \in \R^d,
\end{equation}
then the solution of \eqref{E:generaletaeq} with $u_0(x) =\langle p_0, x\rangle$, which is given by
\[
	\oline{u}^\eps(x,t) = \langle p_0,x \rangle - \frac{1}{\eps^{2\gamma}} \int_0^t \oline{H}\pars{ p_0, \xi \pars{ \frac{s}{\eps^{2\gamma}} } }ds,
\]
converges, as $\eps \to 0$ and  in distribution,  to $p_0 \cdot x + \sigma(p_0) B(t)$, where $B$ is a standard Brownian motion and
\[
	\sigma(p_0) = \pars{\E\left[ \oline{H}(p_0, \xi(0))^2 \right]}^{1/2}.
\]

Due, however, to the nonlinearity of the map $\xi \mapsto \oline{H}(\cdot,\xi)$ and the difficulties associated with the ``rough'' pathwise solutions, it is not clear how to study the  \eqref{E:generaletaeq} for an arbitrary $u_0 \in UC(\R^d)$. It turns out that the answers are subtle, and, in the multiple path case considered below, depend strongly on the nature of the mixing field $\xi$. 
\smallskip

When $m = 1$, the characterization of $\oline{H}(p,\xi)$ reduces to the study of the two Hamiltonians 
\[
	\oline{H}(p) = \oline{H}(p,1) \quad \text{and} \quad \oline{(-H)}(p) = \oline{H}(p,-1).
\]

Then \eqref{E:generaletaeq}  takes the form
\begin{equation}\label{E:withvariation}
%	\left\{
%	\begin{split}
	\oline{u}^\eps_t +\frac{1}{\eps^\gamma} \oline{H}^1(D \oline{u}^\eps) \xi\pars{ \frac{t}{\eps^{2\gamma}}} + \frac{1}{\eps^\gamma} \oline{H}^2(D \oline{u}^\eps)\abs{ \xi \pars{ \frac{t}{\eps^{2\gamma}}}} = 0 \  \text{in } \ Q_T \quad \oline{u}^\eps(\cdot,0) = u_0, %\R^d \times (0,\oo) \quad \text{and}\\
	%&\oline{u}^\eps(\cdot,0) = u_0 \quad \text{in } \R^d,
%	\end{split}
%	\right.
\end{equation}
where
\[
	\oline{H}^1(p) = \frac{ \oline{H}(p) - \oline{(-H)}(p)}{2} \quad \text{and} \quad \oline{H}^2(p) = \frac{ \oline{H}(p) + \oline{(-H)}(p)}{2}.
\]

Note that $\oline{H}^2=0$ if and only if
\begin{equation}\label{consistentH}
	\oline{(-H)} = -\oline{H},
\end{equation}
and, moreover, that \eqref{centeredeffective} is equivalent to \eqref{consistentH}  when $m = 1$.
\smallskip

Since  \eqref{E:generalcellproblem} is interpreted in the viscosity solution sense, it is not possible to multiply the equation by $-1$, and so \eqref{consistentH} is not only not obvious, but actually false in general.% as can be seen in the following example. %where,  for some $p_0 \in \R^d$, $\oline{(-H)}(p_0) \ne - \oline{H}(p_0)$. 
\smallskip

Indeed, assume that, for some $p_0 \in \R^d$, $\oline{(-H)}(p_0) \ne - \oline{H}(p_0)$. Then  $\oline{u}^\eps$ with $u_0(x) =\langle p_0, x\rangle$ is given by 
\begin{align*}
	\oline{u}^\eps(x,t) = \langle p_0, x\rangle &- \eps^\gamma \frac{ \oline{H}(p_0) - \oline{(-H)}(p_0)}{2} \zeta \pars{ \frac{t}{\eps^{2\gamma}} }
	 - \eps^\gamma \frac{ \oline{H}(p_0) + \oline{(-H)}(p_0)}{2} \int_0^{t/\eps^{2\gamma}} \abs{ \xi\pars{ \frac{s}{\eps^{2\gamma} } }}\;ds,
\end{align*}
and, hence, 
\[
	\eps^\gamma  \oline{u}^\eps(x,t)  \xrightarrow{\eps \to 0} - \frac{ \oline{H}(p_0) + \oline{(-H)}(p_0)}{2} \E\abs{ \xi(0)} t \quad \text{in distribution.}
\]

On the other hand, if \eqref{consistentH} holds, then \eqref{E:generaletaeq} becomes
\begin{equation} \label{E:etaeq}
	\oline{u}^\eps_t + \frac{1}{\eps^\gamma}\oline{H}(D \oline{u}^\eps) \xi\pars{ \frac{t}{\eps^{2\gamma}}} = 0 \ \ \text{in } Q_T \quad \oline{u}^\eps(\cdot,0) = u_0, % \quad \text{in } \R^d,
\end{equation}

and the determination of whether or not $\oline{u}^\eps$ has a limit depends on the properties of the effective Hamiltonian $\oline{H}$, and, in particular, whether or not it is the difference of two convex functions.

\subsection* { The single-noise case.} I state next some results about 
\begin{equation}\label{E:introsinglenoise}
	u^\eps_t + \frac{1}{\eps^\gamma} H\pars{Du^\eps,\frac{x}{\eps}} \xi\pars{\frac{t}{\eps^{2\gamma}}} =0 \ \  \text{in } \ \ Q_T \quad u^\eps(\cdot,0) = u_0.
\end{equation}
As suggested in the previous subsection, the fact that there is only one source of noise simplifies the structure of the problem. Consequently, the results are more comprehensive than in the multiple-path setting.
\smallskip

It is assumed that
\beq\label{paris40}
H \in C(\R^d \times \R^d) \ \text{is convex and coercive in the gradient variable}.
\end{equation}

%The Hamiltonian $H \in C(\R^d \times \R^d)$ is assumed to be convex and coercive in the gradient variable. 
The convexity assumption is important for two reasons. It guarantees that the consistency condition \eqref{consistentH} holds, and it also implies strong path-stability estimates for the solutions. The latter  were already alluded to earlier in the notes, in the section on the comparison principle for equations with convex, spatially-dependent Hamiltonians.
\smallskip

Regarding the spatial environment, the results are general enough to allow for a variety of different assumptions. Here, I  list two well-studied examples. 
\smallskip

The first possible self-averaging assumption is that
\begin{equation}\label{A:Hperiodic}
	y \mapsto H(p,y) \text{ is $\Z^d$-periodic}.
\end{equation}
The periodic homogenization of (time-homogenous) Hamilton-Jacobi equations has a vast literature going back to Lions, Papanicolaou, and Varadhan \cite{LPV} and Evans \cite{E1,E2}.
\smallskip

Another type of averaging dependence, which in general is more physically relevant, is stationary-ergodicity. In this setting, the  Hamiltonians $H = H(p,x,\omega)$ are defined on a probability space $(\mbf \Omega, \mbf F)$ that is independent of the random field $\xi$ and is equipped with %. %Consider the group of translation operators $T_z: \mbf \Omega \to \mbf \Omega$ defined  
a  group of translation operators $T_z: \mbf \Omega \to \mbf \Omega$ such that $H(\cdot,T_z y) = H(\cdot,y+z)$. It is assumed that $\{T_z\}_{z \in \R^d}$ is stationary and ergodic, that is, 
% a family of maps $T_z:\Omega \to \Omega$ such that % is measure-preserving and ergodic, that is,
\begin{equation}\label{A:Hstatergod}
	\left\{
	\begin{split}
		&\mbf P = \mbf P \circ T_z \text{ for all $z \in \R^d$, and} \\
		&\text{if } E \in \mbf F \text{ and }T_z E =  E \text{ for all } z \in \R^d, \text{ then } \mbf P[E] = 1 \text{ or } \mbf P[E] = 0.
	\end{split}
	\right.
\end{equation}
In the time-inhomogenous setting, this homogenization problem was studied by Souganidis \cite{souganidishomo} and Rezakhanlou and Tarver \cite{RT}. %This continues to be an active research field, especially with regards to error estimates, qualitative behavior of the effective equation, and counterexamples.
\smallskip

The first result is stated next.
\begin{thm}\label{T:singlenoise}
	There exists a Brownian motion $B: [0,\oo) \to \R$ such that, as $\eps \to 0$, $(u^\eps,\zeta^\eps)$ converges in distribution to $(\oline{u}, B)$ in $BUC(\R^d \times [0,\oo)) \times C([0,\oo))$, where $\oline{u}$ is the pathwise viscosity solution of
	\begin{equation}\label{E:singlenoiseeq}
		d\oline{u} + \oline{H}(D\oline{u}) \cdot dB =0 \  \text{in } \ Q_\oo \quad \oline{u}(\cdot,0) = u_0. % \quad \text{in } \R^d.
	\end{equation}
\end{thm}
%The topology in $BUC(\R^d \times [0,\oo))$ and $C([0,\oo))$ is that of local uniform convergence. 
%\smallskip

Since $\delta B(t/\delta^2)$ equals $B(t)$ in distribution, it 
%Observe that the singular random field $\xi(t) = dB(t)$ is strongly mixing and 
%\[
%	\zeta(t) = \int_0^t \xi(s)\;ds = \int_0^t dB(s) = B(t),
%\]
%so that $t\mapsto \zeta^\delta(t) = \delta B(t/\delta^2)$ is equal in distribution to $B$ by the scaling properties of Brownian motion. It is therefore 
is also an interesting question to study the limiting behavior of 
\begin{equation}\label{E:singleBrowniannoise}
	du^\eps + H(Du^\eps, x/\eps) \cdot dB = 0 \  \text{in } \ Q_\oo \quad u^\eps(\cdot,0) = u_0. % \quad \text{in } \R^d.
\end{equation}

\begin{thm}\label{T:singleBrowniannoise}
	In addition to the hypotheses of Theorem \ref{T:singlenoise}, assume that the comparison principle holds for \eqref{E:singleBrowniannoise}. Then, with probability one, as $\eps \to 0$, the solution $u^\eps$ of \eqref{E:singleBrowniannoise} converges locally uniformly to the solution of \eqref{E:singlenoiseeq}.
\end{thm}

The final  remark is  that the  theorems above can be applied to a variety of other settings like, for instance, the homogenization of 
\[
	u^\eps_t + H\pars{Du^\eps,x,\frac{x}{\eps}} \dot \zeta^\eps(t)=0 \  \text{in } \ Q_\oo  \quad u^\eps(\cdot,0) = u_0
\]
with $(\zeta^\eps)_{\eps > 0}$ any collection of paths converging locally uniformly and almost surely (or in distribution) to a Brownian motion or other stochastic process, and with the dependence of $H$ on the fast variable being, for instance, periodic, quasi-periodic, or stationary-ergodic.

\subsection*{The multiple-noise case.} Since in this setting the results so far are less general and quite technical, I only present an overview here. Details and more results can be found in a forthcoming work of Seeger \cite{seegerh2}.
\smallskip

The problem is the the behavior of equations like 
\begin{equation}\label{E:multiplepaths}
	u^\eps_t + \frac{1}{\eps^\gamma} \sum_{i=0}^m H^i(Du^\eps, x/\eps) \xi_i(\frac{t}{\eps^{2\gamma}})=0 \ 
	\text{in } \     Q_\oo  \quad  u^\eps(\cdot,0) = u_0, % \quad \text{in } \R^d.
\end{equation}
where, for each $i = 0, \ldots, m$, $\xi^i$ is a mixing field satisfying \eqref{A:xi}. More assumptions on the Hamiltonians and the paths will need to be made later. % although we will be making more specific assumptions about these and the Hamiltonians.
\smallskip

To simplify the presentation, here I only consider the periodic setting \eqref{A:Hperiodic}. It turns out that, under appropriate conditions on the $H^i$'s which are made more specific below, %the following statements will hold.
%\smallskip
for every $p \in \R^d$ and $\xi \in \R^m$, there exists a unique constant $\oline{H}(p,\xi)$ such that the cell problem 
\begin{equation}\label{E:multiplecellproblem}
	\sum_{i=1}^m H^i(p + D_y v, y)\xi^i = \oline{H}(p, \xi)
\end{equation}
admits periodic solutions $v: \T^d \to \R$. Moreover, $\xi \mapsto \oline{H}(p,\xi)$ is positively homogenous, and
\begin{equation}\label{centering}
	\E[\oline{H}(p,\xi(0)] = 0 \ \text{for all } \ p \in \R^d.
\end{equation}

Using error estimates for the theory of periodic homogenization of Hamilton-Jacobi equations, it is possible to show that $u^\eps$ is closely approximated by the solution $\oline{u}^\eps$ of
\begin{equation}\label{E:nonlinearmixing}
	\oline{u}^\eps_t + \frac{1}{\eps^\gamma} \oline{H} \pars{ D\oline{u}^\eps, \xi\pars{ \frac{t}{\eps^{2\gamma}}}} = 0 \  \text{in } \  Q_\oo \quad  \oline{u}^\eps(\cdot,0) = u_0.
\end{equation}
The limiting behavior of  \eqref{E:nonlinearmixing} is well understood if $u_0(x) =\langle p_0,  x\rangle$ for some fixed $p_0 \in \R^d$. Indeed, in view of the mixing properties of $\xi$ and the centering property \eqref{centering}, there exists a Brownian motion $B$ such that, as $\eps \to 0$, $\oline{u}^\eps$ converges locally uniformly in distribution to
\[
	\langle p_0,  x\rangle + \E \left[ \oline{H}(p_0,\xi(0))^2 \right]^{1/2} B(t).
\]
%where $B$ is a Brownian motion. 

I comment next about the limit of $\oline{u}^\eps$ for arbitrary initial data $u_0$. The goal is to show that, under assumptions on the Hamiltonians and mixing fields, there exists $M \ge 1$  and, for each $j = 1, \ldots, M$, an effective Hamiltonian $\oline{H}^j: \R^d \to \R$ which is the difference of two convex functions, and a Brownian motion $B^j$ such that, as $\eps \to 0$ and in distribution, $\oline{u}^\eps$ and, therefore, $u^\eps$ converges in $BUC(Q_T)$ to the pathwise viscosity solution $\oline{u}$ of
\begin{equation}\label{E:effectivemultipleeq}
	d\oline{u} + \sum_{j=1}^M \oline{H}^j(D \oline u) \cdot dB_j = 0 \   \text{in }  \  Q_\oo \quad  \oline{u} = u_0. % \quad \text{in } \R^d.
\end{equation}

Although at first glance, the nature of the problem is similar to the single path case, there are  some fundamental differences. Most importantly, the deterministic effective Hamiltonians $\{\oline{H}^j\}_{j=1}^M$, and even the number $M$, depend on the particular law of the mixing field $\xi$. 
\smallskip

Next I introduce  some further assumptions that give rise to a rich class of examples and results.
\smallskip

As far as the Hamiltonians $(H_i, \ldots, H_m)$ are concerned, it is assumed that
\begin{equation}\label{A:Hs}
	\left\{
	\begin{split}
		&H^i \in C^{0,1}(\R^d \times \T^d), \\
		&p \mapsto H^1(p,\cdot) + \sum_{i=2}^m H^i(p, \cdot)\xi_i \text{ is convex for all } \xi_2, \ldots, \xi_m \in \{-1,1\}, \text{ and}\\
		&\lim_{|p| \to +\oo} \inf_{y \in \T^d} \pars{ H^1(p,y) - \sum_{i=2}^m \abs{ H^i(p,y)}} = +\oo.
	\end{split}
	\right.
\end{equation}
As a consequence, the cell problem \eqref{E:multiplecellproblem} is solvable for all $p \in \R^d$ and $\xi \in \{-1,1\}^m$, and furthermore, $p \mapsto \oline{H}(p,1,\xi)$ is convex and $\xi \mapsto \oline{H}(p,\xi)$ is homogenous, that is, for all $\lambda \in \R$ and $\xi \in \{-1,1\}^m$,
\begin{equation}\label{homogenouseffective}
	\oline{H}(\cdot,\lambda \xi) = \lambda \oline{H}(\cdot,\xi).
\end{equation}
The mixing fields are assumed to be, for $i = 1,\ldots, m$, of the form 
\begin{equation}\label{A:checkerboards}
	\left\{
	\begin{split}
	&\xi_i = \sum_{k=0}^{\oo} X^i_k \ind_{(k,k+1)} \quad \text{where} \\
	&\left(X^i_k\right)_{i=1,2,\ldots, m, \; k = 0,1,\ldots} \quad \text{are independent Rademacher random variables}.
	\end{split}
	\right.
\end{equation}
In particular, if
\begin{equation}\label{zetas}
	\xi^{\eps}_i(t) = \frac{1}{\eps^{\gamma}} \xi_i (t/\eps^{2\gamma}) \quad \text{and} \quad \zeta^{\eps}_i(t) = \int_0^t \xi^{\eps}(s)_i\;ds,
\end{equation}
then each $\zeta^{i,\eps}$ is a scaled, linearly-interpolated, simple random walk on $\Z$, and there exists an $m$-dimensional Brownian motion $(B_1, \ldots, B_m)$, such that, in distribution,
\[
	(\zeta^{1,\eps}, \zeta^{2,\eps}, \ldots, \zeta^{m,\eps}) \xrightarrow{\eps \to 0} (B_1, \ldots, B_m) \ \text{in} \ C([0,\oo); \R^m).
\]
%where $(B_1, \ldots, B_m)$ is an $m$-dimensional Brownian motion. 

Consider the sets of indices 
\[
	\left\{
	\begin{split}
	&\mcl A^m := \{ \mbf j = (j_1, \ldots, j_l) : j_i \in \{1,  \ldots, m\}, \; j_1 < \cdots < j_l \} \ \text{with} \ 
	l=\abs{ \mbf j} = \abs{ (j_1, j_2, \ldots, j_l)} \\[1.5mm]
	&\mcl A^m_O := \{ \mbf j \in \mcl A^m : \abs{ \mbf j} \text{ is odd} \},
	\end{split}
	\right.
\]
noting  that $\#\mcl A^m = 2^m - 1$ and $\#\mcl A^m_0 = 2^{m-1}$.
\smallskip

For any $\mbf j = (j_1,j_2,\ldots, j_l) \in \mcl A^m$, define
\begin{equation}\label{effectiveingredients}
	\left\{
	\begin{split}
		&\xi_{\mbf j} := \xi_{j_1}  \cdots \xi_{j_l} \ \  \text{for } \ \ \xi = (\xi_1, \ldots, \xi_m) \in \{-1,1\}^m,\\
		&\oline{H}^{\mbf j}(p) := \sum_{\xi \in \{-1,1\}^m} 2^{-m} \oline{H}(p,\xi) \xi_{\mbf j},\\
		&X^{\mbf j}_k := X^{j_1}_k X^{j_2}_k \cdots X^{j_l}_k,\\
		&\zeta_{\mbf j}(0) := 0, \quad \dot \zeta_{\mbf j} = \sum_{k=0}^{\oo} X^{\mbf j}_k \ind_{(k,k+1)}, \quad \text{and} \quad \zeta_{\mbf j}^\eps(t) = \eps^\gamma \zeta_{\mbf j}(t/\eps^{2\gamma}),
	\end{split}
	\right.
\end{equation}
and observe that, for each $\mbf j \in \mcl A^m_0$, $\oline{H}^{\mbf j}$ is a difference of convex functions. Note also that,  if $\abs{ \mbf j}$ is even, then the homogeneity property \eqref{homogenouseffective} implies that $\oline{H}^{\mbf j} = 0$.
\smallskip

The following is true.
\begin{thm}\label{T:generalmultpaths}
	Assume that $ \gamma \in (0,1/6)$, $u_0 \in C^{0,1}(\R^d)$, \eqref{A:Hs}, and \eqref{A:checkerboards}, and let $u^\eps$ be the solution of \eqref{E:multiplepaths}. Then there exist $2^{m-1}$ independent Brownian motions $\{B^{\mbf j}\}_{\mbf j \in \mcl A^m_o}$, such that, in distribution,
	\[
		\pars{ u^\eps, \{\zeta^{\mbf j,\eps} \}_{\mbf j \in \mcl A^m_o} } \xrightarrow{\eps \to 0} \pars{ \oline{u}, \{B^{\mbf j} \}_{\mbf j \in \mcl A^m_o} } \  \text{in } \  BUC(Q_T) \times C\pars{[0,T]; \R^{2^{m-1}} },
	\]
	where $\oline{u}$ is the stochastic viscosity solution of
	\begin{equation}\label{E:effectivemultiple}
		d\oline{u} + \sum_{\mbf j \in \mcl A^{m}_0} \oline{H}^{\mbf j}(D \oline{u}) \cdot  dB_{\mbf j} = 0 \  \text{in } \  Q_\oo \quad  \oline{u}(\cdot,0) = u_0.% \quad \text{in } \R^d.
	\end{equation}
\end{thm}

The result relies on the fact that, in view of  the assumptions on the mixing fields $\xi_i$, which take their values only in $\{-1,1\}$, the general effective Hamiltonian $\oline{H}(p,\xi)$ can be decomposed using a combinatorial argument. 
\smallskip

As already mentioned, the above theorem covers only some of the possible homogenization problems that can be studied in the multiple-noise case. In particular, it is shown in \cite{seegerh2} that the limiting equation depends on the law of the mixing field $\xi$. This is in stark contrast to the single-noise case, where the limiting equation is independent of the mild-noise approximation.

\section{stochastically perturbed reaction-diffusion equations and front propagation}
%%%%%%%%

I discuss here a result of Lions and Souganidis \cite{lsallencahn} about  the onset of fronts in the long time  and large space asymptotics of  bistable reaction-diffusion equations which are  additively perturbed by small relatively smooth (mild) stochastic in time forcing. The prototype problem is the so called stochastic Allen-Cahn equation.
%which is  additively perturbed by small relatively smooth (mild) stochastic in time forcing. 
The interfaces evolve with curvature dependent normal velocity which is additively perturbed by time white noise.
No regularity assumptions are made on  the fronts.  
%We also discuss the need to use approximations of the white noise by showing that perturbations with small time  white noise ``destroys'' the stability of the equilibria. 
The results can be extended to more complicated equations with anisotropic diffusion, drift and reaction which may be periodically oscillatory in space. To keep the ideas  simple, in this section I concentrate on the classical Allen-Cahn equation. 
\smallskip

The goal is to study  the behavior, as $\ep\to 0$, of the parabolically rescaled Allen-Cahn equation 
% \begin{equation}\label{allencahn}
% \begin{cases}
%u^\ep_t -\Delta u^\ep + \dfrac{1}{\ep^2} (f(u^\ep) + \ep \dot B^\ep(t,\omega))=0  \ \text{in} \ \R^d\times (0,\infty),\\[1mm]
%u^\ep(\cdot,0)=u_0.
%\end{cases}
%\end{equation} 
\begin{equation}\label{allencahn}
u^\ep_t -\Delta u^\ep + \dfrac{1}{\ep^2} (f(u^\ep) - \ep \dot B^\ep(t,\omega))=0  \ \text{in} \ Q_\oo \quad 
u^\ep(\cdot,0)=u^\ep_0,
\end{equation}
%that $u^\ep(\cdot,0)$ is independent of $\ep$ is made only for simplicity, and
%We assume that $f$ is the derivative of a double well potential with wells of equal depth at, for definiteness, $\pm 1$ and in between maximum at $0$, that is  
%f\in C^2(\R^d;\R)$ has two stable equilibria and one unstable, which for simplicity are taken to be $\pm 1$ and $0$ respectively and is balanced, 
where, $f\in C^2(\R^d;\R)$  is such that
\begin{equation}\label{nonl}
\begin{cases} f(\pm 1)=f(0)=0, f'(\pm 1)>0, \ f'(0)<0 \\[1mm]  f >0 \ \text{in} \ (-1,0), \ f<0 \ \text{in}  \ (0,1), \ \text{and} \ \displaystyle\int_{-1}^{+1} f(u)du=0,
\end{cases}
\end{equation}
that is, $f$ is the derivative of a double well potential with wells of equal depth at, for definiteness, $\pm 1$ and in between maximum at $0$, 
\vskip-.125in
\begin{equation}\label{takis1}
B^\ep (\cdot,\omega) \in C^2([0,\infty);\R) \ \text{is an a.s. mild approximation of the Brownian motion}  \ B(\cdot, \omega),
\end{equation}
that is,  a.s. in $\omega$ and locally uniformly $[0,\infty)$, 
\begin{equation}\label{takis1.1}
\lim_{\ep \to 0}B^\ep(t, \omega)=B, \ B^\ep(0,\omega)=0, %\\[1mm]
 \ \text{and} \  \lim_{\ep \to 0} \ep |\ddot B^\ep (t, \omega)| =0, 
\end{equation}
and there exists an open $\mathcal {O}_0\subset \R^d$ such that 
\begin{equation}\label{takis2}
\begin{cases}
\mathcal {O}_0 =\{x\in \R^d: u^\ep_0(x) >0\},  \ \R^d\setminus \overline {\mathcal {O}_0}=\{x\in \R^d: u^\ep_0(x) <0\}, \ \text{and} \\[1mm] \Gamma_0=\partial \mathcal {O}_0=\partial (\R^d\setminus \overline{ \mathcal {O}_0})=\{x\in \R^d: u^\ep_0(x)=0\}.
\end{cases}
\end{equation}
Although it is not stated explicitly,  it assumed that there exists an underlying probability space, but, for ease of the notation, we omit the dependence on $\omega$ unless  necessary.
\smallskip

Here are two classical examples of mild approximations. The first is  the convolution $B^\ep(t)=B\star \rho^\ep (t)$, where  $\rho^\ep (t)= \ep^{-\gamma} \rho(\ep^{-\gamma} t)$ with $\rho \in C^\infty$, even and  compactly supported in $(-1,1)$, $\int \rho (t) dt=1$ and $\gamma \in (0,1/2)$. The second is $\dot B^\ep(t)=\ep^{-\gamma} \xi (\ep^{-2\gamma}t)$, where $\xi(t)$  is a stationary, strongly mixing,  mean zero stochastic process such that $\max(|\xi|, |\dot \xi |)\leq M$ and $\gamma \in (0,1/3)$. I refer to \cite{kunita} %We refer to Ikeda and Watanabe~\cite{ik} 
for a discussion. 
\smallskip

Next I use the notion of stochastic viscosity solutions and the  level set approach to describe  the  generalized evolution (past singularities) of a set with normal velocity
\begin{equation}\label{takis12}
V=-\text{tr} [Dn] \ dt + d\zeta,
\end{equation}
for some  a continuous path $\zeta \in C_0([0,\infty);\R)$. Here  $n$ is the external normal to the front and, hence, $\text{tr}[Dn]$ is the mean curvature. 
\smallskip

Given a triplet   $(\mathcal {O}_0, \Gamma_0, \R^d\setminus \overline {\mathcal {O}_0})$ with $\mathcal {O}_0\subset \R^d$ open,  we say that the sets $(\Gamma_t)_{t> 0}$ move with normal velocity  \eqref{takis12}, if, for each $t> 0$,  there exists a triplet  $(\mathcal {O}_t, \Gamma_t, \R^d\setminus \overline {\mathcal {O}_t})$, with   $\mathcal {O}_t\subset \R^d$ open, such that 
\begin{equation}\label{takis10}
\begin{cases}
{\mathcal O}_t =\{x\in \R^d: w(x, t) > 0\},  \ \R^d\setminus \overline {\mathcal {O}_t}=\{x\in \R^d: w(x,t) <0\}, \ \text{and} \\  \Gamma_t=\{x\in \R^d: w(x,t)=0\},
\end{cases}
\end{equation}
where  $w \in \text{BUC}(\R^d\times [0,\infty))$ is the unique stochastic (pathwise) solution
of the level-set initial value pde
\begin{equation}\label{ivp1000}
dw=(I - \widehat {Dw} \otimes \widehat {Dw} ): D^2w  -  |Dw|\cdot d\zeta \   \text{in} \  Q_\oo \quad  w(\cdot,0)=w_0,
\end{equation}
 with $\hat p:=p/|p|$ and $w_0\in \text{BUC}(\R^d)$ such that 
 \begin{equation}\label{takis11}
\begin{cases}
\mathcal {O}_0 =\{x\in \R^d: w_0(x) >0\},  \ \R^d\setminus \overline {\mathcal {O}_0}=\{x\in \R^d: w_0(x)<0\},
 \ \text{and} \\
\Gamma_0=\{x\in \R^d: w_0(x)=0\}. 
\end{cases}
\end{equation}
%Above  $\zeta$ can be an arbitrary continuous function, in which case ``$\cdot$'' means multiplication.  When, however, $\zeta$ is a Brownian oath, ``$\circ$'' should be interpreted as the classical Stratonovich differential. % $\circ$.
%\smallskip

The properties of \eqref{ivp1000} are used here to adapt the approach  introduced in Evans, Soner and Souganidis~\cite{evanssonersouganidis}, Barles, Soner and Souganidis~\cite{barlessonersouganidis},  and Barles and Souganidis \cite{barlessouganidis} to study the onset of moving fronts in the asymptotic limit of reaction-diffusion equations and interacting particle systems with long range interactions.  This methodology allows to prove global in time asymptotic results and is not restricted to smoothly evolving fronts. 
\smallskip

The main result of the paper is stated next.
\begin{thm}\label{main}
Assume \eqref{nonl},   \eqref{takis1}, \eqref{takis1.1}, \eqref{takis2}, and  let $u^\ep$ be the solution of   \eqref{allencahn}.
% and consider the sets $({\mathcal O}_0, \R^d\setminus \overline{{\mathcal O}_0}, \Gamma_0)$ defined in \eqref{takis2}. 
There exists $\alpha_0 \in \R$ such that, % with initial value $u_0$ satisfying \eqref{takis2}, 
if $w$ is the solution of \eqref{ivp1000} with $w_0$ satisfying \eqref{takis11} and $\zeta\equiv \alpha_0 B$, where $B$ is a standard Brownian path, 
%and $({\mathcal O}_0, \R^d\setminus \overline{{\mathcal O}_0}, \Gamma_0)$ as in \eqref{takis2}, 
then, as $\ep\to 0$,  a.s. in $\omega$ and locally uniformly in $(x,t)$,  $u^\ep \to 1$ in $\{(x,t)\in \R^d\times (0,\infty): w(x,t) >0\}$ and $u^\ep \to -1$ in $\{(x,t)\in \R^d\times (0,\infty): w(x,t) < 0\}$,  that is, $u^\ep \to 1$ (resp. $u^\ep \to -1$) inside (resp. outside) a front moving with normal velocity $V=-\text{tr} [Dn] \ dt + \alpha_0 dB$.
\end{thm} 

%%\end{document}
%{\it The asymptotic behavior of the solutions to the perturbed Allen-Cahn equation with $c=dB$ and $B$ a space time Brownian motion was conjectured by Otha, Jasnow and Kawasaki \cite{othajasnowkawasaki} , while the unperturbed equation was proposed by Allen and Cahn \cite{allencahn}  as a model to study phase transitions. The rigorous justification of the conjectured behavior for the latter as well as its perturbed version with $c$ smooth was obtained by Evans, Soner and Souganidis \cite{evanssonersouganidis} and Barles and Souganidis \cite{barlessouganidis};  see Souganidis \cite{bardicime, montrealfronts} for a  comprehensive overview of the theory. For the former problem with $c=dB^\ep$ and $B^\ep$ a mild approximation of  a time Brownian motion, Yip \cite{yip} and Funaki \cite{funaki} obtained results for short time and convex initial interfaces. Lions and Souganidis \cite{lionssouganidisbook} proved the full global in time result in this case and they also showed that the behavior conjectured in \cite{othajasnowkawasaki} 
%cannot be correct if $c = dB$ and $B$ a Brownian motion in time--(see Section The explanation for the latter is that the oscillations due to the presence of $dB$ are so strong that they interfere with the stability properties of the equilibria of the potential. 
%}
%
%
%
%
%
%
Theorem~\ref{main} provides a complete characterization of the asymptotic behavior of the  Allen-Cahn equation perturbed by mild approximations of the time white noise. The result holds in all dimensions, it is global in time and does not require any regularity assumptions on the moving interface. 
\vskip.125in

In  \cite{funaki} Funaki studied  the asymptotics of \eqref{allencahn} when $d=2$ assuming  that the initial set is a smooth curve bounding a convex  set. Under these assumptions the evolving curve remains smooth and \eqref{ivp} reduces to a stochastic differential equation in the  arc length variable. Under the assumption that the evolving set is smooth, which is true if the initial set is smooth and for small time, a similar result was announced recently by Alfaro, Antonopoulou, Karali and Matano \cite{aakm}. Assuming convexity at $t=0$, Yip~\cite{yip} showed a similar result for all times using  a variational approach. There have also been several other attempts to study the asymptotics of \eqref{allencahn}  in the graph-like setting and always for small time. 

\smallskip
Reaction-diffusion equations perturbed additively by white noise arise naturally in the study of hydrodynamic limits of interacting particles. %in regimes of external magnetization such that, in the mesoscopic limit, the law of large numbers does not ``absorb'' all the randomness of  the system. 
The relationship between the long time, large space behavior of the Allen-Cahn perturbed additively by space-time white noise and fronts moving by additively perturbed mean curvature  was conjectured by Ohta, Jasnow and Kawasaki~\cite{othajasnowkawasaki}. Funaki  \cite{funaki2} obtained results in this direction when $d=1$ where there is no curvature effect.  A recent observation of Lions and Souganidis \cite{lswhitenoise} shows that the general conjecture cannot be correct. Indeed, it is shown in \cite{lswhitenoise} that   the formally conjectured interfaces, which should move by mean curvature additively perturbed with space-time white noise, are not well defined. 

\smallskip

From the phenomelogical point of view, problems like \eqref{allencahn} arise naturally in the phase-field  theory when modeling double-well potentials with depths (stochastically) oscillating in space-time around a common one. This leads to stable equilibria that are only  formally close to $\pm 1$. As a matter of fact, the locations of the equilibria  may diverge due to the strong effect of the white noise. 
\smallskip

%We note that our aim in this note is to study  a simple problem. In order to study more complicated reaction-diffusion with spatial and possibly behavior it is necessary to extend  in the stochastic setting the full methodology of \cite{bs}.

%\smallskip

The history and literature about the asymptotics of \eqref{allencahn} with or without additive continuous perturbations is rather long. I refer to \cite{barlessouganidis} for an extensive review as well as references.
\smallskip

%The main result  about the wellposedness and stability properties of the stochastic viscosity solutions of \eqref{ivp} is stated next. 
% \begin{thm}\label{takis20}
% For each $w_0\in \text{BUC}(\R^d)$ and $\zeta \in C([0,\infty);\R)$ the initial value problem \eqref{ivp} has a unique pathwise solution in $\text{BUC}(\R^d\times [0,\infty))$. Moreover, if $w_n \in \text{BUC}(\R^d\times [0,\infty))$ is the unique solution of \eqref{ivp} for a path $\zeta _n$ and $w_n(\cdot,0)=w_{n,0}$ such that, as $n\to \infty$ and locally uniformly in time and uniformly in space,  $\zeta_n\to \zeta$ and $w_{n,0} \to w_0$,  then, as $n \to \infty$, locally uniformly in time and uniformly in space, $w_n \to w$, the unique solution of \eqref{ivp} with path $\zeta.$
% \end{thm}
% 
 An important tool in the study of evolving fronts is the signed distance function to the front which is defined as 
 \begin{equation}\label{takis22}
 \rho(x,t)=\begin{cases} \rho(x,  {\{x\in \R^d: w(x,t) \leq 0\}}),\\[1mm]
 -\rho(x,  {\{x\in \R^d: w(x,t) \geq 0\}}),
 \end{cases}
 \end{equation}
 where $\rho (x,A)$ is the usual distance between a point $x$ and a set $A$.
 \smallskip
 
 When there is no interior, that is, 
 $$\partial \{x\in \R^d: w(x,t) < 0\}=\partial \{x\in \R^d: w(x,t) >0\},$$
 then
  \begin{equation*}
 \rho(x,t)=\begin{cases} \rho(x,\Gamma_t) \ \text{if} \  w(x,t) >0,\\[1.3mm]
  - \rho(x,\Gamma_t) \ \text{if} \  w(x,t) <0.
  \end{cases}
  \end{equation*}
  
 The next claim is a direct consequence of the stability properties of the pathwise solutions and the fact that a nondecreasing function of the solution is also a solution. When  $\zeta$ is a smooth path, the claim below is established in \cite{barlessonersouganidis}. The result for the general path follows by the stability of the pathwise viscosity solutions with respect to the local uniform convergence of the paths.
 %The proof follows as the analogous statement in  \cite{bss} using the stability of the pathwise solutions, that is, we consider the equations satisfied by the distance  then pass in the limit in the equations satisfied by the signed distance function with smooth paths. 
 
 \begin{thm}\label{takis21}
 Let $w \in \text{BUC}(\R^d\times [0,\infty))$ be the  solution of  \eqref{ivp} and $\rho$ the signed distance function defined by \eqref{takis22}. Then $\underline \rho= \min (\rho, 0)$ and $\overline  \rho=\max(\rho, 0)$ satisfy respectively
 \begin{equation}\label{takis23}
d \underline \rho \leq \left[\left(I-\dfrac{D \underline\rho\otimes D\underline\rho}{|D\underline \rho|^2}\right): D^2 \underline \rho \right]dt + |D\underline \rho|\circ d\zeta \leq \ \text{in} \ Q_\oo,
 \end{equation}
 and 
 \beq\label{takis24}
 d \overline \rho \geq \left[\left(I-\dfrac{D\overline\rho\otimes D\overline\rho}{|D\overline \rho|^2}\right) : D^2 \overline\rho\right] dt + |D\overline\rho|\circ d\zeta \geq 0 \ \text{in} \  Q_\oo.
\Eq
In addition,
\beq\label{takis25}
 -(D^2\underline \rho D\underline \rho, \underline \rho)\leq 0 \ \ \text{and} \ \  d\underline \rho \leq \Delta \underline \rho - d\zeta  \  \text{in} \ \{\rho < 0\}, 
 \Eq
 and
 \beq\label{takis251}
   -(D^2\overline \rho D\overline \rho, \overline \rho)\geq 0 \ \ \text{and} \ \  d\overline \rho \geq \Delta  \overline \rho - d\zeta  \ \text{in} \ \{ \rho > 0\}.
 \Eq
 \end{thm}
Following the arguments of \cite{barlessonersouganidis}, it is possible to construct global in time subsolutions and supersolutions of \eqref{allencahn} which do not rely on the regularity of the evolving fronts. In view of the stabilities of the solutions, it is then possible to conclude. %Then we use the stability properties of \eqref{ivp} to conclude. 

\smallskip

An important ingredient of the argument is the existence and properties of traveling wave solutions of \eqref{allencahn} and small additive perturbations of it, which we describe next.

\smallskip

It is well known (see, for example, \cite{barlessonersouganidis} for a long list of references) that, if $f$ satisfies \eqref{nonl}, then 
 for every sufficiently small $b$,  there exists a unique strictly increasing traveling wave solution $q=q(x,b)$  and a unique speed $c=c(b)$ of 
\beq\label{takis51}
cq_{\xi}+ q_{\xi\xi}= f(q) -b   \ \text{in} \ \R \quad 
q(\pm \infty, a)=h_{\pm} (b) \quad  q(0,a)=h_0(b),
\Eq
where $h_{-} (b) < h_0(b) < h_{+} (b) )$ are the three solutions of the algebraic equation $f(u) =b$.  Moreover, as $b\to 0$,  
\beq\label{takis51.1}
h_{\pm} (b) \to \pm 1 \ \text{ and} \  h_0(b) \to 0. 
\Eq

%In what follows we normalize $q$ so that it is increasing in $\delta$ and, moreover, $q(0,\delta)=q_0^\star(\delta)$. 

%Finally, we observe that \eqref{nonl} implies that $q_\pm^\star(0)=\pm$, $q_0^\star(0)=0$ and $c(0)=0$. 
%\smallskip
%
The results needed here are summarized  in the next lemma. For a sketch of its proof I refer to \cite{barlessonersouganidis} and the references therein.  In  what follows, $q_\xi$ and $q_{\xi \xi}$ denote first and second derivatives of $q$  in $\xi$ and $q_b$ the derivative with respect to $b$.

\begin{lem}\label{takis76}
Assume \eqref{nonl}. There exist $b_0 >0, C >0, \lambda >0$ such that, for all $|b|<b_0$, there exist a unique  $c(b)\in \R$, a unique strictly increasing $q(\cdot, b):\R\to \R$ satisfying \eqref{takis51},  \eqref{takis51.1}  and $\alpha_0 \in \R$  such that 
\beq\label{takis77}
0 <h_+(b) -q(\xi;b)\leq C e^{-\lambda |\xi|} \  \text{if} \ \xi\geq0 \ \text{and}
\ 0 <q(\xi;b)-h_{-}(b) \leq C e^{-\lambda |\xi|} \  \text{if} \ \xi 
\leq 0,
\Eq
\beq\label{takis771}
0< q_\xi (\xi;b) \leq C e^{-\lambda |\xi|}, \ |q_{\xi \xi}(\xi;b)|\leq C e^{-\lambda |\xi|} \ \text{and} \ |q_b| \leq C,
\Eq
\beq\label{takis772}
c(b)=-\frac{h_+(b)-h_-(b)}{\displaystyle \int_{-\infty}^\infty q_\xi(\xi; b)^2 d\xi },  \quad-\alpha_0:= - \frac{d c}{d b}(0)= \frac{2}{\displaystyle \int_{-1}^{1} q_\xi^2 (\xi,0)d\xi,} \quad \text{and} \quad |\dfrac {c(b)}{b} + \alpha_0 | \leq C|b|. %\ %\text{where} \ F(u):=
%\int_{u}^{m_+}f(z)dz.
\Eq
\end{lem}
\smallskip

In the proof of Theorem~\ref{main} we work with $b=\ep \dot B^\ep(t) -\ep a$ for $a\in (-1,1)$; note that, in view of \eqref{takis1.1}, for $\ep$ sufficiently small, $|b| <b_0$.  To ease the notation, I write 
$$q^\ep(\xi, t, a)=q(\xi, \ep (\dot B^\ep(t) - a)) \ \text{and} \ c^\ep(a)=c(\ep (\dot B^\ep(t) - a)),$$
and I summarize in the next lemma, without a proof, the key properties of $q^\ep$ and $c^\ep$ that we  need later.
\begin{lem}\label{takis76.1}
Assume the hypotheses of  Lemma~\ref{takis76} and \eqref{takis1.1}.  Then, there exists $C> 0$ such that 
\beq\label{takis76.2}
\lim_{\ep\to} \ep |q_t^\ep(\xi, t, a)|=0 \ \text{uniformly on $\xi$ and  $a$ and locally uniformly in  $t\in [0,\infty)$},
\Eq
\beq\label{takis76.3}
\dfrac{1}{\ep} q^\ep_\xi(\xi, t, a) + \dfrac{1}{\ep^2}|q^\ep_{\xi \xi} \xi, t, a)| \leq C e^{-C\eta/\ep} \ \text{for all $|\xi|\geq \eta$ and all $\eta> 0$},
\Eq
%and, for   all $t\geq 0$ and $\ep, |a|$ sufficiently small, 
\beq\label{takis76.4}
q^\ep_\xi \geq 0 \ \text{and} \ q^\ep_a \geq 0 \ \text{for   all $t\geq 0$ and $\ep, |a|$ sufficiently small, 
} \Eq
and
\beq\label{takis76.5}
|\dfrac{c^\ep}{\ep} + \alpha_0\ep ( \dot B^\ep (t) - a)|=\text{o}(1) \ \text{uniformly for bounded $t$ and $a$}.
\Eq
\end{lem}

 Theorem~\ref{main} is proved assuming  that $u^\ep_0$ in \eqref{allencahn} is well prepared, that is, has the form 
\beq\label{takis100.1}
u^\ep_0(x)=q^\ep(\dfrac{\rho (x)}{\ep}, 0),
\Eq
where $\rho$ is the signed distance function to $\Gamma_0$ and $q(\cdot, 0)$ is the standing wave solution of \eqref{takis51}.
\smallskip

Going from \eqref{takis100.1} to a general $u^\ep_0$ as in the statement of the theorem is standard in the theory of front propagation. It amounts to   showing  that, in a conveniently 
small time interval, $u^\ep$ can be ``sandwiched'' between functions like the ones in \eqref{takis100.1}. Since this is only technical, I omit the details and I refer to \cite{barlessouganidis} for the details.  
\smallskip

The proof of the result is a refinement of the analogous results of \cite{evanssonersouganidis} and \cite{barlessonersouganidis}. It is based on using two approximate flows, which evolve with normal velocity $V=-\text{tr}[Dn] + \alpha_0 ( \dot B^\ep(t) - \ep a)$, to construct a subsolution and supesolution \eqref{allencahn}. Since the arguments are similar, here we show the details only for the supersolution construction. 
\smallskip

For fixed $\delta, a >0$ to be chosen below and any $T>0$, consider the solution $w^{a, \delta, \ep}$ of \beq\label{aegean1}
\begin{cases}
w^{a, \delta, \ep}_t - \left(I - \widehat {Dw^{a, \delta, \ep}} \otimes \widehat {Dw^{a, \delta, \ep}} \right) : D^2w^{a, \delta, \ep}  + \alpha_0(  \dot B^\ep - a) |Dw^{a, \delta, \ep}| =0 \  \text{in} \ Q_T, \\[1mm]
 w^{a, \delta, \ep}(\cdot,0)=\rho +\delta.
 \end{cases}
\Eq 

Let  $\rho^{a, \delta, \ep}$ be the signed distance from $\{w^{a, \delta, \ep}=0\}$. It follows from Theorem~\ref{takis21} (see also Theorem~$3.1$ in \cite{barlessonersouganidis}) that 
\beq\label{aegean2}
 \rho^{a, \delta, \ep} - \Delta \rho^{a, \delta, \ep} {-} \alpha_0( \dot B^\ep -a) \geq 0 \ \text{in} \ \{\rho^{a, \delta, \ep}\>0\}.
 \end{equation}

Following the proof of Lemma~$3.1$ of \cite{evanssonersouganidis},  define 
\beq\label{aegean3}
W^{a, \delta, \ep}=\eta_\delta (\rho^{a, \delta, \ep}),
\Eq
where $\eta_\delta:\R\to \R$ is  smooth and such that, for some $C >0$ independent of $\delta$,
\beq\label{aegean4}
\begin{cases}
\eta_\delta \equiv -\delta \ \text{in} \ (-\infty, \delta/4], \quad \eta_\delta \leq -\delta/2   \ \text{in} \ (-\infty, \delta/2], \quad \eta_\delta (z)= z-\delta  \ \text{in} \ [\delta/2, \infty), \ \text{and} \\[1mm]
%\eta_\delta \leq \delta/2   \ \text{in} \ (-\infty, \delta/2], \\[1mm]
%\eta_\delta (z)= z-\delta  \ \text{in} \ [\delta/2, \infty), \\[1mm]
0\leq \eta_\delta'\leq C \ \text{and} \ |\eta_\delta''|\leq C\delta ^{-1} \ \text{on} \ \R.
\end{cases}
\Eq
Let $T^\star$ be the extinction time of $\{w^{a, \delta, \ep}=0\}$.
A straightforward modification of Lemma~$3.1$ of \cite{evanssonersouganidis} leads to the  following claim.
\begin{lem}\label{aegean4} There exists a constant $C >0$, which is independent of $\ep, \delta$ and $a$, such that
\beq\label{aegean5}
W^{a, \delta, \ep}_t - \Delta W^{a, \delta, \ep} - \alpha_0( \dot B^\ep - a) |D W^{a, \delta, \ep}|\geq -\dfrac{C}{\delta} \ \text{in} \ \R^d \times [0, T^\star],
\Eq
%\vskip-.25in
\beq\label{aegean6}
W^{a, \delta, \ep}_t - \Delta W^{a, \delta, \ep} - \alpha_0( \dot B^\ep - a) \geq 0 \ \text{in} 
\ \{\rho^{a, \delta, \ep} > \delta /2\},
\Eq
and  
\beq\label{aegean7}
|D W^{a, \delta, \ep}|=1 \ \text{in} 
\ \{\rho^{a, \delta, \ep} >\delta /2\}.
\Eq
\end{lem}
\smallskip

Finally, set
\begin{equation}\label{aegean8}
U^{a,  \delta, \ep}(x,t)=q^\ep\left(\dfrac{W^{a,  \delta, \ep}(x,t)}{\ep}, t, a\right) \text{ on } \ \R^d\times [0,\infty).
\Eq
\begin{prop}\label{athens1}
Assume \eqref{nonl}, \eqref{takis1.1} and \eqref{takis2}. Then, for every $a\in (0,1)$, $U^\ep$ is a supersolution of \eqref{allencahn} if $\ep\leq \ep_0=\ep_0(\delta,a)$ and $\delta\leq \delta_0=\delta_0(a)$.
\end{prop}
\begin{proof} Since the arguments are similar to the ones used to  prove the analogous result (Proposition~$10.2$) in  \cite{barlessonersouganidis}, here I only sketch the argument.  Note that since everything takes place  at the $\ep>0$ level, there is no reason to be concerned about anything ``rough''. Below, for simplicity, I  argue  as if $w^{\ep,\delta, a}$ had actual derivatives, and is left  up to the reader to argue in the viscosity sense. Note that, throughout the proof, $\text{o}(1)$ stands for a  function such that $\lim_{\ep\to 0} \text{o}(1)=0.$ Finally, throughout the proof  $q^\ep$ and its derivatives are evaluated at $(W^{a,\delta,\ep}/\ep,t,a)$. 
%We  point out, however, that, since we work with  $\ep> 0$, there is no need to be concerned  about anything stochastic. Thus, since the arguments are very similar to the ones used to prove Proposition~$10.2$ in \cite{bss}, here  we present only a  sketch. 
%\smallskip

\smallskip

Using the equation satisfied by $q^\ep$ gives 
\beq\label{athens100}
\begin{split}
U^{a,  \delta, \ep}_t -\Delta U^{a,  \delta, \ep}  + \dfrac{1}{\ep^2}[f(U^{a,  \delta, \ep} ) - \ep \dot B(t))=&J^\ep -\dfrac{1}{\ep^2}q^\ep_{\xi \xi}(|DW^{a,  \delta, \ep} |^2 -1)\\[1mm] 
&+\dfrac{1}{\ep}q^\ep_\xi (DW^{a,  \delta, \ep}_t-\Delta W^{a,  \delta, \ep} +\dfrac{c^\ep}{\ep}) + \dfrac{a}{\ep}, 
\end{split}
\Eq
%where $q^\ep_\xi$ and $q^\ep_{\xi \xi}$ are evaluated at $(W^{a,  \delta, \ep} (x,t)/\ep, t, a)$,
and 
\beq\label{athens101}
J^\ep(x,t)=q_b \left(\dfrac{W^{a,  \delta, \ep} (x,t)}{\ep}, \ep \dot B^\ep(t) -\ep a\right)\ep \ddot B^\ep(t).
\Eq 
In view of its definition, it is immediate that $|DW^{a,  \delta, \ep}|\leq C$ with $C$ as in \eqref{aegean4}, while it follows from Lemma~\ref{takis76}  that, as $\ep \to 0$ and uniformly in $(x,t,
\delta, a)$
\beq\label{athens102} J^\ep=\dfrac{\text{o}(1)}{\ep}.\Eq
%\smallskip

Three different cases, which depend on the relationship bewteen $\rho^{a, \delta, \ep}$ and $\delta$, need to be considered.
\smallskip

If $\delta /2< \rho^{a, \delta, \ep}<2\delta$, then \eqref{aegean6}, \eqref{aegean7}, \eqref{takis76.5} and the form of $\eta_\delta$ allow to rewrite  \eqref{athens100} as 
\beq\label{athens102}
\begin{split}
U^{a,  \delta, \ep}_t -\Delta U^{a,  \delta, \ep}  + \dfrac{1}{\ep^2}[f(U^{a,  \delta, \ep} ) - \ep \dot B^(t))&\geq - \dfrac{1}{\ep}\left[ q^\ep_\xi\right(\frac {c^\ep}{\ep} + \alpha_0 (\ep \dot B^\ep -\ep a)) + a  + \text{o}(1),]\\[1mm]
& \geq - \dfrac{1}{\ep}\left[ q^\ep_\xi \text{o}(1)+ a  + \text{o}(1)\right].
\end{split}
\Eq
%with $q^\ep_\xi$ again evaluated at $(W^{a,  \delta, \ep} (x,t)/\ep, t, a)$. 
It easily now follows that the right side of \eqref{athens102} is positive, if $\ep$ and $\delta$ are small. 

\smallskip

If $d^{a,\delta,\ep} \leq \delta/2$, the choice of $\eta_\delta$ implies that
$W^{a, \delta, \ep} \leq -\delta/2.$ 
Hence, \eqref{takis76.3} yields that, for some $C>0$, 
$$\dfrac{1}{\ep} q^\ep_\xi  +\dfrac{1}{\ep^2} |q^\ep_{\xi \xi}| \leq C e^{-C\delta/\ep}.$$
%\smallskip

Then  $|DW^{a,  \delta, \ep}| \leq C$ and \eqref{aegean5} and \eqref{aegean6} in \eqref{athens102} give 
$$U^{a,  \delta, \ep}_t -\Delta U^{a,  \delta, \ep}  + \dfrac{1}{\ep^2}[f(U^{a,  \delta, \ep} ) - \ep \dot B(t)] \leq -C(\frac{1}{\delta} + 1) e^{-C\delta/\ep} + \text{o}(1) + \dfrac{a}{\ep};$$
note that, for $\ep$ small enough the right hand side of the inequality above is positive.
\smallskip

Finally, if $\rho^{a,\delta,\ep} >\delta$, it is possible to conclude as in the previous  case using \eqref{aegean6} and \eqref{takis76.3}. % to conclude as in the previous case.

\end{proof}

The proof of the main result is sketched next.

\begin{proof}[The proof of Theorem~\ref{main}]
Fix $(x_0,t_0) \in \R^d\times [0,T^\star)$ such that $w(x_0,t_0)=-\beta<0.$ The stability of the pathwise solutions yields that, in the limit $\ep\to 0$, $\delta \to 0$ and $a\to 0$ and uniformly in $(x,t)$, $w^{a,\delta,\ep} \to w$. Thus, for  sufficiently small $\ep, \delta$ and $a$, % so that 
\beq\label{takis110}
w^{a,\delta,\ep} (x_0,t_0) <-\dfrac{\beta}{2} <0.
\Eq 

Then $U^{a,\delta,\ep}$, which is defined in \eqref{aegean7}, is a supersolution of \eqref{allencahn} for sufficiently small $\ep$ and also satisfies, in view of \eqref{takis76.4}, 
$$U^{a,\delta,\ep}(x,0) \geq q^\ep(\dfrac{\rho(x)}{\ep}, 0) \ \text{on} \  \R^d,$$
since 
$$w^{a,\delta,\ep}(x,0)=\eta_\delta(\rho(x) +\delta) \geq \rho(x).$$
The comparison of viscosity solutions of \eqref{allencahn}  then gives 
$$u^\ep \leq U^{a,\delta,\ep} \ \text{in} \ \R^d\times [0,T^\star).$$
Recall that, in view of \eqref{takis110},  $\rho^{a,\delta,\ep}(x_0,t_0)<0$, and, hence,
$$\limsup\limits_{\ep\to 0} u^\ep(x_0,t_0)\leq \limsup\limits_{\ep\to 0} U^{a, \delta, \ep}(x_0,t_0)=-1.$$
For the reverse inequality, observe  that $\hat U(x,t)=-1-\gamma$  is a subsolution of \eqref{allencahn} if $\ep$ and $\gamma >0$ are chosen sufficiently small as can be seen easily from
$$\hat U_t-\Delta \hat U +\dfrac{1}{\ep^2} (f(\hat U) + \ep \dot B^\ep)\leq C + \dfrac{1}{\ep^2}[-\gamma f'(-1) +\text{o}(1)].$$
The maximum principle then gives, for all $(x,t)$ and sufficiently small $\gamma> 0$, 
$$\liminf\limits_{\ep\to 0} u^\ep(x_0,t_0) \geq -1 -\gamma.$$
The conclusion now follows after letting $\gamma \to 0$. 
\smallskip

Finally note that a  simple modification of the argument above yields the local uniform convergence of $u^\ep$ to $-1$ in compact subsets of $\{w<0\}$.

\end{proof}
%\end{document}

%%%%%%%%%%%%%%%%%%%%%%%%%%%%%%%%%%
\section{Pathwise entropy/kinetic solutions for scalar conservation laws with multiplicative rough time signals.}
%%%%%%%%%%%

\subsection*{Introduction} Ideas similar to the ones described up to the previous sections were used by Lions, Perthame and Souganidis \cite{lionsperthamesouganidis1, lionsperthamesouganidis2}, Gess and Souganidis \cite{gesssouganidis1, gesssouganidis2, gesssouganidis3} and Gess, Perthame and Souganidis \cite{gessperthamesouganidis} to study pathwise entropy/kinetic solutions for scalar conservation laws with multiplicative rough time signals as well as their long time behavior, the existence of invariant measures and the convergence of general relaxation schemes with error estimates.
\smallskip

To keep the ideas simple the presentation here is about the simplest possible case, that is the spatially homogeneous  initial value problem 
\begin{equation}\label{scl}
%\begin{cases}
du + \displaystyle \sum_{i=1}^{d} A^i(u)_{x_i} \cdot dB_i = 0 \  \text{ in } \  Q_T  \quad u_0(\cdot,0)=u_0,%\\[2mm]
%u(\cdot,0)=u^0 \quad \text{ on } \quad \R^d,
%\end{cases}
\end{equation}
with 
\begin{equation}\label{flux}
{\bf{A}} =(A_1,...,A_d) \in C^2(\R;\R^d)
\end{equation}
and merely continuous paths
\begin{equation}\label{path1}
{\bf B }=(B_1,...,B_d) \in C([0,\infty);\R^d).
\end{equation}
If, instead of \eqref{path1}, ${\bf B} \in C^1([0,\infty);\R^d)$,  \eqref{scl} is a ``classical'' problem with a well known theory; see, for example, the books by Dafermos \cite{dafermosbook} and Serre \cite{serrebook}. The solution can develop singularities in the form of shocks (discontinuities). Hence it is necessary to consider entropy solutions which, although not regular, satisfy the $L^1$ -contraction property established by Kruzkov \cite{kruzkov}.
% that yields  uniqueness.

\smallskip
Solutions of deterministic non-degenerate conservation laws have  remarkable regularizing effects in Sobolev spaces of low order. It is an interesting question to see if they are still true in the present case. This is certainly possible with different exponents as shown in  \cite{lionsperthamesouganidis2} and  \cite{gesssouganidis2}. 
% in view of the case of kinetic equations where stochastic averaging lemmas are established \cite{}. 

\smallskip

Contrary to the Hamilton-Jacobi  equation, the approach put forward for \eqref{scl} does not work for conservation laws  with semilinear rough path dependence like
\begin{equation}\label{scl2}
%\begin{cases}
du + \displaystyle\sum_{i=1}^{d} (A^i(u))_{x_i} dt= {\bf \Phi}(u) \cdot  d {\bf {\tilde B}} \  \text{ in } \ Q_T \quad 
u(\cdot,0)=u_0,
%end{cases}
\end{equation}
for  ${\bf \Phi}=(\Phi_1,...,\Phi_m) \in C^2(\R;\R^m)$ and  an m-dimensional path ${\bf {\tilde B}}=({\tilde B}_1,...,{\tilde B}_m)$.
%and
%$${\bf \Phi}(u) \circ d {\bf {\tilde W}}=\sum_{j=1}^{m}\Phi_j \circ d{\tilde W}_j .$$

\smallskip

Semilinear stochastic conservation laws in It\^o's form like 
\begin{equation}\label{scl22}
%\begin{cases}
du + \displaystyle\sum_{i=1}^{d} (A^i(u))_{x_i} dt  = {\bf \Phi}(u) d {\bf {\tilde B}} \ \text{ in } \ Q_T\\[2mm]
%u(\cdot,0)=u^0 \quad \text{ on } \quad \R^d,
%\end{cases}
\end{equation}
have been studied  by Debussche and Vovelle \cite{debusschevovelle1, debusschevovelle2, debusschevovelle3}, Feng and Nualart \cite{fengnualard}, Chen, Ding and Karlsen \cite{chenetall},
% Debussche, Hofmanova and Vovelle \cite{debusscheetal}, 
  and Hofmanova \cite{hofmanova1, hofmanova2}).
% put forward a theory of weak entropy solutions of scalar conservation laws with Ito-type semilinear (but no stochastic quasilinear dependence), which in our setting take the form
%\begin{equation}\label{scl22}
%\begin{cases}
%du + \displaystyle\sum_{i=1}^{d} (A^i(u))_{x_i} dt  = {\bf \Phi}(u) d {\bf {\tilde W}} \quad \text{ in } \quad \R^d\times(0,\infty), \\[2mm]
%u=u^0 \quad \text{ on } \quad \R^d\times\{0\},
%\end{cases}
%\end{equation}
%with $\bf Phi$ as above and ${\bf {\tilde W}}$ an $m$-dimensional Brownian motion.
\smallskip

It turns out that pathwise solutions are natural in problems with nonlinear dependence. Indeed, let $u, v$ be solutions
of  the simple one dimensional problems
$$du + A(u)_x\cdot dB=0 \quad \text{and} \quad dv + A(v)_x\cdot dB=0.$$
Then %difference $w=u-v$ satisfies the equation
$$d(u-v)+ (A(u)-A(v))_x\cdot dB=0.$$ 
Multiplying by the $\text{sign}(u-v)$ and integrating over $\R$ formally leads to 
$$d\int_\R |u-v|dx + \int_\R (\text{sign}(u-v)(A(u)-A(v)))_x\cdot dB=0$$
and, hence,
$$d\int_\R |u-v|dx=0.$$
On the other hand, if $du=\Phi(u)\cdot dB$ and $dv=\Phi(v)\cdot dB$, then the previous argument cannot be used 
since the term $\int_R \text{sign}(u-v)(\Phi(u)-\Phi(v))\cdot dB$ is neither $0$ nor has a sign.  More about this is presented in the last subsection. 

\subsection*{The kinetic theory when $B$ is smooth} 
%The basic concepts of the kinetic theory of scalar conservation laws. We are going to show that it allows us to define a global change of variable along the ``kinetic'' characteristics, a very convenient tool for our purpose. Recall that for the conservation laws in the physical space the characteristics are only defined for short times (before crossing) and the method is not so convenient. Such a conclusion was also drawn in \cite{dv} but for a different reason. There the kinetic setting keeps better track of the entropy dissipation (due to the noise).
%
%
%\smallskip
%
To make the connection with the ``non rough'' theory, assume that ${\bf B} \in C^1((0,\infty);\R^d),$  
in which case $du$ stands for the usual derivative and $\cdot$ is the usual multiplication and, hence, should be ignored.
% \noindent Although we use the notation of the Introduction, throughout the discussion in this section 
%we assume that
%\begin{equation}\label{path3}
%{\bf W} \in C^1((0,\infty);\R^d),
%\end{equation}
%in which case $du$ stands for the usual derivative and $\circ$ is the usual multiplication and, hence, should be ignored.
%
\smallskip

\noindent The entropy inequality (see \cite{dafermosbook, serrebook}), which guarantees the uniqueness of the weak solutions, is that 
\begin{equation}\label{eq:rsentr}
%\left\{\begin{array}{l}
dS(u) + \displaystyle\sum_{i=1}^{d} (A^{i,S}(u))_{x_i} \cdot dB_i   \leq 0  \ \text{ in } \ Q_T  \quad S( u(\cdot,0))=S(u_0),
%\\ \\
%S( u)=S(u^0) \quad \text{ on } \quad \R^d\times\{0\},
%\end{array} \right.
%\label{eq:rsentr}
\end{equation}
for all $C^2$ -convex functions $S$ and fluxes ${\bf A^S}$ defined by
$$
\left( {\bf A^S}(u)\right)' = {\bf a}(u) S'(u) \quad \text{ with  } \quad  {\bf a} = {\bf A}'.
$$
%%where
%%$$
%where ${\bf a} = {\bf A}'.$
%
%\end{document}
%\smallskip
%
%In view of these inequalities, it appears that a pathwise theory is more appropriate to study \eqref{scl}. Indeed it easily follows, when the paths are smooth, that \eqref{scl} satisfies an $L^1(\R^d)$-contraction property. If in the stochastic setting, i.e., when ${\bf W}$ is actually a Brownian motion, we wanted a theory involving expectations, then Ito-calculus
%creates terms -- for simplicity here we take $N=1$-- of the form 
%$$
%E(\int_{0}^{t}S''(u)(f'(u))^2(u_x)^2 dt,
%$$  
%which cannot be handled due to the lack of appropriate estimates.
%
%\smallskip
%
It is by now well established that the simplest way to handle conservation laws is through their kinetic formulation developed through a series of papers -- see Perthame and Tadmor \cite{perthametadmor}, Lions, Perthame and Tadmor \cite{lionsperthametadmor}, Perthame \cite{perthame1, perthame2}, and Lions, Perthame and Souganidis \cite{lionsperthamesouganidis6}. The basic idea is to write a linear equation on the nonlinear function
\begin{equation}
\chi(x,\xi, t) = \chi (u(x,t), \xi)=   \left\{\begin{array}{l}
+1 \quad \text{ if } \quad 0 \leq \xi \leq u(x,t),
\\[2mm]
-1  \quad \text{ if } \quad u(x,t) \leq \xi \leq 0,
\\[2mm]
\; 0 \quad \text{ otherwise}.
\end{array} \right.
\label{eq:chi}
\end{equation}
\noindent The kinetic formulation states that using the entropy inequalities (\ref{eq:rsentr}) for all convex entropies $S$ is equivalent  to  $\chi$ solving, in the sense of distributions,
\begin{equation} %\left\{\begin{array}{l}
d \chi + \displaystyle  \sum_{i=1}^d A^i(\xi)\partial_{x_i} \chi \cdot dB_i   = \partial_\xi m dt    \  \text{ in }  \  \R^d \times \R \times (0,\infty) \quad \chi(x,\xi,0) =\chi(u_0(x), \xi), % \  \text{ on }  \  \R^d\times\R\times\{0\},
%\\ \\
%\chi =\chi(u^0(\cdot), \cdot)  \  \text{ on }  \  \R^d\times\R\times\{0\},
%\end{array} \right.
\label{eq:skf}
\end{equation}
where
\begin{equation}\label{measure}
m \  \text{ is a nonnegative bounded measure in } \  \R^d \times \R \times (0,\infty).
\end{equation}

\smallskip

\noindent At least formally, one direction of this equivalence can be seen easily. Indeed since, for all $(x,t)\in \R^d\times(0,\infty)$,
$$
S\big(u(x,t) \big) -S(0)= \int S'(\xi) \chi \big(u(x,t), \xi \big) d\xi,
$$
multiplying \eqref{eq:skf} by $S'(\xi)$ and integrating in $\xi$ leads to  \eqref{eq:rsentr}. %For the converse see \cite{LPTscal,Peuniq,PeKF}.

%\end{document}
\smallskip

The next proposition, which is  stated without proof, summarizes the basic estimates of the kinetic theory,
which hold for smooth paths and are independent of the regularity of the paths. They are  the  $L^p(Q_T)$ and $BV(Q_T)$ bounds (for all $T > 0$) for the solutions, as well as the bounds on the kinetic defect measures $m$, which imply that the latter %defect measures 
are weakly continuous in $\xi$ as measures on $Q_T$.
%-----------------------------------
\begin{prop} Assume \eqref{flux}. %and \eqref{path3}. 
The entropy solutions to \eqref{scl} satisfy, for all $t> 0$,
\begin{equation}\label{reg1}
\| u (\cdot,t) \|_{L^p(\R^d)} \leq \| u_0  \|_{L^p(\R^d)}  \quad \text{ for all } \quad p \in [1,\infty],
\end{equation}
\begin{equation}\label{reg2}
\| Du (\cdot,t) \|_{L^1(\R^d)} \leq \| Du_0  \|_{L^1(\R^d)},
\end{equation}
\begin{equation}\label{reg29}
 \ \{\xi\in \R: |\chi(x,\xi,t) > 0\} \subset  [-|u(x,t)|, |u(x,t)|] \quad \text{for all $(x,t)\in \R\times (0,\infty),$} 
%\xi|\leq \ |u|\leq \|u^0\|_\infty \quad \text{ in } \quad \{(x,\xi,t)\in \R^d\times\R\times(0,\infty):|\chi(x,\xi,t) > 0 \},
\end{equation}
\begin{equation}\label{kf1}
 \int_{0}^\infty  \int_{\R^{d}} \int_{\R} m(x,\xi,t) dx d\xi dt  \leq \frac{1}{2} \| u_0 \|^2_{L^2(\R^d)},
\end{equation}
\begin{equation}\label{kf2}
\int_{0}^\infty  \int_{\R^{d}} m(x,\xi,t)  dx \; dt  \leq \| u_0  \|_{L^1(\R^d)}  \quad \text{ for all } \quad  \xi  \in  \R,
\end{equation}
and, for all smooth test functions $\psi$,
\begin{equation}\label{kf3}
\frac{d }{d\xi} \int_{0}^\infty  \int_{\R^{d}} \psi(x,t) m(x,\xi,t)  dx \; dt  \leq \left[ \| D_{x,t} \psi\|_{L^\infty(\R^{d+1})}+ \| \psi(\cdot,0) \|_{L^\infty(\R^{d})} \right] \| u^0  \|_{L^1(\R^d)}.
\end{equation}
\label{prop:kf}
\end{prop}
%-----------------------------------
%\end{document}
\noindent The next observation is the backbone  of the theory of pathwise entropy/kinetic solutions. The reader will recognize ideas described already in the earlier parts of these notes.
% Its origin goes back to
%\cite{LSsemilinear, LShjscras2000b, LShjscras98,LShjscras98b}, where similar arguments for stochastic Hamilton-Jacobi equations form the basis of the theory of stochastic viscosity solutions.
%
\smallskip

\noindent Since the flux in \eqref{scl} is independent of $x$, it is possible to  use the characteristics associated with \eqref{eq:skf} to derive an identity which is equivalent to solving \eqref{eq:skf} in the sense of distributions. Indeed,
choose
\begin{equation}\label{rho0}
\rho_0 \in C^{\infty}(\R^d) \quad \text{such that} \quad \rho_0 \geq 0 \quad \text{ and } \quad \int_{\R^d} \rho_0(x)dx =1,
\end{equation}
and observe that
\begin{equation}
\rho(y, x,\xi,t)= \rho_0\big(y-x + {\bf a}(\xi){\bf B}(t)\big),
\label{eq:rho}
\end{equation}
where
\begin{equation}\label{rho12}
{\bf a}(\xi) {\bf B}(t):=(a_1(\xi)B_1(t), a_2(\xi)B_2(t),...,a_N(\xi)B_N(t)),
\end{equation}
solves the linear transport equation (recall that in this subsection it is assumed that ${\bf B}$ is smooth)
$$
d \rho  + \sum_{i=1}^{d} A^i(\xi)\partial_{x_i}\rho \cdot dB_i= 0 \  \text{ in } \  \R^d \times \R \times (0,\infty),
$$
and, hence,
\begin{equation}\label{rho1}
d( \rho(y, x,\xi,t) \chi(x,\xi,t))  + \sum_{i=1}^{d} A^i(\xi)\partial_{x_i}(\rho (y,x,\xi,t)\chi(x,\xi,t)) \cdot dB_i  = \rho(y, x,\xi,t) \partial_\xi m(x,\xi,t)dt.
\end{equation}
%
%\smallskip
%
\noindent Integrating \eqref{rho1} with respect to $x$ (recall that $\rho_0$ has compact support) yields that, in the sense of distributions in $\R\times(0,\infty)$,
\begin{equation}\label{rho11}
\frac{d}{dt} \int_{\R^{d}} \chi(x,\xi, t) \rho(y, x,\xi,t) dx  = \int_{\R^{d}}  \rho(y, x,\xi,t) \partial_\xi m(x,\xi,t) dx.
\end{equation}
%
%
%\smallskip
\noindent Observe that, although the regularity of the path was used to derive \eqref{rho11},  the actual conclusion does not need it. In particular, \eqref{rho11} holds for paths which are only continuous.
Moreover, \eqref{rho11} is basically equivalent to the kinetic formulation, if the measure $m$ satisfies \eqref{measure}.

\smallskip

\noindent Finally, note that \eqref{rho11}  makes sense only after integrating with respect to $\xi$ against a test function. This requires  that
${\bf a}'\in C^1(\R;\R^d)$ as long as we only use that $m$ is a measure. Indeed, integrating against a test function $\Psi$, yields
$$
\begin{array}{rl}
\displaystyle\int_{\R^{d+1}} & \Psi(\xi) \rho(y, x,\xi,t) \partial_\xi m(x,\xi,t) \ dx d\xi=
\\[2mm]
&=- \displaystyle \int_{\R^{d+1}} \Psi'(\xi)  \rho(y, x,\xi,t) \ m(x,\xi,t) \ dx d\xi \\[2.75mm]
& + \displaystyle\int_{\R^{d+1}} \Psi(\xi) (\sum_{i=1}^{d} \partial_{x_i} \rho(y,x,\xi,t) (a^i)'(\xi)B_i(t)) \ m(x,\xi,t) \ dx d\xi
\end{array}
$$
and all the terms make sense as continuous functions tested against a measure.

\smallskip

\noindent  Some (new) estimates and identities, needed for the proof of the main results of this section and derived from \eqref{rho11}, are stated next. Here   $\delta$ denotes the Dirac mass at the origin.
\begin{prop}\label{prop:new}
Assume  \eqref{flux} and $ u_0 \in (L^1\cap L^\infty\cap BV)(\R^d)$. Then, for all $t> 0$,
\begin{equation}
\frac{d}{dt} \int_{\R^{d+1}} |\chi(x,\xi, t)| dx \; d\xi =-2 \int_{\R^{d}} m(x,0,t) dx,
\label{eq:unpr1}
\end{equation}
and
\begin{equation}
\begin{array}{rl}
\displaystyle\int_{\R^{d+1} }  \displaystyle \int_{\R^{2d}} \delta(\xi - u(z,t))\ \rho(y, z,\xi,t)  \rho(y, x,\xi,t) \ m(t,x,\xi) dx dy dz d\xi
\\[2.5mm]
=\frac{1}{2} \frac{d}{dt}  \displaystyle\int_{\R^{d+1}}[\left( \displaystyle\int_{\R^{d}} \chi(x,\xi, t) \rho(y, x,\xi,t) dx\right)^2 -|\chi(y,\xi, t)| ] dy d\xi.
\end{array}
\label{eq:unpr2}
\end{equation}
\end{prop}

\begin{proof}
The first identity is classical and  is obtained  from multiplying \eqref{scl} by $\text{sign} (\xi)$ and using that the fact that 
$\text{sign} (\xi) \chi(x, \xi,t)= |\chi(x, \xi,t)|$. Notice that taking the value $\xi=0$ in  $m$ is allowed by the Lipschitz regularity in Proposition \ref{prop:kf}.

\smallskip

\noindent The proof of  \eqref{eq:unpr2} uses the regularization kernel along the characteristics \eqref{eq:rho}. Indeed,   \eqref{rho11} and the fact that $ \chi_\xi(z,\xi, t) = \delta(\xi)-\delta(\xi - u(z,t))$
yield
\begin{equation}\label{eq}
\begin{array}{rl}
\qquad \qquad  &\dfrac{1}{2} \dfrac{d}{dt}  \displaystyle\int_{\R^{d+1}} \left( \displaystyle\int_{\R^{d}} \chi(x,\xi, t) \rho(y, x,\xi,t) dx\right)^2 dy d\xi
\\[2.5mm]
& = \displaystyle \int_{\R^{d+1}} \left[ \displaystyle \int_{\R^{d}} \chi(z,\xi, t) \rho(y, z,\xi,t) dz \;  \displaystyle \int_{\R^{d}}  \rho(y, x,\xi,t) \partial_\xi m(x,\xi,t) \ dx \right]dy  d\xi
\\[2.5mm]
&=-\displaystyle \int_{\R^{d+1} }  \int_{\R^{2d}} [\delta(\xi)-\delta(\xi - u(z,t))] \rho(y, z,\xi,t) \rho(y, x,\xi,t) \  m(x,\xi,t) dz  dx  dy   d\xi
\\[2.5mm]
&=-\displaystyle\int_{\R^{d}} m(x,0,t) dx \\[2.5mm]
& + \displaystyle \int_{\R^{d+1} }  \int_{\R^{2d}}  \delta(\xi-u(z,t)) \rho(y, z,\xi,t) \rho(y, x,\xi,t) \  m(x,\xi,t) dz  dx  dy  d\xi.
\end{array}
\end{equation}
\noindent An important step in the calculation above is that, for all $\xi \in \R$,
$$
\displaystyle \int_{\R^{d} }  \int_{\R^{2d}} \chi(z,\xi, t) [ D_y\rho(y, z,\xi,t)  \rho(y, x,\xi,t)  +\rho(y, z,\xi,t) D_y \rho(y, x,\xi,t)] \  m(t,x,\xi) dz  dx  dy  =0,
$$
which follows from the observation that the integrand is an exact derivative with respect to $y$.

\smallskip

\noindent Using  \eqref{eq:unpr1} in \eqref{eq} gives  \eqref{eq:unpr2}.
\end{proof}
\subsection*{Dissipative solutions} The notion of dissipative solutions, which was studied  by Perthame and Souganidis \cite{perthamesouganidisdissipative}, is  equivalent to that of entropy solutions. The interest in them is twofold. Firstly, the definition resembles and enjoys the same flexibility as the one for viscosity solutions in, of course, the appropriate function space. Secondly, in defining them, it is not necessary to talk at all about entropies, shocks, etc..
\smallskip 

\noindent It is said that $u \in L^\infty((0,T), (L^1 \cap L^\infty)(\R^d))$ is a dissipative solution of \eqref{scl}, if, for all $\Psi \in C([0,\infty);C^\infty_c(\R^d))$ and all $\psi \in C^{\infty}_c(\R; [0,\infty))$, where the subscript $c$ means compactly supported, in the sense of distributions,
\begin{equation*}%\label{dis1}
\frac{d}{dt} \int_{\R^d} \int_{\R} \psi(k) (u-k-\Psi)_+ dxdk \leq
\int_{\R^d} \int_{\R} \psi(k) \text{sign}_+(u-k-\Psi)(-\Psi_t - \sum_{i=1}^{d} \partial_{x_i}(A^i(\Psi)) \cdot dB_i) dxdk.
\end{equation*}
%where $(\cdot)_+$ and $\text{sign}_+$ denote respectively the positive part and its derivative, and
%\smallskip
%\end{document}
\noindent To provide an equivalent definition which will allow to go around the difficulties with inequalities mentioned earlier, it is necessary to take a small detour to recall
the classical fact that, under our regularity assumptions on the flux and paths,
for any $\phi \in C^\infty_c(\R^d)$ and any $t_0> 0$, there exists $h> 0$, which depends on $\phi$, such that the problem
\begin{equation}\label{dis2}
%\begin{cases}
d \bar{\Psi} + \displaystyle \sum_{i=1}^{d} \partial_{x_i}(A^i(\bar{\Psi})) \cdot dB_i = 0 \ \text{ in } \ \R^d\times(t_0-h,t_0+h) \qquad %\[2mm]
\bar{\Psi}(\cdot t_0)=\phi, % \  \text{ on } \  \R^d\times\{t_0\},
%\end{cases}
\end{equation}
has a smooth solution given by the method of characteristics.

\smallskip

\noindent It is left up to the reader to check that the definition of the dissipative solution is equivalent to saying that, for   $\phi \in C^\infty_c(\R^d)$,  $\psi \in C^{\infty}_c(\R;[0,\infty))$  and any $t_0> 0$, there exists $h> 0$, which depends on $\phi$, such that, if $\bar{\Psi}$ and $h> 0$ are as in \eqref{dis2}, then  in the sense of distributions
\begin{equation*}%\label{dis3}
\frac{d}{dt} \int_{\R^d} \int_{\R} \psi(k) (u-k-\bar{\Psi})_+ dxdk \leq 0 \ \text{ in } \ (t_0-h,t_0+h).
\end{equation*}
%where $\bar{\Psi}$ and $h> 0$ are as in \eqref{dis2}.
%\end{document}
%%%%%%%%%%%%%%%%%%%%%%%%%%%%%%%%%%%%%%%%%%%%
\subsection*{Pathwise kinetic/entropy solutions} The following definition is motivated by the theory of pathwise  viscosity solutions.
%\label{sec:wsol}
%-------------------------------------------
%%%%%%%%%%%%%%%%%%%%%%%%%%%%%%%%%%%%%%%%%%%%
%Neither the notions of entropy and dissipative solutions nor the kinetic formulation can be used to study \eqref{scl}, since all involve either inequalities or quantities with sign which do not make sense.  
%for equations/expressions with, in principle, are nowhere differentiable functions;  see Section~2 for an extensive discussion of these and other related problems.
%We refer to \cite{LSsemilinear, LShjscras2000b, LShjscras98,LShjscras98b} for a general discussion about the difficulties encountered when attempting to use the classical weak solution approaches to study fully nonlinear stochastic pde.
%The following definition is motivated by the theory of pathwise  viscosity solutions.
%(\cite{LSsemilinear, LShjscras2000b, LShjscras98,LShjscras98b}) and \eqref{rho11} we introduce next the notion of pathwise stochastic entropy solutions for SSCL. The key fact is the observation in the middle  of the previous section.
%
%\smallskip
%
%\noindent  The basic idea of \cite{LShjscras2000b, LShjscras98,LShjscras98b} is to invert locally the characteristics of the stochastic Hamilton-Jacobi equations to eliminate the stochastic part at the level of the test functions. In our context this is done using the $\rho$'s a fact which leads to \eqref{rho11}, which does not have any stochastic terms.
%
%\smallskip
%
%\noindent  We have:

%----------------------------------
\begin{defn} 
Assume \eqref{flux} and \eqref{path1}. Then $u\in (L^1\cap L^\infty)(Q_T)$ is a pathwise kinetic/entropy solution to \eqref{scl}, if there exists a nonnegative bounded measure $m$ on $\R^d\times \R \times (0,\infty)$ such that, for all test functions $\rho$ given by \eqref{eq:rho} with $\rho_0$ satisfying \eqref{rho0}, in the sense of distributions in $\R \times (0,\infty)$,
\begin{equation}
\frac{d}{dt} \int_{\R^{d}} \chi(x,\xi, t) \rho(y, x,\xi,t) dx  = \int_{\R^{d}}  \rho(y, x,\xi,t) \partial_\xi m(x,\xi,t) dx.
\label{eq:conv}
\end{equation}
\label{def:sol}
\end{defn}
\noindent  The main result is:
%-------------------------------
\begin{thm} 
Assume \eqref{flux}, \eqref{path1} and $ u_0 \in (L^1\cap L^{\infty})(\R^d)$. For all $T >0$  there exists a  unique pathwise  entropy/kinetic  solution $u \in C\big([0,\infty); L^1(\R^d)\big)\cap L^\infty (Q_T) $ to \eqref{scl} and \eqref{reg1}, \eqref{reg2}, \eqref{kf1}, \eqref{kf2} and \eqref{kf3} hold.
In addition, any pathwise entropy solutions  $u_1, u_2\in C\big([0,\infty); L^1(\R^d)\big)$ to \eqref{scl} satisfy, for all $t >0$, the contraction property
\begin{equation}\label{cont}
\| u_2(\cdot,t) -u_1(\cdot,t) \|_{L^1(\R^d)} \leq \| u_2(\cdot,0)-u_1(\cdot,0) \|_{L^1(\R^d)}.
\end{equation}
Moreover, there exists a uniform constant $C> 0$ such that, if, for $i=1,2$, $u_i$ is the pathwise entropy/kinetic  solution to \eqref{scl} with path ${\bf B_i}$ and $u_{i,0}\in BV(\R^d)$, then $u_1$ and $u_2$ satisfy, for all $t>0$, the contraction property
\begin{equation} \label{contraction}
\begin{array}{rl}
& \| u_2(\cdot,t) -u_1(\cdot,t) \|_{L^1(\R^d)} \leq \| u_{2,0}-u_{1,0} \|_{L^1(\R^d)}\\[3mm]
& + C [\|{\bf a}\|(|u_{1,0}|_{BV(\R^d)} + |u_{2,0}|_{BV(\R^d)})|({\bf B_1}-{\bf B_2})(t)| \\ [3mm]
& +(\sup_{s\in (0,t)}|({\bf B_1 -B_2})(s)| \|{\bf a'}\|
[ \| u_{1,0}\|^2_{L^2(\R^{N})} + \| u_{2,0}\|^2_{L^2(\R^{N})}])^{1/2}].
\end{array}
\end{equation}
\label{th:mainindep1}
\end{thm}
%%--------------------------------
%
%\smallskip
%
\noindent  Looking carefully into the proof of \eqref{contraction} for smooth paths, it is possible to establish, after some approximations, an estimate similar to \eqref{contraction}, for non $BV$-data, with a rate that depends on the modulus of continuity in $L^1$ of the initial data. It is also possible to obtain an error estimate for different fluxes. The details for both  are left to the interested reader.

\subsection*{Estimates for regular paths}
%\label{sec:reg}
%%-------------------------------------------
%%%%%%%%%%%%%%%%%%%%%%%%%%%%%%%%%%%%%%%%%%%%%
Following ideas from the earlier parts of the notes,  the solution operator of \eqref{scl} may be thought of as the unique extension of the  solution operators  with regular paths. It is therefore necessary to study first  \eqref{scl}
with smooth  paths and  to obtain estimates that allow  to prove that the solutions corresponding to any regularization of the same path converge to the same limit,  which is a pathwise entropy/kinetic solution. The  intrinsic uniqueness for the latter is proved later.

\smallskip

\noindent  The key step is a new estimate, which depends only on the sup-norm of ${\bf B}$ and yields
compactness with respect to time.
%
%\smallskip
%
%\noindent  We have:
%
%%----------------------------------
\begin{thm} Assume \eqref{flux} and, for $i=1,2$, $u_{i,0}\in (L^1\cap L^\infty \cap BV)(\R^d)$. Consider two smooth paths  ${\bf B_1}$  and ${\bf B_2}$ and the corresponding solutions $u_1$ and $u_2$ to \eqref{scl}.  There exists a uniform constant  $C >0$ such that,
 for all $t > 0$, \eqref{contraction} holds.
\label{th:timeindep}
\end{thm}
%----------------------------------
The proof of Theorem~\ref{th:timeindep}, which  is long and technical, can be found in \cite{lionsperthamesouganidis1}. It combines the uniqueness proof for scalar conservation laws  based on the kinetic formulation of \cite{perthame1, perthame2} and the regularization method along the characteristics introduced for Hamilton-Jacobi equations in \cite{lionssouganidis1, lionssouganidis2, lionssouganidis3, lionssouganidis4, lionssouganidisbook}.
\subsection*{The proof of Theorem \ref{th:mainindep1}}
%\begin{proof}[The proof of Theorem \ref{th:mainindep1}]
The existence of a pathwise kinetic/entropy solution follows easily. Indeed, the estimate of Theorem \ref{th:timeindep} implies that, for every $u_0 \in (L^1 \cap L^\infty \cap BV)(\R^d)$ and for every $T> 0$, the mapping $
{\bf B} \in C([0, T];\R^d) \mapsto u \in C\big([0,T]; L^1(\R^d)\big)
$
is well defined and uniformly continuous with the respect to the norm of  $C([0,T];\R^d)$. Therefore, by density,  it has a unique extension to $C([0, T])$. Passing to the limit gives the contraction properties \eqref{cont} and  \eqref{contraction} as well as \eqref{def:sol}. Once \eqref{cont} is available for initial data in $BV(\R^d)$, the extension to general data is immediate by density.

\smallskip

\noindent The next step is to  show that pathwise kinetic/entropy satisfying \eqref{def:sol} are intrinsically unique in an  intrinsic sense. The contraction property only proves uniqueness of the solution built by the above regularization process. It is, however, possible to prove that  \eqref{eq:conv} implies uniqueness. Indeed, for $BV$-data, the estimates  in the proof of Theorem \ref{th:timeindep} only use the equality of Definition~\ref{def:sol}. From there the only nonlinear manipulation needed is to check that
$$
\begin{array}{rl}
 \dfrac{1}{ 2} \dfrac{d}{dt}&  \displaystyle\int_{\R^{d+1}} \left(\displaystyle \int_{\R^{d}} \chi(x,\xi, t) \rho(y, x,\xi,t) dx\right)^2
\\[3mm]
&=   \displaystyle\int_{\R^{d+1}} \left( \displaystyle \int_{\R^{d}} \chi(x,\xi, t) \rho(y, x,\xi,t) dx\right)
\;  \dfrac{d}{dt}  \displaystyle\int_{\R^{d+1}} \left( \displaystyle\int_{\R^{d}} \chi(x,\xi, t) \rho(y, x,\xi,t) dx\right).
\end{array}
$$
\noindent  This is justified after time regularization by convolution because it  has been assumed that solutions belong to  $C\big([0,T); L^1(\R^d)\big)$ for all $T> 0$. This fact also allows to justify that the right hand side
$$
  \displaystyle\int_{\R^{d+1}} \left(\displaystyle \int_{\R^{d}} \chi(x,\xi, t) \rho(y, x,\xi,t) dx\right) \ \displaystyle\int_{\R^{d+1}} \displaystyle \int_{\R^{d}} \chi(z,\xi, t) \rho(y, z,\xi,t) dz \;  \displaystyle\int_{\R^{d}}  \rho(y, x,\xi,t) \partial_{x_i} \chi(x,\xi,t) \ dx
$$
%\smallskip
can be analyzed by a usual integration by parts, because it is possible to incorporate a convolution in $\xi$ before forming the square. All these technicalities are standard and I omit them. The uniqueness for general data requires one more layer of approximation.
%\end{proof}
%\end{document}
%%%%%%%%%%%%%%%%%%%%%%%%%%%%%%%%%%%%
\subsection*{The semilinear problem }
%\label{sec:semilinear}
Based on the results of Section~4, it is natural to expect that the approach developed earlier will also be applicable to the semilinear problem \eqref{scl2} to yield a pathwise theory of stochastic entropy solutions. It turns out, however, that this not the case.

\smallskip

\noindent To keep things simple, here it is assumed that  $d=1$, ${\bf B} =t$ and $\bf{ \tilde B} \in  C([0,\infty);\R)$ is a single continuous path.  Consider, for $\Phi \in C^2(\R;\R)$, the problem
\begin{equation}\label{scl3}
%\begin{cases}
du + \text{div} {A}(u)dt = \Phi (u) \cdot dB \ \text{ in } \ Q_T \quad % \\[2mm]
u=u_0. % \ \text{ on } \ \R^d\times\{0\}.
%\end{cases}
\end{equation}
%for  %${\bf \Phi}=(\Phi_1,...,\Phi_m) \in C^2(\R;\R^m)$, ${\bf {\tilde W}}=({\tilde W}_1,...,{\tilde W}_m) \in
%$\Phi \in C^2(\R;\R)$ and $W \in C([0,\infty);\R)$.
%\smallskip
\noindent  Following the earlier considerations as well as the analogous problem for Hamilton-Jacobi equations,   it is assumed that, for each $v\in\R$ and $T >0$,  the initial value problem
\begin{equation}\label{ivp}
%\begin{cases}
d\Psi ={ \Phi}(\Psi) \cdot d{\tilde B} \ \text{ in } \ (0,\infty) \quad %\\[2mm]
\Psi(0)=v,
%\end{cases}
\end{equation}
has a unique solution
\begin{equation}\label{ivp1}
\Psi(v;\cdot)\in C([0,T];\R) \ \ \text{ such that, for all $t\in [0,T]$,} \ \   \Psi( \cdot,t)\in C^1(\R;\R).
\end{equation}
\noindent  According to  \cite{lionssouganidis4}, to study \eqref{scl3} it is natural to consider a change of unknown given by the Doss-Sussman-type transformation
\begin{equation}\label{ds}
u(x,t)=\Psi(v(x,t),t).
\end{equation}
%\smallskip
\noindent  Assuming for a moment that $ {\tilde B}$ and, hence, $\Psi$ are smooth with respect to $t$
and \eqref{scl3} and \eqref{ivp} have classical solutions, it follows, after a straightforward calculation, that
\begin{equation}\label{scl5}
%\begin{cases}
v_t + \text{div} { {\tilde A}} (v,t) = 0 \ \text{ in } \ Q_T
\quad % \\[2mm]
v=u_0,
%\end{cases}
\end{equation}
where $ { {\tilde A}} \in C^{0,1}(\R\times [0,T])$ is given by
%\begin{equation}\label{newA}
${ {\tilde A}'}(v,t) = {A'}(\Psi(v,t)).$
%\end{equation}
\smallskip

\noindent  Under the above assumptions on the flux and the forcing term, the theory of entropy solutions of scalar conservation laws applies to \eqref{scl5} and yields the existence of a unique entropy solution.

\smallskip

\noindent  Hence, exactly as in Section~4, it is tempting to define $u\in (L^1\cap L^\infty)(Q_T)$, for all $T> 0$, to be a pathwise entropy/kinetic solution of \eqref{scl3} if $v\in (L^1\cap L^\infty)(Q_T)$ defined, for all $T> 0$, by \eqref{ds} is an entropy solution of \eqref{scl5}.

\smallskip

\noindent  This does not, however, lead to a well-posed theory. The difficulty is best seen
%We have
%----------------------------------
%\begin{definition}
%Assume \eqref{ivp}. $u\in (L^1\cap L^\infty)(\R^d \times (0,T))$ for all $T> 0$ is a stochastic entropy solution to %\eqref{scl3} if $v\in (L^1\cap L^\infty)(\R^d \times (0,T))$ defined, for all $T> 0$, by \eqref{ds} is an entropy %solution of \eqref{scl4}.
%\label{def:sol1}
%\end{definition}
when adding a small viscosity $\nu$ to \eqref{scl3}, and, hence, considering the approximate equation
$$
u_t + \text{div} {\bf A}(u)= {\bf \Phi}(u) \cdot dB + \nu \Delta u ,
$$
and, after the transformation \eqref{ds}, the problem
\begin{align*}
v_t + \langle{\bf a}\big( \Psi(v(x,t),t)\big), Dv\rangle =   \frac{\nu}{\Psi_v(v(x,t),t)} \Delta \Psi(v(x,t),t)
= \nu \Delta v + \nu (\frac{\Psi_{vv}}{\Psi_v})(v(x,t),t) |D v|^2 .
\end{align*}
%{\bf Benoit: We never a stated a theorem like that for the solutions i.e. the limit of viscosity approximations?}

\smallskip

\noindent  If the approach based on \eqref{ds} were correct, one would expect to get, after letting $\nu \to 0$,  \eqref{scl5}. This, however, does not seem to be the case due to the lack of the necessary a priori bounds to pass to the limit.

\smallskip

\noindent  The problem is, however, not just a technicality but something deeper. Indeed  the transformation \eqref{ds} does not, in general, preserve the shocks unless, as an easy calculation shows, the forcing is linear.

\smallskip

\noindent  Assume that $d=1$ and  $B(t)=t$, let $H$ be the Heaviside step  function and consider the semilinear Burgers equation
\begin{equation}\label{sclsemi}
%\begin {cases}
u_t  + \frac{1}{2} (u^2)_x   =  \Phi(u)  \  \text{ in } \ Q_T  \quad  u_0= H, \\[2mm]
%u^0= \begin{cases}
%            1 \qquad \text{for }Êx \langle 0 , \\[2mm]
%            0 \qquad \text{for }Êx > 0 ,\\[2mm]
%        \end{cases}
%u_0= H, %1 \ \text{ if } \ x\langle0 \ \text{ and} \ 0 \ \text{ if } \ x> 0,
%\end{cases}
\end{equation}
with $ \Phi$ such that
\begin{equation}\label{phi}
 \Phi(0)=0, \quad  \Phi(1)=0, \ \text{ and } \   \Phi(u) > 0  \ \text{ for } \ u \in (0,1).
\end{equation}

\smallskip

\noindent  It is easily seen that the entropy solution of  \eqref{sclsemi} is
$$
u(x,t)= \begin{cases}
            1 \ \text{for } \ x < t/2 , \\
            0 \ \text{for } \ x > t/2 .
        \end{cases}
$$
\noindent  Next consider the transformation $u = \Psi(v,t)$ with
$
\dot  \Psi(v;t) = \Phi(\Psi(v;t)), \  \Psi(v;0) =v.
$
\smallskip

\noindent  Since, in view of \eqref{phi},
$
\Psi(0;t)=\Psi(1;t) \equiv 1, \text{ and } \Psi(v;t) >v \ \text{ for } \ v \in (0,1),
$
it follows that the flux for the equation for $v$ is
$$
\widetilde A(v,t) =  \int_0^v \Psi(w;t) dw,
$$
and the entropy solution with initial data $u_0$ is $v(x,t)= H(x- \bar x (t))$ with the Rankine-Hugoniot condition
$$
\dot{\bar x}(t) = \int_0^1 \Psi(w;t) dw >  \int_0^1 w dw = \frac 12,
$$
which shows that the shock waves are not preserved.
\smallskip

%\noindent  It is, therefore, clear that the shock waves are not preserved.
%
%\smallskip
%
\noindent  The final point is that, when $\bf B$ is a Brownian path, it  it is more natural to consider contractions in $L^1(\R^d \times \Omega)$ instead of 
$L^1(\R^d)$ a.s. in $\omega$ for \eqref{scl2}. To fix the ideas take ${\bf A}=0$ and $B$ a Brownian motion and consider the stochastic  initial value problem
\begin{equation}\label{scl10}
%\begin{cases}
du = \Phi(u) \circ dB \ \text{ in } \ (0,\infty) \quad %\\[2mm]
u(\cdot, 0)=u_0.
%\end{cases}
\end{equation}
\noindent  If $u_{1}, u_{2}$ are solutions to \eqref{scl10} with initial data
$u_{1,0}, u_{2,0}$ respectively, then, subtracting the two equations, multiplying by $\text{sign}(u_1-u_2)$, taking expectations and using It\^o's calculus, gives, for some $C> 0$ depending on bounds on $\Phi$ and its derivatives,
$$    
E \int |u_1(x,t)-u_2(x,t)| dx \leq \exp(Ct) E \int |u^0_1(x)-u^0_2(x)| dx,
$$
while it is not possible, in general, to get an almost sure inequality on $\int |u_1(x,t;\omega)-u_2(x,t,\omega)| dx$.

%%%%%%%%%%%%%%%
\appendix
\section{A brief review of the theory of viscosity solutions\\
 in the deterministic setting}
%%%%%%%%%%%%%%%%%%%%%%%%%%%%%%%%%%
%\setcounter{equation}{0}
%\renewcommand{\theequation}{A.\arabic{equation}}
%%%%%%%%%%%%%%%%%%%%%%%%%%%%%%%%%%
This is a summary of several facts  about the theory of viscosity solutions 
of Hamilton-Jacobi equations that are used in these notes. %are presented here.
At several places, an attempt is made to motivate the definitions and the arguments.
This review is very limited in scope. Good references are  the books by 
Bardi and Capuzzo-Dolceta \cite{bardibook},
Barles \cite{barlesbook}, 
Fleming and Soner \cite{flemingsoner},
the CIME notes \cite{bardicime}
and the ``User's Guide'' by Crandall, Ishii and Lions \cite{cil}. 
%for an extensive discussion of the theory and its applications.
\smallskip

\subsection*{Themethod of characteristics} Consider the initial value problem
\begin{equation}\label{eq:A0}
%\begin{cases}
u_t = H(Du,x) \ \text{ in }\  Q_T \quad 
%\noalign{\vskip6pt}
u(\cdot,0)= u_0.% \ \text{ on }\ \R^d. 
%\end{cases}
\end{equation}
The classical method of characteristics yields, 
for smooth $H$ and $u_0$, short time smooth solutions of \eqref{eq:A0}.
%\smallskip
Indeed, assume that $H,u_0\in C^2$. 
The characteristics associated with \eqref{eq:A0} are  the solutions of the system of odes
\begin{equation}\label{eq:odesA}
%\begin{cases}
\dot{X} = - D_p H(P,X), \quad
\dot{P} = D_x H(P,X), \quad
\dot{U} = (H(P,X) - \langle D_p H (P,X), P)\rangle
%\end{cases}
\end{equation}
with initial conditions
\begin{equation}\label{eq:icA}
X(x,0) = x,\quad P(x,0) = Du_0(x)\quad\text{ and }\quad U(x,0) = u_0(x)\ .
\end{equation}
The connection between \eqref{eq:A0} and \eqref{eq:odesA}  is made through the relationship 
$$U(t) = u(X(x,t),t)\quad \text{ and }\quad P(t) = Du (X(x,t),t)\ .$$
The issue is then the invertibility, with respect to $x$, of the map 
$x \mapsto X(x,t)$. 
A simple calculation involving the Jacobian of $X$ shows that 
$x\mapsto X(t,x)$  is a diffeomorphism in $(-T^*,T^*)$ with 
$$T^* = (\|D^2 H\|\, \|D^2 u_0\|)^{-1}\ .$$
\subsection*{Viscosity solutions and comparison principle} 
Passing next  to the issues of the definition and well-posedness of weak solutions, 
to keep the ideas simple,  it is convenient to consider the two simple problems
\begin{equation}\label{ivpA}	%% \label{eq:18}
%\begin{cases} 
u_t = H (Du)\dot B \ \text{in} \ Q_T \quad 
%\noalign{\vskip6pt}
u(\cdot, 0)= u_0,
%\end{cases}
\end{equation} 
and 
\begin{equation}\label{statA}	%% \label{eq:19}
u + H (Du) =f \ \text{in} \  \R^d\ ,
\end{equation}
and to assume that $H\in C(\R^d)$, $B\in C^1(\R)$ and 
$f\in BUC (\R^d)$.
\smallskip

Nonlinear first-order equations do not 
have in general smooth solutions. 
This can be easily seen with explicit examples.
On the other hand, it is natural to expect, in view of the many applications,
like control theory, front propagation, etc., that global, 
not necessarily smooth solutions, must exist for all time 
and must satisfy a comparison principle.
For \eqref{ivpA} this will mean that 
if $u_0\leqq v_0$, then $u(\cdot,t)\leqq v(\cdot,t)$ for 
all $t> 0$ and, for \eqref{statA}, if $f\leqq g$, then $u\leqq v$.
\smallskip

To motivate the definition of the viscosity solutions it is useful to  proceed, in 
a formal way, to prove this comparison principle.
\smallskip

Beginning  with \eqref{statA}, it is assumed that, for $i=1,2$, $u_i$ solves   
 \eqref{statA} with right hand side $f_i$. 
To avoid further technicalities it is  further assumed that the 
$u_i$'s and $f_i$'s are periodic in the unit cube. %$\tee^N$. 
The goal is to show that if $f_1\leqq f_2$,  then $u_1\leqq u_2$.
\smallskip

The ``classical'' proof consists of looking at $\max (u_1-u_2)$ which, 
in view of the assumed periodicity, is attained at some $x_0 \in\R^d$, that is, 
$$(u_1-u_2)(x_0) = \max (u_1 -u_2) \ .$$
%\smallskip
If both $u_1$ and $u_2$ are differentiable at $x_0$, then 
$Du_1 (x_0) = Du_2 (x_0)\ ,$
and then it follows from the equations that 
$$(u_1 - u_2)(x_0) \leqq (f_1 -f_2) (x_0)\ .$$
\smallskip

Observe that to prove that $u_1 \leqq u_2$, it is enough 
to have that 
$$u_1 + H(Du_1) \leqq f_1\quad\text{and}\quad u_2 + H(Du_2) \geqq f_2\ ,$$
that is,  it suffices  for $u_1$ and $u_2$ to be respectively a 
subsolution and  a supersolution.
\smallskip

Turning now to  \eqref{ivpA}, it is again 
assumed the data is periodic in space. 
If, for $i=1,2$, $u_i$ solves \eqref{ivpA} and 
$u_1(\cdot,0) \leqq u_2(\cdot,0)$, the aim is to 
show that, for all $t> 0$,  $u_1(\cdot,t) \leqq u_2(\cdot,t)$.
\smallskip

Fix $\delta > 0$ and let $(x_0,t_0)$ be such that 
$$(u_1- u_2) (x_0,t_0) - \delta t_0 
= \max_{(x,t)\in\R^d\times [0,T]} (u_1(x,t) - u_2(x,t) - \delta t)\ .$$

If $t_0 \in (0,T]$ and $u_1,u_2$ are differentiable at $(x_0,t_0)$, 
then 
$$D u_1(x_0,t_0) = D u_2(x_0,t_0)\quad\text{and}\quad u_{1,t}(x_0,t_0) \geq u_{2,t}(x_0,t_0) +\delta.$$
%while, if $t_0 =T$,
%$$D u_1 = D u_2\quad \text{and}\quad u_{1t} \geqq u_{2t} +\delta\ .$$
Since  evaluating the equations  at $(x_0,t_0)$  yields a contradiction, it must be that 
 $t_0 =0$, and, hence, 
$$\max_{(x,t)\in \R^d\times [0,T]} ((u_1-u_2)(x,t) - \delta t) \leqq 
\max_{\R^d} (u_1(\cdot,0) - u_2(\cdot,0)) \leqq 0\ .$$
Letting $\delta\to0$ leads to the desired conclusion.
\smallskip

The previous arguments use, of course, strongly the fact that $u_1$ and 
$u_2$ are both differentiable at the maximum of $u_1-u_2$, 
which is not the case in general.
This is a major difficulty that is 
overcome using the notion of viscosity solution, 
which relaxes the need to have differentiable solutions.
\smallskip

The definition of the viscosity solutions  for the general problems 
\begin{equation}\label{givpA}	%% \label{eq:20}
u_t = F(D^2 u, Du, u,x,t) \ \ \text{in}\ \  U\times (0,T] \ ,
\end{equation}
and
\begin{equation}\label{gstatA}	%% \label{eq:21}
 F(D^2 u,Du,u,x) = 0\quad\text{in}\quad U\ ,
\end{equation}
where $U$ is an open subset of $\R^d$, is introduced next.  

\begin{defn}
{\rm (i)}
$u\in C(U\times (0,T])$ (resp. $u\in C(U)$) is a viscosity subsolution 
of \eqref{givpA} (resp. \eqref{gstatA}), if, for all smooth test functions 
of $u-\phi$ and all maximum points $(x_0,t_0)\in U\times (0,T]$ (resp. resp. $x_0\in U$) of $u-\phi$
$$
\phi_t + F(D^2 \phi  , D\phi , u, x_0, t_0) \leqq 0\ \quad
\text{(resp. } F (D^2 \phi , D\phi , u, x_0) \leqq 0)\ .
$$
{\rm (ii)}
$u\in C(U\times (0,T])$ (resp. $u\in C(U)$) is a viscosity supersolution 
of \eqref{givpA} (resp. \eqref{gstatA}) if, for all smooth test functions $\phi$
and all minimum points $(x_0,t_0)\in U\times (0,T]$ 
(resp. $x_0\in U$) of $u-\phi$, 
$$
\phi_t + F(D^2\phi , D\phi , u, x_0,t_0) \geqq 0\ \quad
\text{(resp. } F (D^2\phi , D\phi ,u, x_0) \geqq 0)\ .
$$
{\rm (iii)} $u\in C(U\times (0,T))$ (resp. $u\in C(U)$) is a viscosity solution 
of \eqref{givpA} (resp. \eqref{gstatA}) if it is both a sub- and super-solution.
of \eqref{givpA} (resp. \eqref{gstatA}).
\end{defn}
In the definition above, maxima (resp. minima) can be either global or local.
Moreover, $\phi$ may have any regularity, $C^1$ being the least required for 
first-order and $C^{2,1}$ for second-order equations.
\smallskip

Using the definition of viscosity solution, 
it is possible to make the previous heuristic proof 
rigorous and to show the well-posedness of the solutions.
\smallskip

A general comparison result for 
 \eqref{ivpA} and 
\eqref{statA} is stated and proved next.
%The general comparison result is:
\begin{thm}
{\rm (i)} Assume $H\in C(\R^d)$, $f,g\in BUC(\R^d)$ and let 
$u,v\in BUC(\R^d)$ be respectively viscosity subsolution  and supersolution
of \eqref{statA} with right hand side $f$ and $g$ respectively. 
Then $\quad \sup_{\R^d} (u-v)_+ \leqq \sup (f-g)_+$.

{\rm (ii)} Assume $H\in\R^d$, $B\in C^1 (\R)$, $u_0,v_0\in BUC (\R^d)$
and let $u,v\in BUC (\overline Q_T)$ be respectively viscosity 
sub- and super-solutions of \eqref{ivpA} with initial data $u_0$ and $v_0$ 
respectively. 
Then, for all $t\in [0,T]$,
$\quad \sup_{\R^d} (u(\cdot,t) - v(\cdot,t))_+ 
\leqq \sup (u_0- v_0)_+$. 
\end{thm}

\begin{proof} 
To simplify the argument it is assumed throughout the proof  that $f$, $g$, $u_0$, $v_0$, $u$ 
and $v$ are periodic in the unit cube.
This assumption guarantees that all suprema in the statement are actually 
achieved and are therefore maxima. 
The general result is proved by introducing appropriate penalization at 
infinity, i.e., considering, in the case of \eqref{statA} for example, 
$\sup (u(x)- v(x) - \alpha |x|^2)$ and then letting $\alpha\to0$;
see \cite{barlesbook} and \cite{cil} for all the 
arguments and variations. 
\smallskip

Consider first  \eqref{statA}. 
The key technical step is to double the variables by introducing 
the new function
$z(x,y) = u(x) - v(y)$
which solves the doubled equation 
\begin{equation}\label{de}
z + H(D_x z) - H(-D_y z) \leqq f(x) - g(y) \ \text{in} \ \R^d\times \R^d. 
\end{equation}

Indeed if, for a test function $\phi$, $z-\phi$ attains a maximum 
at $(x_0,y_0)$, then $u(x) - \phi (x,y_0)$ and $v(y) +\phi (x_0,y)$
attain respectively a maximum at $x_0$,   and  a minimum at $y_0$. 
Therefore 
$$u(x_0) + H(D_x\phi (x_0,y_0)) \leqq f(x_0)\quad\text{and}\quad 
v(y_0) + H(-D_y\phi (x_0,y_0)) \geqq f(y_0),$$
and the claim follows by 
subtracting these two inequalities.
\smallskip

To prove the comparison result, $z$ is compared with a smooth function, which 
is ``almost'' a solution, that is, in the case at hand, a function of $x-y$.
\smallskip

It turns out that the most convenient choice is, for an appropriate 
$a_\ep$, 
$$\phi_\ep (x,y) = \frac1{2\ep} |x-y|^2 + a_\ep.$$
Indeed 
$$\phi_\ep + H(D_x \phi_\ep)) - H(-D_y\phi_\ep) - (f(x) - g(y)) 
= \frac1{2\ep} |x-y|^2 + a_\ep - (f(x) - g(y)) \geqq 0,$$
if 
$$a_\ep = \max (f-g) + \nu_\ep\quad   \text{ and }\quad  
\nu_\ep = \max \left(g(x) - g(y) - \frac1{2\ep} |x-y|^2\right);$$
note that, since $g$ is uniformly continuous, $\lim_{\ep\to0}\nu_\ep =0$. % as $\ep\to0$.
\smallskip

Let $(x_\ep,y_\ep)$ be such that 
$$z(x_\ep, y_\ep) - \phi (x_\ep,y_\ep) 
= \max_{\R^d\times\R^d} (z-\phi)\ .$$
%The definition of viscosity sub-solution and the fact $\phi$ is a 
%super-solution yields the estimate 
%$$u(x) - v(y) = z(x,y) \leqq \phi (x,y)\ .$$
Then 
$$z(x_\ep, y_\ep)+ H(D_x \phi (x_\ep,y_\ep) - H(-D_y\phi (x_\ep,y_\ep)) 
\leqq f(x_\ep) - g(y_\ep)\ .$$
On the other hand, it is known that 
$$\phi (x_\ep,y_\ep) + H(D_x\phi (x_\ep,y_\ep)) - H(-D_y\phi (x_\ep,y_\ep))
\geqq f(x_\ep) - g(y_\ep)\ .$$
It follows that 
$$z(x_\ep,y_\ep) \leqq \phi (x_\ep,y_\ep)$$
and, hence, 
$$z\leqq \phi\ \text{ in }\ \R^d\times\R^d\ .$$  

Letting $x=y$ in the above inequality yields 
$$u(x) - v(x) \leqq \phi (x,x) = a_\ep = \max (f-g) +\nu_\ep$$
and, after sending $\ep\to0$, 
$$\max (u-v) \leqq \max (f-g)\ .$$
The comparison  for  \eqref{ivpA} is proved  similary.
%The comparison principle is proved as for \eqref{statA}. 
%
%We discuss now several points related to the comparison for the initial 
%value problem \eqref{eq:19} always under the assumption that $W\in C^1$. 
%Of course it is again true that, for all $T> 0$,  
%$$\max_{\R^d\times [0,T]} (u-v) = \max_{\R^d} (u_0-v_0)\ .$$
%
In the course of the proof, however, it is not necessary to double the 
$t$-variable, since the equation is linear in the time derivative. 
This fact plays an important role in the analysis of the 
pathwise pde when $B$ is merely continuous. 
\smallskip

To this end, define the function 
$$z(x,y,t,s) = u(x,t) - v(y,s)\ ,$$
and observe, as before, that 
$$z_t-z_s  \leqq H(D_xz) \dot B (t) - H(-D_y z) \dot B(s) \ \text{ in }\ 
\R^d\times \R^d \times (0,\infty) \times (0,\infty)\ .$$
On the other hand, it is possible to show that 
$$z(x,y,t) = u(x,t) - v(y,t)$$
actually satisfies 
\begin{equation}\label{ntd}	%% \label{eq:22} 
z_t \leqq (H (D_xz) - H(-D_y z))\dot B\ \ \text{ in }\ 
\R^d\times \R^d \times (0,\infty)\ .
\end{equation}
Indeed, fix a  smooth 		%% periodic in $(x,y)$ function 
$\phi$ and let $(x_0,y_0,t_0)$ be a 
(strict) local maximum of $(x,y,t) \to z(x,y,t)  - \phi (x,y,t)$. 
Since all functions are assumed to be periodic with respect to the spatial 
variable, the penalized function 
$$u(x,t) - v(y,s) - \phi (x,y,t) - \frac1{2\theta} (t-s)^2$$
achieves a local maximum at $(x_\theta,y_\theta,t_\theta, s_theta).$ 
It follows that, as $\theta\to0$, $(x_\theta,y_\theta,t_\theta, s_\theta) 
\to (x_0,y_0,t_0,t_0)$.
\smallskip

Applying the definition to the function $u(x,t) - v(y,s)$ gives  
at $(x_\theta,y_\theta,t_\theta, s_\theta)$ 
$$\phi_t + \frac1{\theta} (t_\theta -s_\theta) 
- \frac1{\theta} (t_\theta-s_\theta) 
\leqq H(D_x \phi) \dot B (t_\theta) - H(-D_y \phi) \dot B (s_\theta)$$
and, after letting $\theta\to0$ and using the assumption that $B\in C^1$, 
at $(x_0,y_0,t_0)$,
$$\phi_t \leqq (H (D_x\phi) - H (-D_y\phi))\dot B.$$
Since 
$$\phi_\ep (x,y) = \frac1{2\ep} |x-y|^2$$ 
is a smooth supersolution of \eqref{ntd}, it follows immediately, 
after repeating an earlier argument, that 
$$u(x,t) - v(y,t) \leqq \frac1{2\ep} |x-y|^2 
+ \max_{x,y\in\R^d} (u(x,0) - v(y,0) - \frac1{2\ep} |x-y|^2)$$
and, after letting $\ep\to0$,
$$\max_{\R^d} (u(x,t) - v(x,t)) \leqq \max_{\R^d} (u(x,0)-v(x,0)))\ .$$

\end{proof}

\subsection*{Formulae for solutions} The next item in this review is the control interpretation of Hamilton-Jacobi 
equation. For simplicity here  $\dot B  \equiv 1$. 
\smallskip

Consider the controlled system of ode
$$%\begin{cases}
\dot{x}(t) = b(x(t),\alpha(t)) \quad 
%\noalign{\vskip6pt}
x(0) = x\in \R^d,$$
%\end{cases}$$
where $b: \R^d\times A \to \R^d$ is bounded and Lipschitz continuous 
with respect to $x$ uniformly in $\alpha$, $A$ is a compact  
subset of $\R^M$ for some $M$, the measurable map $t\mapsto\alpha_t \in A$ is the control, 
and $(x_t)_{t\geqq0}$ is the state variable.
\smallskip

The associated cost function is given by 
$$J (x,t,\alpha_\cdot) = \int_0^t f(x(s),\alpha(s)) dt + u_0 (x(t)),$$
where $u_0\in BUC (\R^d)$ is the terminal cost and 
$f:\R^d\times A\to\R$ is the running cost, which is also assumed to 
be bounded and Lipschitz continuous with respect to $x$ uniformly in $\alpha$.
\smallskip

The goal is to minimize 
--- one can, of course, consider maximization ---
the cost function $J$ over all possible controls. 
The value function is
\begin{equation}\label{cf}
u(x,t) = \inf_{\alpha_\cdot} J(x,t,\alpha_\cdot)\ . 
\end{equation}
The key tool to study $u$ is the dynamic programming principle, which is 
nothing more than the semigroup property. 
It states that, for any $\tau \in (0,t)$,
\begin{equation}\label{dpp}	%% \label{eq:23} 
u(x,t) = \inf_{\alpha_\cdot} \bigg[ \int_0^\tau f(x(s),\alpha (s)) ds 
+ u(x({\tau}), t-\tau)\bigg]\ .
\end{equation}
Its proof, which is straightforward, is based on the elementary observation 
that when pieced together, optimal controls and paths in $[a,b]$ and $[b,c]$ 
form an optimal path for $[a,c]$.
\smallskip

The following formal argument, which can be made rigorous 
using viscosity solutions
and test functions shows the connection between the dynamic programming and the Hamilton-Jacobi equation.
%\smallskip
\vskip.05in

Using  the dynamic programming identity, with $\tau =h$ small, yields 
$$u(x,t) \approx \inf_{\alpha_\cdot}
(hf(x,\alpha) + D_x u(x,t) \cdot b(x,\alpha) h) 
+ u(x,t) - u_t (x,t) h 
\ ,$$
and, hence, 
$$u_t + \sup_\alpha [-D_x u\cdot b(x,\alpha) - f(x,\alpha)] =0\ ,$$
that is 
$$%\begin{cases}
u_t + H(Du,x) =0\quad\text{in}\quad \R^d\times (0,\infty)\ ,
$$
where the (convex) Hamiltonian $H$ is given by the formula  
$$H(p,x) = \sup_\alpha [-\langle p, b(x,\alpha)\rangle  -f(x,\alpha)]\ .$$
%The control theoretical representation of the solution can be considerably 
%simplified if $\circW\equiv 1$ and $H(p,x) = H(p)$ is either convex or 
%concave.
%
%Indeed if $u$ is the viscosity solution  of 
%$$\begin{cases}
%u_t + H(Du) =0\quad\text{in}\quad \R^d\times (0,\infty)\ ,\\
%\noalign{\vskip6pt}
%u=u_0\quad\text{on}\quad \R^d\times \{0\}\ ,
%\end{cases}$$
%and $H$ is convex (resp. concave), and, hence, is given by 
%$$H(p) = \sup_{q\in\R^d}  [(p,q) -H^* (q)]\qquad 
%\text{(resp. } H(p) =\inf_{q\in\R^d} [(p,q)+H^*(q)]) \ ,$$
%the Lax-Oleinik formula \cite{Li3} yields the following simple formula 
%$$u(x,t) =\inf_{y\in\R^d} [u_0 (y) + \frac1t H^* (\frac{x-y}t)]\qquad
%\text{(resp. } 
%u(x,t) = \sup_{y\in\R^d} [u_0(y) - \frac1t H^* (\frac{x-y}t)]$$
%which follows from \eqref{eq:23} using Jensen's inequality.
%
%The control formula representation \eqref{cf} can be simplified considerably 
%if $H$ is independent of $x $ and convex.
Recall that, if $H: \R^d\to \R$ is convex, then  
$$H(p) = \sup_{q\in\R^d} (\langle p,q\rangle - H^*(p)),$$
where $H^*$ is the Legendre transform of $H$ defined by 
$$H^* (q) = \sup_{p\in\R^d} (\langle q,p\rangle - H(p))\ .$$
% Recall that the Legendre transform $(H^*)^* = H^{**}$ of $H$ is 
The Legendre transform $H^*$ of any continuous, not necessarily convex, 
$H: \R^d\to \R$ is convex. 
%The Legendre transform $(H^*)^* = H^{**}$ of $H^*$ is 
%the convex envelope of $H$.
%Hence, if $H$ is convex, then $H= H^{**}$. 
\smallskip

When $H:\R^d \times\R^d\to \R$ is convex, the previous discussion provides a formula for the viscosity solution of 
the Hamilton-Jacobi equation 
\begin{equation}\label{ivp100}
%\begin{cases}
u_t + H(Du,x) =0\ \text{ in }\ Q_T, \quad 
%\noalign{\vskip6pt}
u(\cdot,0)= u_0 \ \text{ on }\ \R^d.
%\end{cases}
\end{equation}
Indeed recall that 
$$H(p,x) = \sup_q [\langle p,q\rangle - H^* (q,x)]$$
and consider the controlled system
$$%\begin{cases}
\dot x(t)= q(t) \quad 
%\noalign{\vskip6pt}
x(0) = x,$$ 
%\end{cases}$$
 and the pay-off 
$${\mathcal J} (x,t,q_\cdot) = u_0 (x(t)) + \int_0^t H^* (q(s),x(s)) ds.$$
The theory of viscosity solutions (see \cite{barlesbook}, \cite{lbook}) 
%\cite{Lions book}) 
yields that 
\begin{equation}\label{cf1}
u (x,t) = \inf_{q_\cdot} \left[u_0 (x(t)) + \int_0^t H^* (q(s),x(s))\,ds \right].
\end{equation}
When $H$ does not depend on $x$, then \eqref{cf1} can be simplified 
considerably. Indeed, applying Jensen's inequality  to the representation formula 
\eqref{cf1} of the viscosity solution $u$ of 
\begin{equation}\label{ivp2}
%\begin{cases} 
u_t + H(Du) =0\ \text{ in }\ Q_T \quad 
%\noalign{\vskip6pt}
u(\cdot,0)=u_0\ \text{ on }\ \R^d,
%\end{cases} 
\end{equation}
%where $H(p) = \sup [\langle p,q\rangle- H^* (q)]$ 
yields the Lax-Oleinik formula
\begin{equation}\label{LO1} 
u(x,t) = \inf_{y\in \R^d} \left[u_0 (y) + t H^* (\frac{x-y}t)\right].
\end{equation}
A similar argument, when $H$  is concave, yields 
\begin{equation}\label{LO2} 
u(x,t) = \sup_{y\in \R^d} \left[u_0 (y) - t H^* (\frac{x-y}t)\right].
\end{equation}
The existence of viscosity solution follows either directly using 
Perron's method (see  \cite{cil}), which yields the solution as the 
maximal (resp. minimal) subsolution (resp. supersolution) or indirectly 
by considering regularizations of the equation, the most commonly used 
consisting of ``adding'' $-\ep \Delta u^\ep$ to the equation and passing 
to the limit $\ep\to 0$.  
\smallskip

A  summary follows of some of the key facts about 
viscosity solutions of the initial value problem  \eqref{ivp2}, which are used in the notes,
%\begin{equation}\label{ivp2}		%% \label{eq:24} 
%\begin{cases} 
%u_t = H(Du) \quad\text{in}\quad \R^d \times (0,\infty)\ ,\\
%\noalign{\vskip6pt}
%u= u_0\quad\text{on}\quad \R^d\times \{0\}\ ,
%\end{cases}
%\end{equation}
for $H\in C(\R^d)$ and $u_0 \in BUC (\R^d)$. 
\smallskip

The results discussed earlier yield that there exists a unique solution 
$u\in BUC(Q_T)$.
In particular, 
$u = S_H (t) u_0\ ,$
with the solution operator 
$S_H(t) : BUC (\R^d) \to BUC (\R^d)$  a strongly continuous 
semigroup, that is, for $s,t> 0$,
$$S_H (t + s) = S_H (t) S_H (s).$$
The time homogeneity of the equation also  yields, for $t> 0$,  
the identity 
$$S_H (t) = S_{tH} (1).$$
Moreover, $S_H$ commutes with translations, additions of constants and 
is order-preserving, and, hence, a contraction in the sup-norm, that is, 
$$\|(S_H (t) u - S_H (t)v)_\pm \|_\infty 
\leqq \|(u-v)_\pm\|_\infty.$$
%Recall that 
%$$r_+ = r\vee 0\quad\text{ and }\quad r_- = (-r) \vee 0\ ,$$

If $u_0 \in C^{0,1}(\R^d)$, the space homogeneity of $H$ and the 
contraction property yield, that, for all $t> 0$,  
$S_H (t) u\in C^{0,1} (\R^d)$  and, moreover, 
$$\|DS_H (t) u_0\| \leqq \|Du_0\|. $$

It also follows from the order preserving property that, for all 
$u,v\in BUC (\R^d)$ and $t> 0$, 
$$S_H (t) \max (u,v) \geqq \max( S_H (t) u, S_H (t) v)\quad\text{and}\quad 
S_H (t) \min (u,v) \leqq \min ( S_H (t) u ,S_H (t) v)\ .$$
%where 
%$$a\vee b = \max (a,b)\quad\text{ and }\quad a\wedge b = \min (a,b)\ .$$
Finally, it can be easily seen from the definition of viscosity solutions 
that, if, 
for $i\in I$, $u_i$ is a sub-(resp. super-solution), then 
$\sup_i u_i$ is a subsolution (resp. $\inf_i u_i$ a supersolution).
\smallskip

A natural question is whether there are any other  explicit formulae for the solutions of 
\eqref{ivp2}; recall that for $H$ convex/concave, the solutions  satisfy the Lax-Oleinik formula.
\smallskip

It turns out there exists another formula, known as the Hopf formula, which 
does not require $H$ to have any concavity/convexity property as long as 
the initial datum is convex/concave.
\smallskip

For definiteness, here it is  assumed that $u_0$ is convex,
and denote by $u_0^*$ its Legendre transform.  
\smallskip

It is immediate that, for any $p\in\R^d$, the function 
$u_p (x,t) =\langle p,x\rangle +tH(p)$ 
is a viscosity solution of \eqref{ivp2} and, hence, in view of the previous 
discussion,  
\begin{equation}\label{Hf}	%% \label{eq:25} 
v(x,t) = \sup_{p\in\R^d} [\langle p,x\rangle  +tH(p) - u_0^* (p)]\ ,
\end{equation}
is a subsolution of \eqref{ivp2}.
\smallskip

The claim is that, if $u_0$ is convex, then $v$ is actually a solution.
Since this fact plays an important role in the analysis, it is stated  
as a separate proposition.
\begin{prop}
Let $H\in C(\R^d)$ and assume that $u_0 \in BUC (\R^d)$ is convex. 
The unique viscosity solution $u\in BUC (Q_T)$ of 
\eqref{ivp2} is given by \eqref{Hf}.
\end{prop}
\begin{proof} 
If $H$ is either convex or concave, the claim follows using the 
Lax-Oleinik formula. 
Assume, for example,  that $H$ is convex. 
Then 
\begin{equation*}
\begin{split}
\sup_{p\in\R^d} [\langle p,x\rangle+ tH(p) - u_0^* (p)]
& = \sup_{p\in\R^d}[\langle p,x\rangle +t\sup_{q\in\R^d}((p,q) -H^*(q))- u_0^* (p)]\\
\noalign{\vskip4pt}
&=\sup_{p\in\R^d} \sup_{q\in\R^d}[\langle p,x\rangle +t(p,q) -tH^* (q) - u_0^* (p)]\\
\noalign{\vskip4pt}
& = \sup_{q\in \R^d} \sup_{p\in\R^d} [\langle p,x+tq\rangle -u_0^* (p) - tH^* (q)]\\
\noalign{\vskip4pt}
& = \sup_{q\in\R^d} [u_0 (x+tq) - tH^* (q)]
= \sup_{y\in \R^d} [u_0(y) - tH^* (\frac{y-x}t)]\ .
\end{split}
\end{equation*}
If $H$ is concave, the argument is similar, provided the min-max theorem
is used to interchange the $\sup$ and $\inf$ that appear in the  formula.
\smallskip

For the general case the first step is that the map 
$F(t) :BUC (\R^d)\to BUC(\R^d)$ defined by 
$$F(t) u_0 (x) = \sup_{p\in\R^d} [\langle p,x\rangle +tH(p) - u_0^* (p)]$$
has the semi-group property, that is, 
$$F(t+s) = F(t) F(s).$$
If $u_0$ is convex, then $F(t) u_0$ is also convex, 
since it is the sup of linear functions, and, moreover,  
$$F(t) u_0 = (u_0^* - tH)^*.$$
In view of this observation and the fact that, if $w$ is convex 
then $w=w^{**}$, the semigroup identity follows if it  shown that 
$$(u_0^* - (t+s)H)^{**} = ((u_0^* - sH)^{**} - tH)^{**}\ .$$
On the other hand, the definition of the Legendre transform, 
the min-max theorems and the fact that 
$$\sup_{x\in\R^d} \langle z,x\rangle = \begin{cases} +\infty&\text{if }\ z\ne0\ ,\\
0&\text{if }\ z=0\ ,
\end{cases}$$
yield, for $\tau> 0$, the following sequence of equalities: 
\begin{equation*}
\begin{split} 
&\\
\noalign{\vskip-12pt}
(u_0^* - \tau H)^{**} (y) 
& = \sup_{x\in\R^d} [\langle y,x\rangle  - (u_0^* - \tau H)^* (x)] 
=\sup_{x\in\R^d}[\langle y,x\rangle -\sup_{p\in\R^d}[\langle x,p\rangle +\tau H (p) - u_0^* (p)]]\\
\noalign{\vskip4pt}
& = \sup_{x\in\R^d} \inf_{p\in\R^d} [\langle y-p,x\rangle -\tau H (p) +u_0^* (p)]\\
\noalign{\vskip4pt}
& = \inf_{p\in\R^d} \sup_{x\in\R^d} [\langle y-p,x\rangle-\tau H (p) +u_0^* (p)] 
= u_0^* (y) - \tau H(p)\ .
\end{split}
\end{equation*}
It follows that
$$(u_0^* - (t+s) H)^{**} = u_0^* - (t+s) H 
= u_0^* - sH - tH = (u_0^* -sH)^{**} - tH 
= ((u_0^* -sH)^{**} - tH)^{**}\ .$$
Next it is shown that actually \eqref{Hf}  is a viscosity solution. In view of the previous discussion, 
it is only needed 
to check the super-solution property.
\smallskip

Assume that, for some smooth $\phi$, $v-\phi$ attains a minimum at 
$(x_0,t_0)$ with $t_0> 0$. 
Let $p = D\phi (x_0,t_0)$ and $\lambda = \phi_t (x_0,t_0)$.
The convexity of $v$ yields that, for all $(x,t)$ and $h\in (0,t_0)$, 
$$v (x,t_0 - h) \geqq v(x_0,t_0) + \langle p,x-x_0\rangle - \lambda h + o(h)\ .$$
Since
$$v(x_0,t_0) = F(h) v(\cdot,t_0-h) (x_0),$$
it follows that 
$$v(x_0,t_0) = F(h) (v(x_0,t_0) + \langle p,\cdot -x_0\rangle)  (x_0) 
- \lambda h + o(h)\ ,$$
and, finally, 
$$\lambda h \geqq h H(p) + o(h)\ .$$
Dividing by $h$ and letting $h\to0$ gives 
$\lambda \geqq H(p).$
\end{proof}
The above proof is a typical argument in the theory of 
viscosity solutions which has been used by Lions \cite{lionssemigroup} to give a 
characterization of viscosity solutions and Souganidis \cite{souganidis} and 
Barles and Souganidis \cite{barlessouganidis2} to prove convergence of approximations to 
viscosity solutions.
Similar arguments were also used by Lions \cite{lions}  in image processing
and Barles and Souganidis \cite{barlessouganidis} to study front propagation.
\smallskip

It is a natural question
%, which may not have anything to do with stochastic pde, 
to investigate whether the Hopf formula can be used for 
more general Hamilton-Jacobi equations with possible dependence on $(u,x)$. 
\smallskip

A first requirement for such formula to hold is that the equation must 
preserve convexity, that is,  if $u_0$ is convex, then 
$u(\cdot,t)$ must be  convex for all $t> 0$.
\smallskip

It turns out  that the general form of Hamiltonian's satisfying this 
latter property is 
$$H(p,u,x) = \sum_{j=1}^d x_j H_j (Du) + uH_0 (Du) + G(Du)\ .$$
To establish a Hopf-type formula, it is necessary to look at solutions 
starting with linear initial data, that is, for some $p\in \R^d$ and $a\in \R$,, 
$$u_0 (x) = \langle p,x\rangle +a\ .$$
If there is a Hopf-type formula, the solution $u$ starting with $u_0$ 
as above must be of the form 
$$u(x,t) = P(t)x +A(t)\quad\text{with}\quad A(0)=a\quad\text{and}\quad 
P(0)=p.$$
A straightforward computation yields that $P$  and $A$ must satisfy, for   $H = (H_1,\ldots, H_N)$,  the ode 
$$\dot P = H(P) + H_0(P)P\quad\text{and}\quad \dot A= H_0(P)A + G(P).$$
%where $H = (H_1,\ldots, H_N)$.
Whether the function 
$$\sup_{p\in\R^d} [\langle P(t),x\rangle + A(t)]$$
with $A(0) = -u_0^* (p)$ is a solution of the Hamilton-Jacobi equation 
is an open question in general.
Some special cases can be analyzed under additional assumptions on
the $H^i$'s, $H_0$, etc..

\section*{Aknowlegment}
I would like to thank Ben Seeger for his help in preparing these notes.

\bibliographystyle{plain}
\bibliography{cimenotes_articles.bib}

\end{document}